\documentclass[12pt]{amsart}
\usepackage{graphicx}
\usepackage{amssymb,amstext,amsfonts,amscd}
\usepackage{dp}
\author{Valery Alexeev and Viacheslav V. Nikulin}

\address{University of Georgia, Athens, GA 30602, USA}
\email{valery@math.uga.edu}
\address{Dept. of Pure Math., Univer. of Liverpool,XW Liverpool
  L69 3BX,~UK;\ \
Steklov Math. Inst., ul. Gubkina 8, Moscow 117966, GSP-1, Russia}
\email{vnikulin@liv.ac.uk,\ vvnikulin@list.ru}

\title[log del Pezzo surfaces of index $\le 2$]%
{Classification of log del Pezzo surfaces of index $\le 2$}
\begin{document}
\bibliographystyle{amsalpha}
\begin{abstract} 

 This is an expanded version of our work \cite{AN1}, 1988, in Russian. 
  
  We classify del Pezzo surfaces over $\bC$ with log terminal
  singularities (equivalently, quotient singularities) of index $\le
  2$.  By classification, we understand a description of the
  intersection graph of all exceptional curves on an appropriate
  (so-called ``right'') resolution of singularities
  together with the subgraph of the curves which are contracted to
  singular points.
  
  The final results are similar to classical results about
  classification of non-singular del Pezzo surfaces and use the usual
  finite root systems. However, the intermediate considerations use
  the theory of K3 surfaces (especially of K3 surfaces with
  non-symplectic involutions) and the theory of reflection groups in
  hyperbolic spaces (especially of reflection groups of hyperbolic
  lattices).

As an ``elementary'' application, our results permit one to classify
sextics in $\bP^2$ with simple singularities and a component of geometric 
genus $\ge 2$.   
\end{abstract}
\maketitle
\tableofcontents

\section*{Introduction}
\label{introduction}

This work is based on our paper \cite{AN1} in Russian. See also
\cite{AN2} for a short exposition of the same results.

The main purpose of this work is to classify del Pezzo surfaces
with log-terminal singularities of index $\le 2$. For this
classification, it is important to use K3 surfaces and K3 surfaces
theory. So, the most part of this work is also devoted to K3
surfaces and can also serve as an introduction to K3 surfaces.

\subsection{Historical remarks and our main principle of
classification of log del Pezzo surfaces of index $\le 2$}
\label{subsec:intro,history}
In this work, we {\it consider algebraic varieties over
the field $\bC$ of complex numbers.}
Further we usually don't mention that again.

A complete algebraic surface $Z$ (over $\bC$) with log terminal
singularities is a del Pezzo surface if its anticanonical divisor
$-K_Z$ is ample.
A 2-dimensional log terminal singularity over $\bC$ is a singularity
which is analytically equivalent to a quotient singularity $\bC^2/G$, where
$G\subset \GL(2,\bC)$ is a finite subgroup. The index $i$ of $z\in
Z$ is the minimal positive integer for which the divisor $iK_Z$ is
a Cartier divisor in a neighbourhood of $z$.

The aim of this work is to classify del Pezzo surfaces with log
terminal singularities (or simply \emph{log del Pezzo surfaces}) of
index $\le2$.

Log del Pezzo surfaces of index $\le 2$ include classical cases of
non-singular del Pezzo surfaces and log del Pezzo surfaces of index
$1$, i. e. Gorenstein log del Pezzo surfaces. First let us consider
classical results about these del Pezzo surfaces.

In 1849, Cayley \cite{Cayley1} and Salmon \cite{Salmon} discovered
$27$ lines on a non-singular cubic surface $Z$. Now we know that
they are all exceptional curves on a non-singular del Pezzo
surface $Z$ of the degree 3, and they are crucial for its
geometry. Here the degree $d$ of a del Pezzo surface $Z$ is
$d=(K_Z)^2$.

Classification of nonsingular del Pezzo surfaces is well known,
and they are classical examples of rational surfaces (see, e.g.
\cite{28,6,7}). A connection between nonsingular del Pezzo
surfaces and reflection groups was noticed a long time ago.
Schoutte \cite{Schoutte} noted that there is an
incidence-preserving bijection between 27 lines on a smooth cubic
and vertices of a certain polytope in $\bR^6$. In modern
terminology, this polytope is the convex hull of an orbit of
reflection group $W(E_6)$. Coxeter \cite{Coxeter28} and Du Val
\cite{DuVal33} noted a similar correspondence between
$(-1)$-curves on del Pezzo surfaces of degree 2 and 1 and
reflection polytopes for groups $W(E_7)$ and $W(E_8)$.

Du Val was the first to investigate the relationship between
reflection groups and singular surfaces. In \cite{20a} he
introduced Du Val singularities.  Possible singularities of cubic
surfaces $Z_3\subset \bP^3$ were classified by Schl\"afli
\cite{Schlafli} and Cayley \cite{Cayley2}. In \cite{20} Du Val
found all possible configurations of Du Val singularities on the
``surfaces of del Pezzo series'' of degree 2 and 1, i.e.  double
covers $Z_2\to \bP^2$ ramified in a quartic and double covers
$Z_1\to Q$ over a quadratic cone ramified in an intersection of
$Q$ with a cubic.  As was proved much later
\cite{18,HidakaWatanabe81}, these are precisely the Gorenstein log
del Pezzo surfaces of degree 2 and 1.

Du Val observed the following amazing fact: the configurations of
singularities on del Pezzo surfaces $Z_d$ of the degree $d$ with
Du Val singularities are in a one-to-one correspondence with
subgroups generated by reflections (i.e. root subsystems) of a
reflection group of type $E_{9-d}$, i.e.
$E_8,E_7,E_6,D_5,A_4,A_2+A_1$ respectively for $d=1,\dots,6$, with
four exceptions: $8A_1$, $7A_1$, $D_44A_1$ for $d=1$ and $7A_1$
for $d=2$.  (These days, we know that the prohibited cases do
appear in characteristic 2.)  He also noted that in some cases
(for example $4A_1$ in $E_8$) there are two non-conjugate ways to
embed a subgroup and, on the other hand, there are two distinct
deformation types of surfaces.

The proof was by comparing two long lists. The reflection
subgroups were conveniently classified by Coxeter \cite{Coxeter34}
in the same 1934 volume of Proceedings of Cambridge Philosophical
Society. Du Val went through all possibilities for quartics on
$\bP^2$ and sextic curves on the quadratic cone $Q$ and computed
the singularities of the corresponding double covers $Z_d$,
$d=1,2$. The modern explanation for the fact that configurations
of singularities correspond to some reflection subgroups is
simple: $(-2)$-curves on the minimal resolution $Y$ of a
Gorenstein del Pezzo $Z$ lie in the lattice $(K_{Y})^{\perp}$
which is a root lattice of type $E_{9-d}$.

In the 1970s, Gorenstein del Pezzo surfaces attracted new attention in
connection with deformations of elliptic singularities, see
\cite{Looijenga77,Pinkham77,BruceWall79}. The list of possible
singularities was rediscovered and reproved using modern methods, see
\cite{HidakaWatanabe81,30,15,Furushima86}.

In addition, Demazure \cite{18} and Hidaka-Watanabe
\cite{HidakaWatanabe81} established a fact which Du Val
intuitively understood but did not prove, lacking modern
definitions and tools: the minimal resolutions $Y_d$ of Gorenstein
log del Pezzo surfaces $Z_d\ne \bP^1\times \bP^1$ are precisely the
blowups of $9-d$ points on $\bP^2$ in ``almost general position'',
and $Z_d$ is obtained from such blowup by contracting all
$(-2)$-curves.

\vskip0.5cm

In addition to clarifying, unifying and providing new results for
the index 1 case, our methods are general enough to obtain
similar results in the much more general case of log del Pezzo
surfaces $Z$ of index $\le 2$. Thus, we admit log terminal
singularities of index 1 and index 2 as well. Classification
of the much larger class of log del Pezzo surfaces of index $\le
2$ (together with the described above classical index-1 case) is
the subject of our work.

{\it By classification, we understand a description of the graphs
of all exceptional curves (i. e. irreducible with negative
self-intersection) on an appropriate resolution of singularities
$\sigma:Y\to Z$, together with the subset of curves contracted by
$\sigma$.} We call $\sigma$ the {\it right resolution.} See
Section \ref{subsec:introduction,classification} below for the
precise definition. For Gorenstein, i. e. of index 1
singularities, $\sigma$ is simply the minimal resolution.

Thus, in principle of classification we follow the classical
discovery by Cayley \cite{Cayley1} and Salmon \cite{Salmon} of
$27$ lines on a non-singular cubic surface which we have mentioned
above.

The graph of exceptional curves provides complete information about
the surface. Indeed, knowing the dual graph of exceptional
curves on $Y$, we can describe all the ways to obtain $Y$ and $Z$ by
blowing up $Y\to \overline{Y}$ from the relatively minimal rational
surfaces $\overline{Y}=\bP^2$ or $\bF_n$, $n=0,2,3\dots$.  Images of
exceptional curves on $Y$ give then a configuration of curves on
$\overline{Y}$ related with these blow ups.  Vice versa, if one starts
with a ``similar'' configuration of curves on $\overline{Y}$ and
performs ``similar'' blowups then the resulting surface $Y$ is
guaranteed to be the right resolution of a log del Pezzo surface $Z$
of index $\le2$, by Theorem \ref{thm3.5.2.a}.

In the singular case of index 1, we add to the classical results
which were described above a description of all graphs of exceptional
curves on the minimal resolution of singularities. In the case of
$\Pic Z = \bZ$ this was done by Bindschadler, Brenton and Drucker
\cite{15}.

In Section \ref{subsec:introduction,K3} we give more detailed
information about our methods and results, and in Section
\ref{subsec:introduction,classification} we formulate our final
classification results about log del Pezzo surfaces of index $\le
2$.

\medskip

\subsection{Classification of
log del Pezzo surfaces of index $\le 2$ and K3 surfaces}
\label{subsec:introduction,K3} The main method for obtaining our
classification of log del Pezzo surfaces of index $\le 2$ is to
reduce it to a classification of K3 surfaces with non-symplectic
involution and to K3 surfaces theory. The main points of the
latter are contained in \cite{8,N1,10,11,31}.

In Chapter~\ref{sec1}, we show that for log del Pezzo surfaces $Z$ of
index $\le2$ the linear system $|-2K_Z|$ contains a nonsingular curve,
and that there exists an appropriate (``right'') resolution of
singularities $\sigma:Y\to Z$ for which the linear system $|-2K_Y|$
contains a nonsingular divisor $C$ (i.e. $Y$ is a right DPN surface)
such that the component of $C$ that belongs to $\sigma^*|-2K_Z|$ has
genus $\ge 2$ (i.e. the DPN surface $Y$ is of elliptic type).

In Chapter~\ref{sec2}, following \cite{8,N1,10,11,31}, we build a
general theory of DPN surfaces $Y$. Here, we use the fact that the
double cover $X$ of $Y$ branched along $C$ is a K3 surface with a
non-symplectic involution $\theta$. In this way, the
classification of DPN surfaces $Y$ and DPN pairs $(Y,C)$ is
equivalent to the classification of K3 surfaces with
non-symplectic involution $(X,\theta)$. The switch to K3 surfaces
is important because it is easy to describe exceptional curves on
them and there are powerful tools available: the global Torelli
Theorem \cite{13} due to Piatetsky-Shapiro and Shafarevich, and
surjectivity of the period map \cite{5} due to Vik. Kulikov.

In Chapter~\ref{sec3}, we extend this theory to the classification
of DPN surfaces $Y$ of elliptic type, i.e. when one of the
components of $C$ has genus $\ge2$, by describing dual diagrams of
exceptional curves on $Y$. See Theorems \ref{thm3.5.1},
\ref{thm3.5.2} and \ref{thm3.5.2.a}. In Section \ref{subsec3.6},
we give an application of this classification to a classification
of curves $D$ of degree 6 on $\bP^2$ (and $D\in |-2K_{\bF_n}|$ as
well) with simple singularities in the case when one of components
of $D$ has geometric genus $\ge2$.

In obtaining results of Chapters~\ref{sec2} and \ref{sec3}, a big
role is played by the arithmetic of quadratic forms and by reflection
groups in hyperbolic spaces which are very important in the theory of K3
surfaces. From this point of view,
the success of our classification hinges mainly on the fact that we
explicitly describe some hyperbolic quadratic forms and their
subgroups generated by all reflections
(2-elementary even hyperbolic lattices of small rank, see
Theorem~\ref{thm3.1.1}). These computations are also important by
themselves for the arithmetic of quadratic forms.

In Chapter~\ref{sec4}, the results of
Chapters~\ref{sec1}---\ref{sec3} are applied to the classification
of log del Pezzo surfaces of index $\le2$. In particular, we show
that there are exactly 18 log del Pezzo surfaces of index 2 with
Picard number 1. For completeness, we also included the list of the
isomorphism classes in the index 1 Picard number 1 case.
This list, which for the most difficult degree 1 case can be
deduced from \cite{MirandaPersson}, is skipped or given with some
inaccuracies in other references.

In Section~\ref{subsec4.2}, following \cite{15}, we give an
application of our classification to describe some rational
compactifications of certain affine surfaces. In
Section~\ref{subsec4.3.1}, we give formulae for the dimension of
moduli spaces of log del Pezzo surfaces of index $\le 2$.

In Section~\ref{subsec2.1a} we review results about K3 surfaces over
$\bC$ which we use. In Appendix (Chapter \ref{sec:appendix}), for
readers' convenience, we review known results about lattices,
discriminant forms of lattices, non-symplectic involutions on K3
which we use (see Sections \ref{subsec:discrforms}---
\ref{subsec:Witt'stheorem}). For instance, in Section
\ref{subsec:maininv} we review classification of main invariants
$(r,a,\delta)$ (see below) of non-symplectic involutions on $K3$
and their geometric interpretation which are very important in
this work. In Section \ref{fundchambers} we give details of
calculations of fundamental chambers of hyperbolic reflection
groups which were skipped in the main part of the work. They are
very important by themselves. Thus, except for some standard
results from Algebraic Geometry (mainly about algebraic surfaces),
and reflection groups and root systems, our work is more or less
self-contained.

\medskip

\subsection{Final classification results about log del Pezzo
surfaces of index $\le 2$}
\label{subsec:introduction,classification}
Below, we try to give an explicit and as elementary exposition as possible
of our final results on classification of log del Pezzo
surfaces of index $\le 2$. In spite of importance of K3 surfaces,
in final classification results, K3 surfaces disappear,
and it is possible to formulate all results in terms of only del Pezzo
surfaces and appropriate their non-singular models
which are DPN surfaces.

Let $Z$ be a log del Pezzo surface of index $\le 2$.  Its
singularities of index 1 are Du Val singularities classified by their
minimal resolution of singularities. They are described by Dynkin
diagrams $\bA_n$, $\bD_n$ or $\bE_n$, with each vertex having weight $-2$.
Singularities of $Z$ of index $2$ are singularities $K_n$ which have
minimal resolutions with dual graphs shown below:

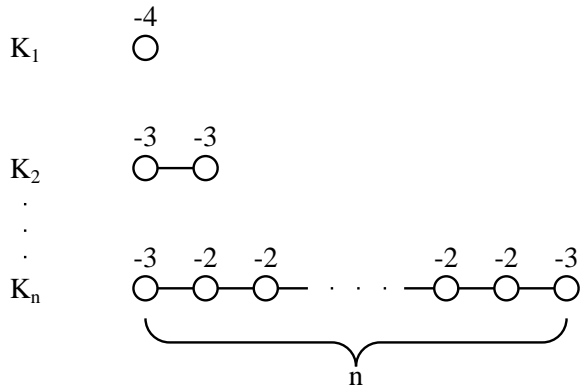
\centerline{\includegraphics[width=8cm]{pics/c-p55.eps}}

\noindent
To get the right resolution of singularities $\sigma:Y\to Z$,
additionally one has to blow up all points of intersection of components
in preimages of singular points
$K_n$. Then the right resolution of a singular point $K_n$ is described
by the graph


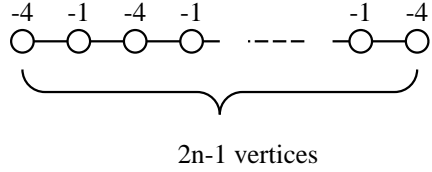
\centerline{\includegraphics[width=6cm]{pics/c-p60.eps}}

\noindent
In these graphs every vertex corresponds to an irreducible
non-singular rational curve $F_i$ with $F_i^2$ equal to the weight of
the vertex. Two vertices are connected by an edge if $F_i\cdot F_j=1$,
and are not connected if $F_i\cdot F_j=0$. Thus, the right
resolution of singularities $\sigma:Y\to Z$ of a log del Pezzo surface
of index $\le 2$ consists of minimal resolutions of singular points of
index 1 and the right resolutions shown on the graphs above of
singular points $K_n$ of index 2.

{\it Our classification of log del Pezzo surfaces $Z$ of index $\le 2$
  and the corresponding DPN surfaces of elliptic type which are
  right resolutions of singularities of $Z$ (i. e. they are appropriate
  non-singular models of the del Pezzo surfaces) is contained in Table
  3 (see Section \ref{subsec3.5}).}

All cases of Table 3 are labelled by a number $1\le N \le 50$.  For
$N=7,\,8,\,9,\,10,\,20$ we add some letters and get cases: 7a,b,
8a--c, 9a--f, 10a--m, 20a--d. Thus, altogether, Table 3 contains
$$
50+(2-1)+(3-1)+(6-1)+(13-1)+(4-1)=73
$$
cases.

The labels $N=1, \dots, 50$ enumerate the so-called {\it main
  invariants} of log del Pezzo surfaces $Z$.  They are triplets
$(r,a,\delta)$ (equivalently $(k=(r-a)/2,g=(22-r-a)/2,\delta)$) where
$r$, $a$, $\delta$ are integers: $r\ge 1$, $a\ge 0$,
$\delta\in\{0,1\}$, $g\ge 2$, $k\ge 0$.  Thus, there exist exactly 50
possibilities for the main invariants $(r,a,\delta)$ (equivalently,
$(k,g,\delta)$) of log del Pezzo surfaces of index $\le 2$.

The main invariants have a very important geometric meaning. Any log
del Pezzo surface $Z$ of index $\le 2$ and its right resolution of
singularities $Y$ are rational. The number
$$
r=\rk \Pic Y
$$
is the Picard number of $Y$, i. e.  $\Pic Y=\bZ^r$. We prove that
$|-2K_Z|$ contains a non-singular irreducible curve $C_g$ of genus
$g\ge 2$ which shows the geometric meaning of $g$. This is equivalent
to saying that there is a curve
$$
C=C_g+E_1+\cdots +E_k\in |-2K_Y|,
$$
where $E_i$ are all exceptional curves on $Y$ with $(E_i)^2=-4$.
(The inequality
$g\ge 2$ means that $Y$ is of {\it elliptic type}).  All these curves
$E_i$ come from the right resolution of singularities of $Z$ described
above.  Thus, the invariant $k$ equals the number of exceptional
curves on $Y$ with square $-4$. All of them are nonsingular and
rational.  E.g., $k=0$ if and only if $Z$ is Gorenstein and all of its
singularities are Du Val.  See Chapter \ref{sec1}.

Let us describe the invariant $\delta\in\{0,1\}$.  The components
$C_g$, $E_1,\dots,$ $E_k$ are disjoint. Since $C$ is divisible by 2 in
$\Pic Y$, it defines a double cover $\pi:X\to Y$ ramified in $C$. Let
$\theta$ be the involution of the double cover. Then the set of fixed
points $X^\theta=C$. Here, $X$ is a K3 surface and
\begin{equation}
\delta=0 \iff X^{\theta} \sim 0\mod 2\ \text{ in } H_2(X,\bZ)\
\iff
\label{intdeltaX}
\end{equation}
there exist signs $(\pm)_i$ for which
\begin{equation}
\dfrac{1}{4}\sum_{i}{(\pm)_i cl(C^{(i)})}\in \Pic Y,
\label{intdeltaY}
\end{equation}
where $C^{(i)}$ are all irreducible components (i. e. $C_g$,
$E_1,\dots,E_k$) of $C$. 
The connection to K3 surfaces with non-symplectic
involution is the main tool of our classification, see Section 
\ref{subsec:introduction,K3} above.  

As promised, our classification describes all intersection (or dual)
graphs $\Gamma (Y)$ of exceptional curves on $Y$ and also shows
exceptional curves which must be contracted by $\sigma:Y\to Z$ to get the
log del Pezzo surface $Z$ of index $\le 2$ from $Y$.
All these graphs can be obtained
from graphs $\Gamma$ in the right column of Table 3 of the same
main invariants $(r,a,\delta)$. Let us describe this in more details.

All exceptional curves $E$ on $Y$ are irreducible non-singular and
rational. They are of three types:

\begin{enumerate}
\item $E^2=-4$, equivalently $E$ is a component of genus 0 of $C\in
  |-2K_Y|$.  In graphs of Table 3 these correspond to double
  transparent vertices;

\item $E^2=-2$. In graphs of Table 3 these correspond to black
  vertices;

\item $E^2=-1$ (the 1st kind). In graphs of Table 3 these correspond to
  simple transparent vertices.
\end{enumerate}

All exceptional curves $E_i$ with $(E_i)^2=-4$, $i=1,\dots k$,
together with all exceptional curves $F$ of the 1st kind such that there
exist two different curves
$E_i$, $E_j$, $i\not=j$, with $(E_i)^2=(E_j)^2=-4$ and
$F\cdot E_i=F\cdot E_j=1$ define the {\it logarithmic
part $\Log \Gamma (Y)\subset \Gamma(Y)$ of $Y$.}
Since $F\cdot (-2K_Y)=F\cdot C=2$,
the curves $F$ are characterized by the property $C_g\cdot F=0$.
The logarithmic part $\Log \Gamma (Y)$ can be easily seen on graphs $\Gamma$
of Table 3: curves $E_i$, $i=1,\dots, k$, are shown as double transparent,
the curves $F$ of the first kind of $\Log \Gamma (Y)$
are shown as simple transparent vertices connected by two
edges with (always two) double transparent vertices. This part of $\Gamma$
is denoted by $\Log \Gamma$ and is also called {\it the logarithmic part of
$\Gamma$.} Thus, we have:
\begin{equation}
\Log (\Gamma(Y))=\Log \Gamma
\label{intlog}
\end{equation}
(with the same main invariants $(r,a,\delta)$).  {\it The logarithmic
  part $\Log \Gamma (Y)$ gives precisely the preimage of singular points
  of $Z$ of index two.}

All exceptional curves $E$ on $Y$ with $E^2=-2$ define the {\it Du Val
  part} $\Duv \Gamma(Y)\subset \Gamma(Y)$ of $\Gamma (Y)$.  Its
connected components are Dynkin graphs $A_n$, $D_n$ or $E_n$ and
they correspond to all Du Val singularities of $Z$. Thus {\it the Du
  Val part $\Duv \Gamma (Y)\subset \Gamma(Y)$ gives precisely the
  preimage of all Du Val (i. e. of index one) singular points of $Z$.}
The Du Val part $\Duv \Gamma$ of a graphs $\Gamma$ of Table 3 is
defined by all its black vertices. We have:
\begin{equation}
D=\Duv \Gamma(Y)\subset \Duv \Gamma
\label{intduv}
\end{equation}
(for the same main invariants $(r,a,\delta)$). \emph{Any subgraph $D$
of $\Duv\Gamma$ can be taken.}

Let us describe the remaining part of $\Gamma (Y)$. Each graph
$\Gamma$ of Table 3 defines a lattice $S_Y$ in the usual way. It is
$$
S_Y=\left(\bigoplus_{v\in V(\Gamma)}{\bZ e_v}\right)/\Ker
$$
defined by the intersection pairing: $e_v^2=-1$, if $v$ is simple
transparent, $e_v^2=-2$, if $v$ is black, $e_v^2=-4$, if $v$ is double
transparent, $e_v\cdot e_{v^\prime}=m$ if the vertices
$v\not=v^\prime$ are connected by $m$ edges. Here $\oplus$ means the direct 
sum of $\bZ$-modules, and  ``$\Ker$'' denotes the
kernel of this pairing. We denote $E_v=e_v\mod \Ker$.  In all cases
except trivial cases $N=1$ when $Y=\bP^2$, $N=2$ when $Y=\bF_0$ or $\bF_2$,
$N=3$ when $Y=\bF_1$,  $N=11$ when $Y=\bF_4$,
the lattice $S_Y$ gives the Picard lattice of $Y$.

Thus, $\Log \Gamma(Y)=\Log \Gamma$ and $D=\Duv
\Gamma(Y)\subset \Duv \Gamma$ define divisor classes $E_v$, $v\in
V(\Log \Gamma(Y)\cup \Duv \Gamma(Y))$, of the corresponding
exceptional curves on $Y$. Each exceptional curve $E$ is evidently
defined by its divisor class.

Black vertices $v\in V(\Duv \Gamma)$ define roots $E_v\in S_Y$
with $E_v^2=-2$ and define reflections $s_{E_v}$ in these roots
which are automorphisms of $S_Y$ such that $s_{E_v}(E_v)=-E_v$ and
$s_{E_v}$ gives identity on the orthogonal complement $E_v^\perp$
to $E_v$ in $S_Y$. These reflections $s_{E_v}$, $v\in V(\Duv
\Gamma)$, generate a finite Weyl group $W\subset O(S_Y)$.

The remaining part
$$
\Var \Gamma (Y)=\Gamma (Y)-(\Log \Gamma (Y)\cup \Duv \Gamma (Y))
$$
(it is called the {\it varying part of $\Gamma (Y)$})
is defined by the {\it varying part}
$$
\Var \Gamma=\Gamma-(\Duv \Gamma\cup \Log \Gamma)
$$
of the graph $\Gamma$ of Table 3. Further, we identify exceptional
curves $v\in V(\Gamma (Y))$ with their divisor classes $E_v\in S_Y$. We
have
\begin{equation}
V(\Var \Gamma (Y))=\{E\in W\left(\{E_v\ |\ v\in V(\Var \Gamma)\}\right)\
|\ E\cdot D\ge 0\}\subset S_Y.
\label{intvar}
\end{equation}
Here $E\cdot D\ge 0$ means $E\cdot E_i\ge 0$ for any $E_i\in D$.
The intersection pairing on $S_Y$ then defines the full graph
$\Gamma (Y)$ of $Y$. This completes description of possible graphs
$\Gamma (Y)$ of exceptional curves of log del Pezzo surfaces $Z$
of index $\le2$.

Thus, to find all possible graphs $\Gamma (Y)$ of exceptional
curves of $\sigma:Y\to Z$, one has to choose one of the graphs
$\Gamma$ of Table 3 (this also defines main invariants
$(r,a,\delta)$ of $Y$ and $Z$), then one has to choose a subgraph
$D=\Duv \Gamma (Y)\subset \Duv \Gamma$. Then $\Gamma (Y)$ consists
of $D$, $\Log \Gamma (Y)=\Log \Gamma$ and the remaining part $\Var
\Gamma (Y)$ defined by \eqref{intvar}, elements in the $W$-orbits
of $\Var\Gamma$ that have non-negative intersection with the Du
Val part. See Theorems \ref{thm3.5.1}, \ref{thm3.5.2},
\ref{thm3.5.2.a}  and \ref{thm4.1.1}. See Section
\ref{subsecN=20en} about such type of calculations in the most
non-trivial case $N=20$.

We note two important opposite cases.

{\it Extremal case.} This is the case when $D=\Duv \Gamma (Y)=\Duv \Gamma$.
Then $\Gamma (Y)=\Gamma$ is completely calculated in Table 3.
This case is called {\it extremal} and gives log del Pezzo surfaces $Z$
with Du Val singularities of the highest rank,
respectively $\rk \Pic Z=r-\#V(\Log \Gamma(Y))-\#V(\Duv \Gamma (Y))$ is
minimal for the fixed main invariants.
In particular, this case delivers all cases of minimal log del Pezzo
surfaces of index $\le 2$ with $\rk \Pic Z=1$.
See Theorems \ref{thm3.5.1}, \ref{thm4.1.2b}, \ref{thm4.1.2a}.

{\it No Du Val singularities.} This is the case when
$D=\Duv \Gamma(Y)=\emptyset$. Equivalently, all singularities of $Z$
have index $2$, if they exist. Then
$\Gamma (Y)=\Log \Gamma \cup \Var \Gamma (Y)$ where
\begin{equation}
V(\Var \Gamma (Y))=W\left(\{E_v\ |\ v\in V(\Var \Gamma)\}\right).
\label{intvarnd}
\end{equation}
Here, all the multiple cases 7a,b, 8a--c, 9a--f, 10a--m, 20a--d
give the same graphs (because they have the same, equal to zero,
root invariant, see below), and one can always take cases 7a,
8a, 9a, 10a, 20a for the main invariants. This case is very
similar to and includes the classical case of non-singular del
Pezzo surfaces corresponding to the cases 1---10. See Theorem
\ref{thm4.1.4} about this (without Du Val singularities) case. Log
del Pezzo surfaces of this case are defined by their main
invariants $(r,a,\delta)$ up to deformation. The Du Val parts
$\Duv \Gamma$ of graphs $\Gamma$ of Table 3 can be considered (for
this case) as analogs of root systems (or Dynkin diagrams) which
one usually associates to non-singular del Pezzo surfaces. Its
true meaning is to give the type of the Weyl group $W$ that
describes the varying part $\Var (\Gamma (Y))$ from $\Var \Gamma$
by \eqref{intvarnd}. In cases 7 --- 10, 20, one can take graphs
$\Gamma$ of cases 7a
--- 10a, 20a.

{\it The Root invariant.} It is possible that two different subgraphs
$D\subset \Duv \Gamma$, $D\subset \Duv \Gamma^\prime$ of graphs of Table 3
(with the same main invariants $(r,a,\delta)$) give isomorphic graphs
$\Gamma(Y)$ and $\Gamma(Y^\prime)$
for the corresponding right resolutions, and then they give similar
log del Pezzo surfaces $Z$ and $Z^\prime$ of index $\le 2$, according
to our classification. The {\it root invariant}
\begin{equation}
([D],\xi)
\label{introotinv}
\end{equation}
gives the necessary and sufficient condition for this to happen.

To define the root invariant \eqref{introotinv}, we first remark that
the main invariants $(r,a,\delta)$ define a unique hyperbolic (i. e. with
one  positive square) even 2-elementary lattice $S$ with these
invariants. Here $r=\rk S$, $S^\ast/S\cong (\bZ/2)^a$, and $\delta=0$, if
and only if $(x^\ast)^2\in \bZ$ for any $x\in S^\ast$. In \eqref{introotinv},
 $[D]$ is the root lattice generated by $D$, and
$\xi:[D]/2[D]\to S^\ast/S$ a homomorphism preserving finite forms
$(x^2)/2\mod 2$, $x\in [D]$, and $y^2\mod 2$, $y\in S^\ast$. The
construction of the root invariant \eqref{introotinv}
uses the double cover
$\pi:X\to Y$ by a K3 surface $X$ (see above) with the non-symplectic
involution $\theta$. Then $S=H^2(X,\bZ)^\theta$ is the
sublattice where $\theta^\ast$ is identical. The root invariant
\eqref{introotinv} is considered up to automorphisms of $S$
and the root lattice
$[D]$.  See Sections \ref{subsec2.4} and \ref{subsec3.2}
about this construction and a very easy criterion (the kernel $H$
of $\xi$ is almost equivalent to $\xi$) about isomorphism of
root invariants. The root invariant was first introduced and
used in \cite{11} and \cite{31}.

In practice, to calculate the root invariant of a log del Pezzo surface
of index $\le 2$, one should just go from the graphs $\Gamma$
of Table 3 to the equivalent graphs $\Gamma (P(\cM^{(2,4)}))$ or
$\Gamma(P(X)_+)$ of Tables 1 or 2 of exceptional curves for the
K3 pairs $(X,\theta)$ (see Sections \ref{subsec3.2}, \ref{subsec3.5}).

Thus, two Du Val subgraphs
$D\subset \Duv \Gamma$, $D\subset \Duv \Gamma^\prime$ of graphs of Table 3
give isomorphic full graphs $\Gamma(Y)$ and $\Gamma(Y^\prime)$ of
their log del Pezzo surfaces if and only if their root invariants
\eqref{introotinv} are isomorphic (see Theorem \ref{thm3.3.1}).
Moreover, we constantly use the root
invariant to prove existence of the corresponding K3 pairs $(X,\theta)$
and log del Pezzo surfaces $Z$. The main invariants $(r,a,\delta)$
and the root invariants \eqref{introotinv} are the main tools
in our classification. They are equivalent to the full graphs $\Gamma(Y)$
of exceptional curves on $Y$, but they are much more convenient to
work with. For non-singular del Pezzo surfaces and
log del Pezzo surfaces of index $\le 2$ without Du Val singularities
the root invariant is zero. This is why, in these cases,
we have such a simple classification as above.

See Section \ref{subsecN=20en} about enumeration of root invariants
(equivalently graphs of exceptional curves) in the most
non-trivial case $N=20$.

{\bf Acknowledgement.} We are grateful to V.A. Iskovskikh for
useful discussions. The first author was supported by NSF for part
of this work. The second author is grateful to Steklov
Mathematical Institute, Max-Planck-Insitut f\"ur Mathematik and
the University of Liverpool for hospitality.

\section{Log del Pezzo surfaces of index
  $\le2$ and \\ Smooth Divisor Theorem.}
\label{sec1}

In this work we always consider {\it algebraic varieties which are
defined over the field $\bC$ of complex numbers.} Further we
don't mention that.

\section{Basic definitions and notation}
\label{subsec1.1}
Let $Z$ be a normal algebraic surface, and $K_Z$ be a canonical Weil divisor
on it. The surface $Z$ is called
$bQ$-Gorenstein
\index{singularity!$\bQ$-Gorenstein}
if a
certain positive multiple of $K_Z$ is Cartier, and
$\bQ$-factorial
\index{singularity!$\bQ$-factorial} if this is true for any Weil
divisor $D$. These properties are local: one has to require all
singularities to be $\bQ$-Gorenstein, respectively
$\bQ$-factorial.

Let us denote by $Z^1(Z)$ and $\Div(Z)$ the groups of Weil and Cartier
divisors on $Z$. Assume that $Z$ is $\bQ$-factorial.  Then the groups
$Z^1(Z)\otimes\bQ$ and $\Div(Z)\otimes\bQ$ of $\bQ$-Cartier divisors and
$\bQ$-Weil divisors coincide. The intersection form defines natural pairings
$$
\Div(Z)\otimes\bQ \times \Div(Z)\otimes\bQ \to \bQ,
$$
$$
\Div(Z)\otimes\bR \times \Div(Z)\otimes\bR \to \bR .
$$
Quotient groups modulo kernels of these pairings are denoted
$N_{\bQ}(Z)$ and $N_{\bR}(Z)$ respectively; if the surface $Z$ is
projective, they are finite-dimensio\-nal linear spaces. 
The Kleiman--Mori cone \index{Kleiman--Mori cone} 
is a convex cone $\NE(Z)$ in
$N_{\bR}(Z)$, the closure of the cone generated by the classes of
effective curves. \label{Moricone}

Let $D$ be a $\bQ$-Cartier divisor on $Z$. We will say that $D$ is
ample if some positive multiple is an ample Cartier divisor in the
usual sense. By {\it Kleiman's criterion} \index{Kleiman's
criterion} \cite{Kleiman66}, for this to hold it is necessary and
sufficient that $D$ defines a strictly positive linear function on
$\NE(Z)-\{0\}$.

One says that the surface $Z$ has only log terminal
singularities \index{singularity!log terminal} if it is
$\bQ$-Gorenstein and for one (and then any) resolution of
singularities $\pi:Y\to Z$, in a natural formula
$K_Y=\pi^*K_Z+\sum\alpha_iF_i$, where $F_i$ are irreducible
divisors and $\alpha_i\in \bQ$, one has $\alpha_i>-1$. The least
common multiple of denominators of $\alpha_i$ is called the
index \index{singularity!index of}
 of $Z$.

It is known that two-dimensional log terminal singularities in
characteristic zero are exactly the quotient singularities \cite{25}.
A self-contained and characteristic-free classification in terms 
of dual graphs of resolutions is given in \cite{Alexeev_LogCanSings}.
Log terminal singularities are rational and $\bQ$-factorial.
We can now formulate the following:
\begin{definition} 
  A normal complete surface $Z$ is called a
log del Pezzo surface \index{log del Pezzo surface}
if it
  has only log terminal singularities and the anticanonical divisor $-K_Z$
  is ample. It has index $\le k$ if all of its singularities are of index $\le
  k$.
\end{definition}

We will use the following notation. If $D$ is a $\bQ$-Weil divisor,
$D=\sum c_iC_i$, $c_i\in\bQ$, then $\ulcorner D \urcorner$ will denote
the round-up $\sum\ulcorner c_i \urcorner D_i$, and
$\{D\}=\sum\{c_i\}C_i$ the fractional part. A divisor $D$ is
nef
\index{divisor!nef}
if for any curve $C$ one has $D\cdot C\ge0$; $D$ is
big and nef
\index{divisor!big and nef}
if in addition $D^2>0$.

Below we will frequently use the following generalization of Kodaira's
vanishing theorem. The two-dimensional case is due to Miyaoka
\cite{Miyaoka_Vanishing} \index{Miyaoka} 
and does not require the normal-crossing condition. The higher-dimensional
case is due to Kawamata \cite{24} \index{Kawamata} and Viehweg \cite{33}.
\index{Viehweg} 

\begin{theorem}[Generalized Kodaira's Vanishing theorem]
  \label{thm1.1.2}
  Let $Y$ be a \linebreak smooth surface and let $D$ be a $\bQ$-divisor on $Y$
  such that
  \begin{enumerate}
  \item $\supp\{D\}$ is a divisor with normal crossings;
  \item $D$ is big and nef.
  \end{enumerate}
  Then $H^i(K_Y+\ulcorner D \urcorner)=0$ for $i>0$.
\end{theorem}

\section{Log terminal singularities of index 2}
\label{subsec1.2}
\index{singularity!log terminal of index 2}
Let $(Z,p)$ be a two-dimensional log terminal singularity of index
$\le2$, and $\pi:\wZ\to Z$ be its minimal resolution. We have
$K_{\wZ}=\pi^*K_Z + \sum\alpha_iF_i$, where $-1<\alpha_i\le0$ and
$F_i^2\le-2$. Therefore, for each $i$ one has $\alpha_i=-1/2$ or
$0$. One can rewrite the set of equations $K_{\wZ}\cdot F_i=-F_i^2-2$
in a matrix form:
$$
M\cdot (\alpha_1,\dots,\alpha_n)^t = (-F_1^2-2,\dots,-F_n^2-2)^t,
$$
where $M=(F_i\cdot F_j)$ is the intersection matrix. By a basic theorem of
Mumford \cite{Mumford_SurfaceSingularities}, \index{Mumford} 
$M$ is negative definite and, in 
particular, nondegenerate. All the entries of the inverse matrix $M^{-1}$ are
strictly negative \cite{Artin_NumericalCriteria}.

Let us give some easy consequences of this formula.
\begin{enumerate}
\item If for some $i_0$, $\alpha_{i_0}=0$ then all $\alpha_i=0$, and
  the singularity $(Z,p)$ is Du Val, of type $A_n$, $D_n$, $E_6$,
  $E_7$ or $E_8$.
\item If all $\alpha_i=-1/2$ then we get the following list of
  singularities:

\medskip

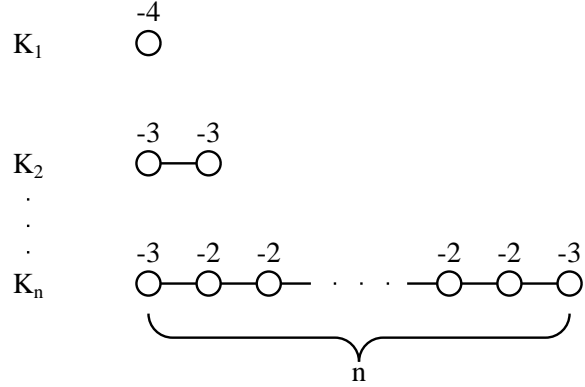
\centerline{\includegraphics[width=8cm]{pics/c-p55.eps}}

\medskip

\noindent
In these graphs every curve $F_i$ corresponds to a vertex with
  weight $F_i^2$, two vertices are connected by an edge if $F_i\cdot
  F_j=1$ and are not connected if $F_i\cdot F_j=0$.
\end{enumerate}

\section{Basic facts about log del Pezzo surfaces}
\label{subsec1.3}
\begin{lemma}\label{lemma1.3.1}
  All log del Pezzo surfaces $Z$ are rational.
\end{lemma}
\begin{proof}
  Let $\pi:\wZ\to Z$ be the minimal resolution of singularities,
  $K_{\wZ}= \pi^*K_Z+\sum\alpha_iF_i$, $-1<\alpha_i\le0$. Then
  $\ul -\pi^*K_Z\ur = -K_{\wZ}$, and so
  \begin{displaymath}
    h^1(O_{\wZ}) =
    h^1(K_{\wZ} + \ul -\pi^*K_Z\ur) =0
  \end{displaymath}
by Theorem \ref{thm1.1.2}.

Also, $h^0(nK_{\wZ})=0$ for any positive integer $n$ since $-\pi_*
K_{\wZ} = -K_Z$ is an effective nonzero $\bQ$-Weil divisor. Therefore,
by Castelnuovo criterion the surface $\wZ$, and hence also $Z$, are
rational.
\end{proof}

\index{exceptional curves of!DPN surface} 

\begin{lemma}\label{lemma1.3.2}
  In the above notation, if $\wZ\ne\bP^2$ or $\bF_n$ then the 
Kleiman--Mori cone
  of the surface $\wZ$ is generated by the curves $F_i$ and
  exceptional curves of the 1st kind. The number of these curves is
  finite. There are no other irreducible curves with negative
  self-intersection number (i. e. exceptional curves) on $\wZ$.
  \index{exceptional curve} 

  Moreover, in this statement the minimal resolution $\wZ$ can be
  replaced by any resolution of singularities $\pi:Z'\to Z$ such that
  $\alpha_i \le 0$, where $K_{Z^\prime} = \pi^*K_Z + \sum \alpha_i F_i$ (for
  example, by the right resolution of $Z$, see Section \ref{subsec1.5} below).
\end{lemma}
\begin{proof}
  Let us show that on the surface $\wZ$ (or $Z'$) there exists a
  $\bQ$-divisor $\Delta$ with $\Delta\ge0$, $[\Delta]=0$ and such that
  the divisor $-(K_{\wZ} +\Delta)$ is ample.

  Choose $E=\sum\beta_iF_i$ so that $Z\cdot F_i>0$. Since the matrix
  $(F_i\cdot F_j)$ is negative definite and $F_i\cdot F_j\ge 0$ for
  $i\ne j$, one has $\beta_i<0$. Let us show that for a small
  $0<\epsilon\ll1$ the divisor $T=-\pi^*K_Z+\epsilon E$ is ample.
  Since $(-K_Z)^2>0$,
  we may assume that $T^2>0$. Now, let us check that, if the positive number
  $\epsilon$ is sufficiently small, then $C\cdot T>0$
for any irreducible curve $C$ on $Z$.

  If $C^2\ge0$, this follows from the fact that the intersection form
  on $N_{\bR}(\wZ)$ is hyperbolic.
  If $C=F_i$ then $F_i\cdot T = \epsilon F_i\cdot E >0$.

  If $C^2<0$ and $C\ne F_i$ then
$$
C\cdot K_{\wZ} = C\cdot (\pi^*K_{\wZ} + \sum \alpha_iF_i) <0.
$$
On the other hand, $p_a(C)= \dfrac{C^2+C\cdot K_{\wZ}}{2} +1\ge0$.
So, $C^2<0$ and $C\cdot K_{\wZ}<0$ imply that $C^2=-1$ and
$p_a(C)=0$, i.e. $C$ is an exceptional curve of the 1st
kind. \index{exceptional curve!1st kind}

If $n$ is the index of $Z$ then $C\cdot (-\pi^*K_Z)\in (1/n)\bZ$.
On the other hand, $0<-C\cdot\pi^*K_Z = 1 +\sum\alpha_iF_i\cdot
C\le 1$. Hence, there are only finitely many possibilities for
$(-\pi^*K_Z)\cdot C$ and $\sum\alpha_iF_i\cdot C$, and for
$\epsilon$ small enough, $C\cdot T>0$.

By Kleiman's criterion, this implies that $T$ is ample. Since the degree of
the $(-1)$-curves with respect to $T$ is bounded, there are only finitely many
of them.

One has $-\pi^*K_Z+\epsilon E
= - \left(K_{\wZ}+ \sum (-\alpha_i-\epsilon\beta_i)F_i\right)$.
Therefore, $\Delta\ge0$, and for  $\epsilon\ll1$,
we have $[\Delta]=0$, since $-\alpha_i<1$.

Now, by Cone theorem \cite[Thm.4.5]{25}, $\NE(\wZ) = \sum R_j$, where
$R_j$ are ``good extremal rays''. The rays generated by the curves $F_i$ and
exceptional curves of the 1st kind are obviously extremal. On the other
hand, let $R_j$ be a ``good extremal ray'', generated by an irreducible curve
$C$. If $C\notin \{F_1,\dots,F_k\}$ then $C\cdot K_{\wZ}= C\cdot (\pi^*K_Z +
\sum \alpha_iF_i)<0$. Hence, by \cite{27} the curve $C$ is an exceptional
curve of the 1st kind, unless $Z\simeq\bP^2$ or $\bF_n$.
\end{proof}

\section{Smooth Divisor Theorem}
\label{subsec1.4}
\begin{theorem}\label{thm1.4.1}
  Let $Z$ be a log del Pezzo surface of index $\le2$. Then the linear
  system $|-2K_Z|$ is nonempty, has no fixed components and contains a
  nonsingular element $D\in |-2K_Z|$.
\end{theorem}

\index{Smooth divisor Theorem}

\begin{proof}
  Let $\pi:\wZ\to Z$ be the minimal resolution of singularities. It is
  sufficient to prove the statement for the linear system
  $|-\pi^*(2K_Z)|$ on $\wZ$. We have $2K_{\wZ}=\pi^*(2K_Z)-\sum a_iF_i$,
  and all $a_i=0$ or 1 ($a_i=-2\alpha_i$).

  1. \emph{Nonemptiness.}
$$
-\pi^*(2K_Z) = K_{\wZ}+ (-3K_{\wZ} -\sum a_iF_i) = K_{\wZ}+ \ul D \ur,
$$
where $D=\dfrac{3}{2}(-2K_{\wZ}-\sum a_iF_i)=-\pi^* (3K_Z)$
is big and nef. Hence, by Vanishing
Theorem \ref{thm1.1.2},  $H^i(-\pi^*(2K_Z)) = 0$ for $i>0$ and
  $h^0(-\pi^*(2K_Z)) = \chi(-\pi^*(2K_Z)) = 3K_Z^2+1>0$.

  2. \emph{Nonexistence of fixed components.} Let $E$ be the fixed part, so
  that $|-\pi^* (2K_Z)-E|$ is a movable linear system. Then

\begin{eqnarray*}
&&  h^0(-\pi^*(2K_Z)) =   h^0(-\pi^*(2K_Z)-E),  \\
&&  -\pi^*(2K_Z)-E = K_{\wZ} + (-3K_{\wZ} - \sum a_iF_i - E)
  = K_{\wZ} + \ul D \ur, \\
&&  D = \frac{3}{2} (-2K_{\wZ} -\sum a_iF_i) -E =
  (-\pi^*(2K_Z)-E) + (-\pi^*K_Z).
\end{eqnarray*}
The first of these divisors is movable and the second is big and nef, so the
sum is big and nef.  By Vanishing Theorem \ref{thm1.1.2}, we have
$$
h^i(-\pi^*(2K_Z)-E) = 0, \quad i>0,
$$
$$
\chi(-\pi^*(2K_Z))=  \chi(-\pi^*(2K_Z)-E),
$$
\begin{equation}
2\chi(-\pi^*(2K_Z))-  2\chi(-\pi^*(2K_Z)-E)  =
  E\cdot (-2\pi^*(2K_Z) - K_{\wZ} - E) .
\label{1.4.1}
\end{equation}
Let us show that this expression \eqref{1.4.1} is not equal to zero. Suppose
$$-\pi^*K_Z\cdot (K_{\wZ}+E) = -\pi^*K_Z\cdot E - K_Z^2 <0.$$
Then the divisor
$K_{\wZ}+E$ cannot be effective. Therefore,
$$
\chi(-E) = h^0(-E)-h^1(-E) + h^0(K_{\wZ}+E) \le 0.
$$
Hence, $E\cdot (K_{\wZ}+E)=2\chi(-E)-2<0$, and the expression \eqref{1.4.1}
is strictly positive.
So, we can assume that $-\pi^*K_Z\cdot E \ge K_Z^2$. Let us write
$$
E=\beta(-\pi^*K_Z)+F, \quad F\in (\pi^*K_Z)^{\perp}
$$
in $N_{\bQ}(\wZ)$. One has $\beta\le2$ since
$-\pi^*K_Z\cdot (-\pi^*(2K_Z)-E)\ge0$. Then
$$
E\cdot (-2\pi^*(2K_Z) - K_{\wZ} - E) = (5-\beta)\beta K_Z^2 -
F\cdot (\sum \alpha_iF_i + F).
$$
The first term in this sum is $\ge3$ since $\beta K_Z^2= -\pi^*K_Z
\cdot E\ge K_Z^2$ and $K_Z^2 = \chi(\pi^*(2K_Z)) -1$ is a positive
integer. The second term achieves the minimum for $F=-\dfrac{1}{2}
\sum \alpha_iF_i$ and equals $-\dfrac{m}{4}$, where $m$ is the number
of non-Du Val singularities. Therefore, all that remains to be shown
is that the surface $Z$ has fewer than 12 non-Du Val singularities.

By Lemma~\ref{lemma1.3.1}, the surface $\wZ$ is rational. By
Noether's formula, $(K_{\wZ})^2+\rk \Pic\wZ =10$. By
Lemma~\ref{prop1.4.2} below, $(K_{\wZ})^2\ge 0$. Hence $Z$ has no
more than $\rk \Pic\wZ - 1\le 9$ singular points.

\begin{lemma}\label{prop1.4.2}
  Let $Z$ be a log del Pezzo surface of index $\le 2$ and $\pi:\wZ\to Z$
  be its minimal resolution of singularities. Then $K_{\wZ}^2\ge 0$.
\end{lemma}
\begin{proof}
  One has $K_{\wZ} = \pi^*K_Z +\sum \alpha_iF_i$.
Denote $\oK = \pi^*K_Z -\sum \alpha_iF_i$. Let us show that $-\oK$
is nef. By Lemma~\ref{lemma1.3.2}, one has to show that $-\oK\cdot
F_i\ge0$ and $-\oK\cdot C\ge0$ if $C$ is an exceptional curve of the
1st kind. We have $-\oK\cdot F_i = K_{\wZ}\cdot F_i = -F_i^2 -2\ge0$
since the resolution $\pi$ is minimal. Next,
$$
-\pi^*K_Z\cdot C = -K_{\wZ}\cdot C + \sum \alpha_iF_i\cdot C =
1+ \sum \alpha_iF_i\cdot C >0.
$$
Since this number is a half-integer,
$$ -\oK\cdot C = 1+ 2\sum \alpha_iF_i\cdot C \ge0. $$
So, $-\oK$ is nef and $K_{\wZ}^2= \oK^2\ge0$. Finally, if $\wZ=\bP^2$
or $\bF_n$ then $K_{\wZ}^2=9$ or $8$ respectively.
\end{proof}

3. \emph{Existence of a smooth element.} Assume that all divisors in
the linear system
$|-\pi^*(2K_Z)|$ are singular. Then there exists a base point $P$, and
for a general element $D\in |-\pi^*(2K_Z)|$ the multiplicity of $D$ at
$P$ is $k\ge2$. This point does not lie on $F_i$ since
$-\pi^*(2K_Z)\cdot F_i=0$. Let $\epsilon:Y\to \wZ$ be the blowup at $P$,
$f=\pi\epsilon:Y\to Z$, and let $L$ be the exceptional divisor of
$\epsilon$. We have: $h^0(-f^*(2K_Z))= h^0(-f^*(2K_Z)-L)$, the linear
system $|-f^*(2K_Z)-kL|$ is movable, and
\begin{eqnarray*}
&&  2K_Y = f^*(2K_Z) - \sum a_iF_i+ 2L, \\
&&  -f^*(2K_Z) = K_Y + (-3K_Y -\sum a_iF_i + 2L) = K_Y + \ul D\ur,\\
&& D = \frac{3}{2}(-2K_Y - \sum a_iF_i + L) =
\frac{3}{2}(-f^*(2K_Z)-L) =\\
&&\frac{3}{2}[ (-f^*(2K_Z) - kL) + (k-1)L].
\end{eqnarray*}
The divisor $D$ is nef since for any irreducible curve $C\ne L$,
$C\cdot D\ge0$, and also $D\cdot L=3/2$. It is big since $(-f^*(2K_Z)-L)^2
= 4K_Z^2-1>0$. Now,
$$
-f^*(2K_Z) -L = K_Y+ (-3K_Y -\sum a_iF_i +L) = K_Y + \ul D \ur, 
$$
$$
D = \frac{3}{2}(-2K_Y - \sum a_iF_i + \frac{2}{3}L) =
  \frac{3}{2}\big((-f^*(2K_Z)-kL) + (k-\frac{4}{3})L\big).
$$
The latter divisor $D$ is nef since for any irreducible curve
$C\ne L$, $C\cdot D\ge0$, and also $D\cdot L=2$; $D$ is big since
$$
(-f^*(2K_Z) - \frac{4}{3}L)^2 = 4K_Z^2 - \frac{16}{9} >0.$$
Now,
again by Vanishing Theorem \ref{thm1.1.2},
$$h^i(-f^*(2K_Z)) =h^i(-f^*(2K_Z) -L) =0\quad \text{for }i>0,$$
and one must have $\chi(-f^*(2K_Z)) = \chi(-f^*(2K_Z) -L)$.
But
\begin{align*}
\chi(-f^*(2K_Z)) &- \chi(-f^*(2K_Z) -L) =
\\
=\frac{1}{2}L\cdot
(-2f^*(2K_Z) - K_Y - L) &=
1+L\cdot (-f^*(2K_Z)) >0.
\end{align*}
The contradiction thus obtained completes the proof of the theorem.
\end{proof}

\begin{remark}\label{rem1.4.3}
  In the same way, parts 1 and 3 can be proved for a log del Pezzo
  surface of arbitrary index $n$ and the linear system $-\pi^* (nK_Z)$.
  Part 2 is easy to prove under the assumption that $\pi(E)$ passes
  only through (some of) the Du Val singularities.
\end{remark}

\subsection{Reduction to DPN surfaces of elliptic type}
\label{subsec1.5}

Let $Z$ be a log del Pezzo surface of index $\le2$. Consider a
resolution of singularities $f:Y\to Z$ for which every Du Val
singularity is resolved by inserting the usual tree of $(-2)$-curves,
and the singularity $K_n$ by the following chain:


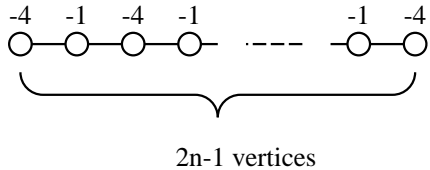
\begin{equation}
\centerline{\includegraphics[width=6cm]{pics/c-p60.eps}}
\label{1.5.1}
\end{equation}

\noindent
The latter resolution is obtained by blowing up all intersection points
of exceptional curves on the minimal resolution of $K_n$ points, see
their diagrams in Section~\ref{subsec1.2}. In contrast to the minimal
resolution, we will call this the {\it right} resolution of
singularities.  Consider a smooth element $C_g\in |-2K_Z|$. It does
not pass through singularities of the surface $Z$. If we identify the
curve $C_g$ with its image under the morphism $f$, then it is easy to
see from the formulae of Section~\ref{subsec1.2} that $-f^*2K_Z$
is linearly equivalent to $C_g$, and $-2K_Y$ with the disjoint union of $C_g$
and curves in the above diagrams which have self-intersection $-4$.
Moreover, it is easy to compute the genus of the curve $C_g$, and it
equals $g=K_Z^2+1\ge 2$. This shows that the surface $Y$ is a right
DPN surface of elliptic type in the sense of the next Chapter (see
Sections~\ref{subsec2.1} and \ref{subsec2.7}).

Vice versa, the results of Chapters~\ref{sec2} and \ref{sec3} will
imply (see Chapter~\ref{sec4}) that a right DPN surface $Y$ of elliptic type
admits a unique contraction of exceptional curves $f:Y\to Z$ to a log
del Pezzo surface of index $\le2$.

In this way, the classification of log del Pezzo surfaces of index $\le2$
is reduced to classification of right DPN surfaces of elliptic type.

\section{General Theory of DPN surfaces and K3 surfaces with non-symplectic
involution}
\label{sec2}
\subsection{General remarks}
\label{subsec2.1}
As it was shown in Chapter~\ref{sec1}, a description of log del
Pezzo surfaces of index $\le2$ is reduced to a description of
rational surfaces $Y$ containing a nonsingular curve $C\in |-~2K_Y|$
and a certain configuration of exceptional curves. Such surfaces $Y$
and exceptional curves on them were studied in the papers \cite{N1,
  10, 11, 31} of the second author. They are one of possible
generalizations of del Pezzo surfaces.

Many other generalizations of del Pezzo surfaces were proposed, see
e.g. \cite{18,21,22,26}, and most authors call their surfaces
``generalized del Pezzo surfaces''. Therefore, we decided following
\cite{31} to call our generalization DPN surfaces. One can consider DPN
surfaces as appropriate non-singular models of log del Pezzo
surfaces of index $\le 2$ and some their natural generalizations.

\begin{definition}
  A nonsingular projective algebraic surface $Y$ is called a DPN
  surface if its irregularity $q(Y)=0$, $K_Y\ne 0$ and there exists an
  effective divisor $C\in|-2K_Y|$ with only simple rational,
  i.e. $A,D,E$-singularities. Such a pair $(Y,C)$ is called a
  {\it DPN pair.} A DPN surface $Y$ is called \emph{right} if there exists a
  nonsingular divisor $C\in |-2K_Y|$; in this case the pair $(Y,C)$ is
  called {\it right DPN pair} or nonsingular DPN pair.
\end{definition}

The classification of algebraic surfaces implies that if
$C=\emptyset$ then a DPN surface $Y$ is an Enriques surface
($\varkappa=p=q=0$). If $C\ne\emptyset$ then $Y$ is a rational
surface ($\varkappa=-1$, $p=q=0$), e.g. see \cite{Shaf65}.

Using the well-known properties of blowups, the following results are
easy to prove. Let $(Y,C)$ be a DPN pair, $E\subset Y$ be an
exceptional curve of the 1st kind on $Y$ and $\sigma:Y\to Y'$ the
contraction of $E$. Then $(Y', \sigma(C))$ is also a DPN pair. In this
way, by contracting exceptional curves of the 1st kind, one can always
arrive at a DPN pair $(Y',C')$ where $Y'$ is a relatively minimal
(i.e. without $(-1)$-curves)
rational surface. In this case, the only possibilities for $Y'$ are
$\bP^2$, $\bF_0$, $\bF_2$, $\bF_3$ and $\bF_4$, since only for them
$|-2K_{Y^\prime}|$ contains a reduced divisor.
If $Y'=\bP^2$ then $C'$ is a curve of
degree 6; if $Y'=\bF_0=\bP^1\times \bP^1$ then $C'$ is a curve of
bidegree $(4,4)$; if $Y'=\bF_2$ then $C'\in |8f+4s_2|$; if $Y'=\bF_3$
then $D'=C_1+s_3$, where $C_1\in |10f+3s_3|$; if $Y'=\bF_4$ then
$C'=D_1+ s_4$, where $D_1\in |12f + 3s_4|$.
Here, the linear system $|f|$ is a pencil of rational curves on
surface $\bF_n$ with a section $s_n$, $s_n^2=-n$. Vice versa, if
$(Y',C')$ is a DPN pair, $P$ is a singular point of $C'$ and
$\sigma:Y\to Y'$ is a blowup of $P$ with an exceptional $(-1)$-curve
$E$ then $(Y,C)$ is a DPN pair, where
$$
C= \left\{
    \begin{array}{ll}
      \sigma\inv_*(C^\prime) &
      \text{if } P \text{ has multiplicity } 2 \text{ on } C^\prime \\
      \sigma\inv_*(C^\prime)+E &
      \text{if } P \text{ has multiplicity } 3 \text{ on } C^\prime
    \end{array}
    \right.\ .
    $$
Here $\sigma\inv_*(C^\prime)$ denotes the proper preimage (or the strict transform) of $C^\prime$, i. e. $\sigma\inv_*(C^\prime)$ is the closure of
$\sigma^{-1}(C^\prime -\{P\})$ in $Y$ where
$\sigma^{-1}(C^\prime -\{P\})$ is the set-theoretic preimage of
$C^\prime-\{P\}$.

    In this way, by blowups, from an arbitrary DPN pair $(Y',C')$
    one can always pass to a right DPN pair $(Y,C)$, i.e. with a
    nonsingular $C$. A description of arbitrary DPN pairs and surfaces
    is thus reduced to a description of right (or nonsingular) DPN pairs
$(Y,C)$ and
    right DPN surfaces $Y$, and exceptional curves on $Y$, where a curve
    $E\subset Y$ is called exceptional if $E$ is irreducible and
    $E^2<0$.

    We shall need a small, elementary, and well-known fact about
    ramified double covers.  Let $\pi:X\to Y$ be a finite morphism of
    degree 2 between smooth algebraic varieties. Then $\pi$ is Galois
    with group $\bZ/2$. Therefore, the $\cO_Y$-algebra $\pi_*\cO_X$
    splits into $(\pm 1)$-eigenspaces as $\cO_Y \oplus L$. Since $\pi$
    is flat, $L$ is flat and hence invertible. The algebra structure
    is given by a homomorphism $L^2\to \cO_Y$, i.e. by a section $s\in
    H^0(Y,L^{-2})$. Locally, $X$ is isomorphic to $y^2=s(x)$. Since
    $X$ is smooth, the ramification divisor $C=(s)$ must be smooth.

    Vice versa, let $L^{-1}$ be a sheaf dividing by two the sheaf
    $\cO_Y(C)$ for an effective divisor $C$ in $\Pic Y$ and let $s$ be
    a section of $\cO_Y(C)$ with $(s)=C$.  Then $s$ defines an algebra
    structure on $\cA=\cO_Y \oplus L$, and $\pi:X:=\Spec \cA \to Y$ is
    a double cover ramified in $C$. The representation of $\cA$ as a
    quotient of $\oplus_{d\ge0} L^d$ gives an embedding of $X$
    into a total space of the line bundle $L\inv$ and a section of
    $\pi^*L\inv$ ramified along $\pi\inv(C)$ with multiplicity one.
    Hence $\pi\inv(C) \sim \pi^*L\inv$.

    Let $(Y,C)$ be a right DPN pair. Since $C\in |-2K_Y|$, there
    exists a double cover $\pi:X\to Y$ defined by $L\inv=-K_Y$, branched
    along $C$. By the above, we have $\pi^\ast (-K_Y)\sim \pi^{-1}(C)$.

    Let $\omega_Y$ be a rational 2-dimensional differential form on
    $Y$ with the divisor $(\omega_Y)$ whose components
do not contain components of
    $C$. Then $(\omega_Y)\sim K_Y$, and the divisor
    $(\pi^\ast\omega_Y)= \pi^\ast (\omega_Y)+\pi^{-1}(C)\sim
    \pi^\ast (\omega_Y)+\pi^\ast(-K_Y)\sim 0$.
    Thus, $K_X=0$. Then $X$ is either a K3 surface (i. e. $q(X)=0$) or an
    Abelian surface (i. e. $q(X)=2$), e.g. see \cite{Shaf65}.
    Let $X$ be an Abelian surface. Then $C$ is
    not empty (otherwise, $Y$ is an Enriques surface and then $X$ is
    a K3 surface, \cite{Shaf65}), and $Y$ is rational.
    It follows that there exists a non-zero regular 1-dimensional
    differential form $\omega_1$ on $X$ such that $\theta^\ast
    (\omega_1)=\omega_1$ for the involution $\theta$ of the
    double cover $\pi$. Then $\omega_1=\pi^\ast
    \widetilde{\omega}_1$ where $\widetilde{\omega}_1$ is a
    regular 1-dimensional differential form on $Y$. This
    contradicts $q(Y)=0$. It proves that $X$ is a K3 surface.

    Let $\omega_X$ be a non-zero regular 2-dimensional
    differential form on $X$. If $\theta^\ast(\omega_X)=\omega_X$,
    then $\omega_X=\pi^\ast(\omega_Y)$ where $\omega_Y$ is a
    regular 2-dimensional differential form on $Y$. This
    contradicts the fact that $Y$ is an Enriques or rational
    surface (e. g. see \cite{Shaf65}). Thus, $\theta^\ast
    (\omega_X)=-\omega_X$, and then $\theta$ is a {\it non-symplectic
    involution} of the K3 surface $X$.

    Vice versa, assume that $(X,\theta)$ is a K3 surface
    with a non-symplectic involution.  Then the set $X^{\theta}$ of
    fixed points of the involution is a nonsingular curve (otherwise, $\theta$
    is symplectic, i. e. $\theta^\ast(\omega_X)=\omega_X$ for any regular
    2-dimensional differential form on $X$). It follows (reversing
    arguments above) that the pair
    $\left(Y=X/\{1,\theta\},\ C = \pi(X^{\theta})\right)$ is a right
    DPN pair where $\pi:X\to Y$ is the quotient morphism.

    Thus, a description of right DPN pairs $(Y,C)$ is reduced to a
    description of K3 surfaces with a non-symplectic involution
    $(X,\theta)$.

\subsection{Reminder of basic facts about K3 surfaces}
\label{subsec2.1a} Here we remind basic results about K3 surfaces
that we use. We follow \cite{Shaf65}, \cite{13}, \cite{5} and
also \cite{8,10,12}. Of course, all these results are well-known.

Let $X$ be an algebraic K3 surface. We recall that this
means that $X$ is a projective non-singular algebraic surface, the
canonical class $K_X=0$ (i. e. there exists a non-zero regular
2-dimensional differential form $\omega_X$ on $X$ with zero
divisor), and $q(X)=\dim \Gamma(X,\Omega^1)=0$ (i. e. $X$ has no non-zero
regular 1-dimensional differential forms). From definition,
$\omega_X$ is unique up to multiplication by $\lambda \in \bC$,
$\lambda \not=0$.

Let $F\subset X$ be an irreducible algebraic curve. By genus formula,
\begin{equation}
p_a(F)=\frac{F^2+(F\cdot K_X)}{2}+1=\frac{F^2}{2}+1\ge 0.
\label{genusK3}
\end{equation}
It follows that $F^2\equiv 0\mod 2$, $F^2\ge -2$, and $F$ is
non-singular rational if $F^2=-2$. In particular, any exceptional
curve $F$ on $X$ (i. e. $F$ is irreducible and $F^2<0$) is
non-singular rational with $F^2=-2$.

By Riemann-Roch Theorem, for any divisor $D\subset X$ we have
$$
l(D)+l(K_X-D)=h^1(D)+\frac{D\cdot(D-K_X)}{2}+\chi(\cO_X),
$$
which gives for a K3 surface $X$ that
\begin{equation}
l(D)+l(-D)=h^1(D)+\frac{D^2}{2}+2\ge \frac{D^2}{2}+2.
\label{R-RK3}
\end{equation}
It follows that one of $\pm D$ is effective if $D^2\ge -2$.

All algebraic curves on $X$ up to linear equivalence generate the
Picard lattice $S_X$ of $X$. For K3 surfaces linear equivalence is
equivalent to numerical, and $S_X\subset H^2(X,\bZ)$ where
$H^2(X,\bZ)$ is an even unimodular lattice of the signature
$(3,19)$. Here, ``even'' means that $x^2$ is even for any $x\in
H^2(X,\bZ)$. Unimodular means that for a basis $\{e_i\}$ of
$H^2(X,\bZ)$ the determinant $\det(e_i\cdot e_j)=\pm 1$.
Such even unimodular lattice is unique up to an
isomorphism, see e.g. \cite{Serre}. By Hodge Index Theorem, the
Picard lattice $S_X$ is hyperbolic, i. e. it has signature
$(1,\rho-1)$ where $\rho=\rk S_X$.  Let
\begin{equation}
V(S_X) = \{ x\in S_X\otimes\bR \suchthat x^2>0 \}.
\label{coneSXK3}
\end{equation}
Since $S_X$ is hyperbolic, $V(S_X)$ is an open cone which has two
convex halves. One of these halves $V^+(X)$ is distinguished by
the fact that it contains the ray $\bR^+h$ of a polarization $h$
(i. e. a hyperplane section) of $X$ where $\bR^+$ denotes the set
of all non-negative real numbers.

Let
\begin{equation}
\NEF(X)=\{x\in S_X\otimes \bR\ |\ x\cdot C\ge 0\  \ \forall\
{effective\ curve}\ C\subset X\}
\label{NEFK31}
\end{equation}
be the $nef$ cone of $X$. Since $S_X$ is hyperbolic, for any
irreducible curve $C$ with $C^2\ge 0$ we have that $C\in
\overline{V^+(X)}$, and $C\cdot V^+(X)>0$. It follows that
\begin{equation}
\NEF(X)=\{x\in \overline{V^+(X)} \suchthat x\cdot P(X)\ge 0\}
\label{NEFK32}
\end{equation}
where $P(X)\subset S_X$ denote the set of all divisor classes of
irreducible non-singular rational (i. e. all exceptional) curves
on $X$.

Let $h\in \NEF(X)$ be a hyperplane section. By Riemann-Roch Theorem
above, $f\in S_X$ with $f^2=-2$ is effective if and only if
$h\cdot f>0$. It follows that $\NEF(X)$ is a fundamental chamber
(in $V^+(X)$) for the group $W^{(2)}(S_X)$ generated by
reflections in all elements $f\in S_X$ with $f^2=-2$. Each such
$f$ gives a reflection $s_f\in O(S_X)$ where
\begin{equation}
s_f(x)=x+(x\cdot f)f,
\label{reflectionK3}
\end{equation}
in particular, $s_f(f)=-f$ and $s_f$ is identical on $f^\perp$.

Since all $F\in P(X)$ have $F^2=-2$, the nef cone $\NEF(X)$ is
locally finite in $V^+(X)$, all its faces of codimension one are
orthogonal to elements of $P(X)$. {\it This gives a one-to-one
correspondence between the faces of codimension one of $\NEF(X)$ and
the elements of $P(X)$.} Indeed, let $\gamma$ be a codimension one
face of $\NEF(X)$. Assume $F\in P(X)$ is orthogonal to $\gamma$, i.
e. $\gamma$ belongs to the edge of the half-space $F\cdot x\ge 0$,
$x\in S_X\otimes \bR$, containing $\NEF(X)$. Such $F\in S_X$ with
$F^2=-2$ is obviously unique because any element $f\in S_X$ which
is orthogonal to $\gamma$ is evidently $\lambda F$, $\lambda \in
\bR^+$. We have $(\lambda F)^2=\lambda^2F^2=\lambda^2(-2)$, and
$F$ is distinguished by the condition $F^2=-2$. In such a way, all
faces of codimension one of $\NEF(X)$ give a subset
$P(X)^\prime\subset P(X)$ of elements of $P(X)$ which are
orthogonal to codimension one faces of $\NEF(X)$. Let us show that
$P(X)^\prime=P(X)$. Obviously it will be enough to show that for
any $E\in P(X)$, the orthogonal projection of $\NEF(X)$ into
hyperplane $E^\perp$ belongs to $\NEF(X)$. The projection is given
by the formula $H\mapsto \widetilde{H}=H+(H\cdot E)E/2$ for $H\in
\NEF(X)$. Let us show that $\widetilde{H}\in \NEF(X)$. Let $C$ be an
irreducible curve on $X$. If $C\not=E$, then $C\cdot
\widetilde{H}=C\cdot H+(H\cdot E)(C\cdot E)/2\ge 0$ because $H$ is
nef and $C$ is different from $E$. If $C=E$, then $C\cdot
\widetilde{H}=E\cdot \widetilde{H}=E\cdot H+(H\cdot E)(E^2)/2=0\ge
0$. Thus, $\widetilde{H}\in \NEF(X)$.

Therefore, we obtain a group-theoretic description of the nef cone
of $X$ and all exceptional curves of $X$: The $\NEF(X)$ is the
fundamental chamber for the reflection group $W^{(2)}(S_X)$ acting
on $V^+(X)$, this chamber is distinguished by the condition that
it contains a hyperplane section of $X$. The set $P(X)$ of all
exceptional curves on $X$ consists of all elements $f\in S_X$
which have $f^2=-2$ and which are orthogonal to codimension one
faces of $\NEF(X)$ and directed outwards (i. e. $f\cdot \NEF(X)\ge
0$).

It is more convenient to work with the corresponding hyperbolic
space
\begin{equation}
\cL(X) = V^+(X)/\bR^+.
\label{LobachevskyK3}
\end{equation}
Elements of this space are rays $\bR^+x$, where $x\in
S_X\otimes\bR$, $x^2>0$ and $x\cdot h>0$. Each element $\beta\in
S_X\otimes \bR$ with square $\beta^2<0$ defines a half-space
\begin{equation}
\cH_{\beta}^+ = \{ \bR^+x\in \cL(X) \suchthat \beta\cdot x \ge 0 \},
\label{halfspaceK3}
\end{equation}
so that $\beta$ is perpendicular to the bounding hyperplane
\begin{equation}
\cH_{\beta} = \{\bR^+x\in \cL(X) \suchthat \beta\cdot x = 0 \},
\label{hyperplaneK3}
\end{equation}
and faces outward. The set
\begin{equation}
\cM(X) = \bigcap_{f\in S_X, f^2=-2 \atop
  f\ \rm{is\  effective}} \cH_{f}^+ =
\bigcap_{f\in P(X)} \cH^+_{f}
\label{M(X)K3}
\end{equation}
is a locally finite convex polytope in $\cL(X)$. The set
$P(\cM(X))$ of vectors with square $-2$, perpendicular to the
facets of $\cM(X)$ and directed outward, is exactly the set $P(X)$
of divisor classes of exceptional curves on $X$. Moreover,
$\cM(X)$ admits a description in terms of groups. Let $O'(S_X)$ be
the subgroup of index two of the full automorphism group $O(S_X)$
of the lattice $S_X$ which consists of automorphisms which
preserve the half-cone $V^+(X)$. Let $W^{(2)}(S_X) \subset
O^\prime(S_X)$ be the subgroup of $O^\prime(S_X)$ generated by
reflections $s_{f}$ with respect to all elements $f\in S_X$ with
square $(-2)$. The action of the group $W^{(2)}(S_X)$, as well as
$O'(S_X)$, on $\cL(X)$ is discrete. $W^{(2)}(S_X)$ is the group
generated by reflections in all hyperplanes $\cH_{f}$, $f\in S_X$
and $f^2=-2$. The set $\cM(X)$ is a fundamental chamber for this
group, i.e. $W^{(2)}(S_X)(\cM(X)) $ defines a decomposition of
$\cL(X)$ into polytopes which are congruent to $\cM(X)$, and
$W^{(2)}(S_X)$ acts transitively and without fixed elements on
this decomposition (cf. \cite{13,3}). The fundamental chamber
$\cM(X)$ is distinguished from other fundamental chambers by the
fact that it contains the ray $\bR^+h$ of polarization.

By Hodge decomposition, we have the direct sum
\begin{equation}
H^2(X,\bZ)\otimes \bC=H^2(X,\bC)=H^{2,0}(X)+H^{1,1}(X)+H^{0,2}(X)
\label{HodgeK3}
\end{equation}
where $H^{2,0}(X)=\bC\omega_X$, $H^{0,2}(X)=\overline{H^{2,0}(X)}$
and $H^{1,1}(X)=(H^{2,0}(X)+H^{0,2}(X))^\perp$. Then the Picard lattice of
$X$ is
\begin{equation}
S_X=H^2(X,\bZ)\cap H^{1,1}(X)= \{x\in H^2(X,\bZ)\suchthat x\cdot
H^{2,0}(X)=0\}.
\label{PicardK3}
\end{equation}
The triplet
\begin{equation}
\left(H^2(X,\bZ),H^{2,0}(X), \cM(X)\right)
\label{periodsK3}
\end{equation}
is called {\it the periods of $X$.}

An {\it isomorphism
\begin{equation}
\phi:\left(H^2(X,\bZ),H^{2,0}(X), \cM(X)\right) \to
\left(H^2(Y,\bZ),H^{2,0}(Y), \cM(Y)\right)
\label{globaltorelliK3}
\end{equation}
of periods} of two K3 surfaces means an isomorphism
$\phi:H^2(X,\bZ) \to H^2(Y,\bZ)$ of cohomology lattices (i. e. modules with
pairing) such that
$\phi(H^{2,0}(X))=H^{2,0}(Y)$, $\phi(\cM(X))=\cM(Y)$ for the
corresponding induced maps which we denote by the same letter
$\phi$. By {\it Global Torelli Theorem for K3 surfaces \cite{13},}
$\phi$ is defined by a unique  isomorphism $f:Y\to X$ of the K3
surfaces: $\phi=f^\ast$.

As an application of the Global Torelli Theorem, let us consider the
description of $\Aut(X)$ from \cite{13}. By Serre duality,
$h^0(\cT_X)=h^2(\Omega^1_X)=h^{1,2}=0$. Thus, $X$ has no regular
vector-fields. It follows, that $\Aut(X)$ acts on $S_X$ with only
a finite kernel. Let
\begin{equation}
Sym(\cM(X))=\{\phi\in O^\prime (S_X)\suchthat
\phi(\cM(X))=\cM(X)\}
\label{symmetryM(X)K3}
\end{equation}
be the symmetry group of the fundamental chamber $\cM(X)$. Let us
denote by $Sym(\cM(X))^0$ its subgroup of finite index which
consists of all symmetries which are identical on the discriminant
group $(S_X)^\ast/S_X$. Elements $\phi\in Sym(\cM(X))^0$ can be
extended to automorphisms of $H^2(X,\bZ)$ which are identical on
the transcendental lattice $T_X=(S_X)^\perp\subset H^2(X,\bZ)$
(see Propositions \ref{overlattice}, \ref{primembedd1} in
Appendix). We denote this extension by the same letter $\phi$
since it is unique. We have $H^{2,0}(X)\subset T_X\otimes \bC$
since $H^{2,0}(X)\cdot S_X=0$. Thus,
$\phi(H^{2,0}(X))=H^{2,0}(X)$, and $\phi$ is an automorphism of
periods of $X$. Thus, $\phi=f^\ast$ where $f\in \Aut(X)$. Thus,
the natural contragradient representation
\begin{equation}
\Aut(X) \to Sym(\cM(X))
\label{autsymK3}
\end{equation}
has a finite kernel and a finite cokernel. It follows that the
groups $\Aut(X)\approx Sym(\cM(X))$ are naturally isomorphic up to
finite groups. Since we have a na\-tural isomorphism \newline
$Sym(\cM(X))\cong $ $O^\prime(S_X)/W^{(2)}(S_X)$, we also obtain
that
\begin{equation}
\Aut(X)\approx O^\prime (S_X)/W^{(2)}(S_X).
\label{autWK3}
\end{equation}
In particular, $\Aut(X)$ is finite if and only if
$[O(S_X):W^{(2)}(S_X)]<\infty$. See \cite{10}, \cite{11},
\cite{31} about the enumeration of all these cases.

 Periods $\left(H^2(X,\bZ),H^{2,0}(X), \cM(X)\right)$ of a K3
surface $X$ satisfy {\it the Riemann relation:} $H^{2,0}(X)\cdot
H^{2,0}(X)=0$ and $\omega_X\cdot \overline{\omega_X}>0$ for
$0\not=\omega_X\in H^{2,0}(X)$ (shortly we will be able to  write
$H^{2,0}(X)\cdot \overline{H^{2,0}(X)}>0$).

{\it Abstract K3 periods} is a triplet
\begin{equation}
(L_{K3},H^{2,0},\cM)
\label{abstractperiodsK3}
\end{equation}
where $L_{K3}$ is an even unimodular lattice of signature
$(3,19)$, $H^{2,0}\subset L_{K3}\otimes \bC$ is a one dimensional
complex linear subspace satisfying $H^{2,0}\cdot H^{2,0}=0$,
$H^{2,0}\cdot \overline{H^{2,0}}>0$, and $\cM$ is a fundamental
chamber of $W^{(2)}(M)\subset \cL(M)$ where $M=\{x\in
L_{K3}\suchthat x\cdot H^{2,0}=0\}$ is an abstract Picard lattice.
By {\it the surjectivity of Torelli map for K3 surfaces \cite{5},} any
abstract K3 periods are isomorphic to periods of an algebraic K3
surface.

As an application of Global Torelli Theorem and Surjectivity of
Torelli map for K3 surfaces, let us describe {\it moduli spaces of
K3 surfaces with conditions on Picard  lattice}. For details see
\cite{8} and for real case \cite{12}.

Let $M$ be an even (i. e. $x^2$ is even for any $x\in M$)
hyperbolic (i. e. of signature $(1,\rk M-1)$) lattice. Like for
$S_X$ above, we consider the light cone
\begin{equation}
V(M)=\{x\in M\otimes \bR  \suchthat x^2>0\}
\label{coneMK3}
\end{equation}
of $M$, and we choose one of its half $V^+(M)$ defining the
corresponding hyperbolic space $\cL(M)=V^+(M)/\bR^+$. We choose a
fundamental chamber $\cM(M)\subset \cL(M)$ for the reflection group
$W^{(2)}(M)$ generated by reflections in all elements $f\in M$
with $f^2=-2$. Remark that the group $\pm W^{(2)}(M)$ acts
transitively on all these additional data $(V^+(M),\cM(M))$ which
shows that they are defined by the lattice $M$ itself (i. e. by
its isomorphism class), and we can fix these additional data
$(V^+(M),\cM(M))$ without loss of generality.

We consider K3 surfaces $X$ such that a primitive sublattice
$M\subset S_X$ is fixed, $V^+(X)\cap (M\otimes \bR)=V^+(M)$,
$\cM(X)\cap \cL(M)\not=\emptyset$, and $\cM(X)\cap \cL(M)\subset
\cM(M)$. (This is one of the weakest possible conditions of
degeneration.) Such a K3 surface $X$ is called a {\it K3 surface
with the condition $M$ on Picard lattice.} A general such a K3 surface $X$ (i. e. with moduli or periods general enough) has $S_X=M$. We will show this till later.
Then $S_ X=M$, $V^+(X)=V^+(M)$, and $\cM(X)=\cM(M)$. Then one can consider
this condition as a marking of elements of the Picard lattice
$S_X$ by elements of the standard lattice $M$.

Let $(X,M\subset S_X)$ be a K3 surface with the condition $M$ on
the Picard lattice. Then $M\subset S_X\subset H^2(X,\bZ)$ defines
a primitive sublattice $M\subset H^2(X,\bZ)$. Depending on the
isomorphism class of this primitive sublattice, we obtain
different connected components of moduli of K3 surfaces with the
condition $M$ on Picard lattice.

We fix a primitive embedding $M\subset L_{K3}$. We consider marked
K3 surfaces $(X,M\subset S_X)$ with the condition $M$ on Picard
lattice and the class $M
\subset L_{K3}$ of the condition $M$ on cohomology.
Here marking means an isomorphism $\xi:H^2(X,\bZ)\cong
L_{K3}$ of lattices such that $\xi|M$ is identity. Taking
\begin{equation}
\left(L_{K3},H^{2,0}=\xi(H^{2,0}(X)),\cM=\xi(\cM(X))\right)
\label{periodmap}
\end{equation}
we obtain periods of a marked K3 surface $(X,M\subset S_X,\xi)$ with
condition $M$ on Picard lattice. By the surjectivity of Torelli map,
any abstract periods
$$
(L_{K3},H^{2,0},\cM)
$$
where $H^{2,0}\cdot M=0$, $\cM\cap \cL(M)\not=\emptyset$,  and
$\cM\cap \cL(M)\subset \cM(M)$ correspond to a marked K3 surface
with the condition $M$ on Picard lattice. Let us denote by
$\widetilde{\Omega}_{M\subset L_{K3}}$ the space of all these
abstract periods. It is called {\it the period domain of K3
surfaces $(X,M\subset S_X)$ with the condition $M$ on Picard
lattice and with the type $M\subset L_{K3}$ of this condition on
cohomology.} Let
\begin{equation} \Omega_{M\subset
L_{K3}}=\{H^{2,0}=\bC\omega \subset L_{K3}\otimes \bC \suchthat
\omega\cdot M=0,\ \omega^2=0\ and\ \omega\cdot
\overline{\omega}>0\}.
\label{Omega1}
\end{equation} We have the natural
projection $p:\widetilde{\Omega}_{M\subset L_{K3}}\to
\Omega_{M\subset L_{K3}}$ of forgetting $\cM$. The space
$\Omega_{M\subset L_{K3}}$ is an open subset of projective quadric
of the dimension $\rk L_{K3}-\rk M-2=20-\rk M$. It follows that
for a general K3 surface $X$ with the condition $M$ on Picard
lattice we have $S_X=M$. Indeed, if $\rk S_X>\rk M$ for all
K3 surfaces $X$ with the condition $M\subset L_{K3}$, then, since
$H^{2,0}\cdot \xi(S_X)=0$, periods $H^{2,0}$ define a quadric of
smaller dimension $20-\rk S_X<20-\rk M$ which leads to a
contradiction. It also follows that the forgetting map
$p:\widetilde{\Omega}_{M\subset L_{K3}}\to \Omega_{M\subset
L_{K3}}$ is an isomorphism in general points: e. g. it
is isomorphism in all points with $S_X=M$. Actually, $p$ gives an
\'etale covering which makes $\widetilde{\Omega}_{M\subset
L_{K3}}$ non-Hausdorff in special points (see
\cite{BurnsRapoport75} about construction and using of this
covering).

Considerations above also show that an even hyperbolic lattice $M$
is isomorphic to a Picard lattice $S_X$ of some K3 surface $X$ if
and only if $M$ has a primitive embedding $M\subset L_{K3}$. In
particular, this is valid if $\rk M\le \rk L_{K3}/2=11$ (see
\cite{9}): any even hyperbolic lattice $M$ of $\rk M\le 11$ is
Picard lattice of some K3 surface. See other sufficient and
necessary conditions in Theorems \ref{primembedd2} and Corollary
\ref{primembedd3} of Appendix.

The period space $\Omega_{M\subset L_{K3}}$ is a Hermitian
symmetric domain of type IV in the classification by \'E. Cartan.
The domains
$\Omega_{M\subset L_{K3}}$ and then $\widetilde{\Omega}_{M\subset
L_{K3}}$ have two connected components which are complex
conjugate. Indeed, $H^{2,0}\subset L_{K3}\otimes\bC$ is equivalent
to an oriented positive definite real subspace
$(H^{2,0}+\overline{H^{2,0}})\cap (L_{K3}\otimes \bR)\subset
L_{K3}\otimes \bR$ which is orthogonal to $M\subset L_{K3}$. Let
us consider the orthogonal complement $T=M^\perp$ in $L_{K3}$ and
the automorphism group $O(2,20-\rk M)$ of $T\otimes \bR$. Then
\begin{equation}
\Omega_{M\subset L_{K3}}=O(2,20-\rk M)/\left(SO(2)\times O(20-\rk
M)\right)
\label{periodsomegaK31}
\end{equation}
has two connected components since $SO(2)\times O(20-\rk M)$ has
index two in the maximal compact subgroup $O(2)\times O(20-\rk M)$
of $O(2,20-\rk M)$. The number of connected components of $O(2,20-\rk M)$ and
$O(2)\times O(20-\rk M)$ coincide.

Let
\begin{equation}
O(M\subset L_{K3})=\{\phi\in O(L_{K3})\ \suchthat\
\phi|M=\text{identity}\}
\label{periodsomegaK32}
\end{equation}
be the automorphism group of the period domain
$\widetilde{\Omega}_{M\subset L_{K3}}$. By Global Torelli Theorem,
the corresponding K3 surfaces are isomorphic if and only if their
periods are conjugate by this group. This group is discrete. Thus
\begin{equation}
Mod_{M\subset L_{K3}}=\widetilde{\Omega}_{M\subset
L_{K3}}/O(M\subset L_{K3})
\label{moduli1}
\end{equation}
gives the coarse {\it moduli space of K3 surfaces with the
condition $M$ on Picard lattice and with the type $M\subset
L_{K3}$ of the embedding in cohomology.} Usually $O(M\subset
L_{K3})$ contains an automorphism which permutes two connected
components of periods (equivalently it has the spinor norm $-1$,
i. e. it does not belong to a connected component of $SO(2)\times
O(20-\rk M)$ of $O(2,20-\rk M)$ above). Then the moduli space
\eqref{moduli1} is connected.

Two primitive embeddings $a:M\subset L_{K3}$, $b:M\subset L_{K3}$
give the same moduli space \eqref{moduli1}, if they are conjugate
by an automorphism of the lattice $L_{K3}$, i. e. they are
equivalent. Taking disjoint union of moduli spaces $Mod_{M\subset
L_{K3}}$ for all equivalence classes $M\subset L_{K3}$ of
primitive embeddings of lattices, we obtain the {\it moduli space
\begin{equation}
Mod_M=\bigsqcup_{class\ of\ M\subset L_{K3}}{Mod_{M \subset L_{K3}}}
\label{moduli2}
\end{equation}
of K3 surfaces with the condition $M$ on Picard lattice.} If the
primitive embedding $M\subset L_{K3}$ is unique up to
isomorphisms, and if $O(M\subset L_{K3})$ has an automorphism of
spinor norm $-1$, then the moduli space $Mod_M$ is connected. We
remark that the same results about connectedness of moduli of K3
surfaces with conditions on Picard lattice can be obtained using
only Global Torelli Theorem and local surjectivity of Torelli map
for K3 surfaces (see the paper \cite{8} which had been written
before the surjectivity of Torelli map for K3 was established).

\subsection{The lattice $S$, and the main invariants $(r,a,\delta)$
(equivalently $(k,g,\delta)$) of pairs $(X,\theta)$ and $(Y,C)$}
\label{subsec2.2} All results of this Section were obtained in
\cite{8, 9, N1, 10} (see also \cite{31}). Here we omit
some technical proofs. They will be given in Section
\ref{subsec:maininv} of Appendix.

Let $X$ be an algebraic K3 surface with a
non-symplectic involution $\theta$. (We remark that existence of a
non-symplectic involution on a K\"ahler K3 surface $X$ implies
that $X$ is algebraic (see \cite{8})).

For a module $Q$ with action of $\theta$ we denote by $Q_\pm$ the
$\pm 1$ eigenspaces of $\theta$.

The lattice (i. e. a free $\bZ$-module with a non-degenerate
symmetric bilinear form)
$$
S=H^2(X,\bZ)_+
$$
considered up to isomorphisms is called {\it the main invariant of
$(X,\theta)$.} Since $\theta$ is non-symplectic, we have
$H^{2,0}(X)\subset H^2(X,\bZ)_-\otimes \bC$. It follows that
$S\cdot H^{2,0}(X)=0$. Thus, $S\subset S_X$ is a sublattice of the
Picard lattice $S_X$ of $X$. Let $h\in S_X$ be a polarization of
$X$. Then $h_1=h+\theta^\ast h\in S$ is also a polarization of
$X$, and $h_1^2>0$. It follows that $S$ is hyperbolic like the
Picard lattice $S_X$. The rank $r=\rk S$ is one of main invariants
of $S$.

The following property of the sublattice $S\subset S_X$ is very
important: {\it The lattice $(S_X)_-$ (i. e. the orthogonal
complement to $S$ in $S_X$) has no elements $f$ with $f^2=-2$}.
Indeed, by Riemann-Roch Theorem for K3, then $\pm f$ is effective
and $\theta^\ast(\pm f)=\mp f$ which is impossible.

Let $T=S^\perp$ be the orthogonal complement to $S$ in
$H^2(X,\bZ)$. Canonical epimorphisms $H^2(X,\bZ)\to S^\ast$ and
$H^2(X,\bZ)\to T^\ast$ defined by intersection pairing give
canonical $\theta$-equivariant epimorphisms
$$
S^\ast/S\cong H^2(X,\bZ)/(S\oplus T)\cong T^\ast/T
$$
because $H^2(X,\bZ)$ is an unimodular lattice. The involution
$\theta$ is $+1$ on $S^\ast/S$, and it is $-1$ on $T^\ast/T$. It
follows that the groups $S^\ast/S\cong T^\ast/T\cong (\bZ/2\bZ)^a$
are $2$-elementary. Only in this case multiplications by $\pm 1$
coincide. Thus, the lattice $S$ is $2$-elementary, which means
that its discriminant group ${\gA}_S=S^\ast/S\cong (\bZ/2\bZ)^a$
is $2$-elementary where $a$ gives another important invariant of
$S$.

There is one more invariant $\delta$ of $S$ which takes values in
$\{0, 1\}$. One has $\delta=0 \iff (x^*)^2\in \bZ$ for every
$x^*\in S^* \iff $ {\it the discriminant quadratic form of $S$}
$$
q_S:{\gA}_S=S^\ast/S\to \bQ/2\bZ,\ \ q_S(x^*+S)=(x^*)^2+2\bZ
$$
is even: it takes values in $(\bZ/2\bZ)$. See Appendix, Section
\ref{subsec:discrforms} about discriminant forms of lattices.

{\it The invariants $(r,a,\delta)$ of $S$ define the isomorphism
class of a 2-elementary even hyperbolic lattice $S$.} See more
general statement and the proof in Appendix, Section
\ref{subsec:maininv} and Theorem \ref{2elementary}.

Thus, any two even
hyperbolic 2-elementary lattices with the same invariants $(r,a,\delta)$
are isomorphic. The invariants $(r,a,\delta)$ of $S$ are
equivalent to the main invariant $S$, and we later call them as
{\it the main invariants of a K3 surface $X$ with non-symplectic
involution $\theta$.}

Vice versa, let $S$ be a hyperbolic even 2-elementary lattice
having a primitive embedding to $L_{K3}$. Let $S\subset L_{K3}$ be
one of primitive embeddings. Considering $T=S^\perp$ in $L_{K3}$
and the diagram similar to above,
$$
S^\ast/S\cong L_{K3}/(S\oplus T)\cong T^\ast/T\ ,
$$
we obtain that there exists an involution $\theta^\ast$ of
$L_{K3}$ which is $+1$ on $S$, and $-1$ on $T$. Let us consider
the moduli space
\begin{equation}
 Mod_S^{\,\prime}\subset Mod_S
\label{moduliS}
\end{equation}
of K3 surfaces $(X,S\subset S_X)$ with condition $S$ on Picard
lattice (see \eqref{moduli2}) where for $(X,S\subset S_X)$ from
$Mod_S^{\,\prime}$ we additionally assume that the orthogonal
complement $S^\perp$ to $S$ in $S_X$ has no elements with square
$(-2)$. One can easily see that $Mod_S^{\,\prime}$ is Zariski open
subset in $Mod_S$. Any general $(X,S\subset S_X)$ (i. e. when
$S=S_X$) belongs to $Mod_S^{\,\prime}$. Thus, the difference between
$Mod_S^{\,\prime}$ and $Mod_S$ is in complex codimension one, and
they have the same connected components. By Global Torelli
Theorem, the action of $\theta^\ast$ in $L_{K3}$ can be lifted to
a non-symplectic involution $\theta$ on $X$ with $H^2(X,\bZ)_+=S$.
Thus, the moduli space $Mod_S^{\,\prime}$ in \eqref{moduliS} can
be considered as {\it the moduli space of K3 surfaces with
non-symplectic involution and the main invariant $S$.} Since $S$
is defined by the main invariants $(r,a,\delta)$, it can also be
denoted as
$$
\cM_{(r,a,\delta)}=Mod_S^{\,\prime}
$$
and can be considered as moduli space of K3 surfaces with
non-symplec\-tic involution and the main invariants
$(r,a,\delta)$. Any even hyperbolic 2-elementary lattice $S$ has a
unique primitive embedding to $L_{K3}$ (up to isomorphisms) if it
exists. Then the group $O(S\subset L_{K3})$ always has an
automorphism of spinor norm $-1$. Thus, {\it the moduli space
$\cM_{(r,a,\delta)}$ is connected.}

Evidently, to classify all possible main invariants $S$
(equivalently $(r,a,\delta)$) one just needs to classify all even
hyperbolic 2-elementary lattices $S$ having a primitive embedding
$S\subset L_{K3}$. All such possibilities for $(r,a,\delta)$
(equivalently, $(k=(r-a)/2,g=(22-r-a)/2,\delta)$, see below) are
known and are shown on Figure~\ref{fig1}.


\begin{figure}
\centerline{\includegraphics[width=10cm]{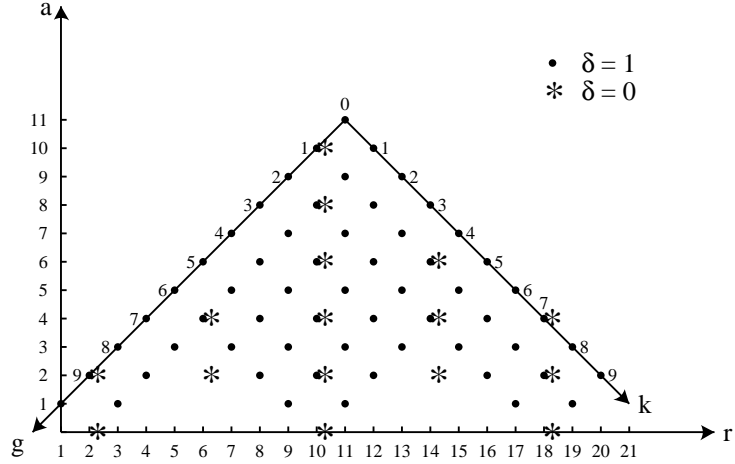}} \caption{All
possible main invariants $(r,a,\delta)$} \label{fig1}
\end{figure}

The triple $(r,a,\delta)$ admits an interpretation in terms of
$X^{\theta} = C$. If $(r,a,\delta) \ne (10,8,0)$ or $(10,10,0)$
then
$$
X^{\theta}=C = C_g + E_1 + \dots + E_k,
$$
where $C_g$ is a nonsingular irreducible curve of genus $g$, and $E_1,
\dots E_k$ are nonsingular irreducible rational curves, the curves are
disjoint to each other,
$$
g = (22-r-a)/2,\ \ k=(r-a)/2
$$
(we shall formally use the same formulae for $g$ and $k$ even in cases
$(r,a,\delta)=(10,8,0)$ or $(10,10,0)$).
If $(r,a,\delta)=(10,8,0)$ then $X^{\theta}\simeq C =
C_1^{(1)} +  C_1^{(2)}$, where $C_1^{(i)}$ are elliptic (genus 1)
curves. If $(r,a,\delta) = (10,10,0)$ then $X^{\theta}= C
=\emptyset$, i.e. in this case $Y$ is an Enriques surface. One has
\begin{equation}
\delta=0 \iff X^{\theta} \sim 0\mod 2\ \text{ in } S_X\
\text{(equivalently in $H_2(X,\bZ)$)}\
\iff
\label{deltaX}
\end{equation}
there exist signs $(\pm)_i$ for which
\begin{equation}
\dfrac{1}{4}\sum_{i}{(\pm)_i cl(C^{(i)})}\in S_Y=H^2(Y,\bZ),
\label{deltaY}
\end{equation}
where $C^{(i)}$ go over all irreducible components of $C$. Signs
$(\pm)_i$ for $\delta=0$ are defined uniquely up to a simultaneous
change.  They define a new natural orientation (different from the complex one)
of the components of $C$; a positive sign gives the complex
orientation and a negative sign the opposite orientation.

The main invariants $S$, equivalently $(r,a,\delta)$ (or
$(k,g,\delta)$) of K3 surfaces with non-symplectic involution, and
the corresponding DPN pairs and DPN surfaces play a crucial role
in our classification.

\subsection{Exceptional curves on $(X,\theta)$ and $Y$}
\label{subsec2.3} A description of exceptional curves on a DPN
surface $Y$ can also be reduced to the K3 surface $X$ with a
non-symplectic involution $\theta$ considered above.

\begin{figure}
\centerline{\includegraphics[width=13cm]{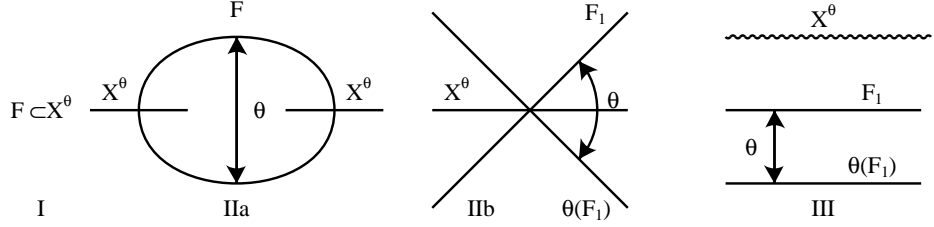}}
\caption{Pictures of exceptional curves.}
\label{fig2}
\end{figure}

Let $(X,\theta)$ be a K3 surface with a non-symplectic involution
and $(Y = X/\{1,\theta\}, C=\pi(X^{\theta}))$ the corresponding
DPN pair. If $E\subset Y$ is an exceptional curve, then the curve
$F=\pi^\ast (E)_{red}$ is either an irreducible curve with negative square
on the K3 surface $X$, or $F=F_1 + \theta(F_1)$, where
$$
F^2=F_1^2 + \theta(F_1)^2 + 2 F_1\cdot \theta(F_1) = 2E^2<0.
$$
The curves $F_1$ and $\theta(F_1)$ are irreducible and have
negative square (i.e. equal to $(-2)$, see Section
\ref{subsec2.1a}).  Using this, in an obvious way we get the
following four possibilities for $E$ and $F$ (see Fig. \ref{fig2})

\begin{enumerate}
\item[I] $E^2 = -4$, $E$ is a component of $C$; respectively $F$ is a
  component of $X^{\theta}$, and $F^2=-2$.
\item[IIa] $E^2=-1$, $E\cdot C =2$ and $E$ intersects $C$
transversally at two points; respectively $F$ is irreducible and $F^2=-2$.
\item[IIb] $E^2=-1$, $E\cdot C=2$ and $E$ is tangent to $C$ at one
  point; respectively $F=F_1+ \theta(F_1)$, where $(F_1)^2=-2$ and
$F_1\cdot \theta(F_1) =1$.
\item[III] $E^2=-2$, $E\cap C =\emptyset$; respectively $F=F_1+\theta(F_1)$,
where $(F_1)^2=-2$ and $F_1\cdot \theta(F_1) = 0$.
\end{enumerate}


If $Y$ is an
Enriques surface, we let $S_Y$ be the Picard lattice of $Y$ modulo torsion.
Let $P(Y)\subset S_Y$ be the subset of divisor
classes of all exceptional curves $E$ on $Y$, and $P(X)_+ \subset
S = (S_X)_+$ be the set of divisor classes of all $F= \pi\inv(E)$.
We call them {\it exceptional classes of the pair $(X,\theta)$.}
By what we said above, $P(Y)$ and $P(X)_+$ are divided into
subsets:
\begin{eqnarray*}
  P(Y) &=&
  P(Y)_{I} \bigsqcup P(Y)_{IIa} \bigsqcup P(Y)_{IIb} \bigsqcup P(Y)_{III}. \\
  P(X)_+ &=&
  P(X)_{+I} \bigsqcup P(X)_{+IIa} \bigsqcup P(X)_{+IIb} \bigsqcup P(X)_{+III}.
\end{eqnarray*}
By projection formula,
$$
\pi^\ast(\NEF(Y))=\NEF(X)\cap
(S\otimes\bR)=\NEF(X)_+\,.
$$
In the same way as for K3 surfaces $X$ in Section
\ref{subsec2.1a}, we have
\begin{equation}
\NEF(Y)=\{y\in \overline{V^+(Y)} \suchthat y\cdot P(Y)\ge 0\},
\label{NEFY}
\end{equation}
a locally finite polyhedron in $V^+(Y)$ whose facets are
orthogonal and numerated by elements of $P(Y)$. Since
$\pi^\ast(\NEF(Y))=\NEF(X)\cap (S\otimes\bR)=\NEF(X)_+$, we obtain
that
\begin{equation}
\NEF(X)_+=\NEF(X)\cap (S\otimes\bR)=\{x\in \overline{V^+(S)}
\suchthat x\cdot P(X)_+\ge 0\} \label{NEFX+}
\end{equation}
is a locally finite polyhedron whose facets are orthogonal and
numerated by elements of $P(X)_+$. Here $V^+(S)=V^+(X)\cap
(S\otimes \bR)$.

As for K3 surfaces in Section \ref{subsec2.1a}, we can interpret
the above results using hyperbolic spaces. Since the lattice
$S=(S_X)_+$ is hyperbolic, and $S\subset S_X$, we have embeddings
of cones
$$
V(S) \subset V(S_X), \qquad V^+(S) \subset V^+(S_X)=V^+(X)
$$
and hyperbolic spaces $\cL(S)=V^+(S)/\bR^+ \subset
V^+(X)/\bR^+=\cL(X)$.

If $h$ is a polarization of $X$, the set $\cL(S)$ contains the
polarization ray $\bR^+(h+\theta^*h)$ of $X$. Therefore, $\cM(X)_+
=\NEF(X)_+/\bR^+= \cM(X) \cap \cL(S)$ is a non-degenerate (i.e.
containing a nonempty open subset of $\cL(S)$) convex locally
finite polytope in $\cL(S)$. Since $S$ is even, it is easy to see
that $P(X)_+$ is precisely the set of primitive elements of $S$,
perpendicular to facets of $\cM(X)_+$ and directed outward. One
has
\begin{eqnarray*}
  P^{(2)}(X)_+ &\overset{{\rm Def}}{=}&
  \{ f\in P(X)_+ \suchthat f^2=-2 \} \\
  &=& P(X)_{+I} \bigsqcup P(X)_{+IIa} \bigsqcup P(X)_{+IIb} , \\
  P^{(4)}(X)_+ &\overset{{\rm Def}}{=}&
  \{ f\in P(X)_+ \suchthat f^2=-4 \} = P(X)_{+III}.
\end{eqnarray*}
Moreover, $\cM(X)_+$, like $\cM(X)$ for K3 surfaces in Section
\ref{subsec2.1a}, admits a description in terms of groups.

Indeed, by Section \ref{subsec2.1a}
\begin{equation}
\cM(X)_+=\{\bR^+x\in \cL(S)\suchthat  x\cdot f\ge 0\}
\label{CM(X)_+2}
\end{equation}
for any effective $f\in S_X$ with $f^2=-2$. Let us write
$f=f_+^\ast + f_-^\ast$ where $f_+^\ast \in S^\ast$ and $f_-^\ast
\in (S_X)_-^\ast$. We have $2f_+^\ast=f+\theta^\ast(f)\in S$ and
$2f_-=f-\theta^\ast(f)\in (S_X)_-$. It follows $f=(f_++f_-)/2$
where $f_+\in S$ and $f_-\in (S_X)_-$. If $f_+^2\ge 0$, then
$f_+\cdot V^+(S)\ge 0$, and $f$ does not influence in
\eqref{CM(X)_+2}. Thus, in \eqref{CM(X)_+2} we can assume that
$f_+^2<0$. Since $(S_X)_-$ is negative definite and the lattice
$S_X$ is even, we then obtain that $f_+\in \Delta^{(2)}_{+}\cup
\Delta^{(4)}_{+}$ defined below.

Let
\begin{eqnarray*}
  \Delta^{(4)}_{\pm} &=&
\{f_{\pm}\in (S_X)_{\pm} \suchthat f^2_{\pm} =-4,
\text{ and } \exists f_{\mp} \in (S_X)_{\mp},\\
&&\text{ for which } f_{\mp}^2=-4
\text{ and } (f_{\pm} + f_{\mp})/2 \in S_X \}; \\
  \Delta^{(2)}_{+} &=& \Delta^{(2)}(S) =
\{f_{+}\in S \suchthat f^2_{+} =-2\}; \\
  \Delta^{(2)}_{+t} &=&
\{f_{+}\in \Delta^{(2)}(S) \suchthat \exists f_{-} \in (S_X)_{-},\\
&&\text{ for which } f_{-}^2=-6
\text{ and } (f_{+} + f_{-})/2 \in S_X \}; \\
  \Delta^{(6)}_{-} &=&
\{f_{-}\in (S_X)_- \suchthat f^2_{-} =-6
\text{ and } \exists f_{+} \in \Delta^{(2)}_{+t}, \\
&&\text{ for which } (f_++f_-)/2\in S_X \}.
\end{eqnarray*}

If $f_{\pm}\in \Delta_{\pm}^{(4)}$, then $f_{\pm}\cdot (S_X)_{\pm}
\equiv 0 \mod 2$. Hence, $f_{\pm}\in \Delta_{\pm}^{(4)}$ are roots
of $(S_X)_\pm$, and there exists a reflection $s_{f_{\pm}} \in
O'((S_X)_{\pm})$ with respect to element $f_{\pm}$:
$$
s_{f_{\pm}} (x) = x+ \frac{(x\cdot f_{\pm})}{2} f_{\pm}, \quad
x\in (S_X)_{\pm}.
$$
One has a very important property:
\begin{equation}
s_{f_{\pm}}(\Delta_{\pm}^{(2)}\cup \Delta_{\pm}^{(4)}) =
\Delta_{\pm}^{(2)}\cup \Delta_{\pm}^{(4)}\ \ \forall f_{\pm}\in
\Delta_{\pm}^{(2)}\cup \Delta_{\pm}^{(4)}
\label{ref1}
\end{equation}
where we formally put $\Delta_-^{(2)}=\emptyset$ because the
lattice $(S_X)_-$ has no elements $f_-$ with $f_-^2=-2$
(see the previous section).

Let us prove \eqref{ref1}. Assume $f_+\in \Delta_{+}^{(2)}$. The
reflection $s_{f_+}\in O(S_X)$ and $s_{f_+}((S_X)_\pm)=(S_X)_\pm$.
It follows \eqref{ref1} for such $s_{f_+}$. Assume $f_+\in
\Delta_{+}^{(4)}$. Then there exists $f_-\in \Delta_{-}^{(4)}$
such that $\alpha_1=(f_+ + f_-)/2\in S_X$. The element
$\alpha_2=(f_+ + f_-)/2-f_-=(f_+ - f_-)/2$ also belongs to $S_X$.
We have $\alpha_1^2=\alpha_2^2=-2$. Thus, the reflections
$s_{\alpha_1}$ and $s_{\alpha_2}$ belong to $O(S_X)$. It follows
that $s=s_{\alpha_2}s_{\alpha_1}\in O(S_X)$. On the other hand,
simple calculation shows that $s((S_X)_\pm)=(S_X)_\pm$, and $s$ in
$(S_X)_\pm$ coincides with the reflection $s_{f_\pm}$. It follows
that
$$
s_{f_{+}}(\Delta_{+}^{(2)}\cup \Delta_{+}^{(4)})=
s(\Delta_{+}^{(2)}\cup \Delta_{+}^{(4)}) = \Delta_{+}^{(2)}\cup
\Delta_{+}^{(4)}.
$$
For $f_-\in \Delta_{-}^{(2)}\cup \Delta_{-}^{(4)}$ the arguments
are the same. This simple but very important trick had been first used
by Dolgachev \cite{19} for Enriques surfaces.

By \eqref{ref1}, reflections with respect to all the elements of
$\Delta^{(2)}(S) \cup \Delta_+^{(4)} = \Delta_+^{(2,4)}(S)$
generate a group $W_+^{(2,4)} \subset O'(S)$ which geometrically
is the group generated by reflections in the hyperplanes of $\cL(S)$
which are orthogonal to $\Delta_+^{(2,4)}(S)$, any reflection in a
hyperplane of $\cL(S)$ from this group is reflection in an element
of $\Delta_+^{(2,4)}(S)$. It follows (by exactly the same
considerations as for the K3 surface $X$ in Section
\ref{subsec2.1a}) that $\cM(X)_+$ is a fundamental chamber for
$W_+^{(2,4)}$. Thus,
$$
P(X)_+ = P(\cM(X)_+)
$$
is the set of primitive elements of $S$, which are orthogonal to
facets of $\cM(X)_+$ and directed outward. We obtain the
description of $P(X)_+$ and $P(Y)$ using the reflection group
$W_+^{(2,4)}$.

Denote by
$$
A(X,\theta) = \{\phi\in O'(S) \suchthat \phi(\cM(X)_+) = \cM(X)_+ \}
$$
the subgroup of automorphisms of $\cM(X)_+$ in $O'(S)$ and by
$\Aut(X,\theta)$ the normalizer of the involution $\theta$ in
$\Aut X$. The action of $\Aut(X,\theta)$ on $S_X$ and $S$ defines
a contravariant representation
\begin{equation}
\Aut(X,\theta) \to A(X,\theta).
\label{Aut1}
\end{equation}
Like for K3 surfaces $X$ in Section \ref{subsec2.1a}, Global Torelli
theorem for K3 surfaces \cite{13} implies that this representation
has a finite kernel and cokernel. Therefore, it defines an
isomorphism up to finite groups: $\Aut(X,\theta) \approx
A(X,\theta)$.

\subsubsection{Computing $P(X)_+$}
\label{subsubsec2.3.1}

First, we consider calculation of the fundamental chamber
$\cM^{(2)} \subset \cL(S)$ of $W^{(2)}(S)$.

For that, it is important to consider a larger group $W^{(2,4)}(S)$
generated by reflections
in all elements of $\Delta^{(2)}(S)$ and all elements of
$$
\Delta^{(4)}(S)=\{f\in S \suchthat f^2=-4\ \text{and}\
f\cdot S\equiv 0 \mod 2\}
$$
of all roots with square $(-4)$ of the lattice $S$. Both
sets $\Delta^{(2)}(S)$ and $\Delta^{(4)}(S)$ are invariant with
respect to $W^{(2,4)}(S)$. It follows that every reflection from
$W^{(2,4)}(S)$ gives a hyperplane $\cH_f$ where $f\in
\Delta^{(2)}(S)\cup \Delta^{(4)}(S)$.  The subgroup $W^{(2)}(S)
\triangleleft W^{(2,4)}(S)$ is normal, and any reflection from
$W^{(2)}(S)$ is reflection in an element of $\Delta^{(2)}(S)$.
Similarly, the subgroup $W^{(4)}(S)\triangleleft W^{(2,4)}(S)$,
generated by reflections in $\Delta^{(4)}(S)$, is normal and any
reflection from $W^{(4)}(S)$ is reflection in an element of
$\Delta^{(4)}(S)$.

This implies the following description of a fundamental chamber
$\cM^{(2)}$ of $W^{(2)}(S)$. Let $\cM^{(2,4)}\subset \cL(S)$ be
a fundamental chamber of $W^{(2,4)}(S)$. It will be extremely
important for our further considerations.
Let $P^{(2)}(\cM^{(2,4)})$ and $P^{(4)}(\cM^{(2,4)})$ be elements of
$\Delta^{(2)}(S)$ and $\Delta^{(4)}(S)$ respectively directed outwards and
orthogonal to $\cM^{(2,4)}$ (i. e. to facets of $\cM^{(2,4)}$).

\begin{proposition} Let $W^{(4)}(\cM^{(2)})$ be the group
generated by reflections in all elements of $P^{(4)}(\cM^{(2,4)})$.

Then the fundamental chamber $\cM^{(2)}$ of $W^{(2)}(S)$
containing $\cM^{(2,4)}$ is equal to
$$
\cM^{(2)}=W^{(4)}(\cM^{(2)})(\cM^{(2,4)}),
$$
$$
P(\cM^{(2)})=W^{(4)}(\cM^{(2)})(P^{(2)}(\cM^{(2,4)})).
$$
Moreover,
$$
W^{(4)}(\cM^{(2)})=\{w\in W^{(4)}(S)\suchthat w(\cM^{(2)})=\cM^{(2)}\}.
$$
Reflections which are contained in $W^{(4)}(\cM^{(2)})$ are
exactly the reflections in elements
$$
\Delta^{(4)}(\cM^{(2)})=\{f\in \Delta^{(4)}(S)\suchthat
\cH_f\ \text{intersects the interior of}\  \cM^{(2)}\}.
$$
\label{prop2.3.1}
\end{proposition}

\begin{proof} This easily follows from the facts that
  $W^{(2)}(S)\triangleleft W^{(2,4)}(S)$ and $W^{(4)}(S)\triangleleft
  W^{(2,4)}(S)$ are normal subgroups, and $\Delta^{(2)}(S)\cap
  \Delta^{(4)}(S)=\emptyset$. We leave details to the reader.
\end{proof}

\begin{example} Let us consider the  hyperbolic 2-elementary lattice
$S=\langle 2 \rangle \oplus 5\langle -2 \rangle$ with the invariants
$(r,a,\delta)=(6,6,1)$. Here and in what follows we denote by
$\langle A \rangle$ the integral lattice given by an integral symmetric
matrix $A$ in some its basis. We denote by $\oplus$ the orthogonal
sum of lattices.

The Dynkin diagram of $W^{(2,4)}(S)$ (which is equivalent to
the Gram matrix of all elements of $P(\cM^{(2,4)}$)) is

\medskip

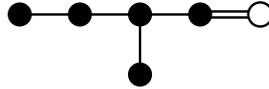
\centerline{\includegraphics[width=4cm]{pics/c-p70.eps}}


\medskip

\noindent (see \cite{3}). Here black vertices correspond to
elements from \linebreak $P^{(4)}(\cM^{(2,4)})$ and white vertices
correspond to elements from \linebreak $P^{(2)}(\cM^{(2,4)})$ (see
Section \ref{subsec3.1} below about edges).  From the diagram, one
can see that $W^{(4)}(\cM^{(2)})$ is the Weyl group of the root
system $\bD_5$, the $\Delta^{(4)}(\cM^{(2)})$ is the root system
$\bD_5$, the set $P(\cM^{(2)})=P^{(2)}(\cM^{(2)})$ is the orbit of
the Weyl group of $\bD_5$ applied to the unique element of
$P^{(2)}(\cM^{(2,4)})$ which corresponds to the white vertex.
Calculations show that the set $P(\cM^{(2)})$ consists of $16$
elements, and it is not easy to draw their Dynkin (or Gram)
diagram, but it is completely defined by the diagram above which
has only $6$ vertices. \label{2.3.2}
\end{example}

Now let us consider a subset $\Delta_+^{(4)}\subset \Delta^{(4)}(S)$
and the subgroup $W_+^{(2,4)}$ of reflections generated by this subset
and by the set $\Delta^{(2)}(S)$.  As in our case \eqref{ref1}, we
shall assume that the set $\Delta_+^{(4)}$ is $W_+^{(2,4)}$-invariant.
Then each reflection from $W_+^{(2,4)}$ is a reflection in a
hyperplane $\cH_f$, $f\in \Delta^{(2)}(S)\cup\Delta_+^{(4)}$. As
before, $W^{(2)}(S)\triangleleft W_+^{(2,4)}$ is a normal subgroup. We
denote by $W^{(4)}_+$ the group generated by reflections in
$\Delta_+^{(4)}$, it is normal in $W_+^{(2,4)}$ as well. Thus, for a
fundamental chamber $\cM_+^{(2,4)}\subset \cM^{(2)}$ of $W_+^{(2,4)}$
we can similarly define $P^{(4)}(\cM_+^{(2,4)})$, $P^{(2)}(\cM_+^{(2,4)})$
(which
are the sets of all elements in $\Delta_+^{(4)}$ and
$\Delta^{(2)}(S)$ respectively which are orthogonal to $\cM_+^{(2,4)}$),
the group
$W_+^{(4)}(\cM^{(2)})$ generated by reflections in $P^{(4)}(\cM_+^{(2,4)})$,
the set
$$
\Delta^{(4)}_+(\cM^{(2)})=W_+^{(4)}(\cM^{(2)})\left(P^{(4)}
(\cM_+^{(2,4)})\right).
$$
We get similar statements to Proposition \ref{prop2.3.1}:
$$
W_+^{(4)}(\cM^{(2)})=\{w\in W_+^{(4)}\suchthat w(\cM^{(2)})=\cM^{(2)}\},
$$
$$
\Delta_+^{(4)}(\cM^{(2)})=\{f\in \Delta_+^{(4)}\suchthat
\cH_f\ \text{intersects the interior of}\ \cM^{(2)}\},
$$
the group $W_+^{(4)}(\cM^{(2)})$ is generated by
reflections in the $\Delta_+^{(4)}(\cM^{(2)})$, and every
reflection from $W_+^{(4)}(\cM^{(2)})$ is a reflection in a hyperplane
$\cH_f$, $f\in \Delta_+^{(4)}(\cM^{(2)})$.

Obviously, the fundamental chamber $\cM_+^{(2,4)}\subset \cM^{(2)}$
for $W_+^{(2,4)}$ is the fundamental chamber of the group
$W_+^{(4)}(\cM^{(2)})$ considered as a group acting on
$\cM^{(2)}$.

Let us show how one can calculate a fundamental chamber
$\cM_+^{(2,4)}$ of $W_+^{(2,4)}$ contained in the fixed fundamental
chamber $\cM^{(2)}$ of $W^{(2)}(S)$.
\begin{proposition} We have:
$$
P(\cM^{(2)})=W_+^{(4)}(\cM^{(2)})(P^{(2)}(\cM_+^{(2,4)}))
$$
and
$$
P^{(2)}(\cM_+^{(2,4)})=\{f\in P(\cM^{(2)})\suchthat
f\cdot P^{(4)}(\cM_+^{(2,4)})\ge 0\}.
$$
\label{prop2.3.3}
\end{proposition}

\begin{proof} The first statement is analogous to Proposition \ref{prop2.3.1}.
We denote the right hand side of the proving second equality as $P^{(2)}$.
Since $P^{(2)}\subset P(\cM^{(2)})$ and $\cM^{(2)}$ has acute angles,
$f\cdot f^\prime\ge 0$ for any two different elements
$f,\,f^\prime \in P^{(2)}\cup P^{(4)}(\cM_+^{(2,4)})$.
It follows that $P(\cM_+^{(2,4)})\subset P^{(2)}\cup P^{(4)}(\cM_+^{(2,4)})$
because the fundamental chamber $\cM_+^{(2,4)}$ must have acute angles. Then
$$
\bigcap_{f\in P^{(2)}\cup P^{(4)}(\cM_+^{(2,4)})}
{\cH_f^+} \subset \cM_+^{(2,4)}
$$
where the left hand side is not empty. Indeed, it contains
the non-empty subset
$$
\bigcap_{f\in P(\cM^{(2)})\cup P^{(4)}(\cM_+^{(2,4)})}
{\cH_f^+}=
\cM^{(2)}\bigcap \left(\bigcap_{f\in P^{(4)}(\cM_+^{(2,4)})}
{\cH_f^+} \right)\supset\cM_+^{(2,4)}.
$$
It follows (see Proposition 3.1 in  \cite{3}) that
$$
P(\cM_+^{(2,4)})=P^{(2)}\cup P^{(4)}(\cM_+^{(2,4)})
$$
because for all $f\not=f^\prime\in P^{(2)}\cup P^{(4)}(\cM_+^{(2,4)})$
we have $f\cdot f^\prime\ge 0$.
\end{proof}

Propositions \ref{prop2.3.1} and \ref{prop2.3.3} imply the
result which will be very important in further considerations.

\begin{theorem}
Let $\cM^{(2,4)}$ be a fundamental chamber of $W^{(2,4)}(S)$
in $\cL(S)$,
and $W^{(4)}(\cM^{(2)})$ the group generated by reflections in
all elements of $P^{(4)}(\cM^{(2,4)})$, and
$\Delta^{(4)}(\cM^{(2)})=W^{(4)}(\cM^{(2)})(P^{(4)}(\cM^{(2,4)}))$.

Then

(1) Subsets $\Delta_+^{(4)}\subset \Delta^{(4)}(S)$ which are invariant
for the group $W_+^{(2,4)}$ generated by reflections in all elements of
$\Delta^{(2)}(S)\cup \Delta_+^{(4)}$ are in one-to-one correspondence with
subsets $\Delta_+^{(4)}(\cM^{(2)})\subset \Delta^{(4)}(\cM^{(2)})$ which are
invariant for the group $W_+^{(4)}(\cM^{(2)})$ generated by reflections in
all elements of $\Delta_+^{(4)}(\cM^{(2)})$. This correspondence is given
by
$$
\Delta_+^{(4)}(\cM^{(2)})=\Delta_+^{(4)}\cap \Delta^{(4)}(\cM^{(2)});\ \
\Delta_+^{(4)}=W^{(2)}(S)(\Delta_+^{(4)}(\cM^{(2)})).
$$

(2) The fundamental chamber $\cM^{(2)}$ of $W^{(2)}(S)$ containing
$\cM_+^{(2,4)}$ is $\cM^{(2)}=W_+^{(4)}(\cM^{(2)})(\cM_+^{(2,4)})$.
Moreover,
$$
P(\cM^{(2)})=W_+^{(4)}(\cM^{(2)})(P^{(2)}(\cM_+^{(2,4)})).
$$

(3) Under the one-to-one correspondence in (1), any fundamental
chamber $\cM_+^{(2,4)}\subset \cM^{(2)}$ of $W_+^{(2,4)}$ can be
obtained as follows: Let $\cM_+^{(4)}(\cM^{(2)})$ be a fundamental
chamber for $W_+^{(4)}(\cM^{(2)})$.  Then
$$
\cM_+^{(2,4)}=\cM^{(2)}\cap \cM_+^{(4)}(\cM^{(2)})\ \text{and}\
P^{(4)}(\cM_+^{(2,4)})=P^{(4)}(\cM_+^{(4)}(\cM^{(2)})),
$$
$$
P^{(2)}(\cM_+^{(2,4)})=\{f\in W_+^{(4)}(\cM^{(2)})(P^{(2)}(\cM_+^{(2,4)}))
\suchthat f\cdot P^{(4)}(\cM_+^{(2,4)})\ge 0\}.
$$
\label{thm2.3.4}
\end{theorem}

\begin{proof} Only the statement (1) requires some clarification.
  Assume that $\Delta_+^{(4)}(\cM^{(2)})\subset
  \Delta^{(4)}(\cM^{(2)})$ is invariant with respect to the subgroup
  $W_+^{(4)}(\cM^{(2)})$ generated by reflections in
  $\Delta_+^{(4)}(\cM^{(2)})$. Note that $\cM^{(2)}$ is invariant
  with respect to $W_+^{(4)}(\cM^{(2)})$.  It follows that the
  fundamental chamber $\cM_+^{(2,4)}$ for $W_+^{(4)}(\cM^{(2)})$
  acting on $\cM^{(2)}$ will be the fundamental chamber for the group
  $W_+^{(2,4)}$ generated by reflections in all elements of
  $\Delta^{(2)}(S)$ and $\Delta_+^{(4)}(\cM^{(2)})$. It follows that
  $\Delta_+^{(4)}=W^{(2)}(S)(\Delta_+^{(4)}(\cM^{(2)}))$ is invariant
  with respect to $W_+^{(2,4)}$. It follows that
  $\Delta_+^{(4)}(\cM^{(2)})=\Delta_+^{(4)}\cap
  \Delta^{(4)}(\cM^{(2)})$.

The remaining statements are obvious.
\end{proof}

In Chapter \ref{sec3} we apply this theorem to describe DPN surfaces of
elliptic type.

\subsection{The root invariant of a pair $(X,\theta)$}
\label{subsec2.4}
To describe $W_+^{(2,4)}$ and sets $P(X)_{+III}$, and $P(X)_{+IIa}$,
$P(X)_{+IIb}$, one should add to the main invariants $(r,a,\delta)$
(equivalently $(k,g,\delta)$) of $(X,\theta)$ the so-called {\it root
  invariants}. We describe them below.  The root
invariants for DPN surfaces had been introduced and considered in
\cite{11}  and \cite{31}.

Everywhere below we follow Appendix, Section \ref{subsec:discrforms}
about lattices and discriminant forms of lattices.

Let $M$ be a lattice (i. e. a non-degenerate integral symmetric
bilinear form). Following \cite{9}, $M(k)$ denotes a lattice
obtained from $M$ by multiplication of the form of $M$ by $k\in
\bQ$.

Let $K(2)$ be the sublattice of $(S_X)_-$ generated by
$\Delta_-^{(4)}\subset (S_X)_-$. Since
$\Delta_-^{(4)}\cdot(S_X)_-\equiv 0\mod 2$, the lattice $K=K(2)(1/2)$
is integral and is generated by its subset $\Delta_-^{(4)}\subset K$
of elements with square $(-2)$ defining in $K$ a root system, since
reflections with respect to elements of $\Delta_-^{(4)}$ send
$\Delta_-^{(4)}$ to itself. If follows that the lattice $K$ is
isomorphic to the orthogonal sum of root lattices $A_n$, $D_m$ and
$E_k$ corresponding to the root systems $\bA_n$, $\bD_m$, $\bE_k$ (or
their Dynkin diagrams), and $\Delta_-^{(4)}=\Delta^{(2)}(K)$ is the
set of all elements of $K$ with square $(-2)$. Equivalently,
$\Delta_-^{(4)}=\Delta^{(4)}(K(2))$ is the set of all elements with
square $(-4)$ of $K(2)$. Moreover, we have a natural homomorphism of
groups
\begin{equation}
\xi:Q=\frac{1}{2}K(2)/K(2)\to {\gA}_S=S^\ast/S
\label{xi}
\end{equation}
such that $\xi(\frac{1}{2}f_-+K(2))=\frac{1}{2}f_+ + S$, if
$f_{\mp}\in \Delta_{\mp}^{(4)}$ and $(f_- + f_+)/2\in S_X$. This
defines a homomorphism of finite quadratic forms $\xi:q_{K(2)}|Q\to
-q_S$ with values in $\frac{1}{2}\bZ/2\bZ\subset \bQ/2\bZ$. (Here
$q_M:{\gA}_M=M^\ast/M\to \bQ/2\bZ$ denotes the discriminant quadratic
form of an even lattice $M$.) The homomorphism $\xi$ is equivalent to
the homomorphism (which we denote by the same letter $\xi$) of the
finite quadratic forms
\begin{equation}
\xi:K\mod 2\to -q_S
\label{xi1},
\end{equation}
by the natural isomorphism $\frac{1}{2}K(2)/K(2)\cong K/2K$, where
$q_{K(2)}|Q$ is replaced by the finite quadratic form
$\frac{1}{2}x^2\mod 2$ for
$x\in K$.

We define the {\it root invariant of the pair $(X,\theta)$ or
the corresponding DPN pair $(Y,C)$} as
the equivalence class of the triplet
\begin{equation}
R(X,\theta)=(K(2),\Delta_-^{(4)},\xi)\cong (K(2),K(2)^{(4)},\xi)\cong
(K,\Delta^{(2)}(K),\xi),
\label{rootinvariant}
\end{equation}
up to isomorphisms of lattices $K$ and automorphisms of the
lattice $S$. Clearly, similarly we can introduce
{\it abstract root invariants,} without any
relation to K3 surfaces with involutions and DPN pairs; see
beginning of Section \ref{subsec2.6} below.

We have the following statement from \cite{9}.

\begin{lemma} Let $S$ be an even hyperbolic $2$-elementary lattice.

Then the natural homomorphism $O(S)\to O(q_S)$ is surjective.
\label{lemma2.4.1}
\end{lemma}

\begin{proof} We remind the proof from \cite{9}.
If $\rk S\ge 3$, this follows from the general theorem 1.14.2 in
\cite{9}. If $\text{rk} S=2$, then $S\cong U=
\left\langle\begin{array}{cc}
0 & 1\\
1 & 0
\end{array}\right\rangle
$, $U(2)$, $\langle 2 \rangle\oplus \langle -2\rangle$. If $\rk S=1$,
then $S=\langle 2\rangle$. For all these lattices one can
check the statement directly. See Appendix Theorems
\ref{uniqueness1}, \ref{2elementary} for more details.
\end{proof}

\begin{lemma}Let $S$ be an even hyperbolic $2$-elementary lattice and
$\rk S\ge 2$. Then every $x\in {\gA}_S$ with
$q_S(x)=n/2\mod 2$, $n\in \bZ$, can be represented as $x=u/2\mod S$ where
$u\in S$ and $u^2=2n$.
\label{lemma2.4.2}
\end{lemma}

\begin{proof} If
$\rk S=2$, then $S\cong U$, $U(2)$, $\langle 2\rangle \oplus
\langle -2 \rangle$, and the statement can be checked directly.
Assume that the statement is valid for $\rk S\le k$ where $k \ge
2$. Let $\rk S=k+1$. By Theorem 1.12.2 from \cite{9} about
existence of an even lattice with a given discriminant quadratic
form (see Appendix, Theorem \ref{discrexistence}), we get that
$S=S^\prime\oplus T$ where $S^\prime$ is a hyperbolic
$2$-elementary lattice of the rank $\rk S^\prime \ge 2$, and $T$
is a negative definite $2$-elementary lattice of the $\rk T\ge 1$.
Let $x=y\oplus z$$, x\in {\gA}_S$, $y\in {\gA}_T$, and assume
$z=u/2\mod T$ where $u\in T$ and $u^2=2m$, $m\in \bZ$ . By the
induction assumption, there exists $v\in S^\prime$ with $y=v/2\mod
S^\prime$ and $v^2=2n-2m$ since
$q_{S^\prime}(y)=q(x)-q(z)=(n-m)/2\mod 2$.
\end{proof}

\begin{lemma} Let $q:A\to \bQ/2\bZ$ be a non-degenerate quadratic form
  on a finite $2$-elementary group $A$ and $\phi:H_1\cong H_2$ be an
  isomorphism of two subgroups in $A$ which preserves $q|H_1$ and
  $q|H_2$. Assume that the characteristic element $a_q$ of $q$ on $A$
  either does not belong to both subgroups $H_1$ and $H_2$ or belongs
  to both of them. In the second case we additionally assume that
  $\phi(a_q)=a_q$.

Then $\phi$ can be extended to an automorphism of $q$.
\label{lemma2.4.3}
\end{lemma}

\begin{proof} See Proposition 1.9.1 in \cite{12} (we repeated the proof in
Appendix, Proposition \ref{propWitt2}). We remind that $a_q\in A$ is
the {\it characteristic element of $q$,} if $q(x)\equiv
b_q(x,a_q)\mod 1$ for any $x\in A$. Here $b_q$ is the bilinear
form of $q$. This defines the characteristic element $a_q$
uniquely.
\end{proof}

Lemmas \ref{lemma2.4.1} --- \ref{lemma2.4.3} imply

\begin{proposition}
The root invariant $R(X,\theta)$ of $(X,\theta)$ (or $(Y,C)$) is
equivalent to the triplet
$$
R(X,\theta)=(K(2);H;\alpha,\overline{a})\cong(K;H;\alpha,\overline{a}).
$$
Here $H=\text{Ker\ }\xi$ is an isotropic for $q_{K(2)}$ subgroup in $Q$
(equivalently in $K\mod 2$); $\alpha=0$, if $\xi(Q)=\xi(K\mod 2)$
contains the characteristic element $a_{q_S}$ of the quadratic form $q_S$,
and $\alpha=1$ otherwise; if $\alpha=0$, the element
$\overline{a}=\xi^{-1}(a_{q_S})+H\in Q/H$ ; if $\alpha=1$, the element
$\overline{a}$ is not defined.
\label{prop2.4.4}
\end{proposition}

The root invariant $R(X,\theta)$ is important because it defines
$$
\Delta_+^{(4)}=\{f_+\in S \suchthat \  f_+^2=-4,\
f_+/2\mod S\in
$$
$$
\xi(\frac{1}{2}\Delta^{(4)}(K(2))\mod K(2))
=\xi(\Delta^{(2)}(K)\mod 2K)\},
$$
and $W_+^{(2,4)}$, $\cM(X)_+$,
$P(\cM(X)_+)=P(X)_+$, up to the action of $O(S)$.
Moreover, for $f_+\in P(X)_+^{(4)}$,
the root invariant defines the decomposition $f_+=f+\theta^\ast(f)$,
$f\in P(X)$, uniquely up to permutation of $f$ and $\theta^\ast(f)$.
More precisely, we have the
following. Let $f_-\in \Delta^{(4)}(K(2))$ and $\xi(f_-/2\mod
K(2))=f_+/2\mod S$. Then $f=(f_+ \pm f_-)/2$, $\theta(f)=(f_+\mp
f_-)/2$. Indeed, if ${f_-}^\prime \in \Delta^{(4)}(K(2))$ satisfies
the same conditions, then $(f_-+{f_-}^\prime)/2\in K(2)$. In $K(2)$,
if ${f_-}^\prime\not=\pm f_-$ then either $f_-\cdot{f_-}^\prime=0$ or
$f_-\cdot{f_-^\prime}=\pm 2$.  If $f_-\cdot{f_-}^\prime=0$ then
$\left((f_-+{f_-}^\prime)/2\right)^2=-2$, and we get a contradiction
since $(S_X)_-$ does not have elements with the square $(-2)$.  If
$f_-\cdot{f_-}^\prime=\pm 2$ then $f_-\cdot(f_-+{f_-}^\prime)/2=-2\pm
1$, and we get a contradiction since $f_-\cdot (S_X)_-\equiv 0 \mod
2$.  Thus, ${f_-}^\prime=\pm f_-$, and the pair of elements $f$ and
$\theta^\ast(f)$ is defined uniquely.

Similarly one can define a {\it generalized root invariant}
$$
R_{gen}(X,\theta)=(K_{gen}(2),\ \Delta_-^{(4)}\cup
\Delta_-^{(6)},\ \xi_{gen})\cong
$$
$$
(K_{gen},\Delta^{(2)}_-\cup \Delta^{(3)}_-,\xi_{gen}),
$$
where for $f_-\in \Delta^{(6)}_-$ one has
$$
\xi_{gen}(f_-/2\mod K_{gen}(2))=f_+/2\mod S
$$
where $f_+\in \Delta_{+t}^{(2)}(S)$ and $(f_++f_-)/2\in S_X$. Here
$K_{gen}(2)\subset (S_X)_-$
is generated by $\Delta_-^{(4)}\cup\Delta_-^{(6)}$.

Using Lemmas \ref{lemma2.4.1} --- \ref{lemma2.4.3}, one can similarly prove
that it is equivalent to
the tuple
$$
R_{gen}(X,\theta)=(K_{gen}(2),\Delta_-^{(4)}\cup\Delta_-^{(6)};
\,H_{gen};\,\alpha_{gen},\,\overline{a}_{gen})\cong
$$
$$
(K_{gen},\Delta^{(2)}_-\cup   \Delta^{(3)}_-;\,H_{gen};\,
\alpha_{gen},\,\overline{a}_{gen}).
$$
It is defined similarly to the root invariant.

Importance of the generalized root invariant is that it contains
the root invariant $R(X,\theta)$. Thus,  it defines $W_+^{(2,4)}$,
$\cM(X)_+$ and also \linebreak$P(\cM(X)_+) = P(X)_+$, up to the
action of $O(S)$. But, it also defines
$$
\Delta_{+t}^{(2)}=\{f_+\in S\ \suchthat\ (f_+)^2=-2,\ f_+/2\mod S\in
$$
$$
\xi_{gen}(\frac{1}{2}\Delta^{(6)}_-\mod K_{gen}(2))=\xi(\Delta^{(3)}_-
\mod 2K_{gen})\},
$$
and then it defines
$$
P(X)_{+IIb}=P^{(2)}(X)_{+t}=\{f_+\in P(X)_+\ |\ f_+\in \Delta_{+t}^{(2)}\}.
$$

Thus, using root invariants, we know how to find \newline
$P(X)_{+III}  =P^{(4)}(X)_+$, $P(X)_{+IIb}=P^{(2)}(X)_{+t}$,
and hence, we know $P(X)_{+I}\cup P(X)_{+IIa}=P^{(2)}(X)_+\,-\,P^{(2)}
(X)_{+t}$. To distinguish
$P(X)_{+I}$ and $P(X)_{+IIa}$, it is sufficient to know
$P(X)_{+I}$.

\subsection{Finding  the locus $X^\theta$}
\label{subsec2.5} Here we show how one can
find $P(X)_{+I}$. This is based on the following considerations
(similar to \cite{10}):

1) Since $W^{(2)}(S)\triangleleft W_+^{(2,4)}$ is a normal
subgroup, the fundamental chamber $\cM_+^{(2,4)}$ is contained in
one fundamental chamber $\cM^{(2)}$ of $W^{(2)}(S)$; we have
$\cM_+^{(2,4)}\subset \cM^{(2)}$. One can consider replacing
$\cM_+^{(2,4)}$ by $\cM^{(2)}$ as a deformation of a pair
$(X,\theta)$ to a general pair
$(\widetilde{X},\widetilde{\theta})$ having $S_{\widetilde{X}}=S$,
$\cM(\widetilde{X})=\cM^{(2)}$ and
$P(\widetilde{X})=P(\cM^{(2)})$. See Section \ref{subsec2.2} on
corresponding results about moduli.

The divisor classes of fixed points of the involution do not change
under this deformation, thus $P(X)_{+I}=P(\widetilde{X})_{+I}$. In
particular, $P(X)_{+I}$ does not change when a root invariant changes
(with fixed main invariants $(r,a,\delta)$ equivalent to the lattice
$S$).

2) Let $\delta_1$,  $\delta_2$ belong to $P^{(2)}(X)_+$ and
$\delta_1\cdot \delta_2=1$, i. e. the curves $D_1$, $D_2$ corresponding
to them intersect transversally. Then one of $\delta_1$,
$\delta_2$ belongs to $P(X)_{+I}$, and another to $P(X)_{+II}=
P(X)_{+IIa}\cup P(X)_{+IIb}=P(X)_{+IIa}$ for the general case we consider.
See the diagrams below where an element of $P(X)_{+I}$ is denoted by a double
transparent vertex, and an element of $P(X)_{+II}$ by a single transparent
vertex.

\medskip

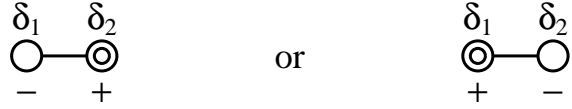
\centerline{\includegraphics[width=8cm]{pics/c-p77.eps}}

\medskip

Indeed, the intersection point $D_1\cap D_2$ is a fixed point of
$\theta$, tangent directions of $D_1$ and $D_2$ at this point are
eigenvectors of $\theta_\ast$.  We know that they have eigenvalues
$+1$ and $-1$.

3) If the Gram diagram of elements $\delta_0$, $\delta_1$, $\delta_2$
and $\delta_3$ from $P^{(2)}(X)_+$ has the form as shown

\medskip

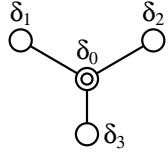
\centerline{\includegraphics[width=2.5cm]{pics/c-p78.eps}}


\medskip

\noindent
then $\delta_0\in P(X)_{+I}$, and $\delta_1,\,\delta_2,\,\delta_3\in
P(X)_{+II}=P(X)_{+IIa}$ (for the general case we consider).  Indeed,
the rational curve corresponding to $\delta_0$ has three different
fixed points of $\theta$, and hence belongs to $X^\theta$.

4) If $\delta\in P(X)_{+III}=P^{(4)}(X)_+$ and $\delta_1\in P(X)_{+I}$,
then $\delta_1\cdot \delta=0$. This is obvious from the definition
of $P(X)_{+III}$.

Considering all possible lattices $S$, it is not difficult to see that
statements 1) --- 3) are sufficient for finding $P(X)_{+I}$ and the divisor
class of the irreducible component $C_g$ of the curve $X^\theta$
of fixed points. The statement 4) simplifies these considerations, if
some elements of $P^{(4)}(X)_+$ are known.

\subsection{Conditions for the existence of root invariants}
\label{subsec2.6}  Assume that the main invariants $(r,a,\delta)$
(equivalently $(k,g,\delta)$) of $(X,\theta)$ are known and fixed.
Here we want to give conditions which are necessary and sufficient for
the existence of a pair $(X,\theta)$ with a given root or
generalized root invariant. We consider the root invariant.
Similarly, one can consider the generalized root invariant.

Assume that $(K,\Delta^{(2)}(K),\xi)$ is the root invariant of a pair
$(X,\theta)$.

Then the conditions 1 and 2 below must be satisfied:

\medskip

{\bf Condition 1.} {\it The lattice
$$
K_H=[K;x/2\ \text{where\ }x+2K\in H]
$$
does not have elements with the square $(-1)$. Equivalently, the lattice
$$
K_H(2)=[K(2);x/2\ \text{where\ }x/2+K(2)\in H]$$
does not have elements with square $(-2)$. We remind that $H=\Ker \xi$.}

\medskip

Indeed, the lattice $K_H(2)\subset (S_X)_-$, but the lattice $(S_X)_-$
does not have elements with square $(-2)$.

\medskip

{\bf Condition 2.} $\rk S+\rk K=r+\rk K\le 20$.

\medskip
Indeed, $S\oplus K(2)\subset S_X$ and $\rk S_X\le 20$.

\medskip

A pair $(X,\theta)$ (or the corresponding DPN pair $(Y,C)$) is called
{\it standard,} if $K_H(2)$ is a primitive
sublattice of $(S_X)_-$, and the primitive sublattice
$[S\oplus K(2)]_{\pr}$ in $S_X$
generated by $S\oplus K(2)$ is defined by the homomorphism $\xi$,
i. e. it is equal to
\begin{equation}
\begin{split}
M=&[S\oplus K(2); \{ a+b\suchthat  \forall a\in S^\ast,\
\forall b\in K(2)/2,\\ &\text{such that\ }
\xi(b+K(2))=a+S\}].
\end{split}
\label{2.6.1}
\end{equation}
Clearly, $M\subset [S\oplus K(2)]_{\pr}$ is always a sublattice of
finite index.

\medskip

Let {\it $l(\gA)$ be the minimal number of generators of a finite
Abelian group $\gA$.} Let {\it $\gA_M=M^\ast/M$ be the
discriminant group of a lattice $M$.}

Let us consider an {\it abstract root invariant $(K(2), \xi)$.}
This means that $K$ is a negative definite lattice generated by its
elements with square $(-2)$, and $K(2)$ is obtained by
multiplication of the form of $K$ by $2$. The map
$$
\xi: q_{K(2)}|Q=\frac{1}{2}K(2)/K(2)\to -q_S
$$
is a homomorphism of finite quadratic forms. We assume that for
each $f_-\in \Delta^{(4)}(K(2))$ there exists $f_+\in
\Delta^{(4)}(S)$ such that $\xi (f_-/2 + K(2))= f_+/2+S$ (by Lemma
\ref{lemma2.4.2}, this condition is always satisfied). As above, we denote
$H=\Ker \xi$.

\begin{proposition}\label{prop2.6.1}
A standard pair $(X,\theta)$ with a given root invariant $(K,\xi)$
satisfying Conditions 1 and 2 does exist, if additionally
\begin{equation}\label{2.6.2}
r+a+2l(H)<22
\end{equation}
and
\begin{equation}\label{2.6.3}
r+\rk K+l({\gA}_{K_p})<22
\end{equation}
for any prime $p>2$. Here $K_p=K\otimes \bZ_p$ where $\bZ_p$ is
the ring of $p$-adic integers.
\end{proposition}

\begin{proof} By Global Torelli Theorem \cite{13} and surjectivity
of Torelli Map \cite{5} for K3 (see Section \ref{subsec2.1a}), the
pair $(X,\theta)$ does exist, if there exists a primitive
embedding of the lattice $M$ described in \eqref{2.6.1} into an
even unimodular lattice $L_{K3}\cong H^2(X,\bZ)$ of the signature
$(3,19)$ (see the proof of Proposition \ref{prop2.6.2} below). By
Corollary 1.12.3 in \cite{9} (see Appendix, Corollary
\ref{primembedd3}), such a primitive embedding does exist, if $\rk
M+l(\gA_{M_p})<22$ for all prime $p\ge 2$.

If $p>2$, then $\rk M+l(\gA_{M_p})=r+rk K+l(\gA_{K_p})<22$ by \eqref{2.6.3}
(here we remember that  the lattice $S$ is 2-elementary).

Assume that $p=2$. Let $\Gamma_\xi$ be the graph of $\xi$. Then
$\gA_M=(\Gamma_\xi)^\perp/\Gamma_\xi$ for the discriminant
form $q_S\oplus q_{K(2)}$. Let $Q=H\oplus Q^\prime$ (see \eqref{xi})
where $Q^\prime$ is a complementary subgroup, and $\xi^\prime=\xi|Q^\prime$.
Then $\Gamma_{\xi^\prime}\subset \Gamma_\xi$,
moreover $\Gamma_{\xi^\prime}\subset \gA_S\oplus \gA_{K(2)}=\gA$ is a
2-elementary subgroup, and $\Gamma_{\xi^\prime}\cap 2\gA=\{0\}$ since
$2\gA=\{0\}\oplus2\gA_{K(2)}$ and $\xi^\prime$ is injective. Let
$\gA^{(2)}$ be the kernel of multiplication by $2$ in $\gA$, and
$q=(q_S\oplus q_{K(2)})|\gA^{(2)}$. It is easy to see that the
kernel $\Ker q=\gA^{(2)}\cap 2\gA$.
Since $\Gamma_{\xi^\prime}\cap 2\gA=\{0\}$,
then $\Gamma_{\xi^\prime}\cap \Ker q=\{0\}$.
Let $\gA^{(2)}_1$ be a subgroup in $\gA^{(2)}$ which is
complementary to $\Ker q$ and contains $\Gamma_{\xi^\prime}$.
Then $q_S\oplus q_{K(2)}=q_1\oplus q_2$ where
$q_1=q_S\oplus q_{K(2)}|\gA_1^{(2)}$ and $q_2$ is the orthogonal complement
to $q_1$ (since $q_1$ is non-degenerate). The subgroup $\Gamma_{\xi^\prime}$
is isotropic for the non-degenerate 2-elementary form $q_1$ and has
rank $\rk K-\rk H$. It follows that
$$
l((\Gamma_{\xi^\prime})^\perp_{q_1}/\Gamma_{\xi^\prime})=
l({\gA}_1^{(2)})-2l(\Gamma_{\xi^\prime}),
$$
and then
$$
l({\gA}_{M_2})\le l(\gA)-2l(\Gamma_{\xi^\prime})=a+\rk K-2(\rk K-l(H)).
$$
This implies that
$$
\rk M+l({\gA}_{M_2})\le r+a+2l(H)<22
$$
by \eqref{2.6.2}.
\end{proof}

Finally, in general, by Global Torelli Theorem \cite{13} and
surjectivity of Torelli Map \cite{5} for K3 (see Section
\ref{subsec2.1a}), we have the following necessary and sufficient
conditions of existence of a pair $(X,\theta)$ with a root
invariant $(K(2),\xi)$. It even takes under consideration the more
delicate invariant which is the isomorphism class of embedding of
lattices $M\subset H^2(X,\bZ)\cong L_{K3}$.

\begin{proposition}\label{prop2.6.2}
  There exists a K3 pair $(X,\theta)$ with a root invariant
  $(K(2),\xi)$ and the isomorphism class of embedding $\phi:M\subset
  L_{K3}$ of lattices (see \eqref{2.6.1}), if and only if

1) $\phi(S)\subset L_{K3}$ is a primitive sublattice;

2) the primitive sublattice $\phi(K(2))_{\pr}\subset L_{K3}$
generated by $\phi(K(2))$ in $L_{K3}$ does not have elements with
square $(-2)$;

3) we have:
\begin{equation}
\begin{split}
\phi(\Delta^{(4)}(K(2)))=&\{f_-\in \phi(K(2))_{\pr}\suchthat
f_-^2=-4,\ \text{and}\ \exists f_+\in S\\ &\text{such that\ }
f_+^2=-4\ \text{and}\ (\phi(f_+) + f_-)/2\in L_{K3}\}.
\end{split}
\label{2.6.3.a}
\end{equation}
Here, the right hand side always contains the left hand side.
We remind that $\Delta^{(4)}(K(2))$ is the set of all elements in $K(2)$
with square $(-4)$.
\end{proposition}

\begin{proof}
Using \eqref{2.6.1}, we construct an even lattice $M$
which contains $S\oplus K(2)\subset M$ as a sublattice of finite index.
It contains $S\subset M$ as a primitive
sublattice, and its primitive sublattice generated by $K(2)$
is $K(2)_H$ where $H=\Ker \xi$.

Let $\phi:M\to L_{K3}$ be an embedding of lattices. If $\phi$ corresponds
to a K3 pair $(X,\theta)$, with the root invariant $(K(2),\xi)$,
then conditions 1), 2) and 3) must be satisfied.
Now we assume that they are valid for the abstract embedding
$\phi:M\to L_{K3}$ of lattices we consider.

Then $\phi(S)\subset L_{K3}$ is a primitive sublattice. To simplify
notation, we identify $S=\phi(S)\subset L_{K3}$ and $K(2)=\phi (K(2))$.
Since $S$ is 2-elementary, there exists an involution $\alpha$ on $L_{K3}$
with $(L_{K3})_+=S$ and $(L_{K3})_-=S^\perp$. Then $\alpha=-id$ on $K(2)$.
We denote by $\widetilde{M}$ the primitive sublattice in $L_{K3}$
generated by $\phi(M)=M$.

Assume that $f\in \widetilde{M}$ satisfies $f^2=-2$, $f=f_-^\ast+f_+^\ast$
where $f_-^\ast \in (K(2)_{\pr})^\ast$, $f_+^\ast \in S^\ast$ and
$(f_+^\ast)^2<0$.
Since $2f_-^\ast=f_-=f-\alpha(f) \in K(2)_{\pr}$, $2f_+^\ast=f_+=
f+\alpha(f)\in S$, $K(2)$ is negative definite and satisfies 2), it follows
that either $f=(f_-+f_+)/2$ where $f_-=0$ and $f=f_+/2\in \Delta^{(2)}(S)$,
or $f=(f_-+f_+)/2$ where $f_-\in K(2)^{(4)}$, $f_+\in \Delta^{(4)}(S)$, or
$f=(f_-+f_+)/2$ where $(f_-)^2=-6$ and $f_+\in \Delta^{(2)}(S)$.

It follows that there exists $h_+\in S$ with $(h_+)^2>0$ such that
$h_+\cdot f\not=0$ for any $f\in \Delta^{(2)}(\widetilde{M})$.

By surjectivity of Torelli map for K3 surfaces \cite{5}, we can
assume that there exists a K3 surface $X$ with
$H^2(X,\bZ)=L_{K3}$, $S_X=\widetilde{M}$ and a polarization $h_+$.
The involution $\alpha$ preserves periods of $X$. By Global
Torelli Theorem for K3 \cite{13}, $\alpha=\theta^\ast$ corresponds
to an automorphism $\theta$ of $X$. The automorphism $\theta$ is
non-symplectic because $H^2(X,\bZ)_+=(S_X)_+=\widetilde{M}_+=S$ is
hyperbolic. By 3), the root invariant of $(X,\theta)$ is
$(K(2),\xi)$. See Sections \ref{subsec2.1a} and \ref{subsec2.2} about
the used results on K3 surfaces.
\end{proof}

We remark that from the proof above we can even describe the moduli
$Mod_{(S,K(2),\xi,\phi)}$ of K3 surfaces with a non-symplectic
involution $\theta$ having the main invariant $S$, the root
invariant $(K(2),\xi)$ and the embedding $\phi:M\to L_{K3}$ of the
corresponding lattice $M$ which satisfies conditions of
Proposition \ref{prop2.6.2}. As in the proof we denote by
$\widetilde{M}\supset M$ the overlattice of $M$ of finite index
such that $\phi(\widetilde{M})\subset L_{K3}$ is the primitive
sublattice in $L_{K3}$ generated by $\phi(M)$.

We consider a fundamental chamber $\cM(\widetilde{M})$ for
$W^{(2)}(\widetilde{M})$ such
that $\cM(\widetilde{M})\cap \cL(S)\not=\emptyset$. Then
$\cM(\widetilde{M})\cap \cL(S)$ defines a unique $\cM(S)$ containing
$\cM(\widetilde{M})\cap \cL(S)$.
Up to isomorphisms of
the pair $S\subset \widetilde{M}$ there exists only finite number of
such $\cM(\widetilde{M})$. We have
\begin{equation}
Mod_{(S,K(2),\xi,\phi)}=\cup_{class\ of\ \cM(\widetilde{M})}
{Mod_{(S,K(2),\xi,\phi,\cM(\widetilde{M}))}}
\label{moduliSrootinv1}
\end{equation}
where
\begin{equation}
Mod_{(S,K(2),\xi,\phi,\cM(\widetilde{M}))}\subset
Mod_{\phi: \widetilde{M}\subset L_{K3}}\cap
Mod_{\phi:S\subset L_{K3}}^{\,\prime}
\label{moduliSrootinv2}
\end{equation}
consists of K3 surfaces $(X,\widetilde{M}\subset S_X)$ with the
condition $\widetilde{M}$ and $\cM(\widetilde{M})$ on the Picard
lattice and the class $\phi:\widetilde{M}\subset L_{K3}$ of the
embedding on cohomology; moreover $X$ has a non-symplectic
involution $\theta$ with the main invariant $S$ (i. e.
$(X,\theta)\in Mod_{\phi:S\subset L_{K3}}^{\,\prime})$ and
$(X,\theta)$ has the root invariant $(K(2),\xi)$. See
\eqref{moduli1}, \eqref{moduliS}. A general such a K3 surface $X$
has $S_X=\widetilde{M}$, and the dimension of moduli is equal to
\begin{equation}
\dim Mod_{(S,K(2),\xi,\phi)}= 20-\rk S-\rk K(2).
\label{dimmoduliSrootinv}
\end{equation}
Taking union by different classes of embeddings $\phi:M\subset L_{K3}$
(their number is obviously finite),
we obtain the moduli space of K3 surfaces $X$ with a non-symplectic involution
$\theta$, and the main invariant $S$, and the root invariant $(K(2),\xi)$.

Proposition \ref{prop2.6.2} implies the following result important for us.

\begin{corollary}\label{cor2.6.3}
Let $(K(2),\xi)$ be the root invariant of a pair $(X,\theta)$ and
$K^\prime(2)\subset K(2)$ a primitive sublattice of $K(2)$ generated by
its elements $\Delta^{(4)}(K^\prime(2))$ with the square $(-4)$.

Then the pair $(K^\prime(2),\xi^\prime=\xi|Q^\prime=
\frac{1}{2}K^\prime(2)/K^\prime(2))$ is also
the root invariant of some K3 pair $(X^\prime,\theta^\prime)$.

If the pair $(X,\theta)$ was standard, the pair $(X^\prime,\theta^\prime)$
also can be taken standard.
\end{corollary}

Corollary \ref{cor2.6.3} shows that to describe all possible root
invariants of
pairs $(X,\theta)$, it is enough to describe all possible root invariants
of extremal pairs. Here a {\it pair} $(X^\prime,\theta^\prime)$ is called
{\it extremal,} if its root invariant
$R(X^\prime,\theta^\prime)=(K^\prime(2),\xi^\prime)$ cannot be obtained
using Corollary \ref{cor2.6.3} from the root invariant
$R(X,\theta)=(K(2),\xi)$ of any other pair $(X,\theta^\prime)$ with
$\rk K(2)>\rk K^\prime(2)$.

\medskip

\subsection{Three types of non-symplectic involutions of K3 surfaces
and DPN surfaces.}
\label{subsec2.7} It is natural to divide non-symplectic
involutions $(X,\theta)$ of K3 and the corresponding DPN surfaces in three
types:

{\it Elliptic type:} $X^\theta\cong C\cong C_g+E_1+\cdots +E_k$
where $C_g$ is an irreducible curve of genus $g\ge 2$ (equivalently,
$(C_g)^2>0$), and $E_1,\dots,E_k$ are non-singular irreducible
rational curves. By Section \ref{subsec2.2}, this is equivalent to
$r+a\le 18$ and $(r,a,\delta)\not=(10,8,0)$. Then $\Aut (X,\theta)$ is finite
because $(C_g)^2>0$, see \cite{N1} , \cite{10} and Section \ref{subsec3.1} below.

{\it Parabolic type:} Either $X^\theta\cong C\cong C_1+E_1+\cdots
+E_k$ (using the same notation), or $X^\theta\cong C\cong
C_1^{(1)}+C_1^{(2)}$ is sum of two elliptic (i. e. of genus 1)
curves. By Section \ref{subsec2.2}, this is equivalent to either
$r+a=20$ and $(r,a,\delta)\not=(10,10,0)$, or
$(r,a,\delta)=(10,8,0)$. Then $\Aut (X,\theta)$ is Abelian up to
finite index and usually non-finite, see \cite{N1},  \cite{10}.
Here $(C_1)^2=0$.

{\it Hyperbolic type:} Either $X^\theta\cong C\cong E_0+\cdots
+E_k$ is sum of non-singular irreducible rational curves, or
$X^\theta=\emptyset$. By Section \ref{subsec2.2}, this is equivalent
to either $r+a=22$, or $(r,a,\delta)=(10,10,0)$. Then $\Aut
(X,\theta)$ is usually non-Abelian up to finite index, see
\cite{N1}, \cite{10}. Here $C_1=E_0$ has $C_1^2=-2$, if
$X^\theta\not=\emptyset$.

Thus, pairs $(X,\theta)$ of elliptic type are the most simple, and
we describe them completely in Chapter \ref{sec3}. On the other hand,
for classification of log del Pezzo surfaces of index $\le 2$
we need only these pairs.

\section{Classification of DPN surfaces of elliptic type}
\label{sec3}

\subsection{Fundamental chambers of $W^{(2,4)}(S)$ for elliptic type.}
\label{subsec3.1} The most important property of the lattices $S$ for
elliptic type is that the subgroup $W^{(2)}(S)\subset O(S)$ has
finite index.
We remark that it is parallel to Lemma \ref{lemma1.3.2}, and it
is an important step to prove that log del Pezzo surfaces of index $\le 2$
are equivalent to DPN surfaces of elliptic type.

This finiteness was first observed and used for classification of
hyperbolic lattices $M$ with finite index $[O(M):W^{(2)}(M)]$ in
\cite{N1}, \cite{10}. We repeat arguments of \cite{N1}, \cite{10}.
Let us take a general pair $(X,\theta)$ with $(S_X)_+=S$. Then
$S_X=S$, and the involution $\theta$ of $X$ is unique by the
condition that it is identical on $S_X=S$ and is $-1$ on the
orthogonal complement to $S_X$ in $H^2(X,\bZ)$. Thus, $\Aut X=\Aut
(X,\theta)$. By Global Torelli Theorem for K3 (see \cite{13}), the
action of $\Aut X$ on $S_X$ gives that $\Aut X$ and
$O(S_X)/W^{(2)}(S_X)$ are isomorphic up to finite groups. In
particular, they are finite simultaneously. Thus,
$[O(S):W^{(2)}(S)]$ is finite, if and only if $\Aut (X,\theta)$ is
finite. If $(X,\theta)$ has elliptic type, then $\Aut (X,\theta)$
preserves $X^\theta$ and its component $C_g$ with $(C_g)^2>0$.
Since $S_X$ is hyperbolic, it follows that the action of $\Aut
(X,\theta)$ in $S_X$ is finite. But it is known for K3 (see
\cite{13}) that the kernel of this action is also finite. It
follows that $\Aut (X,\theta)$ and $[O(S):W^{(2)}(S)]$ are finite.
See more details on used results about K3 in Section
\ref{subsec2.1a}.

Since $O(S)$ is arithmetic, $W^{(2)}(S)$ has a fundamental chamber
$\cM^{(2)}$ in $\cL(S)$ of finite volume and with a finite number of faces
(e. g. see \cite{3}). Since $W^{(2)}(S)\subset W^{(2,4)}(S)\subset O(S)$,
the same is valid for $W^{(2,4)}(S)$.

Let $\cM^{(2,4)}\subset \cL(S)$
be a fundamental chamber of $W^{(2,4)}(S)$, and
$\Gamma(P(\cM^{(2,4)}))$ its Dynkin diagram (see \cite{3}).
Vertices corresponding to different elements
$f_1,\ f_2\in P(\cM^{(2,4)})$ {\it are not connected by any edge,}
if $f_1\cdot f_2=0$. They are connected by a {\it simple edge of the
weight $m$} (equivalently, by {\it $m-2$ simple edges,} if $m>2$ is small), if
$$
{2\,f_1\cdot f_2\over \sqrt{f_1^2f_2^2}}=2\cos{\frac{\pi}{m}},  \ \ m\in \bN.
$$
They are connected by a {\it thick edge,} if
$$
{2\,f_1\cdot f_2\over \sqrt{f_1^2f_2^2}}=2.
$$
They are connected by a {\it broken edge of the weight $t$,} if
$$
{2\,f_1\cdot f_2\over \sqrt{f_1^2f_2^2}}=t>2.
$$
Moreover, a vertex corresponding to $f\in P^{(4)}(\cM^{(2,4)})$
is {\it black.} It is {\it transparent,} if $f\in P^{(2)}(\cM^{(2,4)})$.
It is {\it double transparent,} if $f\in P(X)_{+I}$ (i. e.
it corresponds to the class of a rational component of
$X^\theta$), otherwise, it is {\it simple transparent.}
Of course, here we assume that $\cM^{(2,4)}\subset \cM(X)_+$
for a K3 surface with involution $(X,\theta)$ and
$(S_X)_+=S$.

Classification of DPN surfaces of elliptic type is based on the
purely arithmetic calculations of the fundamental chambers $\cM^{(2,4)}$
(equivalently, of the graphs $\Gamma(P(\cM^{(2,4)})$) of the
reflection groups $W^{(2,4)}(S)$ of the lattices $S$ of elliptic type. Since
$S$ is 2-elementary and even, $W^{(2,4)}(S)=W(S)$ is the full reflection
group of the lattice $S$, any root $f\in S$ has $f^2=-2$ or $-4$. We have

\begin{theorem}\label{thm3.1.1} 2-elementary even hyperbolic lattices $S$
of elliptic type have fundamental chambers $\cM^{(2,4)}$ for their
reflection groups $W^{(2,4)}(S)$
(it is the full reflection group of $S$),
equivalently the corresponding Dynkin
diagrams $\Gamma(P(\cM^{(2,4)}))$, which are given in Table 1 below,
where the lattice $S$ is defined by its invariants
$(r,a,\delta)$ (equivalently,
$(k,g,\delta)$), see Section \ref{subsec2.2}.
\end{theorem}

\noindent
\begin{table}
\label{table1}
\caption{Fundamental chambers $\cM^{(2,4)}$ of reflection
groups $W^{(2,4)}(S)$ for 2-elementary even hyperbolic lattices $S$ of
elliptic type.}

\index{Table 1}

\addtocontents{toc}{\contentsline {section}{\tocsection {}{T.1}{Table 1}}
{\pageref{table1}}}

\begin{tabular}{|r||r|r|r|r|r|c|c|}
\hline
 $N$& $r$& $a$&$\delta$& $k$&$g$&$l(H)$& $\Gamma(P(\cM^{(2,4)}))$\\
\hline
  1 &  1 &  1 &  1     &  0 &10 &   0  &
\includegraphics[width=2cm]{pics/c-p83-1.eps}\\
\hline
  2 &  2 &  2 &  0     &  0 & 9 &   0  &
\includegraphics[width=1cm]{pics/c-p83-2.eps}\\
\hline
  3 &  2 &  2 &  1     &  0 & 9 &   0  &
\includegraphics[width=1cm]{pics/c-p83-3.eps}\\
\hline
  4 &  3 &  3 &  1     &  0 & 8 &   0  &
\includegraphics[width=2.5cm]{pics/c-p83-4.eps}\\
\hline
  5 &  4 &  4 &  1     &  0 & 7 &   0  &
\includegraphics[width=2.5cm]{pics/c-p83-5.eps}\\
\hline
  6 &  5 &  5 &  1     &  0 & 6 &   0  &
\includegraphics[width=2.5cm]{pics/c-p83-6.eps}\\
\hline
  7 &  6 &  6 &  1     &  0 & 5 &   0  &
\includegraphics[width=3cm]{pics/c-p83-7.eps}\\
\hline
  8 &  7 &  7 &  1     &  0 & 4 &   0  &
\includegraphics[width=3.5cm]{pics/c-p83-8.eps}\\
\hline
  9 &  8 &  8 &  1     &  0 & 3 &   0  &
\includegraphics[width=4cm]{pics/c-p83-9.eps}\\
\hline
  10&  9 &  9 &  1     &  0 & 2 &   0  &
\includegraphics[width=4.5cm]{pics/c-p83-10.eps}\\
\hline\hline
  11&  2 &  0 &  0     &  1 & 10&   0  &
\includegraphics[width=1cm]{pics/c-p83-11.eps}\\
\hline
  12&  3 &  1 &  1     &  1 & 9&   0  &
\includegraphics[width=2.5cm]{pics/c-p83-12.eps}\\
\hline
  13&  4 &  2 &  1     &  1 & 8&   0  &
\includegraphics[width=2.5cm]{pics/c-p83-13.eps}\\
\hline
  14&  5 &  3 &  1     &  1 & 7&   0  &
\includegraphics[width=3cm]{pics/c-p84-1.eps}\\
\hline
  15&  6 &  4 &  0     &  1 & 6&   0  &
\includegraphics[width=3.5cm]{pics/c-p84-2.eps}\\
\hline
  16&  6 &  4 &  1     &  1 & 6&   0  &
\includegraphics[width=4cm]{pics/c-p84-3.eps}\\
\hline
  17&  7 &  5 &  1     &  1 & 5&   0  &
\includegraphics[width=4.7cm]{pics/c-p84-4.eps}\\
\hline
  18&  8 &  6 &  1     &  1 & 4&   1  &
\includegraphics[width=5.2cm]{pics/c-p84-5.eps}\\
\hline
\end{tabular}

\end{table}

\begin{tabular}{|r||r|r|r|r|r|c|c|}
\hline
 $N$& $r$& $a$&$\delta$& $k$&$g$&$l(H)$& $\Gamma(P(\cM^{(2,4)}))$\\
\hline
  19&  9 &  7 &  1     &  1 & 3&   1  &
\includegraphics[width=5.2cm]{pics/c-p84-6.eps}\\
\hline
  20&  10&  8 &  1     &  1 & 2&   1  &
\includegraphics[width=5.8cm]{pics/c-p84-7.eps}\\
\hline\hline
  21&  6&  2 &  0     &  2 & 7&   0  &
\includegraphics[width=3.5cm]{pics/c-p84-8.eps}\\
\hline
  22&  7&  3 &  1     &  2 & 6&   0  &
\includegraphics[width=4.7cm]{pics/c-p84-9.eps}\\
\hline
  23&  8&  4 &  1     &  2 & 5&   0  &
\includegraphics[width=5.2cm]{pics/c-p84-10.eps}\\
\hline
  24&  9&  5 &  1     &  2 & 4&   0  &
\includegraphics[width=5.2cm]{pics/c-p84-11.eps}\\
\hline
  25& 10&  6 &  0     &  2 & 3&   1  &
\includegraphics[width=6cm]{pics/c-p84-12.eps}\\
\hline
  26& 10&  6 &  1     &  2 & 3&   1  &
\includegraphics[width=5.2cm]{pics/c-p84-13.eps}\\
\hline
  27& 11&  7 &  1     &  2 & 2&   1  &
\includegraphics[width=3cm]{pics/c-p85-1.eps}\\
\hline\hline
  28& 8&   2 &  1     &  3 & 6&   0  &
\includegraphics[width=4cm]{pics/c-p85-2.eps}\\
\hline
  29& 9&   3 &  1     &  3 & 5&   0  &
\includegraphics[width=5.2cm]{pics/c-p85-3.eps}\\
\hline
  30&10&   4 &  0     &  3 & 4&   0  &
\includegraphics[width=6cm]{pics/c-p85-4.eps}\\
\hline
  31&10&   4 &  1     &  3 & 4&   0  &
\includegraphics[width=4cm]{pics/c-p85-5.eps}\\
\hline
  32&11&   5 &  1     &  3 & 3&   0  &
\includegraphics[width=4cm]{pics/new-p85-6.eps}\\
\hline
  33&12&   6 &  1     &  3 & 2&   1  &
\includegraphics[width=3cm]{pics/c-p85-7.eps}\\
\hline\hline
\end{tabular}

\begin{tabular}{|r||r|r|r|r|r|c|c|}
\hline
 $N$& $r$& $a$&$\delta$& $k$&$g$&$l(H)$& $\Gamma(P(\cM^{(2,4)}))$\\
\hline
  34& 9&   1 &  1     &  4 & 6&   0  &
\includegraphics[width=4.7cm]{pics/c-p85-8.eps}\\
\hline
  35&10&   2 &  0     &  4 & 5&   0  &
\includegraphics[width=5.2cm]{pics/c-p86-1.eps}\\
\hline
  36&10&   2 &  1     &  4 & 5&   0  &
\includegraphics[width=5.8cm]{pics/c-p86-2.eps}\\
\hline
  37&11&   3 &  1     &  4 & 4&   0  &
\includegraphics[width=5.2cm]{pics/new-p86-3.eps}\\
\hline
  38&12&   4 &  1     &  4 & 3&   0  &
\includegraphics[width=5.2cm]{pics/c-p86-4.eps}\\
\hline
  39&13&   5 &  1     &  4 & 2&   0  &
\includegraphics[width=5.2cm]{pics/c-p86-5.eps}\\
\hline\hline
  40&10&   0 &  0     &  5 & 6&   0  &
\includegraphics[width=5.2cm]{pics/c-p86-6.eps}\\
\hline
  41&11&   1 &  1     &  5 & 5&   0  &
\includegraphics[width=6cm]{pics/c-p86-7.eps}\\
\hline
  42&12&   2 &  1     &  5 & 4&   0  &
\includegraphics[width=6.5cm]{pics/c-p86-8.eps}\\
\hline
  43&13&   3 &  1     &  5 & 3&   0  &
\includegraphics[width=5.2cm]{pics/c-p87-1.eps}\\
\hline
  44&14&   4 &  0     &  5 & 2&   0  &
\includegraphics[width=5.2cm]{pics/c-p87-2.eps}\\
\hline
\end{tabular}

\begin{tabular}{|r||r|r|r|r|r|c|c|}
\hline
 $N$& $r$& $a$&$\delta$& $k$&$g$&$l(H)$& $\Gamma(P(\cM^{(2,4)}))$\\
\hline
45&14&   4 &  1     &  5 & 2&   0  &
\includegraphics[width=5.2cm]{pics/c-p87-3.eps}\\
\hline\hline
46&14&   2 &  0     &  6 & 3&   0  &
\includegraphics[width=6.7cm]{pics/c-p87-4.eps}\\
\hline
47&15&   3 &  1     &  6 & 2&   0  &
\includegraphics[width=6cm]{pics/c-p87-5.eps}\\
\hline\hline
48&16&   2 &  1     &  7 & 2&   0  &
\includegraphics[width=6.5cm]{pics/c-p87-6.eps}\\
\hline\hline
49&17&   1 &  1     &  8 & 2&   0  &
\includegraphics[width=6.5cm]{pics/c-p87-7.eps}\\
\hline\hline
50&18&   0 &  0     &  9 & 2&   0  &
\includegraphics[width=6.5cm]{pics/c-p87-8.eps}\\
\hline
\end{tabular}

\medskip

\begin{proof} When $S$ is unimodular (i.e. $a=0$) or $r=a$ (then
  $S(1/2)$ is unimodular), i. e. for cases 1---11, 40, 50, these
  calculations were done by Vinberg \cite{2}. In all other cases they can be
  done using Vinberg's algorithm for calculation of the fundamental
  chamber of a hyperbolic reflection group. See \cite{2} and also \cite{3}.
  These technical calculations take too much space and will be
  presented in Appendix,
  Section \ref{fundchambersofS}.

  To describe elements of $P(X)_{+I}$ (i. e. double transparent
  vertices), we use the results of Section \ref{subsec2.5} and
  that their number $k$ is known by Section \ref{subsec2.2}.
\end{proof}

\begin{remark}\label{rem3.1.2} Using diagrams of Theorem
  \ref{thm3.1.1}, one can easily find the class in $S$ of the
  component $C_g$ of $X^\theta$ as an element $C_g\in S$
  such that $C_g\cdot x=0$, if $x$ corresponds to a black or a double
  transparent vertex, and $C_g\cdot x=2-s$ if $x$ corresponds to a
  simple transparent vertex which has $s$ edges to double transparent
  vertices.
\end{remark}

\subsection{Root invariants, and subsystems of roots in
$\Delta^{(4)}(\cM^{(2)})$ for elliptic case}
\label{subsec3.2} We use the notation and results of
Section \ref{subsubsec2.3.1}.  Let $\cM^{(2)}\supset \cM^{(2,4)}$ be the
fundamental chamber of $W^{(2)}(S)$ containing $\cM^{(2,4)}$.  Dynkin
diagram of $P^{(4)}(\cM^{(2,4)})$ (i. e. black vertices) consists of
components of types $\bA$, $\bD$ or $\bE$ (see Table 1). Thus, the
group $W^{(4)}(\cM^{(2)})$ generated by reflections in all elements of
$P^{(4)}(\cM^{(2,4)})$ is a finite Weyl group. It has to be finite because
$W^{(4)}(\cM^{(2)})(\cM^{(2,4)})=\cM^{(2)}$ has finite volume, and
$\cM^{(2,4)}$ is the fundamental chamber for the action of
$W^{(4)}(\cM^{(2)})$ in $\cM^{(2)}$.  Thus,
$$
\Delta^{(4)}(\cM^{(2)})=W^{(4)}(\cM^{(2)})P^{(4)}(\cM^{(2,4)})
$$
is a finite root system of the corresponding type with the
negative definite root sublattice
$$
R(2)=[P^{(4)}(\cM^{(2,4)})]\subset S.
$$

Let $(X,\theta)$ be a K3 surface with a non-symplectic involution,
and \linebreak $(S_X)_+=S$. Let $\Delta_+^{(4)}\subset
\Delta^{(4)}(S)$ be the subset defined by $(X,\theta)$ which is
invariant with respect to $W_+^{(2,4)}$ (we remind that it is
generated by reflections in $\Delta^{(2)}(S)$ and
$\Delta_+^{(4)}$). By Theorem \ref{thm2.3.4},
$\Delta_+^{(4)}=W^{(2)}(S)\Delta_+^{(4)}(\cM^{(2)})$ where
$\Delta_+^{(4)}(\cM^{(2)})=\Delta_+^{(4)}\cap
\Delta^{(4)}(\cM^{(2)})$ is a root subsystem in
$\Delta^{(4)}(\cM^{(2)})$. Let
\begin{equation}
K^+(2)=[\Delta_+^{(4)}(\cM^{(2)})]\subset R(2)\subset S
\label{3.2.1}
\end{equation}
be its negative definite root sublattice in $S$, and
\begin{equation}
Q=\frac{1}{2}K^+(2)/K^+(2), \ \ \ \xi^+:q_{K^+(2)}|Q\to q_S
\label{3.2.2}
\end{equation}
a homomorphism such that $\xi^+(x/2+K^+(2))=x/2+S$, $x\in K^+(2)$.
We obtain a pair $(K^+(2),\xi^+)$ which is similar to a root invariant,
and it is equivalent to the root invariant for elliptic type.

\begin{proposition}\label{prop3.2.1} Let $(X,\theta)$ be a K3 surface
with a non-symplectic involution of elliptic type, and $S=(S_X)_+$.

In this case, the root invariant $R(X,\theta)$ is equivalent
to the root subsystem
$\Delta_+^{(4)}(\cM^{(2)})\subset \Delta^{(4)}(\cM^{(2)})$,
considered up to the action of $O(S)$ (i. e. two root subsystems
$\Delta_+^{(4)}(\cM^{(2)})\subset \Delta^{(4)}(\cM^{(2)})$ and
$\Delta_+^{(4)}(\cM^{(2)})^\prime\subset \Delta^{(4)}(\cM^{(2)})$
are equivalent, if
$\Delta_+^{(4)}(\cM^{(2)})^\prime =\phi(\Delta_+^{(4)}(\cM^{(2)}))$
for some $\phi\in O(S)$):

The root invariant $R(X,\theta)\cong (K^+(2),\xi^+)$ is defined by
\eqref{3.2.1} and \eqref{3.2.2}.

The fundamental chamber $\cM(X)_+$ is defined by the
root subsystem $\Delta_+^{(4)}(\cM^{(2)})\subset \Delta^{(4)}(\cM^{(2)})$
(up to above equivalence), by Theorem \ref{thm2.3.4}. Moreover,
$P^{(4)}(\cM(X)_+)$ coincides with a basis of the root subsystem
$\Delta_+^{(4)}(\cM^{(2)})$.
\end{proposition}

\begin{proof} Let $E_i$, $i\in I$, be all non-singular rational curves
on $X$ such that $E_i\cdot\theta (E_i)=0$, i. e.
$$
\cl(E)+\cl(\theta(E))=\delta\in P^{(4)}(\cM(X)_+)=P^{(4)}(X)_+=
P(X)_{+III}.
$$
Since $E_i\cdot C_g=0$ and $C_g^2=2g-2>0$, the curves $E_i$, $i\in I$,
generate in $S_X$ a negative definite sublattice. Thus, their components
define a Dynkin diagram $\Gamma$ which consists of several connected
components $\bA_n$, $\bD_m$ or $\bE_k$. The involution $\theta$
acts on these diagrams and corresponding curves without fixed points.
Thus it necessarily changes connected components of $\Gamma$. Let
$\Gamma=\Gamma_1\bigsqcup \Gamma_2$ where $\theta(\Gamma_1)=\Gamma_2$, and
$I=I_1\bigsqcup I_2$ the corresponding subdivision of vertices of
$\Gamma$. Then
$$
\delta_i^+=\cl(E_i)+\cl(\theta(E_i)),\ \ i\in I_1,
$$
and
$$
\delta_i^-=\cl(E_i)-\cl(\theta(E_i)),\ \ i\in I_1
$$
give bases of root systems
$\Delta_+^{(4)}(\cM^2)$ and $\Delta_-^{(4)}=\Delta^{(4)}(K(2))$
respectively. The map
$$
\delta_i^-=\cl(E_i)-\cl(\theta(E_i))\mapsto \delta_i^+=
\cl(E_i)+\cl(\theta(E_i)),\ i\in I_1,
$$
defines an isomorphism $\Delta_-^{(4)}\cong\Delta_+^{(4)}(\cM^{(2)})$
of root systems, since it evidently preserves the intersection pairing.
The homomorphism $\xi$ of the root invariant
$R(X,\theta)=(K(2),\xi)$ of the pair $(X,\theta)$ then goes to
$(K^+(2),\xi^+)$.

In the opposite direction, the root invariant $R(X,\theta)$ defines
$\Delta_+^{(4)}$ and $\Delta_+^{(4)}(\cM^{(2)})=
\Delta^{(4)}(\cM^{(2)})\cap \Delta_+^{(4)}$.

The last statement follows from Section \ref{subsubsec2.3.1}.
\end{proof}

By Proposition \ref{prop3.2.1}, in the elliptic case instead of root
invariants one can consider root subsystems
$\Delta_+^{(4)}(\cM^{(2)})$ (in $\Delta^{(4)}(\cM^{(2)})$). We say
that a {\it root subsystem $\Delta_+^{(4)}(\cM^{(2)})$ ``is
  contained'' (respectively ``is primitively contained'') in a root
  subsystem $\Delta_+^{(4)}(\cM^{(2)})^\prime$, if
  $\phi(\Delta_+^{(4)}(\cM^{(2)}))\subset
  \Delta_+^{(4)}(\cM^{(2)})^\prime$ (respectively
  $[\phi(\Delta_+^{(4)}(\cM^{(2)}))]\subset
  [\Delta_+^{(4)}(\cM^{(2)})^\prime]$ is a primitive embedding of
  lattices) for some $\phi\in O(S)$.} By Corollary \ref{cor2.6.3}, we
obtain

\begin{proposition}\label{prop3.2.2} If a root subsystem
$\Delta_+^{(4)}(\cM^{(2)})$ in $\Delta^{(4)}(\cM^{(2)})$
corresponds to a K3 surface with non-symplectic
involution $(X,\theta)$, then any primitive root subsystem in
$\Delta_+^{(4)}(\cM^{(2)})$ corresponds to a K3 surface with non-symplectic
involution.

Thus, it is enough to describe {\bf extremal pairs $(X,\theta)$}
such that their root subsystems $\Delta_+^{(4)}(\cM^{(2)})$ in
$\Delta^{(4)}(\cM^{(2)})$ are not contained as primitive root subsystems
of strictly smaller rank in a root subsystem
$\Delta_+^{(4)}(\cM^{(2)})^\prime$ in $\Delta^{(4)}(\cM^{(2)})$
corresponding to another pair $(X^\prime,\theta^\prime)$.
\end{proposition}

\subsection{Classification of non-symplectic involutions $(X,\theta)$
of elliptic type of K3 surfaces}\label{subsec3.3}
We have

\begin{theorem}\label{thm3.3.1} Let $(X,\theta)$ and
  $(X^\prime,\theta^\prime)$ be two non-symplectic involutions of
  elliptic type of K3 surfaces.

Then the following three conditions are equivalent:

(i) Their main invariants $(r,a,\delta)$ (equivalently, $(k,g,\delta)$)
coincide, and their root invariants are isomorphic.

(ii) Their main invariants $(r,a,\delta)$
coincide, and the root subsystems $\Delta_+^{(4)}(\cM^{(2)})$ are
equivalent.

(iii) Dynkin diagrams $\Gamma(P(X)_+)$ and $\Gamma(P(X^\prime)_+)$
of their exceptional curves are isomorphic, and additionally the
genera $g$ are equal, if these diagrams are empty.  The diagram
$\Gamma(P(X)_+)$ is empty if and only if either $(r,a,\delta)=(1,1,1)$
(then $g=10$), or $(r,a,\delta)=(2,2,0)$ (then $g=9$)
and the root invariant is zero.
The corresponding DPN surfaces are  $\bP^2$ or $\bF_0$ respectively.
\end{theorem}

\begin{proof} By Sections \ref{subsec3.2} and \ref{subsec2.4}, the
conditions (i) and (ii) are equivalent, and they imply (iii).

Let us show that (iii) implies (i).

Assume that $r=\rk S\ge 3$.

First, let us show that $S$ is
generated by $\Delta^{(2)}(S)$, if $r=\rk S\ge 3$. If $r\ge a+2$, then
it is easy to see that either $S\cong U\oplus T$ or $S\cong U(2)\oplus T$
where $T$ is orthogonal sum of $A_1$, $D_{2m}$, $E_7$, $E_8$ (one can
get all possible invariants $(r, a,\delta)$ of $S$ taking these
orthogonal sums). We have $U=[c_1,c_2]$ where $c_1^2=c_2^2=0$ and
$c_1\cdot c_2=1$ (the same for $U(2)$, only $c_1\cdot c_2=2$). Then
$S$ is generated by elements with square $-2$ which are
$$
\Delta^{(2)}(T)\cup (c_1\oplus \Delta^{(2)}(T))\cup (c_2\oplus
\Delta^{(2)}(T)).
$$
If $r=a$ then $S\cong \langle 2 \rangle\oplus tA_1$. Let $h,e_1,\dots,e_t$
be the corresponding orthogonal basis of $S$ where $h^2=2$ and $e_i^2=-2$,
$i=1,\dots,t$. Then $S$ is generated by elements with square $(-2)$ which
are $e_1,\dots,e_t$ and $h-e_1-e_2$.

Now, let us show that $P(X)_+$ generates $S$.  Indeed, every element of
$\Delta^{(2)}(S)\cup \Delta_+^{(4)}$ can be obtained by composition of
reflections in elements of $P(X)_+$ from some element of $P(X)_+$.  It
follows, that it is an integral linear combination of elements of
$P(X)_+$. Since we can get in this way all elements of
$\Delta^{(2)}(S)$ and they generate $S$, it follows that $P(X)_+$
generates $S$.

It follows that the lattice $S$ with its elements $P(X)_+$ is defined
by the Dynkin diagram $\Gamma(P(X)_+)$. From $S$, we can find
invariants $(r,a,\delta)$ of $S$, and they define invariants
$(k,g,\delta)$.

Let $K^+(2)\subset S$ be a sublattice generated by $P^{(4)}(X)_+$ (i. e.
by the black vertices), and $\xi^+:Q=(1/2)K^+(2)/K^+(2)\to q_S$ the
homomorphism with $\xi^+(x/2+K^+(2))=x/2+S$. By Proposition \ref{prop3.2.1},
the pair $(K^+(2),\xi^+)$ coincides with the root invariant $R(X,\theta)$.

Now assume that $r=\rk S=1,\,2$ for the pair $(X,\theta)$. Then
$S\cong\langle 2\rangle$, $U(2)$, $U$ or
$\langle 2 \rangle \oplus \langle -2 \rangle$.

In the first two cases $\Delta^{(2)}(S)=\emptyset$ and then
$P^{(2)}(X)_+=\emptyset$. In the last two cases $\Delta^{(2)}(S)$ and
$P^{(2)}(X)_+$ are not empty.

Thus, only the first two cases give an empty diagram $P^{(2)}(X)_+$.
This distinguishes these two cases from all others.  In the case
$S=\langle 2 \rangle$, the invariant $g=10$, and the root invariant is
always zero because $S$ has no elements with square $-4$.  Thus, in
this case, the diagram $P(X)_+$ is always empty. This case gives
$Y=X/\{1,\theta\}\cong \bP^2$. In the case $S=U(2)$, the diagram
$P^{(2)}(X)_+$ is empty, but $P^{(4)}(X)_+=\emptyset$, if the root
invariant is zero, and $P^{(4)}(X)_+$ consists of one black vertex, if
the root invariant is not zero (see Table 1 for this case).  First
case gives $Y=\bF_0$. Second case gives $Y=\bF_2$.  In both these
cases $g=9$. Thus difference between two cases when the diagram is
empty ($\bP^2$ or $\bF_1$) is in genus: $g=10$ for the first case, and
$g=9$ for the second.

The difference of $S=U(2)$ with a non-empty diagram $\Gamma(P(X)_+)$
from all other cases is that this diagram consists of only one black
vertex.  All cases with $\rk S\ge 3$ must have at least 3 different
vertices to generate $S$. In cases $S=U$ and $S=\langle 2
\rangle\oplus \langle -2 \rangle$, the diagram $\Gamma(P(X)_+)$ also
consists of one vertex, but it is respectively double transparent and
simple transparent (see Table 1).  Moreover, this consideration also
shows the difference between cases $S=U$ and $S=\langle 2 \rangle\oplus
\langle -2 \rangle$ and with all other cases.
\end{proof}

Theorem \ref{thm3.3.1} shows that to classify pairs $(X,\theta)$ of
elliptic type, we can use any of the following invariants: either the
root invariant, or the root subsystem (together with the main
invariants $(k,g,\delta)$ or $(r,a,\delta)$), or the Dynkin diagram of
exceptional curves.

It seems that the most natural and geometric is the classification by
the Dynkin diagram.  Using this diagram, on the one hand, it easy to
calculate all other invariants. On the other hand, considering the
corresponding DPN surface, we get the Gram diagram of all exceptional
curves on it and all possibilities to get the DPN surface by blow-ups
from relatively minimal rational surfaces.

However, the statements (i) and (ii) of Theorem \ref{thm3.3.1} are
also very important since they give a simple way to find out if two
pairs $(X,\theta)$ and $(X^\prime, \theta^\prime)$ (equivalently, the
corresponding DPN surfaces) have isomorphic Dynkin diagrams of
exceptional curves. Moreover, the classification in terms of root
invariants and root subsystems is much more compact, since the full
Gram diagram of exceptional curves can be very large (e. g. recall the
classical non-singular del Pezzo surface corresponding to $\bE_8$).

We have the following

\begin{theorem}[Classification Theorem in the extremal case of
elliptic type] \label{thm3.3.2} A
  K3 surface with a non-symplectic involution $(X,\theta)$ of elliptic
  type is {\bf extremal,} if and only if the number of its exceptional
  curves with the square $(-4)$, i. e. $\#P^{(4)}(X)_+$, is equal to
  $\#P^{(4)}(\cM^{(2,4)})$ (see Theorem \ref{thm3.1.1}) where
  $\cM^{(2,4)}$ is a fundamental chamber of $W^{(2,4)}(S)$,
  $S=(S_X)_+$.  Equivalently, numbers of black vertices of Dynkin
  diagrams $\Gamma(P(X)_+)$ and $\Gamma(P(\cM^{(2,4)}))$ with the same
  invariants $(r,a,\delta)$ are equal.

Moreover, the diagram $\Gamma(P(X)_+)$ is isomorphic to
(i. e. coincides with) $\Gamma(P(\cM^{(2,4)}))$ (see Table 1)
in all cases of Theorem \ref{thm3.1.1} except cases 7, 8, 9, 10
and 20 of Table 1.
In the last five cases, all possible diagrams $\Gamma(P(X)_+)$ are given
in Table~2. All diagrams of Tables~1 and 2 correspond to some
extremal standard K3 pairs $(X,\theta)$.
\end{theorem}

\begin{proof} See Section \ref{subsec3.4} below.
\end{proof}

Now let us consider a description of non-extremal pairs $(X,\theta)$.
The worst way to describe them is using full diagrams
$\Gamma(P(X)_+)$, since the number of non-extremal pairs $(X,\theta)$
is very large and diagrams $\Gamma(P(X)_+)$ can be huge.  It is better
to describe them using Proposition \ref{prop3.2.2} and Theorem
\ref{thm3.3.1}. It is better to describe them by primitive root
subsystems $\Delta_+^{(4)}(\cM^{(2)})^\prime$ in the root subsystems
$\Delta_+^{(4)}(\cM^{(2)})\subset \Delta^{(4)}(\cM^{(2)})$ of extremal
pairs $(\widetilde{X}, \widetilde{\theta})$.

Let us choose $\cM^{(2)}$ in such a way that $\cM^{(2)}\supset
\cM(\widetilde{X})_+$. By Section \ref{subsubsec2.3.1}, then
$\Delta_+^{(4)}(\cM^{(2)})=\Delta^{(4)}([P^{(4)}(\widetilde{X})_+])$
is the subsystem of roots with the basis $P^{(4)}(\widetilde{X})_+$,
i. e.
$\Delta_+^{(4)}(\cM^{(2)})=\Delta^{(4)}([P^{(4)}(\widetilde{X})_+])$
is the set of all elements with the square $(-4)$ in the sublattice
$[P^{(4)}(\widetilde{X})_+]$ generated by $P^{(4)}(\widetilde{X})_+$
in $S=(S_{\widetilde{X}})_+$. Equivalently,
$\Delta^{(4)}([P^{(4)}(\widetilde{X})_+])=
W_+^{(4)}(\widetilde{X})(P^{(4)}(\widetilde{X})_+)$, where
$W_+^{(4)}(\widetilde{X})$ is the finite Weyl group generated by
reflections in all elements of $P^{(4)}(\widetilde{X})_+$.

Replacing a primitive root subsystem
$\Delta_+^{(4)}(\cM^{(2)})^\prime\subset \Delta^{(4)}([P^{(4)}
(\widetilde{X})_+])$
for a non-extremal pair $(X,\theta)$
by an equivalent root subsystem \newline  $\phi(\Delta_+^{(4)}
(\cM^{(2)})^\prime)$,
$\phi\in W_+^{(4)}(\widetilde{X})$, we can assume (by primitivity)
that a basis of $\Delta_+^{(4)}(\cM^{(2)})^\prime$ is a part of the basis
$P^{(4)}(\widetilde{X})_+$ of the root system $\Delta^{(4)}([P^{(4)}
(\widetilde{X})_+])$. Thus, we can assume that
the root subsystem $\Delta_+^{(4)}(\cM^{(2)})^\prime$
is defined by a subdiagram
$$
D\subset \Gamma(P^{(4)}(\widetilde{X})_+)
$$ where
$\Gamma(P^{(4)}(\widetilde{X})_+)$ is the subdiagram of the full diagram
$\Gamma(P(\widetilde{X})_+)$ generated by all its black vertices.
The $D$ is a basis of $\Delta_+^{(4)}(\cM^{(2)})^\prime$.

By Propositions \ref{prop2.3.1}, \ref{prop2.3.3} and Theorem \ref{thm2.3.4},
the subdiagram $D\subset \Gamma(P^{(4)}(\widetilde{X})_+)$ defines the
full Dynkin diagram $\Gamma(P(X)_+)$ of the pair $(X,\theta)$ with the
root subsystem $\Delta_+^{(4)}(\cM^{(2)})^\prime$: We have
\begin{equation}
P^{(2)}(X)_+=\{f\in W_+^{(4)}(\widetilde{X})(P^{(2)}
(\widetilde{X})_+)\ |\ f\cdot D\ge 0\}.
\label{3.3.1}
\end{equation}
The subdiagram of $\Gamma(P(X)_+)$ defined by all its black vertices
coincides with $D$. It is called {\it Du Val's part of $\Gamma(P(X)_+)$,}
and it is denoted by $\Duv{\Gamma(P(X)_+)}$. Thus,
$$
\Duv \Gamma(P(X)_+)=D\subset \Duv \Gamma(P(\widetilde{X})_+).
$$
Double transparent vertices of $\Gamma(P(X)_+)$ are identified with
double transparent vertices of $\Gamma(P(\widetilde{X})_+)$ (see
Section \ref{subsec2.5}), and single transparent vertices of $P(X)_+$
which are connected by two edges with double transparent vertices of
$\Gamma(P(X)_+)$ are identified with such vertices of
$\Gamma(P(\widetilde{X})_+)$. Indeed, they are orthogonal to the set
$P^{(4)}(\widetilde{X})_+$ which defines the reflection group
$W_+^{(4)}(\widetilde{X})$ as the group generated by reflections in all
elements of $P^{(4)}(\widetilde{X})_+$. Thus, the group
$W_+^{(4)}(\widetilde{X})$ acts identically on all these vertices, and all
of them satisfy \eqref{3.3.1}. All double transparent vertices and all single
transparent vertices connected by two edges with double transparent
vertices of $\Gamma(P(X)_+)$ define the {\it logarithmic part of
$\Gamma(P(X)_+)$,} and it is denoted by
$\Log \Gamma(P({X})_+)$. Thus, we have
$$
\Log \Gamma(P({X})_+)=\Log \Gamma(P(\widetilde{X})_+),
$$
logarithmic parts of $X$ and $\widetilde{X}$ are identified. Moreover,
the Du Val part  $\Duv \Gamma(P(X)_+)$ and the logarithmic part
$\Log \Gamma(P({X})_+)$ are {\it disjoint in $\Gamma(P(X)_+)$}
because they are orthogonal to each other. Thus, the logarithmic part of
$\Gamma(P(X)_+)$ is stable, it is the same for all pairs
$(X,\theta)$ with the same main invariants
$(r,a,\delta)$. On the Du Val part of $\Gamma(P(X)_+)$
we have only a restriction: it is
a subdiagram of Du Val part of one of extremal pairs
$(\widetilde{X},\widetilde{\theta})$ described in Theorems
\ref{thm3.1.1} and \ref{thm3.3.2}
(with the same main invariants $(r,a,\delta)$).

All vertices of $\Gamma(P({X})_+)$ which do not belong to
$\Duv \Gamma(P({X})_+)\cup \Log \Gamma(P({X})_+)$ define a subdiagram
$\Var \Gamma(P({X})_+)$ which is called the {\it varying part of
$\Gamma(P({X})_+)$.} By \eqref{3.3.1}, we have
$$
\Var P({X})_+=
\{f\in W_+^{(4)}(\widetilde{X})(\Var P(\widetilde{X})_+)\suchthat
f\cdot D\ge 0\}
$$
(we skip $\Gamma$ when we consider only vertices). It describes
$\Var \Gamma(P({X})_+)$ by the intersection pairing in $S$.

Of course, two  Dynkin subdiagrams
$D\subset \Gamma(P^{(4)}(\widetilde{X})_+)$ and $D^\prime \subset
\Gamma(P^{(4)}(\widetilde{X}^\prime)_+)$, with isomorphic Dynkin
diagrams $D\cong D^\prime$, of two extremal pairs
$(\widetilde{X},\widetilde{\theta})$ and
$(\widetilde{X}^\prime,\widetilde{\theta}^\prime)$ with the same
main invariants can give isomorphic Dynkin diagrams
$\Gamma(P({X})_+)$ and $\Gamma(P({X}^\prime)_+)$ for defining by
them K3 pairs $(X,\theta)$ and $(X^\prime, \theta^\prime)$.
To have that, it is necessary and sufficient that
root invariants $([D],\xi^+)$ and $([D^\prime],(\xi^\prime)^+)$
defined by them are isomorphic. We remind that they can be obtained
by restriction on $[D]$ and $[D^\prime]$ of the root invariants of pairs
$(\widetilde{X},\widetilde{\theta})$ and
$(\widetilde{X}^\prime,\widetilde{\theta}^\prime)$ respectively,
and they can be easily computed. We remind that to have
$([D],\xi^+)$ and $([D^\prime],(\xi^\prime)^+)$ isomorphic, there
must exist an isomorphism $\gamma:[D]\to [D^\prime]$ of the root
lattices and an automorphism $\overline{\phi}\in O(q_S)$ of the
discriminant quadratic form of the lattice $S$ which send $\xi^+$
for $(\xi^\prime)^+$. Section \ref{subsec2.4} gives the very simple
and effective method for that. Thus, we have a very simple and
effective method to find out when different subdiagrams $D$ above
give K3 pairs with isomorphic diagrams.

Note that we have used all equivalent conditions (i), (ii) and (iii)
of Theorem \ref{thm3.3.1} which shows their importance. Finally, we
get

\begin{theorem}[Classification Theorem in the non-extremal, i. e. arbitrary,
case of elliptic type] \label{thm3.3.3} Dynkin diagrams $\Gamma(P(X)_+)$ of
exceptional curves of non-extremal (i. e. arbitrary) non-symplectic
involutions $(X,\theta)$ of elliptic type of K3 surfaces are described
by arbitrary (without restrictions) Dynkin subdiagrams
$D\subset \Duv \Gamma(P(\widetilde{X})_+)$ of extremal pairs
$(\widetilde{X}, \widetilde{\theta})$ (see Theorem \ref{thm3.3.2})
with the same main invariants $(r,a,\delta)$ (equivalently $(k,g,\delta)$).
Moreover,
$$
\Duv \Gamma(P(X)_+)=D,\ \ \Log \Gamma(P(X)_+)=\Log \Gamma(P(\widetilde{X})_+),
$$
and they are disjoint to each other,
$$
\Var P({X})_+=
\{f\in W_+^{(4)}(\widetilde{X})(\Var P(\widetilde{X})_+)\suchthat
f\cdot D\ge 0\}
$$
where $W_+^{(4)}(\widetilde{X})$ is generated by reflections in all elements
of $\Duv P(\widetilde{X})_+$ $=P^{(4)}(\widetilde{X})_+$.

Dynkin subdiagrams $D\subset \Duv \Gamma(P(\widetilde{X})_+)$,
$D^\prime \subset \Duv \Gamma(P(\widetilde{X^\prime})_+)$
(with the same main invariants)
give K3 pairs $(X,\theta)$, $(X^\prime,\theta^\prime)$
with isomorphic Dynkin diagrams $\Gamma(P(X)_+)\cong \Gamma(P(X^\prime)_+)$,
if and only if the root invariants $([D],\xi^+)$,
$([D^\prime],(\xi^\prime)^+)$ defined by
$D\subset \Duv \Gamma(P(\widetilde{X})_+)$,
$D^\prime \subset \Duv \Gamma(P(\widetilde{X^\prime})_+)$ are isomorphic.

\end{theorem}

\noindent
\begin{table}\label{table2}
\caption{Diagrams $\Gamma(P(X)_+)$ of extremal K3 surfaces
$(X,\theta)$ of
elliptic type which are different from Table 1}

\index{Table 2}

\addtocontents{toc}{\contentsline {section}{\tocsection {}{T.2}{Table 2}}
{\pageref{table2}}}

(In (a) we repeat the corresponding case of Table 1)
\medskip

\begin{tabular}{|r||r|r|r|r|r|c|c|}
\hline
 $N$& $r$& $a$&$\delta$& $k$&$g$&$l(H)$& $\Gamma(P(X)_+)$\\
\hline
  7 &  6 &  6 &  1     &  0 &5  &     & \\
a &    &    &        &    &   &   0  &
\includegraphics[width=3cm]{pics/c-p94-1.eps}\\
b &    &    &        &    &   &   1  &
\includegraphics[width=4cm]{pics/c-p94-2.eps}\\
\hline
  8 &  7 &  7 &  1     &  0 &4  &     & \\
a &    &    &        &    &   &   0  &
\includegraphics[width=3.5cm]{pics/c-p94-3.eps}\\
b &    &    &        &    &   &   1  &
\includegraphics[width=4cm]{pics/c-06p94-4.eps}\\
c &    &    &        &    &   &   0  &
\includegraphics[width=3cm]{pics/c-p94-5.eps}\\
\hline
  9 &  8 &  8 &  1     &  0 &3  &     & \\
a &    &    &        &    &   &   0  &
\includegraphics[width=4cm]{pics/c-p95-1.eps}\\
b &    &    &        &    &   &   1  &
\includegraphics[width=4cm]{pics/c-p95-2.eps}\\
c &    &    &        &    &   &   0  &
\includegraphics[width=3cm]{pics/c-p95-3.eps}\\
d &    &    &        &    &   &   1  &
\includegraphics[width=3cm]{pics/c-p95-4.eps}\\
e &    &    &        &    &   &   1  &
\includegraphics[width=4.5cm]{pics/c-p95-5.eps}\\
\hline
\end{tabular}

\end{table}

\begin{tabular}{|r||r|r|r|r|r|c|c|}
\hline
 $N$& $r$& $a$&$\delta$& $k$&$g$&$l(H)$& $\Gamma(P(X)_+)$\\
\hline
  9 &  8 &  8 &  1     &  0 &3  &     & \\
f &    &    &        &    &   &   2  &
\includegraphics[width=3cm]{pics/c-06p95-6.eps}\\
\hline
 10 &  9 &  9 &  1     &  0 &2  &     & \\
a &    &    &        &    &   &   0  &
\includegraphics[width=4.5cm]{pics/c-p96-1.eps}\\
b &    &    &        &    &   &   0  &
\includegraphics[width=4.5cm]{pics/c-p96-2.eps}\\
c &    &    &        &    &   &   1  &
\includegraphics[width=3.5cm]{pics/c-p96-3.eps}\\
d &    &    &        &    &   &   1  &
\includegraphics[width=6cm]{pics/c-p96-4.eps}\\
e &    &    &        &    &   &   0  &
\includegraphics[width=3.5cm]{pics/c-p96-5.eps}\\
\hline
\end{tabular}

\begin{tabular}{|r||r|r|r|r|r|c|c|}
\hline
 $N$& $r$& $a$&$\delta$& $k$&$g$&$l(H)$& $\Gamma(P(X)_+)$\\
\hline
 10 &  9 &  9 &  1     &  0 &2  &     & \\
\hline
f &    &    &        &    &   &   1  &
\includegraphics[width=5cm]{pics/c-p97-1.eps}\\
g &    &    &        &    &   &   1  &
\includegraphics[width=6cm]{pics/c-p97-2.eps}\\
h &    &    &        &    &   &   0  &
\includegraphics[width=6cm]{pics/c-p97-3.eps}\\
i &    &    &        &    &   &   1  &
\includegraphics[width=5cm]{pics/c-p97-4.eps}\\
j &    &    &        &    &   &   2  &
\includegraphics[width=7cm]{pics/c-p97-5.eps}\\
k &    &    &        &    &   &   2  &
\includegraphics[width=5.5cm]{pics/c-p98-1.eps}\\
\hline
\end{tabular}

\begin{tabular}{|r||r|r|r|r|r|c|c|}
\hline
 $N$& $r$& $a$&$\delta$& $k$&$g$&$l(H)$& $\Gamma(P(X)_+)$\\
\hline
 10 &  9 &  9 &  1     &  0 &2  &     & \\
    &    &    &        &    &   &     &
\includegraphics[width=6cm]{pics/c-p98-2.eps}\\
l &    &    &        &    &   &   2  &
\includegraphics[width=4cm]{pics/new-p98-3.eps}\\
    &    &    &        &    &   &     &
\includegraphics[width=6cm]{pics/c-p98-4.eps}\\
m &    &    &        &    &   &   0  &
\includegraphics[width=4cm]{pics/c-p98-5.eps}\\
\hline\hline
\end{tabular}

\begin{tabular}{|r||r|r|r|r|r|c|c|}
\hline
 $N$& $r$& $a$&$\delta$& $k$&$g$&$l(H)$& $\Gamma(P(X)_+)$\\
\hline
 20 &  10&  8 &  1     &  1 &2  &     & \\
 a  &    &    &        &    &   &   1 &
\includegraphics[width=6cm]{pics/c-p99-1.eps}\\
 b  &    &    &        &    &   &   2 &
\includegraphics[width=4.5cm]{pics/c-p99-2.eps}\\
 c  &    &    &        &    &   &   1 &
\includegraphics[width=4.5cm]{pics/c-p99-3.eps}\\
 d  &    &    &        &    &   &   2 &
\includegraphics[width=4.5cm]{pics/c-p99-4.eps}\\
\hline\hline

\end{tabular}
\bigskip


\subsection{Proof of Theorem \ref{thm3.3.2}}\label{subsec3.4}
Let $(X,\theta)$ be a non-symplectic involution of elliptic type
of a K3 surface, with the main invariants $(r,a,\delta)$,
and $(X,\theta)$ is an extremal pair.

By Theorem \ref{thm3.3.1}, the $\Gamma(P(X)_+)$ is
defined by the root subsystem
$\Delta_+^{(4)}(\cM^{(2)})\subset \Delta^{(4)}(\cM^{(2)})$
corresponding to $(X,\theta)$ where
$\cM^{(2)}$ is a fundamental chamber of $W^{(2)}(S)$, and $S=(S_X)_+$ has
the invariants $(r,a,\delta)$. We can assume that
$\cM^{(2)}\supset \cM(X)_+\supset\cM^{(2,4)}$ where
$\cM^{(2,4)}$ is a fundamental chamber of $W^{(2,4)}(S)$ defined by
a choice of a basis $P^{(4)}(\cM^{(2,4)})$ of the root system
$\Delta^{(4)}(\cM^{(2)})$ (see Section \ref{subsubsec2.3.1}).

Let $\Gamma(P^{(4)}(\cM^{(2,4)}))$ be the Dynkin diagram of the root
system \linebreak
$\Delta^{(4)}(\cM^{(2)})$ and $W^{(4)}(\cM^{(2)})$ the Weyl
group of the root system \linebreak
$\Delta^{(4)}(\cM^{(2)})$. We use the
following description of a root subsystem
$\Delta_+^{(4)}(\cM^{(2)})\subset \Delta^{(4)}(\cM^{(2)})$.

\subsubsection{}\label{subsubsec3.4.1} Let $T\subset R$ be a root subsystem
of a root system $R$ and all components of $R$ have types $\bA$, $\bD$ or
$\bE$. We consider two particular cases of root subsystems.

Let $B$ be a basis of $R$. Let $T\subset R$ be a primitive root subsystem.
Then $T$ can be replaced
by an equivalent root subsystem $\phi(T)$, $\phi\in W(R)$, such that
a part of the basis $B$ gives a basis of $T$ (see \cite{1}).
Thus (up to equivalence defined by the Weyl group $W(R)$),
primitive root subsystems $T\subset R$ can be described by Dynkin
subdiagrams $\Gamma\subset \Gamma(B)$.

Now let $T\subset R$ be a root subsystem of a finite index.
Let $R_i$ be a component of $R$. Let $r_j$, $j\in J$, be a basis of $R_i$.
Let  $r_{\max}=\sum_{j\in J}{k_jr_j}$ be the maximal root of $R_i$
corresponding to this basis. Dynkin diagram of the set of roots
$$
\{r_j\suchthat j\in J\}\cup \{-r_{\max} \}
$$
is an extended Dynkin diagram of the Dynkin diagram
$\Gamma(\{r_j\suchthat j\in J\})$. Let us replace the component $R_i$
of the root system $R$ by the root subsystem $R_i^\prime\subset R_i$
having by its basis the set
$(\{r_j\suchthat j\in J\}\cup \{-r_{\max} \})-\{r_t\}$ where
$t\in J$ is some fixed element. We get a root subsystem
$R^\prime\subset R$ of finite index $k_t$. It can be shown
\cite{4} that iterations of this procedure give any root subsystem
of finite index of $R$ up to the action of $W(R)$.

Description of an arbitrary root subsystem $T\subset R$ can be reduced
to these two particular cases, moreover it can be done in two ways.

Firstly, any root subsystem $T\subset R$ is a subsystem of finite index
$T\subset T_{\pr}$ where $T_{\pr}\subset R$ is a primitive root subsystem
generated by $T$.

Secondly, any root subsystem $T\subset R$ can be considered as a
primitive root subsystem $T\subset R_1$ where $R_1\subset R$ is
root subsystem of finite index. One can take
$R_1$ generated by $T$ and by any $u=\rk R-\rk T$ roots
$r_1,\dots, r_u$ such that $\rk [T,r_1,\dots,r_u]=\rk R$.

\subsubsection{}\label{subsubsec3.4.2} Here we show that the root
subsystems $\Delta_+(\cM^{(2)})$ which coincide with the full
root systems $\Delta^{(4)}(\cM^{(2)})$ can be realized by K3
pairs $(X,\theta)$. Obviously, they are extremal. For them
$\cM(X)_+=\cM^{(2,4)}$, and the Dynkin diagrams
$\Gamma(P(X)_+)=\Gamma(P(\cM^{(2,4)}))$ coincide. All these diagrams
are described in Table 1 of Theorem \ref{thm3.1.1}. It is natural to
call such pairs $(X,\theta)$ as super-extremal. Thus, a non-symplectic
involution $(X,\theta)$ of elliptic type of K3 (equivalently,
the corresponding DPN pair $(Y,C)$ or DPN surface) is called
{\it super-extremal,} if for the corresponding root subsystem
$\Delta_+^{(4)}(\cM^{(2)})\subset \Delta^{(4)}(\cM^{(2)})$ we have
$\Delta_+^{(4)}(\cM^{(2)})=\Delta^{(4)}(\cM^{(2)})$ (equivalently,
$\Delta_+^{(4)}=\Delta^{(4)}(S)$). We have

\begin{proposition}\label{prop3.4.1} For any possible elliptic triplet
of main invariants $(r,a,\delta)$ there exists
a super-extremal, i. e.
$$
\Gamma(P(X)_+)=\Gamma(P(\cM^{(2,4)})),
$$
and standard (see Section \ref{subsec2.6}) K3 pair $(X,\theta)$.

See the description of their graphs
$\Gamma(P(X)_+)=\Gamma(P(\cM^{(2,4)}))$ in Table 1 of Theorem \ref{thm3.1.1}.
\end{proposition}

\begin{proof} Let us consider an elliptic triplet of main
invariants $(r,a,\delta)$
and the corresponding Dynkin diagram $\Gamma(P(\cM^{(2,4)}))$
which is described in Theorem \ref{thm3.1.1}. Denote
$K^+(2)=[P^{(4)}(\cM^{(2,4)})]$, i. e. it is the sublattice generated
by all black vertices of $\Gamma(P(\cM^{(2,4)}))$.
Consider the corresponding root invariant
$(K^+(2),\xi^+)$, see \eqref{3.2.1} and \eqref{3.2.2}. Consider
$H=\Ker \xi^+$. By Propositions \ref{prop3.2.1} and \ref{prop2.6.1},
there exists a super-extremal standard pair $(X,\theta)$, if the
inequalities
$$
r+\rk K^++l({\gA}_{(K^+)_p})<22\ \text{for all prime}\ p>2,
$$
$$
r+a+2l(H)<22
$$
are valid together with Conditions 1 and 2 from Section \ref{subsec2.6}.

By trivial inspection of all cases in Table 1, we can see that first
inequality is valid. To prove second inequality, it is enough to show that
$l(H)\le 1$ since $r+a\le 18$ in elliptic case. The inequality $l(H)\le 1$
can be proved by direct calculation of $l(H)$ in all cases of Table 1
of Theorem 3.1.1. These calculations are simplified by the
general statement.

\begin{lemma}\label{lemma3.4.2} In elliptic super-extremal case,
$$
l(H)=\#P^{(4)}(\cM^{(2,4)})-l({\gA}_S^{(1)})
$$
where ${\gA}_S^{(1)}\subset {\gA}_S$ is the subgroup generated by all elements
$x\in\gA_S$ such that $q_S(x)=1\mod 2$. Moreover, we have:

\medskip

\noindent
If $\delta=0$ then $l({\gA}_S^{(1)})=a$ except ($a=2$ and
$\sign S=2-r\equiv 0\mod 8)$. In the last case $l({\gA}_S^{(1)})=a-1$.

\medskip

\noindent
If $\delta=1$, then $l({\gA}_S^{(1)})=a-1$ except cases
($a=2$ and $\sign S\equiv 0\mod 8$), ($a=3$ and $\sign S\equiv\pm 1\mod 8$),
and  ($a=4$ and $\sign S\equiv 0\mod 8$). In these cases
$l({\gA}_S^{(1)})=a-2$.
\end{lemma}

\begin{proof} We know (see Section \ref{subsubsec2.3.1}) that
\newline
$\Delta^{(4)}(S)=W^{(2)}(S)(\Delta^{(4)}(\cM^{(2,4)}))$. The group
$W^{(2)}(S)$ acts identically on ${\gA}_S$. Therefore,
$$
\IM \xi^+=[\{\xi^+(f/2+K^+(2))\suchthat f\in \Delta^{(4)}(\cM^{(2,4)})\}]=
$$
$$
[\{f/2+S\suchthat f\in \Delta^{(4)}(S)\}={\gA}_S^{(1)}.
$$
In the last equality, we use Lemma \ref{lemma2.4.2}. For
$Q=(K^+(2)/2)/K^+(2)$, we have $l(Q)=\rk K^+=\#P^{(4)}(\cM^{(2,4)})$.
Thus, $l(H)=l(Q)-l({\gA}_S^{(1)})=\#P^{(4)}(\cM^{(2,4)})-l({\gA}_S^{(1)})$.

The rest of statements of Lemma can be proved by direct calculations using
a decomposition of a 2-elementary non-degenerate finite quadratic form as
sum of elementary ones:
$q_{\pm 1}^{(2)}(2)$, $u_+^{(2)}(2)$ and $v_+^{(2)}(2)$ (in notation of
\cite{9}).
\end{proof}

One can easily check Condition 2 of Section \ref{subsec2.6}.

To check Condition 1 of Section \ref{subsec2.6}, note that if the
lattice $K^+_H(2)$ has elements with the square $(-2)$, then the
sublattice $[P^{(4)}(\cM^{(2,4)})]_{\pr}$ of
$S$ also has elements with the square $(-2)$.
Let us show that this is not the case.

Let us consider the subspace
$$
\gamma=\bigcap_{f\in P^{(4)}(\cM^{(2,4)})}{\cH_f}
$$
of $\cL(S)$ which is orthogonal to $[P^{(4)}(\cM^{(2,4)})]$ (equivalently,
we consider the corresponding face $\gamma\cap \cM^{(2,4)}$ of
$\cM^{(2,4)}$). If the
sublattice $[P^{(4)}(\cM^{(2,4)})]_{\pr}\subset S$ has elements with
square $(-2)$, then some hyperplanes $\cH_e$, $e\in \Delta^{(2)}(S)$,
also contain $\gamma$ and give reflections from $W^{(2,4)}(S)$.
On the other hand (e. g. see \cite{3}), all hyperplanes of reflections
from $W^{(2,4)}(S)$ containing $\gamma$ must be obtained from the hyperplanes
$\cH_f$, $f\in P^{(4)}(\cM^{(2,4)})$, by the group generated by reflections
in $P^{(4)}(\cM^{(2,4)})$. All these hyperplanes are then also orthogonal
to elements with square $(-4)$ from $S$. They cannot be orthogonal
to elements with square $(-2)$ from $S$ too.

This finishes the proof of Proposition \ref{prop3.4.1}.

\subsubsection{}\label{subsubsect3.4.3} Let us prove Theorem \ref{thm3.3.2}
in all cases except 7 --- 10 and 20 of Table 1. These cases
(i. e. different from 7 --- 10 and 20 of Table 1) are characterized by
the property that Dynkin diagram $\Gamma(P^{(4)}(\cM^{(2,4)}))$
consists of components of type $\bA$ only. By Section \ref{subsubsec3.4.1},
any root subsystem
$\Delta_+^{(4)}(\cM^{(2)})\subset \Delta^{(4)}(\cM^{(2)})$ is then primitive.
In particular, any root subsystem
$\Delta_+^{(4)}(\cM^{(2)})\subset \Delta^{(4)}(\cM^{(2)})$ of finite index is
$\Delta_+^{(4)}(\cM^{(2)})=\Delta^{(4)}(\cM^{(2)})$. By Proposition
\ref{prop3.4.1}, we then obtain

\begin{proposition}\label{prop3.4.3} For any elliptic triplet $(r,a,\delta)$
of main invariants which is different from (6,6,1), (7,7,1), (8,8,1),
(9,9,1) and (10,8,1), any extremal K3 pair $(X,\theta)$ is super-extremal,
i. e. $\Gamma(P(X)_+)=\Gamma(P(\cM^{(2,4)}))$ (see their description in
Table 1 of Theorem \ref{thm3.1.1}).
\end{proposition}

Above, we have proved that the primitive sublattice
$[P^{(4)}(\cM^{(2,4)})]_{\pr}$ in $S$ generated by $P^{(4)}(\cM^{(2,4)})$
has no elements with square $-2$. The lattice $[P^{(4)}(\cM^{(2,4)})]$
coincides with the root lattice $[\Delta^{(4)}(\cM^{(2)})]$.
Thus, its primitive sublattice $[\Delta^{(4)}(\cM^{(2)})]_{\pr}$ in
$S$ also has no elements with square $-2$. This fact is very important. Using
\eqref{3.2.1} and \eqref{3.2.2}, we can define the {\it root invariant}
$(K^+(2),\xi^+)$ for any root subsystem
$\Delta^{(4)}_+(\cM^{(2)})\subset \Delta^{(4)}(\cM^{(2)})$. Like
for root subsystems of K3 pairs $(X,\theta)$, we then have

\begin{lemma}\label{lemma3.4.2a} Two root subsystems
$\Delta^{(4)}_+(\cM^{(2)})\subset \Delta^{(4)}(\cM^{(2)})$ and
$\Delta^{(4)}_+(\cM^{(2)})^\prime \subset \Delta^{(4)}(\cM^{(2)})$
are $O(S)$ equivalent, if and only if their root invariants are
isomorphic.
\end{lemma}

\begin{proof} Assume that the root invariants are isomorphic. Since
$\pm 1$ and $W^{(2)}(S)$ act identically on the discriminant form
$q_S$, there exists an automorphism $\phi\in O(S)$
such that $\phi(\Delta^{(4)}(\cM^{(2)}))=\Delta^{(4)}(\cM^{(2)})$
and, identifying by $\phi$ the root subsystem
$\Delta^{(4)}_+(\cM^{(2)})\subset \Delta^{(4)}(\cM^{(2)})$
with $\phi\left(\Delta^{(4)}_+(\cM^{(2)})\right)\subset
\Delta^{(4)}(\cM^{(2)})$,
we have the following.
There exists an isomorphism
$\alpha:\Delta^{(4)}_+(\cM^{(2)})\cong \Delta^{(4)}_+(\cM^{(2)})^\prime$
of root systems such that
$\alpha(f)/2+S=f/2+S$ for any $f\in \Delta^{(4)}_+(\cM^{(2)})$.
Equivalently, $(\alpha(f)+f)/2\in S$.

Assume that $\alpha(f)\not=\pm f$. Then,
since $\alpha(f)$ and $f$ are two elements of a finite root system
$\Delta^{(2)}(\cM^{(2)})$ which is a sum of $\bA_n$, $\bD_m$, $\bE_k$,
it follows that either $\alpha(f)\cdot f=\pm 2$, or $\alpha(f)\cdot f=0$.
First case gives $f\cdot (\alpha(f)+f)/2\equiv 1\mod 2$ which is impossible
because $f\in S$ is a root. Second case gives that $\beta=(\alpha(f)+f)/2$
has $\beta^2=-2$ which is impossible because
$[\Delta^{(4)}(\cM^{(2)})]_{\pr}$ has no elements with square $-2$.
Thus, $\alpha(f)=\pm f$. It follows that
$\Delta^{(4)}_+(\cM^{(2)})=\Delta^{(4)}_+(\cM^{(2)})^\prime$ are
identically the same root subsystems of $\Delta^{(4)}_+(\cM^{(2)})$.
\end{proof}

\subsubsection{}\label{subsubsec3.4.4} Now let us consider cases
7---10 and 20 of Table 1. In these cases, the root system
$R=\Delta^{(4)}(\cM^{(2)})$ is
$\bD_5$ in the case 7, $\bE_6$ in the case 8,
$\bE_7$ in the case 9, $\bE_8$ in the case 10, and $\bD_8$ in the case 20.

We have

\begin{lemma}\label{lemma3.4.4} If $R$ is a root system of one of types
$\bD_5$, $\bE_6$, $\bE_7$, $\bE_8$ or $\bD_8$, then its root subsystem
$T\subset R$ of finite index is determined by the isomorphism type of
the root system $T$ itself, up to the action of $W(R)$.
Moreover, the type of $T$ can be the following and only the following
which is given in Table of Lemma \ref{lemma3.4.4} below (we identify the
type with the isomorphism class of the corresponding root lattice).

Moreover, in the corresponding cases labelled by $N=7, 8, 9, 10$ and $20$ of
Table 1 the above
statement is equivalent to the fact that the root invariant of
the corresponding root subsystem
$T\subset R=\Delta^{(4)}(\cM^{(2)})$ of finite index is defined by its type.
The root invariants $(T,\xi^+)$ for them are given below
by showing the kernel $H=\Ker \xi^+$
and the invariants $\alpha$ and $\overline{a}$, if $\alpha=0$
(we use Proposition \ref{prop2.4.4}).
\end{lemma}

\newpage
\noindent
{\bf Table of Lemma \ref{lemma3.4.4}.}

\medskip

\begin{tabular}{|r||c|l|}
\hline
 $N$& $R$& \hskip3cm $T$\\
\hline
  7 &  $D_5$ & a) $D_5$, b) $A_3\oplus 2A_1$\\
\hline
  8 &  $E_6$ & a) $E_6$, b) $A_5\oplus A_1$, c) $3A_2$\\
\hline
    &        & a) $E_7$, b) $A_7$, c) $A_5\oplus A_2$, d) $2A_3\oplus A_1$,
               e) $D_6\oplus A_1$,\\
  9 &  $E_7$ & f) $D_4\oplus 3A_1$, g) $7A_1$\\
\hline
    &        & a) $E_8$, b) $A_8$,   c) $A_7\oplus A_1$,
               d) $A_5\oplus A_2\oplus A_1$, e) $2A_4$,\\
    &        & f) $D_8$,  g) $D_5\oplus A_3$, h) $E_6\oplus A_2$,
               i) $E_7\oplus A_1$, j) $D_6\oplus 2A_1$,\\
 10 &  $E_8$ & k) $2D_4$, l) $2A_3\oplus 2A_1$, m) $4A_2$,
               n) $D_4\oplus 4A_1$, o) $8A_1$\\
\hline\hline
    &        & a) $D_8$, b) $D_6\oplus 2A_1$, c) $D_5\oplus A_3$, d) $2D_4$,
               e) $2A_3\oplus 2A_1$,\\
 20 &  $D_8$ & f) $D_4\oplus 4A_1$, g) $8A_1$\\
\hline

\end{tabular}

\medskip

\noindent
{\bf The root invariants of $T\subset R$:}

\medskip

7a, $D_5\subset D_5$: with the basis in $T$

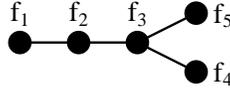
\medskip

\centerline{\includegraphics[width=3.5cm]{pics/c-d5.eps}}

\medskip

\noindent
$H=0\mod T$, $\overline{a}=(f_4+f_5)/2\mod H$ (since $\overline{a}$ is
defined, the invariant $\alpha=0$).

7b, $A_3\oplus A_1\subset D_5$: with the basis (in $T$)


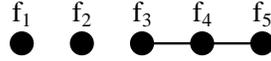
\medskip

\centerline{\includegraphics[width=4cm]{pics/c-p115-1.eps}}

\medskip

\noindent
$H=[(f_1+f_2+f_3+f_5)/2]\mod T$, $\overline{a}=(f_3+f_5)/2\mod H$.

8a, $E_6\subset E_6$: Then $H=0\mod T$ and $\alpha=1$ (it follows that
$\alpha=1$ and $\overline{a}$ is not defined for all cases 8a---c below).

8b, $A_1\oplus A_5\subset E_6$:  with the basis


\medskip

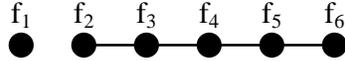
\centerline{\includegraphics[width=5cm]{pics/c-p115-2.eps}}

\medskip

\noindent
$H=[(f_1+f_2+f_4+f_6)/2]\mod T$ and $\alpha=1$.

8c, $3A_2\subset E_6$: Then $H=0\mod T$ and $\alpha=1$.

9a, $E_7\subset E_7$: with the basis

\medskip

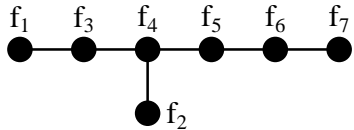
\centerline{\includegraphics[width=5.5cm]{pics/c-e7.eps}}

\medskip

\noindent
$H=0\mod T$ and $\overline{a}=(f_2+f_5+f_7)/2\mod T$.

9b, $A_7\subset E_7$: with the basis


\medskip

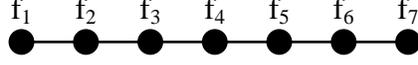
\centerline{\includegraphics[width=6cm]{pics/c-p115-3.eps}}

\medskip

\noindent
$H=[(f_1+f_3+f_5+f_7)/2]\mod T$ and $\alpha=1$.

9c, $A_5\oplus A_2\subset E_7$:  with the basis


\medskip

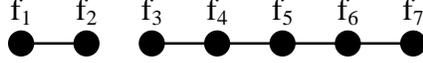
\centerline{\includegraphics[width=6.5cm]{pics/c-a2a5.eps}}

\medskip

\noindent
$H=0\mod T$ and $\overline{a}=(f_3+f_5+f_7)/2\mod H$.

9d, $2A_3\oplus A_1\subset E_7$:  with the basis


\medskip

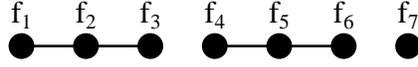
\centerline{\includegraphics[width=6cm]{pics/c-p116-1.eps}}

\medskip

\noindent
$H=[(f_1+f_3+f_4+f_6)/2]\mod T$ and $\overline{a}=(f_1+f_3+f_7)/2\mod H$.

9e, $D_6\oplus A_1\subset E_7$: with the basis


\medskip

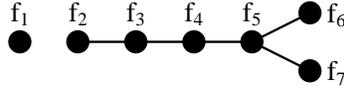
\centerline{\includegraphics[width=5cm]{pics/c-p116-2.eps}}

\medskip

\noindent
$H=[(f_1+f_2+f_4+f_6)/2]\mod T$ and $\overline{a}=(f_1+f_6+f_7)/2\mod H$.

9f, $D_4\oplus 3A_1\subset E_7$: with the basis


\medskip

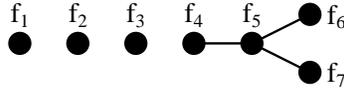
\centerline{\includegraphics[width=5cm]{pics/c-p116-3.eps}}

\medskip

\noindent
the $H=[(f_1+f_2+f_4+f_6)/2, (f_2+f_3+f_6+f_7)/2]\mod T$ and
$\overline{a}=(f_1+f_2+f_3)/2\mod H$.

9g, $7A_1\subset E_7$: with the basis $f_v$, $v\in \bP^2(F_2)$ where
$\bP^2(F_2)$ is the projective plane over the field $F_2$ with two elements,
the group $H$ is generated by $\left(\sum_{v\in \bP^2(F_2)-l}{f_v}\right)/2$
where $l$
is any line in $\bP^2(F_2)$. The element
$\overline{a}=\left(\sum_{v\in l}{f_v}\right)/2$ where $l$ is any line in
$\bP^2(F_2)$.

10a, $E_8\subset E_8$: Then $H=0\mod T$ and $\alpha=1$
(it follows that $\alpha=1$
and the element $\overline{a}$ is not defined for all cases 10a---o).

10b, $A_8\subset E_8$: Then $H=0\mod T$ and $\alpha=1$.

10c, $A_7\oplus A_1\subset E_8$:  with the basis


\medskip

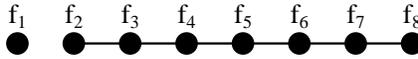
\centerline{\includegraphics[width=6cm]{pics/c-p116-4.eps}}

\medskip

\noindent
$H=[(f_2+f_4+f_6+f_8)/2]\mod T$ and $\alpha=1$.

10d, $A_5\oplus A_2\oplus A_1\subset E_8$: with the basis


\medskip

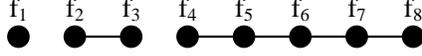
\centerline{\includegraphics[width=6cm]{pics/c-p116-5.eps}}

\medskip

\noindent
$H=[(f_1+f_4+f_6+f_8)/2]\mod T$ and $\alpha=1$.

10e, $2A_4\subset E_8$: Then $H=0\mod T$ and $\alpha=1$.

10f, $D_8\subset E_8$: with the basis


\medskip

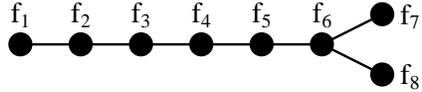
\centerline{\includegraphics[width=6cm]{pics/c-p116-6.eps}}

\medskip

\noindent
$H=[(f_1+f_3+f_5+f_7)/2]\mod T$ and $\alpha=1$.

10g, $D_5\oplus A_3\subset E_8$:  with the basis


\medskip

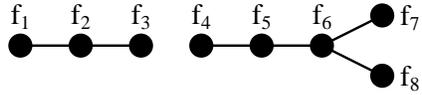
\centerline{\includegraphics[width=6cm]{pics/c-p117-1.eps}}

\medskip

\noindent
$H=[(f_1+f_3+f_7+f_8)/2]\mod T$ and $\alpha=1$.

10h, $E_6\oplus A_2\subset E_8$:  Then $H=0$ and $\alpha=1$.

10i, $E_7\oplus A_1\subset E_8$: with the basis


\medskip

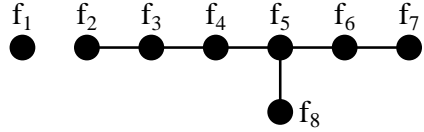
\centerline{\includegraphics[width=6cm]{pics/c-p117-2.eps}}

\medskip

\noindent
$H=[(f_1+f_2+f_4+f_8)/2]\mod T$ and $\alpha=1$.

10j, $D_6\oplus 2A_1\subset E_8$: with the basis


\medskip

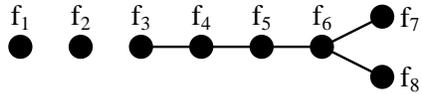
\centerline{\includegraphics[width=6cm]{pics/c-p117-3.eps}}

\medskip

\noindent
$H=[(f_1+f_3+f_5+f_7)/2,\ (f_2+f_3+f_5+f_8)/2]\mod T$ and $\alpha=1$.

10k, $2D_4\subset E_8$: with the basis


\medskip

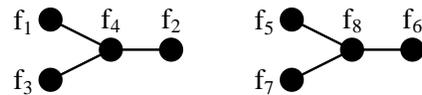
\centerline{\includegraphics[width=6cm]{pics/c-p117-4.eps}}

\medskip

\noindent
$H=[(f_1+f_2+f_5+f_6)/2,\ (f_2+f_3+f_6+f_7)/2]\mod T$ and $\alpha=1$.

10l, $2A_3\oplus 2A_1\subset E_8$: with the basis


\medskip

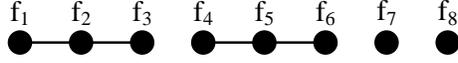
\centerline{\includegraphics[width=6.5cm]{pics/c-p117-5.eps}}

\medskip

\noindent
$H=[(f_1+f_3+f_7+f_8)/2,\ (f_4+f_6+f_7+f_8)/2]\mod T$ and $\alpha=1$.

10m, $4A_2\subset E_8$: Then $H=0\mod T$ and $\alpha=1$.

10n, $D_4\oplus 4A_1\subset E_8$:

\medskip

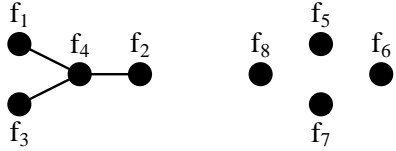
\centerline{\includegraphics[width=6cm]{pics/c-d44a1.eps}}

\medskip

\noindent
$H=[(f_1+f_2+f_5+f_6)/2,\ (f_2+f_3+f_6+f_7)/2,\ (f_5+f_6+f_7+f_8)/2]\mod T$ and
$\alpha=1$.

10o, $8A_1\subset E_8$: with the basis $f_v$, $v\in V$ and $V$ has
the structure
of 3-dimensional affine space over $F_2$, the group $H$ is generated by
$\left(\sum_{v\in \pi}{f_v}\right)/2$ where $\pi\subset V$ is any
2-dimensional
affine subspace in $V$. The invariant $\alpha=1$.

20a, $D_8\subset D_8$: with the basis $f_1,\dots,f_8$ shown below


\medskip

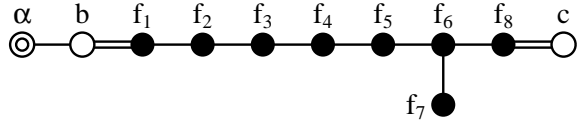
\centerline{\includegraphics[width=8cm]{pics/c-p112.eps}}

\medskip

\noindent
$H=[(f_1+f_3+f_5+f_7)/2]\mod T$, $\overline{a}=(f_7+f_8)/2\mod H$.

20b, $D_6\oplus 2A_1\subset D_8$: with the basis


\medskip

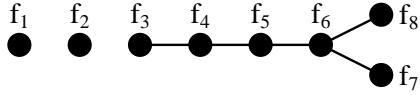
\centerline{\includegraphics[width=6cm]{pics/c-p118-2.eps}}

\medskip

\noindent $H=[(f_1+f_3+f_5+f_7)/2,\ (f_2+f_3+f_5+f_8)/2]\mod T$
and $\overline{a}=(f_7+f_8)/2\mod H$.

20c, $D_5\oplus A_3$:  with the basis


\medskip

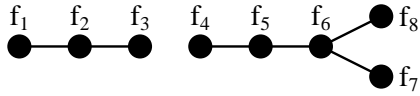
\centerline{\includegraphics[width=6cm]{pics/c-p118-3.eps}}

\medskip

\noindent
$H=[(f_1+f_3+f_7+f_8)/2]\mod T$ and $\overline{a}=(f_7+f_8)/2\mod H$.

20d, $2D_4\subset D_8$:  with the basis


\medskip

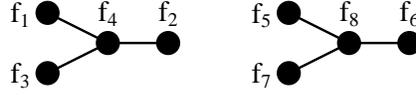
\centerline{\includegraphics[width=6cm]{pics/c-p118-4.eps}}

\medskip

\noindent
$H=[(f_1+f_2+f_5+f_6)/2,\ (f_2+f_3+f_6+f_7)/2]\mod T$ and
$\overline{a}=(f_6+f_7)/2\mod H$.

20e, $2A_3\oplus 2A_1\subset D_8$: with the basis


\medskip

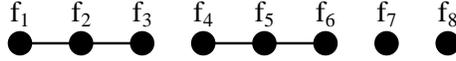
\centerline{\includegraphics[width=6.5cm]{pics/c-p117-5.eps}}

\medskip

\noindent
$H=[(f_1+f_3+f_7+f_8)/2,\ (f_4+f_6+f_7+f_8)/2]\mod T$ and
$\overline{a}=(f_7+f_8)/2\mod H$.

20f, $D_4\oplus 4A_1\subset D_8$:  with the basis


\medskip

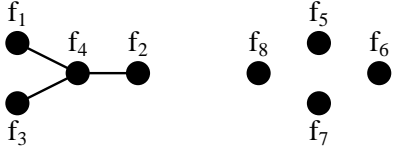
\centerline{\includegraphics[width=6cm]{pics/c-d44a1.eps}}

\medskip

\noindent
$H=[(f_1+f_2+f_5+f_6)/2,\ (f_2+f_3+f_6+f_7)/2,\ (f_5+f_6+f_7+f_8)/2]\mod T$
and $\overline{a}=(f_7+f_8)/2\mod H$.

20g, $8A_1\subset D_8$: with the basis $f_v$, $v\in V$ and $V$
has the structure
of 3-dimensional affine space over $F_2$, the group $H$ is generated by
$\left(\sum_{v\in \pi}{f_v}\right)/2$ where $\pi\subset V$ is any 2-dimensional
affine subspace in $V$. The element $\overline{a}=(f_{v_1}+f_{v_2})/2\mod H$
where $v_1v_2$ is a fixed non-zero vector in $V$. This structure can be seen
in Figure \ref{fig6} below.

\begin{proof} Let us consider cases $N=7$, $8$, $9$, $10$ and $20$
of the main invariants
$S$ in Table 1. By Lemma \ref{lemma2.4.1},
the canonical homomorphism $O(S)\to O(q_S)$ is epimorphic.
Since $\pm 1$ acts identically on the 2-elementary form $q_S$,
it follows that $O^\prime(S)\to O(q_S)$ is
epimorphic. The group $O^\prime(S)$ is the semi-direct product of
$W^{(2,4)}(S)$ and the automorphism group of the diagram
$\Gamma(P(\cM^{(2,4)}))$. The last group is
trivial in all these cases. Thus $W^{(2,4)}(S)\to O(q_S)$ is epimorphic. The
group $W^{(2,4)}(S)$ is the semi-direct product of $W^{(2)}(S)$ and the
symmetry group $W^{(4)}(\cM^{(2)})$ of the fundamental chamber $\cM^{(2)}$.
The group $W^{(2)}(S)$ acts identically on $O(q_S)$. It follows that
the corresponding homomorphism $W^{(4)}(\cM^{(4)})\to O(q_S)$ is epimorphic.
Here $W^{(4)}(\cM^{(2)})$ is exactly the Weyl group of the root system
$R$ defined by black vertices of the diagram  $\Gamma(P(\cM^{(2,4)}))$.

$N=7$: Then $q_S\cong q_1^{(2)}(2)\oplus q_{-1}^{(2)}(2)\oplus
u_+^{(2)}(2)\oplus v_+^{(2)}(2)$ (we use notation of \cite{9}), and $R=D_5$.
By direct calculation (using Lemma \ref{lemma2.4.3}), we get
$\#O(q_S)=5\cdot 3\cdot 2^7$. It is known \cite{1}, that
$\#W(D_5)=5\cdot 3 \cdot 2^7$. Thus we get the canonical isomorphism
$W(D_5)\cong O(q_S)$. By Lemma \ref{lemma3.4.2a}, it follows that any
two root subsystems $T_1\subset D_5$ and $T_2\subset D_5$ are conjugate by
$W(D_5)$, if and only if their root invariants $(T_1,\xi^+_1)$ and
$(T_2,\xi^+_2)$ are isomorphic.

In all other cases considerations are the same.

$N=8$: Then $q_S\cong q_{-1}^{(2)}(2)\oplus v^{(2)}_+(2)\oplus 2u_+^{(2)}(2)$
and
$R=E_6$. We have $\#O(q_S)=\#W(E_6)=5\cdot 3^4\cdot 2^7$. It follows,
$W(E_6)\cong O(q_S)$.

$N=9$: Then $q_S\cong 2q_1^{(2)}(2)\oplus 3u_+^{(2)}(2)$ and $R=E_7$.
We have $\#O(q_S)=\#W(E_7)=7\cdot 5\cdot 3^4\cdot 2^{10}$. It follows,
$W(E_7)\cong O(q_S)$.

$N=10$: Then $q_S\cong q_1^{(2)}(2)\oplus 4u_+^{(2)}(2)$ and $R=E_8$.
We have $\#O(q_S)=7\cdot 5^2\cdot 3^5\cdot 2^{13}$ and
$\#W(E_8)=7\cdot 5^2\cdot 3^5\cdot 2^{14}$. It follows that the
homomorphism $W(E_8)\to O(q_S)$ is epimorphic and has the kernel $\pm 1$.

$N=20$: Then $q_S\cong q_1^{(2)}(2)\oplus q_{-1}^{(2)}(2)\oplus 3u^{(2)}_+(2)$
and $R=D_8$.
We have $\#O(q_S)=7\cdot 5\cdot 3^2\cdot 2^{13}$ and
$\#W(E_8)=7\cdot 5\cdot 3^2\cdot 2^{14}$. It follows that the
homomorphism $W(D_8)\to O(q_S)$ is epimorphic and has the kernel $\pm 1$.

Any root subsystem $T\subset R$ of finite index can be obtained
by the procedure described in Section \ref{subsubsec3.4.1}. In each case $N=7$,
$8$, $9$, $10$ and $20$ of $R$, applying this procedure,
it is very easy to find all root subsystems $T\subset R$
of finite index and calculate their root invariants.
One can see that it is prescribed by the type of the root system $T$ itself.
We leave these routine calculations to a reader. They are
presented above and will be also very important for
further considerations.
\end{proof}

\begin{remark} \label{rem1a}
Like in the proof above,
using the homomorphism $W^{(4)}(\cM^{(2)})\to O(q_S)$,
one can give the direct proof of the important lemma \ref{lemma2.4.1}
in all elliptic cases of main invariants.
Indeed, it is easy to study its kernel and calculate orders of
the groups. This proof uses calculations of $W^{(2,4)}(S)$ and
$O(S)$ of Theorem \ref{thm3.3.1}.
\end{remark}

Consider a root subsystem $T\subset R$ of Lemma \ref{lemma3.4.4}.
By Theorem \ref{thm2.3.4}, the root subsystem $T\subset R$ defines
a subset $\Delta_+^{(4)}(S)\subset \Delta^{(4)}(S)$, the
corresponding reflection group $W_+^{(2,4)}$, and Dynkin diagram
$\Gamma(P(\cM^{(2,4)}_+))$ of its fundamental chamber
$\cM^{(2,4)}_+$. Direct calculation of these diagrams using
Theorem \ref{thm2.3.4} gives diagrams of Table~2 of Theorem
\ref{thm3.3.2} (where $\Gamma(P(\cM^{(2,4)}_+))$ is replaced  by
$\Gamma(P(X)_+)$) in all cases 7a, b; 8a --- c; 9a --- f; 10a ---
m; 20a --- d. In the remaining cases  9g; 10n, o; 20e --- g we get
diagrams $\Gamma(P(\cM_+^{(2,4)}))$ which we describe below.
Details of these calculations are presented in Appendix, Sections
\ref{fundchambersN=7}--\ref{fundchambersN=20}.

In the {\it Case 9g,} it is better to describe $\Gamma(P(\cM_+^{(2,4)}))$
indirectly. Its black vertices correspond to all points of
$\bP^2(F_2)$ which is the projective plane over the field $F_2$ with
two elements. Its transparent vertices correspond to all lines in $\bP^2(F_2)$.
Both sets have seven elements. Black vertices are disjoint;
transparent vertices are also disjoint; a black vertex is connected
with a transparent vertex by the double edge, if the corresponding point
belongs to the corresponding line, otherwise, they are disjoint.

In the {\it Case 10n,} the diagram $\Gamma(P(\cM_+^{(2,4)}))$ is given
below in Figure \ref{fig3}. Since it is quite complicated, we divide
it in three subdiagrams shown. The first one shows all its edges
connecting black and transparent vertices. The second one shows the
edge connecting the transparent vertices numerated by 1 and 2.  The
third one shows edges connecting transparent vertices 3 --- 12. Each
edge of $\Gamma(P(\cM_+^{(2,4)}))$ is shown in one of these diagrams.
All other our similar descriptions of diagrams as unions of their
subdiagrams have the same meaning. In particular, we have used it in
some diagrams of Table~2.


\begin{figure}
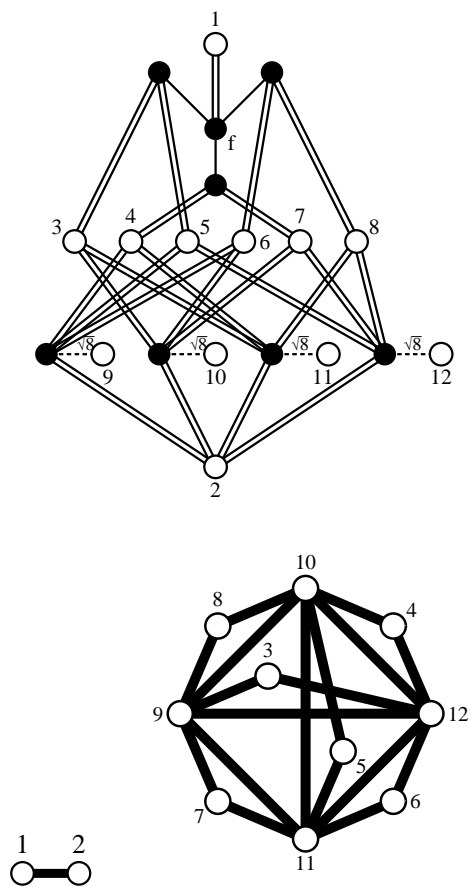

\centerline{\includegraphics[width=6cm]{pics/c-p107-1.eps}}
\includegraphics[width=1.5cm]{pics/c-p107-2.eps}
\includegraphics[width=5cm]{pics/c-p107-3.eps}
\caption{The diagram 10n}
\label{fig3}
\end{figure}

In the {\it Case 10o,} we describe the diagram $\Gamma(P(\cM_+^{(2,4)}))$
indirectly. Its black vertices $f_v$, $v\in V$, correspond to all points
of a three-dimensional affine space $V$ over $F_2$.
Its transparent vertices are of two types. Vertices
$e_v$ of the first type also
correspond to all points $v\in V$. Vertices $e_{\pi}$ of the second type
correspond to all (affine) planes $\pi\subset V$ (there are 14 of them).
Black vertices $f_v$ are disjoint. A black vertex $f_v$
is connected with a transparent vertex $e_{v^\prime}$,
if and only if $v=v^\prime$; the edge has the weight $\sqrt{8}$.
A black vertex $f_v$ is connected with a transparent vertex $e_{\pi}$,
if and only if $v\in \pi$; the edge is double. Transparent
vertices $e_v$, $e_{v^\prime}$ are connected by a thick
edge. A transparent vertex $e_v$ is connected with a transparent
vertex $e_{\pi}$, if and only if $v\notin \pi$; the edge is
thick. Transparent vertices $e_{\pi}$, $e_{\pi^\prime}$ are connected
by edge, if and only if $\pi\|\pi^\prime$; the edge is thick.

In {\it Cases 20e, 20f and 20g} diagrams $\Gamma(P(\cM_+^{(2,4)}))$
are shown in figures \ref{fig4}---\ref{fig6} below.


\begin{figure}
\centerline{\includegraphics[width=6cm]{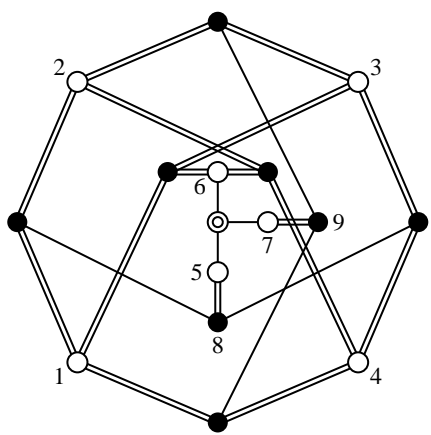}}
\caption{The diagram 20e}
\label{fig4}
\end{figure}


\begin{figure}
\centerline{\includegraphics[width=5cm]{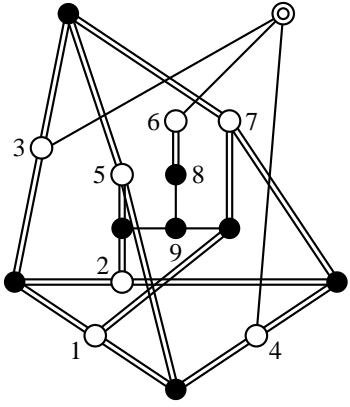}}
\caption{The diagram 20f}
\label{fig5}
\end{figure}


\begin{figure}
\centerline{\includegraphics[width=8cm]{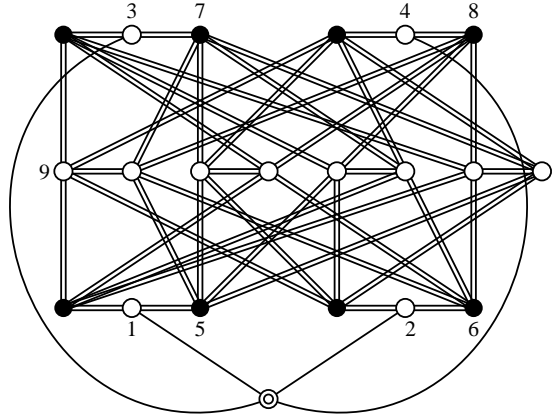}}
\caption{The diagram 20g}
\label{fig6}
\end{figure}

We remark that calculation of $\Gamma(P(\cM_+^{(2,4)}))$ in
cases 7a,b, 8a---c, 9a---g and 10a---o can be obtained from
results of \cite{15} where (in our notation) the dual diagram of
all exceptional curves on the quotient $Y=X/\{1,\theta\}$ is calculated
using completely different method (under the assumption that $Y$ does exist).
By Section \ref{subsec2.4}, both diagrams can be easily obtained from one
another (compare with Section \ref{subsec3.5} below). Therefore, we explain
our method of calculation of $\Gamma(P(\cM_+^{(2,4)}))$ in
more details than it has done in Section \ref{subsubsec2.3.1} only in
the Case 20 (i. e. cases 20a---g).

In the {\it Case 20,} the lattice $S$ has invariants
$(r,a,\delta)=(10,8,1)$, and we can take in $S\otimes \bQ$ an
orthogonal basis $h$, $\alpha$, $v_1,\dots ,v_8$ with $h^2=2$,
$\alpha^2=v_1^2=\cdots =v_8^2=-2$. As $P(\cM^{(2,4)})$, we can
take
\begin{equation}
\begin{split}
P^{(4)}(\cM^{(2,4)})=&\{f_1=v_1-v_2,f_2=v_2-v_3,f_3=v_3-v_4,
f_4=v_4-v_5,\\
&f_5=v_5-v_6,f_6=v_6-v_7,f_7=v_7-v_8,f_8=v_7+v_8,\}
\end{split}
\end{equation}
and
$$
P^{(2)}(\cM^{(2,4)})=
\{\alpha,\,b=\frac{\,h\,}{\,2\,}-\frac{\,\alpha\,}{\,2\,}-v_1,\,
c=h-\frac{\,1\,}{\,2\,}(v_1+v_2+\cdots +v_8)\}.
$$
These elements have Dynkin diagram


\medskip

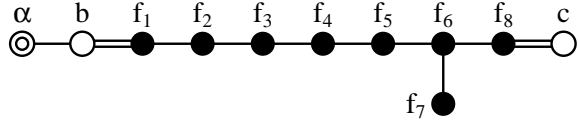
\centerline{\includegraphics[width=8cm]{pics/c-p112.eps}}

\medskip
\noindent
of the case 20a, and they generate and define $S$.

By Section \ref{subsubsec2.3.1}, the set $P^{(2)}(\cM^{(2)})$, where
$\cM^{(2)}\supset \cM^{(2,4)}$, is
$$
W^{(4)}(\cM^{(2)})(\{\alpha,b,c\})
$$
where $W^{(4)}(\cM^{(2)})$ is generated by reflections in $f_1,\dots,f_8$.
It follows that
$$
P(\cM^{(2)})=P^{(2)}(\cM^{(2)})=\{\alpha; b_{\pm i};c_{i_1\dots
i_k}\}
$$
where
$$
b_{\pm i}=\frac{\,h\,}{\,2\,}-\frac{\,\alpha\,}{2}\pm v_i,\
i=1,2,\dots,8;
$$
$$
c_{i_1\dots i_k}=h+\frac{\,1\,}{2}(v_1+v_2+\cdots +v_8)-v_{i_1}-\cdots
-v_{i_k},
$$
where $1\le i_1<i_2<\dots <i_k\le 8$ and $k\equiv 0\mod 2$.
Here all $b_{\pm i}$ give the $W^{(4)}(\cM^{(2)})$-orbit of $b$, and all
$c_{i_1\dots i_k}$ give the $W^{(4)}(\cM^{(2)})$-orbit of $c$.

Elements $f_1,\dots,f_8$ give a basis of the root system $R$ of type $\bD_8$.
If $T\subset R$ is its subsystem of rank $m$, and $t_1,\dots,t_m$ a basis of
$T$, then the fundamental
chamber $\cM_+^{(2,4)}\subset \cM^{(2)}$
defined by $T$ and by its basis $t_1,\dots,t_m$
has
$P(\cM_+^{(2,4)})=P^{(4)}(\cM_+^{(2,4)})\cup P^{(2)}(\cM_+^{(2,4)})$
where
\begin{equation}
\begin{split}
P^{(4)}(\cM_+^{(2,4)})=&\{t_1,\dots,t_m\},\\
P^{(2)}(\cM_+^{(2,4)})=&\{ \alpha \}\cup \{ b_{\pm i}\suchthat
b_{\pm i}\cdot t_s\ge 0,\ 1\le s\le m\}\\
\cup &\{c_{i_1\dots i_k}\suchthat c_{i_1\dots i_k}\cdot t_s\ge 0,\
1\le s\le m \}.
\end{split}
\end{equation}
This describes $\Gamma(P(\cM_+^{(2,4)}))$ completely.

For example, assume that $T\subset R$ has the type $2A_1\oplus D_6$
with the basis $f_1,\,f_9=-v_1-v_2,\,f_3,\,\dots,\,f_8$. Then we get
(after simple calculations)
$$
P^{(2)}(\cM_+^{(2,4)})=\{\alpha,b_{+2},b_{-3},c_{345678},c_{134567}\},
$$
and $\Gamma(P(\cM_+^{(2,4)}))$ is


\medskip

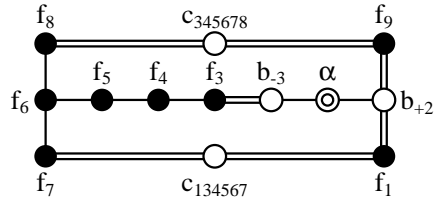
\centerline{\includegraphics[width=6cm]{pics/c-p113.eps}}

\medskip
This gives Case 20b of Table~2.

Exactly the same calculations of $\Gamma(P(\cM_+^{(2,4)}))$ can be
done in all cases 20a---g, and cases 7a,b --- 10a---o of Table of
Lemma \ref{lemma3.4.4} as well. See Appendix, Sections
\ref{fundchambersN=7}--\ref{fundchambersN=20}.

\subsubsection{}\label{subsubsec3.4.5} Here we prove

\begin{proposition}\label{prop3.4.5} Cases 9g, 10n,o and 20 e --- g
of root subsystems $T\subset R$ of Lemma \ref{lemma3.4.4} do not
correspond to non-symplectic involutions $(X,\theta)$ of K3
(in characteristic $0$ and even in characteristic $\ge 3$).
\end{proposition}

\begin{proof} Assume that a root subsystem $T\subset R$ corresponds to
a K3 pair $(X,\theta)$. Then the corresponding Dynkin diagram
$\Gamma(P(\cM_+^{(2,4)}))$ given in Section \ref{subsubsec3.4.4}
coincides with Dynkin diagram $\Gamma(P(X)_+)$ of
exceptional curves of the pair $(X,\theta)$. It
follows the dual diagram of exceptional curves $\Gamma(P(Y))$ on the
corresponding DPN surface $Y=X/\{1,\theta\}$ (see Section \ref{subsec2.3}).
Using this diagram, it is easy to find a sequence of exceptional curves
$E_1,\dots,E_k$ on $Y$ where $k=r-1$ such that their contraction gives
a morphism $\sigma:Y\to \bP^2$. Then other (different from $E_1,\dots,E_k$)
exceptional curves on $Y$ corresponding to Du Val and logarithmic
part of $\Gamma(P(Y))$ give a configuration of rational curves on $\bP^2$
which cannot exist in characteristic 0 and even in characteristic $\ge 3$
(but it exists in characteristic 2). In cases 9g; 10n,o; 20e,f we get
Fano's configuration of seven lines of the finite projective plane over
$F_2$ which can exist only in characteristic 2. In the case 9g one
should contract exceptional curves corresponding to all transparent vertices.
In the case 10n --- corresponding to vertices 1, $f$, 3 --- 8.
In the case  10o --- corresponding to vertices $e_\pi$ where $\pi$ contains
a fixed point $0\in V$ and $e_0$; then curves corresponding to $f_v$,
$v\not=0$, give Fano's configuration. In cases 20e,f --- corresponding to
vertices 1 --- 9. In the case 20g --- corresponding to vertices 1 --- 9,
then we get a conic (corresponding to the double transparent vertex) and four
its tangent lines (corresponding to black vertices different from 5 --- 8)
passing through one point. It is possible only in characteristic 2.
\end{proof}

Another purely arithmetic proof of Proposition \ref{prop3.4.5} (over $\bC$)
can be obtained using Proposition \ref{prop2.6.2}.
This proof is more complicated,
but it can also be done. Here we preferred shorter and geometric
considerations (if diagrams have calculated).
\end{proof}

\subsubsection{} \label{subsubsec3.4.6} Here we prove

\begin{proposition} \label{prop3.4.6} Cases 7a,b; 8a---c, 9a---f; 10a---m
and 20a---d of Table~2 of Theorem \ref{thm3.3.2} correspond to standard
extremal non-symplectic K3 involutions $(X,\theta)$.
\end{proposition}

\begin{proof} Let us calculate root invariants $(K^+(2),\xi^+)$
corresponding to these cases.

Consider the sequence of embeddings of lattices
$$
K^+(2)=[T]\subset [R]\subset S .
$$
It defines the homomorphism
$$
\xi^+:Q=\frac{\,1\,}{2}K^+(2)/K^+(2)\to S^\ast/S\subset \frac{1}{2}S/S
$$
with the kernel $H$. It can be decomposed as
\begin{equation}
\xi^+:Q\stackrel{\widetilde{\xi}^+}{\longrightarrow}\frac{\,1\,}{2}[R]/[R]
\stackrel{\xi^+_R}{\longrightarrow}S^\ast/S\subset
\frac{\,1\,}{2}S/S.
\label{3.4.1}
\end{equation}
Let $H_R=\Ker \xi^+_R$. Then $H=(\widetilde{\xi}^+)^{-1}(H_R)$.
As we know (from our considerations in the super-extremal case),
$H_R=0$ in cases 7, 8, 9, 10.  In the case 20, the $H_R=\bZ/2\bZ$ is
$$
H_R=[\frac{\,1\,}{2}(f_1+f_3+f_5+f_7)+R]/[R]
$$
(see Section \ref{subsubsec3.4.4} about this case). Thus, $H$
can be identified with
$H=(\frac{\,1\,}{2}[T]\cap [R])/[T]$ in cases 7, 8, 9, 10, and with
$$
H=(\frac{\,1\,}{2}[T]\cap[\frac{\,1\,}{2}(f_1+f_3+f_5+f_7)+R])/[T]
$$
in the case 20.

Further details of this calculations in all cases N=7, 8, 9, 10 and 20
are presented in Lemma \ref{lemma3.4.4}.

From these calculations, we get values of $l(H)$ given in
Table~2 of Theorem \ref{thm3.3.2}.

As in Section \ref{subsubsec3.4.2}, using Proposition \ref{prop2.6.1},
one can prove that all these cases when
\begin{equation}
r+a+2l(H)<22
\label{3.4.2}
\end{equation}
correspond to standard extremal non-symplectic K3 involutions
$(X,\theta)$.
Thus, we only need to consider cases when the inequality \eqref{3.4.2}
fails. There are exactly five such cases: 10j,k,l and 20b,d.
Further we consider these cases only.

Below we use some notation and results from \cite{9} about
lattices and their discriminant forms. They are all presented in
Appendix, Section \ref{subsec:discrforms} and Sect
\ref{subsec:maininv}.

 In cases 10j,k,l the discriminant form of $S$ is
$q_S=q_1^{(2)}(2)\oplus 4u_+^{(2)}(2)$. Here, the generator of the
first summand $q_1^{(2)}(2)$ gives the characteristic element
$a_{q_S}$ of the $q_S$, and the second summand $4u_+^{(2)}(2)$
gives the image of $\xi_R^{+}$ from \eqref{3.4.1}, by Lemma
\ref{lemma3.4.2}. Thus, the image of $\xi^+$ belongs to
$4u_+^{(2)}(2)$. The discriminant form of the lattice $M$ (from
\ref{subsec2.6}) is obtained as follows. Let
$$
\Gamma_{\xi^+}\subset Q\oplus {\gA}_S\subset {\gA}_{K^+(2)}\oplus {\gA}_S
$$
be the graph of the homomorphism $\xi^+$ in ${\gA}_{K^+(2)}\oplus {\gA}_S$.
Then
\begin{equation}
q_M=(q_{K^+(2)}\oplus q_S\,|\,(\Gamma_{\xi^+})^
\perp_{q_{K^+(2)}\oplus q_S})/\Gamma_{\xi^+}
\label{3.4.3}
\end{equation}
(here $\Gamma_{\xi^+}$ is an isotropic subgroup).Therefore, $q_{M}\cong
q_1^{(2)}(2)\oplus q^\prime$
since the image of $\xi^+$ belongs to the orthogonal complement of the summand
$q_1^{(2)}(2)$. Considerations in the proof of Proposition \ref{prop2.6.1}
show that
\begin{equation}
\rk M+l({\gA}_{M_2})\le 22
\label{3.4.4}
\end{equation}
since $r+a+2l(H)=22$ in cases 10j,k,l. It is easy to see that
$$
\rk M+l({\gA}_{M_p})<22
$$
for all prime $p>2$. Then, by Theorem 1.12.2 in \cite{9} (see
Appendix,  Theorem \ref{primembedd2}), there exists a primitive
embedding $M\subset L_{K3}$ when either the inequality
\eqref{3.4.4} is strict or $q_{M_2}\cong q_{\pm 1}^{(2)}(2)\oplus
q^\prime$, if it gives the equality. Thus, it always does exist.
It follows that all cases 10j,k,l correspond to standard extremal
non-symplectic K3 involutions $(X,\theta)$ by
Proposition \ref{prop2.6.2} (where we used fundamental Global
Torelli Theorem \cite{13} and surjectivity of Torelli map \cite{5}
for K3).

In cases 20b,d, the proof is exactly the same,
but it is more difficult to prove that
$q_{M_2} \cong q_\theta^{(2)}(2)\oplus q^\prime$ where $\theta=\pm 1$.
In these cases
$$
q_S=3u_+^{(2)}(2)\oplus q_1^{(2)}(2)\oplus q_{-1}^{(2)}(2).
$$
If $\alpha_1$ and $\alpha_2$ are generators of the summands
$q_1^{(2)}(2)$ and $q_{-1}^{(2)}(2)$ respectively,
then $\alpha_{q_S}=\alpha_1+\alpha_2$ is the characteristic element of $q_S$,
and the image of
$\xi^+$ belongs to $3u_+^{(2)}(2)\oplus  [\alpha_{q_S}]$.
In these cases, the lattice
$K^+_{H}(2)$ (see Section \ref{subsec2.6}) is isomorphic to $E_8(2)$.
For example, this is valid because
the subgroups $H$ are the same in cases 10j and 20b, and in cases 10k
and 20d, besides, in
cases 10j and 10k we have $E_8/K^+\cong H$. It follows that
$$
q_{K^+_H(2)}=(q_{K^+(2)}\,|\,(H)^\perp_{q_{K^+(2)}})/H\cong q_{E_ 8(2)}
\cong4u_+^{(2)}(2).
$$
We set $\overline{\Gamma}_{\xi^+}=\Gamma_{\xi^+}/H$. By \eqref{3.4.3}
$$
q_M=(q_{K^+_H(2)}\oplus q_S\,|\,(\overline{\Gamma}_{\xi^+})^\perp_{q_{K^+_H(2)}
\oplus q_S})/\overline{\Gamma}_{\xi^+}.
$$
We have $q_{K^+_H(2)}\oplus q_S=7u_+^{(2)}(2)\oplus q_1^{(2)}(2)\oplus
q_{-1}^{(2)}(2)$. Since
$u_+^{(2)}(2)$ takes values in $\bZ/2\bZ$, the element $\alpha_{q_S}$
(more exactly,
$0\oplus \alpha_{q_S}$) is the characteristic element of
$q_{K^+_H(2)}\oplus q_S$ again.
Moreover, $\alpha_{q_S}\notin \overline{\Gamma}_{\xi^+}$ since
$\Gamma_{\xi^+}$ is the graph
of a homomorphism with the kernel $H$. Therefore
$(\overline{\Gamma}_{\xi^+})^\perp_{q_{K^+_H(2)}
\oplus q_S}$ contains $v$ which is not orthogonal to $\alpha_{q_S}$. Then
$$
(q_{K^+_H(2)}\oplus q_S)(v)=\pm \frac{1}{2}\mod 2
$$
and
$$
[v\mod \overline{\Gamma}_{\xi^+}]\cong q_\theta^{(2)}(2),\ \theta=\pm 1,
$$
is the orthogonal summand of $q_{M_2}$ we were looking for.
\end{proof}

\begin{remark}\label{rem3.4.7} We can give another proof of Proposition
\ref{prop3.4.6} which uses Theorem \ref{thm1.4.1} and considerations which
are inverse to the proof of the previous Proposition \ref{prop3.4.5}.
Indeed, by Theorem \ref{thm1.4.1}, it is enough to prove existence of
rational surfaces with Picard number $r$ and configuration of rational curves
defined by Dynkin diagram of Table~2 of Theorem \ref{thm3.3.2}
(assuming that these Dynkin diagrams correspond to K3 pairs $(X,\theta)$
and considering the quotient by $\theta$). One can prove existence
of these rational surfaces considering
appropriate sequences of blow-ups of  appropriate
relatively minimal rational surfaces $\bP^2$, $\bF_0$,  $\bF_1$,
$\bF_2$, $\bF_3$ or $\bF_4$ with appropriate configurations of rational
curves defined by Dynkin diagrams of Table~2 of Theorem \ref{thm3.3.2}
(see the proof of Proposition \ref{prop3.4.5}). This proof does not
use Global Torelli Theorem and surjectivity of Torelli map for K3.
This gives a hope that results of Chapter \ref{sec2} and
Chapter \ref{sec3} can be generalized to
characteristic $p>0$. Unfortunately, we have proved Theorem \ref{thm1.4.1}
in characteristic $0$ only. Thus, we preferred
the proof of Proposition \ref{prop3.4.6} which is independent of
the results of Chapter \ref{sec1}.
\end{remark}

\subsubsection{}\label{subsubsec3.4.7}
To finish the proof of Theorem \ref{thm3.3.2}, we need to
prove only

\begin{proposition} \label{prop3.4.8} Let $(X,\theta)$ be a non-symplectic
involution of K3 which corresponds to one of cases 7 --- 10 or 20 of
Table 1 of Theorem \ref{thm3.1.1} and
a root subsystem $T\subset R=\Delta^{(4)}(\cM^{(2)})$.

If $(X,\theta)$ is extremal, then $\rk T=\rk R$.
\end{proposition}

\begin{proof}
We can assume (see Section \ref{subsubsec3.4.1}) that $T$ has a basis which
gives a part of a basis of a root subsystem
$\widetilde{T}\subset R=\Delta^{(4)}(\cM^{(2)})$ of the same rank
$\rk \widetilde{T}=\rk R$. Then $\widetilde{T}\subset R$ is one of
root subsystems of Lemma \ref{lemma3.4.4}. If the root subsystem
$\widetilde{T}\subset R$ corresponds to a non-symplectic involution of K3,
i. e. $\widetilde{T}$ gives cases 7a---b, 8a---c, 9a---f, 10a---m and
20a---d, then $T$ is extremal, only if $T=\widetilde{T}$ (by definition).
Then $\rk T=\rk \widetilde{T}=\rk R$ as we want. Thus, it is enough to
consider $\widetilde{T}$ of cases 9g, 10n---o, 20e---g and
$T\subset \widetilde{T}$ to be a primitive root subsystem of a strictly
smaller rank.

Below we consider all these cases.
The following is very important.
In Lemma \ref{lemma3.4.4} we calculated root invariants of
root subsystems $\widetilde{T}\subset R$ of finite index. Restricting
the root invariant of $\widetilde{T}$ on a root subsystem
$T\subset \widetilde{T}$, we get the root invariant of $T\subset R$.
In considerations below, we always consider $T\subset R$ together with
its root invariant. Two root subsystems of $R$ are considered to be
the same, if and only if they are isomorphic root systems together with
their root invariants: then they give equivalent root subsystems
(even with respect to the finite Weyl group $W^{(4)}(\cM^{(2)})$, see
the proof of Lemma \ref{lemma3.4.4}) and isomorphic diagrams.

{\it Case 9g.} Then $\widetilde{T}=7A_1$, and $T=k{A_1}$, $k\le 6$,
is its root subsystem (it is always primitive).
It is easy to see that the same root subsystem $T$
can be obtained as a primitive root subsystem
$T\subset D_4\oplus 3A_1$. Then $T$ is not extremal because
$D_4\oplus 3A_1$ corresponds to K3.

{\it Case 10n.} Then $\widetilde{T}=D_4\oplus 4A_1$ and
$T\subset \widetilde{T}$ is a primitive root subsystem of
the rank $\le 7$. It is easy to see that the same
root subsystem can be obtained as a primitive root subsystem $T$
of $D_6\oplus 2A_1$ or $D_4\oplus D_4$ (then it is not extremal
because $D_6\oplus 2A_1$ and $D_4\oplus D_4$ correspond to K3)
in all cases except when $T=7A_1$.

Let us consider the last case $T=7A_1$ and show (as in Section
\ref{subsubsec3.4.5}) that it does not correspond to K3. As in
Section \ref{subsubsec3.4.4}, one can calculate Dynkin diagram
$\Gamma=\Gamma(P(\cM_+^{(2,4)}))$. See Appendix, Section
\ref{fundchambersN=10}, Case $7A_1\subset E_8$. It is similar to
the case 10o (see Section \ref{subsubsec3.4.4}), but it is more
complicated. We describe it indirectly. One can relate with this
diagram a 3-dimensional linear vector space $V$ over $F_2$.

Black vertices $f_v$ of $\Gamma$ correspond to
$v\in V-\{0\}$ (there are seven of them). Its transparent vertices
(all of them are simple) are

\noindent
$e_v$, $v\in V-\{0\}$; $e_0^{(+)}$, $e_0^{(-)}$;

\noindent
$e_\pi$, $\pi\subset V$ is any affine hyperspace in $V$ which does
not contain $0$;

\noindent
$e_\pi^{(+)}$, $e_\pi^{(-)}$, $\pi\subset V$ is any hyperspace
($0\in \pi$) of $V$.

Edges which connect $f_v$, $e_v$, $e_0^{(+)}$, $e_\pi$,
$e_\pi^{(+)}$ are the same as for the diagram 10o (forget about
$(+)$). The same is valid for $f_v$, $e_v$, $e_0^{(-)}$, $e_\pi$,
$e_\pi^{(-)}$ (forget about $(-)$). Vertices $e_0^{(+)}$ and
$e_0^{(-)}$ are connected by the broken edge of the weight $6$.
Vertices $e_0^{(+)}$ and $e_\pi^{(-)}$ (and $e_0^{(-)}$,
$e_\pi^{(+)}$ as well) are connected by the broken edge of the
weight $4$. This gives all edges of $\Gamma$.

Assume that $\Gamma$ corresponds to a K3 pair $(X,\theta)$.
Consider the corresponding DPN surface and contract exceptional
curves corresponding to $e_{\pi}^{(+)}$ and
$e_0^{(+)}$. Then exceptional curves of $f_v$, $v\in V-\{0\}$,
give Fano's configuration on $\bP^2$ which exists only in
characteristic 2. We get a contradiction.

{\it Case 10o.} This is similar to the previous case.

 {\it Case 20e.} Then $\widetilde{T}=2A_3\oplus 2A_1$ and
$T$ is its primitive root subsystem of the rank $\le 7$. It is
easy to see that the same root subsystem can be obtained as
a primitive root subsystem of $D_6\oplus 2A_1$ or $D_5\oplus A_3$
(and it is not then extremal because $D_6\oplus 2A_1$ and $D_5\oplus A_3$
correspond to K3) in all cases except $T=A_3\oplus 4A_1$.

Let us consider the last case $T=A_3\oplus 4A_1$ and show (as in
Section \ref{subsubsec3.4.5}) that it does not correspond to K3. As
in Section \ref{subsubsec3.4.4}, one can calculate Dynkin diagram
$\Gamma=\Gamma(P(\cM_+^{(2,4)}))$. See Appendix, Section
\ref{fundchambersN=20}, Case $4A_1\oplus A_3\subset D_8$. It has
exactly one transparent double vertex $\alpha$ and eight simple
transparent vertices $c_v$, $v\in V(K)$, where $V(K)$ is the set
of vertices of a 3-dimensional cube $K$ with distinguished two
opposite 2-dimensional faces $\beta,\, \beta^\prime\in \gamma(K)$
where $\gamma(K)$ is the set of all 2-dimensional faces of $K$.
Black vertices of $\Gamma$ are $f_\gamma$, $\gamma\in \gamma(K)$,
and one more black vertex $f_0$. Simple transparent vertices of
$\Gamma$ which are connected by a simple edge with $\alpha$ are
either $b_{\overline{\gamma}}$, $\overline{\gamma}\in
\overline{\gamma(K)}$, where $\overline{\gamma(K)}$ is the set of
pairs of opposite 2-dimensional faces of $K$, or $b_t$, $t\in
\overline{V(K)}$. Here $\overline{V(K)}$ consists of two elements
corresponding to a choice of one vertex from each pair of opposite
vertices of $K$ in such a way that neither of three of them are
contained in a 2-dimensional face $\gamma\in \gamma(K)$ (they
define a regular tetrahedron with edges which are diagonals of
2-dimensional faces of $K$).

Let us describe edges of $\Gamma$ different from above. Thick
edges connect $c_v$ corresponding to opposite vertices $v\in
V(K)$, vertices $b_{t_1}$ and $b_{t_2}$ where
$\{t_1,t_2\}=\overline{V(K)}$, vertices $b_t$ and $c_v$ where
$v\in t$. Simple edges connect $f_0$ with $f_\beta$ and
$f_{\beta^\prime}$. Double simple edges connect $c_v$ with
$f_\gamma$, if $v\in \gamma$, and $b_{\overline{\gamma}}$ with
$f_\gamma$, if $\gamma\in
\overline{\gamma}-\{\beta,\beta^\prime\}$, and the vertex
$b_{\overline{\beta}}$ with $f_0$.

Assume that $\Gamma$ corresponds to a K3 pair $(X,\theta)$. On its
DPN surface, let us contract exceptional curves corresponding to
$c_v$, $v\in t$; $b_{\overline{\gamma}}$, $\overline{\gamma}\in
\overline{\gamma(K)}$; $f_0$ and $b_{t^\prime}$, $t^\prime\not=t$
(here $t\in \overline{V(K)}$ is fixed). Then curves corresponding
to $f_v$, $v\in V(K)$, and the vertex $\alpha$ define Fano's
configuration of lines in $\bP^2$ which can exist only in
characteristic 2.

{\it Cases 20f,g.} In these cases, $\widetilde{T}=D_4\oplus 4A_1$ or
$\widetilde{T}=8A_1$. As for analogous cases 10n,o, everything is
reduced to prove that $T=7A_1$ does not correspond to a K3 pair
$(X,\theta)$.

In this case, $\Gamma=\Gamma(P(\cM_+^{(2,4)}))$ is as follows. See
Appendix, Section \ref{fundchambersN=20}, Case $7A_1\subset D_8$.
Let $I=\{1,2,3,4\}$ and $J=\{1,2\}$. The $\Gamma$ has: exactly one
double transparent vertex $\alpha$; black vertices $f_{ij}$, $i\in
I$, $j\in J$, and $(i,j)\not=(4,2)$; simple transparent vertices
$b_i$, $i=1,2,3$, and $b_{4(+)}$, $b_{4(-)}$ which are connected
by a simple edge with $\alpha$; simple transparent vertices
$c_{j_1j_2j_3j_4}$ where $j_1,j_2,j_3\in J$, $j_4\in \{1,-2,+2\}$
and $j_1+j_2+j_3+j_4\equiv 0\mod 2$ which are disjoint to
$\alpha$.

Edges of $\Gamma$ which are different from above, are as follows.

Double edges connect $b_i$ with $f_{ij}$, if $i=1,2,3$, and
$b_{4(+)}$, $b_{4(-)}$ with $f_{41}$, and $c_{j_1j_2j_3j_4}$ with
$f_{1j_1}$, $f_{2j_2}$, $f_{3j_3}$, and $c_{j_1j_2j_31}$ with
$f_{41}$.

Thick edges connect $b_{4(\pm)}$ with $c_{j_1j_2j_3(\mp 2)}$,
and $c_{j_1j_2j_3j_4}$ with $c_{j_1^\prime j_2^\prime j_3^\prime j_4\prime}$,
if $j_1\not=j_1^\prime$, $j_2\not=j_2^\prime$, $j_3\not=j_3^\prime$,
$|j_4|\not=|j_4^\prime|$, and $c_{j_1j_2j_3(+2)}$ with
$c_{j_1^\prime j_2^\prime j_3^\prime (-2)}$, if
$(j_1,j_2,j_3)\not=(j_1^\prime, j_2^\prime, j_3^\prime)$.

Assume that $\Gamma$ corresponds to a K3 pair $(X,\theta)$. On its
DPN surface, let us contract exceptional curves corresponding to
$b_1$, $b_2$, $b_3$, $b_{4(+)}$, $f_{11}$, $f_{21}$, $f_{31}$,
$f_{41}$, $c_{222(+2)}$. The curve corresponding to $\alpha$ gives
a conic in $\bP^2$. Curves corresponding to $f_{12}$, $f_{22}$, $f_{32}$
give lines touching to the conic and having a common point.
This is possible in characteristic 2 only.

This finishes the proof of Theorem \ref{thm3.3.2}
\end{proof}

\subsection{Classification of DPN surfaces of elliptic type}
\label{subsec3.5}
Each non-symplectic involution  of elliptic type
$(X,\theta)$ of K3 gives rise to the right DPN pair  $(Y,C)$ where
\begin{equation}
Y=X/\{1,\theta\},\ C=\pi(X^\theta)\in |-2K_Y|,
\end{equation}
$\pi:X\to Y$ the quotient morphism; and vice versa.
From Theorems \ref{thm3.3.2}, \ref{thm3.3.3}, we then get
classification of right DPN pairs $(Y,C)$ and DPN surfaces $Y$ of
elliptic type.
See Chapter \ref{sec2} and especially
Sections \ref{subsec2.1} and \ref{subsec2.7}. It is obtained by
the reformulation of Theorems \ref{thm3.3.2} and \ref{thm3.3.3} and by
redrawing of the diagrams. But, for readers' convenience, we do it below.

\begin{theorem}[Classification Theorem for right DPN surfaces of
elliptic type in the extremal case] \label{thm3.5.1}
A right DPN surface $Y$ of elliptic type is {\bf extremal}
if and only if the number of its exceptional curves
with the square $(-2)$  is maximal for the fixed main invariants
$(r,a,\delta)$ (equivalently, $(k,g,\delta)$). (It is equal to
the number of black vertices in the diagram $\Gamma$ of Table 3 below.)

Moreover, the dual diagram $\Gamma(Y)$ of all exceptional curves on
extremal $Y$ is isomorphic to one of diagrams $\Gamma$ given in Table 3.
Vice versa any diagram $\Gamma $ of Table 3 corresponds to some of the $Y$
(the $Y$ can be even taken standard).

In the diagrams $\Gamma$, simple transparent vertices correspond to
curves of the 1st kind (i. e. to non-singular rational irreducible
curves with the square $(-1)$), double transparent vertices correspond to
non-singular rational irreducible curves with the square $(-4)$,
black vertices correspond to non-singular rational irreducible
curves with the square $(-2)$, a $m$-multiple edge (or an edge with
the weight $m$ when $m$ is large) means the intersection index $m$ for
the corresponding curves. Any exceptional curve on $Y$ is one of
these curves.
\end{theorem}

For a not necessarily extremal right DPN surface $Y$ of elliptic
type the dual diagram $\Gamma=\Gamma(Y)$ of all exceptional curves on
$Y$ also consists of simple transparent, double transparent and black
vertices which have exactly the same meaning as in Theorem
\ref{thm3.5.1} above. All black vertices of $\Gamma$ define the{\it Du
  Val part $\Duv \Gamma$} of $\Gamma$. All double transparent vertices
of $\Gamma$, and all simple transparent vertices of $\Gamma$ which are
connected by two edges with double transparent vertices of $\Gamma$
(there are always two of these double transparent vertices) define the
{\it logarithmic part $\Log \Gamma$} of $\Gamma$. The rest of vertices
(different from vertices of $\Duv \Gamma$ and $\Log \Gamma$) define the
{\it varying part $\Var \Gamma$} of $\Gamma$.  In Theorem below we
identify vertices of $\Gamma (Y)$ with elements of Picard lattice
$\Pic \widetilde{Y}$, then weights of edges are equal to the
corresponding intersection pairing in this lattice which makes sense
to the descriptions of the graphs $\Var \Gamma(Y)$ and $\Gamma (Y)$.

\begin{theorem}[Classification Theorem for right DPN surfaces
in the non-extremal, i. e. arbitrary, case of elliptic type]
\label{thm3.5.2}
Dual diagrams $\Gamma (Y)$ of all exceptional curves of not
necessarily extremal right DPN surfaces $Y$ of elliptic type
are described by arbitrary (without any restrictions) subdiagrams
$D\subset \Duv \Gamma$ of extremal DPN surfaces described in  Theorem
\ref{thm3.5.1} above with the
same main invariants $(r,a,\delta)$ (equivalently $(k,g,\delta)$).

Moreover, $\Duv \Gamma (Y)=D$, $\Log \Gamma (Y)=\Log \Gamma$,
and these subdiagrams are disjoint to each other;
$$
\Var \Gamma (Y)=\{f\in W(\Var \Gamma)\ |\
f\cdot D\ge 0\}
$$
where $W$ is the subgroup of automorphisms of the
Picard lattice of the extremal DPN surface (the Picard lattice is
defined by the diagram $\Gamma$), generated by reflections in
elements with square $-2$ corresponding to all vertices of
$\Duv \Gamma$.

Two such subdiagrams $D\subset \Duv \Gamma$ and
$D^\prime \subset \Duv \Gamma^\prime$ (with the same main invariants)
give DPN surfaces $Y$ and $Y^\prime$
with isomorphic diagrams $\Gamma(Y)\cong \Gamma(Y^\prime)$,
if and only if they have
isomorphic root invariants $([D],\xi^+)$  and
$([D^\prime],(\xi^\prime)^+)$ (see Theorem \ref{thm3.3.3}).

To calculate the root invariant
$([D],\xi^+)$ of a DPN surface, one has to go back from the
graph $\Gamma$ of Table 3 to the corresponding graph of Tables 1 or 2.
\end{theorem}

From our point of view, classification above by graphs of
exceptional curves is the best classification of DPN surfaces $Y$.
It shows a sequence (actually all sequences) of $-1$ curves which
should be contracted to get the corresponding relatively minimal
rational surface $\overline{Y}$ isomorphic to $\bP^2$ or $\bF_n$, $n\le 4$
(see Section \ref{subsec3.6} and Table 4 below).
Images of exceptional curves on $Y$ which are not contracted then give
some configuration of rational curves on $\overline{Y}$ which should exist to
get the DPN surface $Y$ back from $\overline{Y}$ by the corresponding
sequence of blow ups. Here the following inverse statement is very
important because it
shows that any surface $Y^\prime$ obtained by a ``similar'' sequence
of blow ups of $\overline{Y}$ which are related with a ``similar''
configuration of rational curves on $\overline{Y}$ will be also a DPN
surface with the graph $\Gamma(Y^\prime)$ of exceptional curves which is
isomorphic to $\Gamma(Y)$.
Here is the exact statement.

\begin{theorem} \label{thm3.5.2.a}
Let $Y$ be a right DPN surface of elliptic type, and the set of
exceptional curves on $Y$
is not empty (i. e. $Y$ is different from $\bP^2$ and $\bF_0$).
Let $Y^\prime$ be a non-singular rational surface such that

1) the Picard number of $Y^\prime$ is equal to the Picard number of $Y$.

2) there exists a set $E_1,\dots, E_m$ of non-singular irreducible rational
exceptional curves on $Y^\prime$ such that
their dual graph is isomorphic to the dual graph $\Gamma(Y)$
of exceptional curves on $Y$.

Then $Y^\prime$ is also a DPN surface and $E_1,\dots, E_m$ are all
exceptional curves on $Y^\prime$
(of course, then $\Gamma(Y^\prime)\cong \Gamma(Y)$).
\end{theorem}

\begin{proof} Let $r$ be the Picard number of $Y$ and $Y^\prime$.
If $r=2$, then obviously $Y\cong Y^\prime\cong \bF_n$, $n>0$.
Further we assume that $r\ge 3$. We denote by $S_Y$ and
$S_{Y^\prime}$ the Picard lattices of $Y$ and $Y^\prime$
respectively. Like for K3 surfaces we shall consider the light
cones $V(Y)\subset S_Y\otimes \bR$, $V(Y^\prime)\subset
S_Y{^\prime}\otimes \bR$ (of elements with positive square) and
their halves $V^+(Y)$ and $V^+(Y^\prime)$ containing
polarizations.

Let $D_1,\dots, D_m$ are all
exceptional curves on $Y$ (corresponding to vertices of $\Gamma(Y)$).
Their number
is finite and they generate $S_Y$ since $r\ge 3$.
We claim that Mori cone
$\NE(Y)=\bR^+D_1+\cdots +\bR^+D_m$ is
generated by $D_1,\dots D_m$. This is equivalent to
\begin{equation}
\overline{V^+(Y)}\subset \bR^+D_1+\cdots +\bR^+D_m
\label{mori1}
\end{equation}
since $D_j$ are all exceptional curves on $Y$ and $V^+(Y)\subset
\NE(Y)$ by Riemann-Roch Theorem on $Y$. The condition \eqref{mori1} is
equivalent to the embedding of dual cones
\begin{equation}
(\bR^+D_1+\cdots +\bR^+D_m)^\ast \subset \overline{V^+(Y)}
\label{mori2}
\end{equation}
because the light cone $V^+(Y)$ is self-dual. By considering the corresponding
K3 double cover $\pi:X\to Y$, the embedding \eqref{mori2} is equivalent to
the embedding
\begin{equation}
(\bR^+\pi^\ast(D_1)+\cdots +\bR^+\pi^\ast(D_m))^\ast \subset \overline{V^+(S)}
\label{mori3}
\end{equation}
which is equivalent to finiteness of volume of $\cM(X)_+\subset \cL(S)$
which we know.

The equivalent conditions \eqref{mori1} and \eqref{mori2} are numerical. Thus,
similar conditions
\begin{equation}
\overline{V^+(Y^\prime)}\subset \bR^+E_1+\cdots +\bR^+E_m
\label{mori4}
\end{equation}
and
\begin{equation}
(\bR^+E_1+\cdots +\bR^+E_m)^\ast \subset \overline{V^+(Y^\prime)}
\label{mori5}
\end{equation}
are valid for $Y^\prime$. This shows that $E_1,\dots E_m$ are the
only exceptional curves on $Y^\prime$.  Indeed, if $E$ is any
other irreducible curve $E$ on $Y^\prime$ satisfying $E\cdot
E_i\ge 0$, then $E^2\ge 0$ by \eqref{mori5} and the curve $E$ is
not exceptional. Thus, $\Gamma(Y^\prime)$ and $\Gamma (Y)$ are
isomorphic. In the same way as for $Y$ above, we then get from
\eqref{mori4} or \eqref{mori5} that the Mori cone
$\NE(Y^\prime)=\bR^+E_1+\cdots +\bR^+E_m$ is generated
by $E_1,\dots, E_m$.

Let us show that $Y^\prime$ is a DPN surface. Definitions of Du Val,
logarithmic
parts of $\Gamma(Y)$ were purely numerical. Since $\Gamma(Y^\prime)$ and
$\Gamma(Y)$
are isomorphic, we can use similar notions for $Y^\prime$.

In Section \ref{subsec4.1} we shall prove
(without using Theorem \ref{thm3.5.2.a}) that there exists a contraction
$p:Y\to Z$ of Du Val and logarithmic parts of exceptional curves of $Y$
which gives the right resolution of singularities of a log del Pezzo
surface $Z$ of index $\le 2$. (Remark that by Lemma \ref{lemma1.3.2} it
also gives another proof of the above statements about Mori cone and
exceptional curves on $Y$ and $Y^\prime$.)
Thus, the element $p^\ast (-2K_Z)\in S_Y$ is defined.
It equals to $-2K_Y$ minus
sum of all exceptional curves on $Y$ with square $-4$. Thus, similar
element can be defined for $Y^\prime$. Let us denote it by
$R\in S_{Y^\prime}$.
In Section \ref{subsec1.4}, we had
proved (for any log del Pezzo surface $Z$ of index $\le 2$)
that the linear system $p^\ast(-2K_Z)$ contains a non-singular curve.
The proof was purely numerical and only used the fact that $-2K_Y-\sum
E_i$ is big and nef.
The same proof for $Y^\prime$ gives that $R$ contains a non-singular curve.
It follows that $Y^\prime$ is a right DPN surface of elliptic type.
\end{proof}

\noindent
\begin{table}\label{table3}
\caption{Dual diagrams $\Gamma$ of all exceptional curves of
extremal right DPN surfaces of elliptic type}

\index{Table 3}

\addtocontents{toc}{\contentsline {section}{\tocsection {}{T.3}{Table 3}}
{\pageref{table3}}}

\begin{tabular}{|r||r|r|r|r|r|c|c|}
\hline
 $N$& $r$& $a$&$\delta$& $k$&$g$&$\widetilde{r}$& $\Gamma$\\
\hline
  1 &  1 &  1 &  1     &  0 &10 &   1  &
$\Gamma=\emptyset$,\ \ $\bP^2$\\
\hline
  2 &  2 &  2 &  0     &  0 & 9 &   1  &
\includegraphics[width=1cm]{pics/c-p126-1.eps}\ \ $\bF_0$ or $\bF_2$\\
\hline
  3 &  2 &  2 &  1     &  0 & 9 &   2  &
\includegraphics[width=1cm]{pics/c-p126-2.eps}\ \ $\bF_1$\\
\hline
  4 &  3 &  3 &  1     &  0 & 8 &   2  &
\includegraphics[width=2.5cm]{pics/c-p126-3.eps}\\
\hline
  5 &  4 &  4 &  1     &  0 & 7 &   1  &
\includegraphics[width=2.5cm]{pics/c-p126-4.eps}\\
\hline
  6 &  5 &  5 &  1     &  0 & 6 &   1  &
\includegraphics[width=2.5cm]{pics/c-p126-5.eps}\\
\hline
  7 &  6 &  6 &  1     &  0 &5  & 1    & \\
a &    &    &          &    &   &      &
\includegraphics[width=3cm]{pics/c-p126-6.eps}\\
b &    &    &          &    &   &      &
\includegraphics[width=4cm]{pics/c-p126-7.eps}\\
\hline
 8 &  7  &  7 &  1     &  0 &4  & 1    & \\
a &    &    &        &    &   &       &
\includegraphics[width=3.5cm]{pics/c-p127-1.eps}\\
b &    &    &        &    &   &     &
\includegraphics[width=4cm]{pics/c-06p127-2.eps}\\
c &    &    &        &    &   &     &
\includegraphics[width=3cm]{pics/c-p127-3.eps}\\
\hline
  9 &  8 &  8 &  1     &  0 &3  &  1 & \\
a &    &    &        &    &   &     &
\includegraphics[width=4cm]{pics/c-p127-4.eps}\\
b &    &    &        &    &   &     &
\includegraphics[width=4cm]{pics/c-p127-5.eps}\\
\hline
\end{tabular}

\end{table}

\begin{tabular}{|r||r|r|r|r|r|c|c|}
\hline
 $N$& $r$& $a$&$\delta$& $k$&$g$&$\widetilde{r}$& $\Gamma$\\
\hline
  9 &  8 &  8 &  1     &  0 &3  &  1 & \\
c &    &    &        &    &   &     &
\includegraphics[width=3cm]{pics/c-p127-6.eps}\\
d &    &    &        &    &   &     &
\includegraphics[width=3cm]{pics/c-p128-1.eps}\\
e &    &    &        &    &   &     &
\includegraphics[width=4.5cm]{pics/c-p128-2.eps}\\
f &    &    &        &    &   &     &
\includegraphics[width=3cm]{pics/c-06p128-3.eps}\\
\hline
 10 &  9 &  9 &  1     &  0 &2  &  1  & \\
a &    &    &        &    &   &     &
\includegraphics[width=4.5cm]{pics/c-p128-4.eps}\\
b &    &    &        &    &   &     &
\includegraphics[width=4.5cm]{pics/c-p128-5.eps}\\
c &    &    &        &    &   &     &
\includegraphics[width=3.5cm]{pics/c-p128-6.eps}\\
\hline
\end{tabular}

\begin{tabular}{|r||r|r|r|r|r|c|c|}
\hline
 $N$& $r$& $a$&$\delta$& $k$&$g$&$\widetilde{r}$& $\Gamma$\\
\hline
 10 &  9 &  9 &  1     &  0 &2  &  1  & \\
d &    &    &        &    &   &     &
\includegraphics[width=6cm]{pics/c-p129-1.eps}\\
e &    &    &        &    &   &     &
\includegraphics[width=3.5cm]{pics/c-p129-2.eps}\\
f &    &    &        &    &   &     &
\includegraphics[width=5cm]{pics/c-p129-3.eps}\\
g &    &    &        &    &   &     &
\includegraphics[width=6cm]{pics/c-p129-4.eps}\\
h &    &    &        &    &   &     &
\includegraphics[width=6cm]{pics/c-p130-1.eps}\\
\hline
\end{tabular}

\begin{tabular}{|r||r|r|r|r|r|c|c|}
\hline
 $N$& $r$& $a$&$\delta$& $k$&$g$&$\widetilde{r}$& $\Gamma$\\
\hline
 10 &  9 &  9 &  1     &  0 &2  &  1  & \\
\hline
i &    &    &        &    &   &     &
\includegraphics[width=5cm]{pics/c-p130-2.eps}\\
j &    &    &        &    &   &     &
\includegraphics[width=7cm]{pics/c-p130-3.eps}\\
k &    &    &        &    &   &     &
\includegraphics[width=5.5cm]{pics/c-p131-1.eps}\\
 &    &    &        &    &   &     &
\includegraphics[width=6cm]{pics/c-p131-2.eps}\\
l &    &    &        &    &   &     &
\includegraphics[width=4cm]{pics/new-p131-3.eps}\\
\hline
\end{tabular}

\begin{tabular}{|r||r|r|r|r|r|c|c|}
\hline
 $N$& $r$& $a$&$\delta$& $k$&$g$&$\widetilde{r}$& $\Gamma$\\
\hline
 10 &  9 &  9 &  1     &  0 &2  &  1  & \\
\hline
    &    &    &        &    &   &     &
\includegraphics[width=6cm]{pics/c-p131-4.eps}\\
m &    &    &        &    &   &     &
\includegraphics[width=4cm]{pics/c-p131-5.eps}\\
\hline\hline
  11&  2 &  0 &  0     &  1 & 10&   1  &
\includegraphics[width=1cm]{pics/c-p132-1.eps}\ \ $\bF_4$\\
\hline
  12&  3 &  1 &  1     &  1 & 9&   2  &
\includegraphics[width=2.5cm]{pics/c-p132-2.eps}\\
\hline
  13&  4 &  2 &  1     &  1 & 8&   2  &
\includegraphics[width=2.5cm]{pics/c-p132-3.eps}\\
\hline
  14&  5 &  3 &  1     &  1 & 7&   2  &
\includegraphics[width=3cm]{pics/c-p132-4.eps}\\
\hline
  15&  6 &  4 &  0     &  1 & 6&   1  &
\includegraphics[width=3.5cm]{pics/c-p132-5.eps}\\
\hline
  16&  6 &  4 &  1     &  1 & 6&   2  &
\includegraphics[width=4cm]{pics/c-p132-6.eps}\\
\hline
  17&  7 &  5 &  1     &  1 & 5&   2  &
\includegraphics[width=4.7cm]{pics/c-p132-7.eps}\\
\hline
  18&  8 &  6 &  1     &  1 & 4&   1  &
\includegraphics[width=5.2cm]{pics/c-p132-8.eps}\\
\hline
  19&  9 &  7 &  1     &  1 & 3&   1  &
\includegraphics[width=5.2cm]{pics/c-p132-9.eps}\\
\hline
\end{tabular}

\begin{tabular}{|r||r|r|r|r|r|c|c|}
\hline
 $N$& $r$& $a$&$\delta$& $k$&$g$&$\widetilde{r}$& $\Gamma$\\
\hline
 20 &  10&  8 &  1     &  1 &2  &   1  & \\
 a  &    &    &        &    &   &    &
\includegraphics[width=6cm]{pics/c-p132-10.eps}\\
 b  &    &    &        &    &   &    &
\includegraphics[width=4.5cm]{pics/c-p132-11.eps}\\
 c  &    &    &        &    &   &    &
\includegraphics[width=4.5cm]{pics/c-p133-1.eps}\\
 d  &    &    &        &    &   &    &
\includegraphics[width=4.5cm]{pics/c-p133-2.eps}\\
\hline\hline
  21&  6&  2 &  0     &  2 & 7&   1  &
\includegraphics[width=3.5cm]{pics/c-p133-3.eps}\\
\hline
  22&  7&  3 &  1     &  2 & 6&   2  &
\includegraphics[width=4.7cm]{pics/c-p133-4.eps}\\
\hline
  23&  8&  4 &  1     &  2 & 5&   2  &
\includegraphics[width=5.2cm]{pics/c-p133-5.eps}\\
\hline
  24&  9&  5 &  1     &  2 & 4&   2  &
\includegraphics[width=5.2cm]{pics/c-p133-6.eps}\\
\hline
  25& 10&  6 &  0     &  2 & 3&   1  &
\includegraphics[width=6cm]{pics/c-p133-7.eps}\\
\hline
  26& 10&  6 &  1     &  2 & 3&   1  &
\includegraphics[width=5.2cm]{pics/c-p133-8.eps}\\
\hline
  27& 11&  7 &  1     &  2 & 2&   1  &
\includegraphics[width=3cm]{pics/c-p134-1.eps}\\
\hline\hline
  28& 8&   2 &  1     &  3 & 6&   2  &
\includegraphics[width=4cm]{pics/c-p134-2.eps}\\
\hline
  29& 9&   3 &  1     &  3 & 5&   3  &
\includegraphics[width=5.2cm]{pics/c-p134-3.eps}\\
\hline
  30&10&   4 &  0     &  3 & 4&   1  &
\includegraphics[width=6cm]{pics/c-p134-4.eps}\\
\hline
\end{tabular}

\begin{tabular}{|r||r|r|r|r|r|c|c|}
\hline
 $N$& $r$& $a$&$\delta$& $k$&$g$&$\widetilde{r}$& $\Gamma$\\
\hline
  31&10&   4 &  1     &  3 & 4&   3  &
\includegraphics[width=4cm]{pics/c-p134-5.eps}\\
\hline
  32&11&   5 &  1     &  3 & 3&   2  &
\includegraphics[width=4cm]{pics/c-p134-6.eps}\\
\hline
  33&12&   6 &  1     &  3 & 2&   1  &
\includegraphics[width=3cm]{pics/c-p134-7.eps}\\
\hline\hline
  34& 9&   1 &  1     &  4 & 6&   2  &
\includegraphics[width=4.7cm]{pics/c-p134-8.eps}\\
\hline
  35&10&   2 &  0     &  4 & 5&   2  &
\includegraphics[width=5.2cm]{pics/c-p135-1.eps}\\
\hline
  36&10&   2 &  1     &  4 & 5&   3  &
\includegraphics[width=5.8cm]{pics/c-p135-2.eps}\\
\hline
  37&11&   3 &  1     &  4 & 4&   3  &
\includegraphics[width=5.2cm]{pics/c-p135-3.eps}\\
\hline
  38&12&   4 &  1     &  4 & 3&   2  &
\includegraphics[width=5.2cm]{pics/c-p135-4.eps}\\
\hline
  39&13&   5 &  1     &  4 & 2&   2  &
\includegraphics[width=5.2cm]{pics/c-p135-5.eps}\\
\hline\hline
  40&10&   0 &  0     &  5 & 6&   1  &
\includegraphics[width=5.2cm]{pics/c-p135-6.eps}\\
\hline
  41&11&   1 &  1     &  5 & 5&   2  &
\includegraphics[width=6cm]{pics/c-p135-7.eps}\\
\hline
\end{tabular}

\begin{tabular}{|r||r|r|r|r|r|c|c|}
\hline
 $N$& $r$& $a$&$\delta$& $k$&$g$&$\widetilde{r}$& $\Gamma$\\
\hline
  42&12&   2 &  1     &  5 & 4&   2  &
\includegraphics[width=6.5cm]{pics/c-p135-8.eps}\\
\hline
  43&13&   3 &  1     &  5 & 3&   2  &
\includegraphics[width=5.2cm]{pics/c-p136-1.eps}\\
\hline
  44&14&   4 &  0     &  5 & 2&   1  &
\includegraphics[width=5.2cm]{pics/c-p136-2.eps}\\
\hline
45&14&   4 &  1     &  5 & 2&   2  &
\includegraphics[width=5.2cm]{pics/c-p136-3.eps}\\
\hline\hline
46&14&   2 &  0     &  6 & 3&   1  &
\includegraphics[width=6.7cm]{pics/c-p136-4.eps}\\
\hline
47&15&   3 &  1     &  6 & 2&   2  &
\includegraphics[width=6cm]{pics/c-p136-5.eps}\\
\hline\hline
48&16&   2 &  1     &  7 & 2&   2  &
\includegraphics[width=6.5cm]{pics/c-p137-1.eps}\\
\hline\hline
49&17&   1 &  1     &  8 & 2&   2  &
\includegraphics[width=6.5cm]{pics/c-p137-2.eps}\\
\hline\hline
50&18&   0 &  0     &  9 & 2&   1  &
\includegraphics[width=6.5cm]{pics/c-p137-3.eps}\\
\hline
\end{tabular}

\medskip

\subsection{Application: On classification of plane sextics
with simple singularities}
\label{subsec3.6} Let $Y$ be a right DPN surface of elliptic
type which were classified in Theorems \ref{thm3.5.1}, \ref{thm3.5.2} and
\ref{thm3.5.2.a}.
Let $\Gamma(Y)$ be the dual diagram of all exceptional curves on $Y$.
By definition of right DPN surfaces, there exists a non-singular curve
\begin{equation}
C=C_g+E_{a_1}+\cdots +E_{a_k}\in |-2K_Y|
\end{equation}
where $E_{a_1},\dots, E_{a_k}$ are exceptional curves with square $(-4)$
corresponding to all double transparent vertices $a_1,\dots,a_k$ of
$\Gamma(Y)$ and $g>1$ the genus of the irreducible non-singular curve $C_g$.
Here $(k,g,\delta)$ (equivalent to $(r,a,\delta)$) are the main invariants
of $Y$.

We denote by $E_v$ the exceptional curve on $Y$ corresponding to a vertex
$v\in V(\Gamma(Y))$.
If $v$ is black, then $C\cdot E_v=C_g\cdot E_v=0$.
If $v$ is simple transparent, then $C\cdot E_v=2$.

If $v$ is simple transparent and $v$ is not
connected by any edge with double transparent vertices of $\Gamma(Y)$
(i. e. $E_v\cdot E_{a_i}=0$, $i=1,\dots, k$) then $C_g\cdot E_v=2$.
This intersection index can be obtained in two ways:
\begin{equation}
C_g\ \text{intersects}\ E_v\ \text{transversally in two points;}
\label{3.6.1}
\end{equation}
\begin{equation}
C_g\ \text{simply touches}\ E_v\ \text{in one point.}
\label{3.6.2}
\end{equation}
(For example, in Case 47 of Table 3 we have two such vertices $v$.)

Up to this ambiguity, we know (from the diagram $\Gamma(Y)$) how components of
$C$ intersect exceptional curves. Which of possibilities \eqref{3.6.1}
or \eqref{3.6.2}
does take place is defined by the generalized root invariant which we
don't consider in this work.

Let $t_1,\dots, t_{r-1}\in V(\Gamma(Y))$ be a sequence of vertices such
that the contraction of exceptional curves $E_{t_1},\dots,E_{t_{r-1}}$
gives a morphism $\sigma:Y\to \bP^2$ which is a sequence of
contractions of curves of the 1st kind. By Section \ref{subsec2.1}, the image
$D=\sigma(C)\subset \bP^2$ is then a sextic (it belongs to
$|-2K_{\bP^2}|$) with simple singularities. What components and what
singularities the curve $D$
does have is defined by the subgraph $\Gamma(t_1,\dots,t_{r-1})$ generated by
vertices $t_1,\dots,t_{r-1}$ in $\Gamma(Y)$. We  formalize that
below.

Let
$$
\widetilde{D}=C_g+\sum_{v_i\in \{a_1,\dots,a_k\}-
\{t_1,\dots,t_{r-1}\}}{E_{v_i}}
$$
be the curve of components of $C$ which are not contracted by $\sigma$.
Then $\sigma:\widetilde{D}\to D$ is the normalization of $D$. In pictures,
we denote $\widetilde{D}$ (or $D$) by the symbol $\otimes$ and evidently
denote the intersection of this curve and its local branches at
the corresponding singular point with the components $E_{t_j}$ which
are contracted to this point.
For connected components of  $\Gamma(t_1,\dots,t_{r-1})$ we then
have possibilities presented in
Table 4 below depending on types of the corresponding singular points of $D$.

By Table 4, the ambiguity \eqref{3.6.1} or \eqref{3.6.2} takes place only
for singularities of the types $A_{2k-1}$ or $A_{2k}$. Thus, we have to
introduce  the notation $\gA_{2k-1}$ for the singularity of the type
$A_{2k-1}$ or $A_{2k}$ of the component
$\sigma(C_g)$ of $D$ of the geometric genus $g>1$.

In the right column of Table 4,
we denote by $\cA_n$, $\cD_n$, $\cE_n$ connected components of
graphs $\Gamma(t_1,\dots,t_{r-1})$ corresponding to singularities
$A_n$, $D_n$ and $E_n$ of the curve $D$ respectively.
Obviously, finding of all possible contractions $\sigma:Y\to \bP^2$
reduces to enumeration of all subgraphs $\Gamma\subset \Gamma(Y)$
with the connected components $\cA_n$, $\cD_n$, $\cE_n$ and with the common
number $r-1$ of vertices. A choice of such a subgraph
$\Gamma \subset \Gamma(Y)$ defines
the sextic $D$ with the corresponding irreducible components and simple
singularities, and the related configuration of rational curves
\begin{equation}
\sigma(E_v),\ v\in V(\Gamma(Y))-\left(\{a_1,\dots,a_k\}\cup
\{t_1,\dots,t_{r-1}\}\right),
\end{equation}
which one can call as {\it the exceptional curves of a sextic $D$ with
simple singularities.}

Thus, the classification in Theorems \ref{thm3.5.1} and \ref{thm3.5.2}
of DPN surfaces of elliptic type implies a quite delicate
classification of sextics $D$ having an irreducible component of the
geometric genus $g\ge 2$.  For this classification, we correspond to a
sextic $D\subset\bP^2$ a subgraph $\Gamma\subset \Gamma(Y)$ up to
isomorphisms of graphs $\Gamma(Y)$ which send the subgraphs $\Gamma$
to one another. The analogous classification can be repeated
to classify curves with simple singularities in $|-2K_{\bF_n}|$,
$n=0,\dots,4$. One should only replace $r-1$ by $r-2$.  We also note
that a choice of different subgraphs $\Gamma\subset \Gamma(Y)$ for the
same curve $C$ defines birational transformations of the corresponding
rational surfaces ($\bP^2$ or $\bF_n$) which transform the curves $D$
to one another. Thus, the graph $\Gamma(Y)$ itself classifies the
corresponding curves $D$ up to some their birational equivalence.

A complete enumeration of all cases has no principal difficulties, and
it is only related to a long enumeration using Theorem
\ref{thm3.5.2} of all possible diagrams $\Gamma(Y)$ and their
subdiagrams $\Gamma$.  Unfortunately, it seems, number of cases is
enormous. But the complete enumeration can be important in some
problems of real algebraic geometry and singularity theory. For
example, it could be important for classification of irreducible
quartics in $\bP^3$ with double rational singularities by the method
of projection from a singular point. To remove the ambiguity
\eqref{3.6.1} or \eqref{3.6.2}, one has to perform similar (to ours)
classification of generalized root invariants.


\noindent
{\bf Table 4.} {\it Correspondence between connected components of
$\Gamma$ and singularities of $D=\sigma(\widetilde{D})$.}

\begin{tabular}{|l|l|l|}
\hline
 Type of &  Equations of       &   Connected Components of $\Gamma$,\\
Singular &  the Singularity      &   and Curve $\widetilde{D}$
(denoted by  $\otimes$)\\
point of $D$& and its Branches     &   \\
\hline
                          &\ \ \ $y^2-x^{2k}=0$     &
\includegraphics[width=5.5cm]{pics/c-p140-1.eps}\\
$A_{2k-1}=$               & I:\ \ $y-x^k=0$      &   \\
${\gA}_{2k-1}$            &II:\ $y+x^k=0$      &   \\
\hline
$A_{2k}={\gA}_{2k-1}$ & $y^2-x^{2k+1}=0$      &
\includegraphics[width=5.5cm]{pics/c-p140-2.eps}\\
\hline
$D_{2k}$                  &\ \ \  $xy^2-x^{2k-1}=0$ &
\includegraphics[width=5.5cm]{pics/c-p140-3.eps}\\
                          & I:\ \ \ $x=0$         & \\
                          &II:\ \  $y-x^{k-1}=0$ &   \\
                          &III:\ $y+x^{k-1}=0$ &   \\
\hline
$D_{2k+1}$                &\ \ \ $xy^2-x^{2k}=0$ &
\includegraphics[width=5.5cm]{pics/c-p140-4.eps}  \\
                          & I:\ \ $x=0$         &   \\
                          &II:\  $y^2-x^{2k-1}=0$ &   \\
\hline $E_{6}$ & $y^3-x^4=0$      &
\includegraphics[width=3.5cm]{pics/c-p140-5.eps}\\
\hline
$E_{7}$                   &\ \ \ $y^3-yx^{3}=0$ &
\includegraphics[width=5.5cm]{pics/c-p140-6.eps}\\
                          & I:\ \ $y=0$         &   \\
                          &II:\  $y^2-x^{3}=0$  &   \\
\hline
$E_{8}$                   & $y^3-x^{5}=0$  &
\includegraphics[width=5.5cm]{pics/c-p140-7.eps}\\
\hline
\end{tabular}

\section{Classification of log del Pezzo surfaces of index $\le2$ and
  applications}
\label{sec4}

\subsection{Classification of log del Pezzo surfaces of index $\le 2$}
\label{subsec4.1} From results of Chapters
\ref{sec1} --- \ref{sec3} we obtain

\begin{theorem}\label{thm4.1.1} For any log del Pezzo
  surface $Z$ of index $\le 2$ there exists a unique resolution of
  singularities $\sigma: Y\to Z$ (it is called {\it right}) such
  that $Y$ is a right DPN surface of elliptic type and $\sigma$
  contracts exactly all exceptional curves of the Du Val and the
  logarithmic part of $\Gamma(Y)$. Vice versa, if $Y$ is a right
  DPN surface of elliptic type, then there exists a unique
  morphism $\sigma:Y\to Z$ of contraction of all exceptional curves
  corresponding to the Du Val and the logarithmic part of $\Gamma(Y)$
  which gives resolution of singularities of log del Pezzo surface $Z$
  of index $\le 2$ (it will be automatically the right resolution).

Thus, classifications of log del Pezzo surfaces of index $\le 2$ and
right DPN surfaces of elliptic type are equivalent, and they are given
by Theorems \ref{thm3.5.1}, \ref{thm3.5.2} and \ref{thm3.5.2.a}.
\end{theorem}

\begin{proof} Let $Z$ be a log del Pezzo surface of index $\le 2$.
  In Chapter \ref{sec1}, a ``canonical'' (i. e. uniquely defined)
  resolution of singularities $\sigma:Y\to Z$ had been suggested such
  that $Y$ is a right DPN surface of elliptic type. First, a
  minimal resolution of singularities $\sigma_1:Y^\prime \to Z$ is
  taken, and second, the blow-up of all intersection points of
  components of curves in preimages of non Du Val singularities $K_n$
  of $Z$ is taken. Let us show that $\sigma$ contracts exactly
  exceptional curves of $\Duv \Gamma(Y)$ and $\Log \Gamma(Y)$.

Let $E$ be an exceptional curve of $Y$ corresponding to a vertex of
$\Duv \Gamma(Y)$ or $\Log \Gamma(Y)$. Let $\widetilde{C}\in |-2K_Z|$ be
a non-singular curve of $Z$ which does not contain singular points of $Z$
(it does exist by Theorem \ref{thm1.4.1}),
and $C_g=\sigma^{-1}(\widetilde{C_g})$. Then (see Sections \ref{subsec1.5} and
\ref{sec2}) $C_g+E_1+\cdots+E_k\in |-2K_Y|$ where $E_i$ are all exceptional
curves on $Y$ with the square $(-4)$ and $C_g$ a non-singular irreducible
curve of genus $g\ge 2$. By Chapter \ref{sec2},
one has $E\cdot C_g=0$. If $\sigma$
does not contract $E$, then for the curve $\sigma(E)$ on $Z$ we have
$\sigma(E)\cdot \widetilde{C_g}=\sigma(E)\cdot (-2K_Z)=0$. Then
$-K_Z$ is not ample. We get a contradiction. Vice versa, by construction,
$\sigma$ contracts only curves from $\Duv \Gamma(Y)$ and $\Log \Gamma(Y)$.
This shows that $\sigma$ is a right resolution.

Now, let $Y$ be a right DPN surface of elliptic type and $\sigma:Y\to Z$
a contraction of all exceptional curves corresponding to vertices
$\Duv \Gamma(Y)$ and $\Log \Gamma(Y)$ (it does exist analytically because
$\Duv \Gamma(Y)\cup \Log \Gamma (Y)$ is negative, and we
show that it does exist algebraically by the direct
construction below).
To prove Theorem, we should prove that $Z$ is a log del Pezzo surface
of index $\le 2$, and $\sigma$ the right
resolution of singularities of $Z$.

Second statement becomes obvious if one decomposes $\sigma$ as the
composition of contractions of all exceptional curves of 1st kind from
\linebreak
$\Log \Gamma(Y)$ (they don't intersect each other) and further
contraction of the remaining exceptional curves.

To prove the first statement, one can use the double covering
$\pi:X\to Y$ with the involution $\theta$ (see Chapter \ref{sec2}), the
relation between exceptional curves of $Y$ and $(X,\theta)$ (see Chapter
\ref{sec2}), and that the contraction of $\bA$, $\bD$ and $\bE$
configurations of $(-2)$ curves on $X$ does exit and gives the
corresponding quotient singularities $\bC/G_i$ where $G_i\subset
SL(2,\bC)$ are finite subgroups.  Using these considerations (i. e.
first we consider the corresponding contraction, and second the
quotient by involution), and Brieskorn's results \cite{16}, we obtain
that all non Du Val singularities of $Z$ are $\bC/\widetilde{G_i}$
where $\widetilde{G_i}\subset GL(2,\bC)$ and $\widetilde{G_i}\cap
SL(2,\bC)=G_i$ have index 2 in $\widetilde{G_i}$, and
$\widetilde{G_i}/G_i=\{1,\theta\}$. It follows that $Z$ is a complete
algebraic surface with log-terminal singularities of index $\le 2$.

Let us show that $-K_Z$ is ample. By Nakai--Moishezon criterion \cite{Nak63},
\cite{Moi67}  (see also Kleiman's criterion \cite{Kleiman66}),
it is enough to
show that $(-K_Z)^2>0$ and $(-K_Z)\cdot D>0$ for any curve $D$ on $Z$.
We have (see Section \ref{subsec1.5})
$$
4(-K_Z)^2=(-2K_Z)^2=(\sigma^\ast(-2K_Z))^2=(C_g)^2>0.
$$
since $Y$ is a DPN surface of elliptic type. Moreover,
$$
-2K_Z\cdot D=-2\sigma^\ast K_Z\cdot \sigma^\ast D=C_g\cdot \sigma^\ast D\ge 0
$$
because $C_g$ is irreducible with $(C_g)^2>0$ and $\sigma^\ast D$ is
effective. Moreover, we get here zero, only if the effective divisor
$\sigma^\ast D$ consists of exceptional curves $F$ on $Y$ with
$C_g\cdot F=0$. But such curves $F$
correspond to vertices of the logarithmic or the Du Val part of $\Gamma (Y)$.
They are contracted by $\sigma$ into points of $Z$ which is impossible for
the divisor $\sigma^\ast D$.
\end{proof}

Using Theorem \ref{thm4.1.1}, we can transfer to log del Pezzo
surfaces $Z$ of index $\le 2$ {\it the main invariants $(r,a,\delta)$,
  equivalently $(k, g,\delta)$, the root invariant, the root
  subsystem, the exceptional curves} which are defined for the surface
$Y$ of the right resolution $\sigma:Y\to Z$. In particular, the
{\it Picard number of $Z$} is
\begin{equation}
\widetilde{r}=\rk \Pic Z=r-\# V(\Duv \Gamma(Y))-\# V(\Log \Gamma(Y)).
\label{4.1.1}
\end{equation}
In Theorem \ref{thm3.5.1} we have shown the Picard number $\widetilde{r}$
in the extremal (for $Y$) case. Obviously, surfaces $Z$ with the extremal
$Y$ are distinguished by the minimal Picard number $\widetilde{r}$ for
the fixed main invariants $(r,a,\delta)$ (equivalently, $(k, g, \delta)$).
Since the $\Log \Gamma(Y)$ is prescribed by the main invariants and is then
fixed, this is equivalent to have the maximal rank (i. e. the number of
$(-2)$-curves for the minimal resolution of singularities) for Du Val
singularities of $Z$.


In Mori Theory, see \cite{27} and \cite{32}, log del Pezzo surfaces
$Z$ with $\rk \Pic Z=1$ are especially important. They give relatively
minimal models in the class of rational surfaces with log-terminal
singularities: any rational surface $X\ne\bP^1\times\bP^1$
with log-terminal singularities
has a contraction morphism onto such a model.  From Theorems
\ref{thm4.1.1} and \ref{thm3.3.2}, we obtain classification of such
models with log-terminal singularities of index $\le 2$.  By Theorem
\ref{thm3.5.1}, they correspond to extremal DPN surfaces of elliptic
type with
$$
\widetilde{r}=r-\# V(\Duv \Gamma(Y))-\# V(\Log \Gamma(Y))=1,
$$
and Theorem \ref{thm3.5.1} gives the classification of the graphs of
exceptional curves on them. This classification can be extended to a
fine classification of the surfaces themselves. Here are results for
the case of $\rk\Pic Z=1$.

\begin{theorem}\label{thm4.1.2b}
  There exist, up to isomorphism, exactly 18 log del Pezzo surfaces
  $Z$ of index $2$ with $\rk\Pic Z=1$.  The DPN surfaces
  $Y$ of their right resolution of singularities are extremal and
  correspond to the following cases of Theorem \ref{thm3.5.1},
  where we show in parentheses the type of singularities of $Z$.
\begin{eqnarray*}
  &&11(K_1),\ 15(K_1A_4),\ 18(K_1A_1A_5),\ 19(K_1A_7),\ 20a(K_1D_8),\\
  &&20b(K_12A_1D_6),\ 20c(K_1A_3D_5),\ 20d(K_12D_4);\
  21(K_2A_2),\\
  &&25(2K_1A_7),\ 26(K_22A_3),\ 27(K_2A_7),\ 30(K_32A_2);\
  33(K_3A_1A_5);\\
  &&40(K_5),\,44(K_5A_4);\ 46(K_6A_2);\ 50(K_9).
\end{eqnarray*}
In particular, the isomorphism class of $Z$ is defined by its
configuration of singularities. The number of non-Du Val
singularities is at most one except when the singularities are
$2K_1A_7$.

In all other cases 11--50 the surface with maximal Du Val part is
also unique.
\end{theorem}
\begin{proof}
  For each of the graphs of Table 3 it is
  straightforward to pick a subgraph such that contracting the
  corresponding curves realizes $Y$ as a sequence of blowups starting
  from $V=\bP^2$ or $\bF_n$, $n\le4$. The images of the remaining curves
  give a configuration of curves on $V$.

  By theorem \ref{thm3.5.2.a} we are guaranteed that, vice versa,
  starting with such a configuration, the corresponding series of
  blowups leads to a right resolution of singularities $\wZ$ of a
  log del Pezzo surface $Z$ of index~$\le 2$.

  So, to compute the number of isomorphism classes, one has to find
  the orbits of the $G$-action on the parameter space for the choices
  of the blowups. Finally, one has to take into account the action of
  the symmetry group of the graph and the (finitely many) choices for
  the contractions to $V$.

  In all the cases this is a straightforward computation which gives
  precisely one orbit.

  A typical case is that of case 48. The configuration of curves can
  be contracted to a ruled surface $\bF_1$ so that the images of
  non-contracted curves are two distinct fibres, the exceptional
  section $s_1$ and an infinite section $s_{\infty} \sim s_1+f$. In
  other words, they are the $(\bC^*)^2$-invariant divisors on the
  toric variety $\bF_1$. The blowups $Y\to \bF_1$ are uniquely
  determined except for the two last blowups corresponding to the two
  white end-vertices. One easily sees that these two blowups
  correspond to a choice of two points $P_1,P_2$ lying on two torus
  orbits $O_1,O_2$ on a toric surface $Y'\to \bF_1$ with $\rk\Pic
  Y'=4$. The surface $Y'$ corresponds to a polytope obtained from the
  polytope of $\bF_1$ by cutting two corners, which adds two new
  sides. These sides are obviously not parallel. Hence, the torus
  $(\bC^*)^2$ acts transitively on $O_1\times O_2$, so the surface $Y$
  is unique.

  The only cases where a similar toric argument does not work are 39,
  45 and 47. In case 39 the surface $Y$ can be contracted to
  $\bP^1\times\bP^1$ with 6 curves, 3 sections and 3 fibres.  This
  configuration is unique and the blowups are uniquely determined, so
  the surface $Y$ is unique. In case 45 the surface $Y$ is similarly
  contracted to $\bP^1\times\bP^1$ with 6 curves, sections $s,s'$,
  fibres $f,f'$ and curves $C\sim C'\sim s+f$ so that $C$ passes
  through $s\cap f$ and $s'\cap f'$ and $C'$ through $s\cap f'$ and
  $s'\cap f$. This configuration is unique as well.

  In the most difficult case 47, $Y$ can be contracted to $\bP^2$ with
  the following configuration:
  \begin{enumerate}
  \item three non-collinear points $P_1,P_2,P_3$,
  \item three lines, $l_1,l_2,l_3$ each passing through two of the three
    points so that $P_i\not\in l_i$.
  \item two conics $q_1,q_2$ such that $q_1$ is tangent to lines $l_2$ and
    $l_3$ respectively at the points $P_3$ and $P_2$; and $q_2$ is tangent to
    lines $l_1$ and $l_3$ respectively at the points $P_3$ and $P_1$.
  \end{enumerate}
  It is easy to see that this configuration is rigid as well.
\end{proof}

The Gorenstein case is ``well known'' to experts but we were
unable to find a complete and accurate description of the isomorphism
classes in the literature. Therefore, we include the following theorem
for completeness. Here we use the {\it degree} $d$ of
a del Pezzo surface $Z$ which is $d=K_Z^2$.

\begin{theorem}\label{thm4.1.2a}
  \begin{enumerate}
  \item There exist 28 configurations of singularities of Gorenstein
    log del Pezzo surfaces of Picard rank 1, and each type
    determines the corresponding surface up to a deformation. The types
    (and the cases $N$ in Table 3) are as follows:
    \begin{enumerate}
    \item $d=9$: $\emptyset$(case 1)
    \item $d=8$:  $A_1$(2)
    \item $d=6$:  $A_2A_1$(5)
    \item $d=5$:  $A_4$(6)
    \item $d=4$:  $D_5$(7a) $A_3 2A_1$(7b)
    \item $d=3$:  $E_6$(8a), $A_5A_1$(8b), $3A_2$(8c)
    \item $d=2$:  $E_7$(9a), $A_7$(9b), $A_5A_2$(9c), $2A_3 A_1$(9d),
      $D_6 A_1$(9e), \newline $D_4 3A_1$(9f)
    \item $d=1$:  $E_8$(10a), $A_8$(10b), $A_7 A_1$(10c),
      $A_5 A_2 A_1$(10d), \newline
      $2A_4$(10e),
      $D_8$(10f), $D_5 A_3$(10g),
      $E_6 A_2$(10h), $E_7 A_1$(10i), \newline
      $D_6 2A_1$(10j),
      $2D_4$(10k), $2A_3 2A_1$(10l), $4A_2$(10m).
    \end{enumerate}
  \item In each type there is exactly one isomorphism class, with the
    following exceptions: in types $E_8$, $E_7A_1$, $E_6A_2$ there are
    two isomorphism classes; and in type $2D_4$ there are infinitely
    many isomorphism classes parameterized by $\bA^1$.
  \item The three extra surfaces of type $E_8$, $E_7A_1$, $E_6A_2$ and
    all surfaces of type $2D_4$ are distinguished by the fact that
    their automorphism groups are 1-dimensional and contain
    $\bC^*$. All other surfaces with $d=1$ have finite automorphism
    groups.
  \end{enumerate}
\end{theorem}
\begin{proof}[1st proof.]
  The first case to consider is $d=1$. Choosing an appropriate
  subgraph in the graph of exceptional curves on $\wZ$, one picks a
  sequence of blowups $\wZ\to \bP^2$.
  These contractions and images of $(-2)$-curves are listed in
  \cite{15}. In addition, one has to compute the images of
  $(-1)$-curves. The result is a configuration of lines, conics and
  cubics on $\bP^2$, and in most cases cubics can be avoided.

  Again, by theorem \ref{thm3.5.2.a} we are guaranteed that, vice versa,
  starting with such a configuration, the corresponding series of
  blowups leads to a minimal resolution of singularities $\wZ$ of a
  Gorenstein log del Pezzo surface $Z$.

  To compute the automorphism groups and the number of isomorphism
  classes, one has to compute the stabilizer $G\subset \PGL(3)$ of a
  projective configuration on $\bP^2$ and the orbits of the $G$-action
  on the parameter space for the configurations and the choices for the
  blowups; and to take into account the action of the symmetry group
  of the graph and the (finitely many) choices for the contractions to
  $\bP^2$.

  In the case $E_8$, the projective configuration is a line and a point on
  it, the group $G$ is the subgroup of upper-triangular matrices, and
  the parameter space is $\bC^*\times \bC^4$ which can be identified
  with the set of power series
  $$ y = \alpha_3 x^3 +\alpha_4x^4 + \alpha_5x^5 + \alpha_6x^6
  +\alpha_7x^7 \mod x^8
  \quad \text{with} \quad \alpha_3\ne 0 .
  $$
  The $G$-action has two orbits: those of $y=x^3$ and of
  $y=x^3+x^7$. The first orbit is in the closure of the second. The
  stabilizer of $y=x^3$ is isomorphic to $\bC^*$ and consists of
  diagonal matrices $(1,c,c^3)$, the second stabilizer is finite. The
  model for the moduli stack is $[\bA^1:\bG_m]$ with $\bC^*$-action
  $\lambda.a = \lambda^4 a$.

  In the case $E_7A_1$, the projective configuration is a line $l_1$, a
  conic $q$ tangent to it, and another line $l_2$. There are two
  cases: when $l_2$ intersects $q$ at 2 distinct points, and when they
  are tangent. One case is a degeneration of another, and in the
  degenerate case the stabilizer of the configuration contains
  $\bC^*$.

  In the case $E_6A_2$, the projective configuration consists of 4
  lines and 3 of them either pass through the same point or they do
  not. Once again, the local model is $[\bA^1:\bG_m]$ with the
  standard action, one configuration degenerates into another, and the
  degenerate configuration has stabilizer~$\bC^*$.

  In the case $2D_4$, the projective configuration consists of 4 lines
  through a point $P$ and the 5th line $l_5\not\ni P$. The parameter
  space for such configuration is $\bP^1\setminus$(3 points). Dividing
  by the symmetry group $\bZ/2\times\bZ/2$ gives $\bA^1$. Every
  configuration has $\bC^*$ as the stabilizer group.

  In all other cases for $d=1$ the computation gives one isomorphism
  class.

  For $d=2$, the surfaces $\wZ_2$ are obtained from the surfaces
  $\wZ_1$ by contracting one $(-1)$-curve. So, the cases where more
  than one isomorphism class is possible are the ones that come from
  the four exceptional cases above.

  The only contraction of the $E_8$-case is the case $E_7$. In this
  case, the group of upper triangular matrices acts on the polynomials
  $y=\alpha_3+\dots+\alpha_6x^6 \mod x^7$ with $\alpha_3\ne0$
  transitively; so there is only one isomorphism class.

  The case $E_7A_1$ produces $D_6A_1$ and $E_7$. In each of these, the
  surface is unique because it can also be obtained by contracting a
  surface of type $D_62A_1$ and $E_8$, respectively.

  The case $E_6A_2$ produces $A_5A_1$, which also comes from the type
  $A_5A_2A_1$ with a unique isomorphism class.  Similarly, the case
  $2D_4$ produces $D_43A_1$, which also comes from the type $D_62A_1$.
  For $d\ge 3$, moreover, there is only one isomorphism class for each
  configuration of singularities.
\end{proof}

\begin{proof}[2nd proof for the $d=1$ case.]
  By Theorem \ref{thm1.4.1} and Remark \ref{rem1.4.3}, a general
  element of the linear system $|-K_Z|$ is smooth.  By Riemann-Roch
  theorem, $h^0(-K_{Z})=2$. Hence, $|-K_Z|$ is a pencil with a unique,
  nonsingular base point $P$. The blowup of $Z$ at $P$ is an elliptic
  surface with a fibration $\pi:Z'\to\bP^1$ and a section. The
  condition $\rk\Pic Z=1$ implies that the minimal resolution of
  singularities $\wZ'\to \bP^1$ is an \emph{extremal} rational
  elliptic surface, as defined in \cite{MirandaPersson}.

  Vice versa, given an extremal relatively minimal (with no
  $(-1)$-curves in fibres of $\pi$) surface $\wZ'$ with a section,
  contracting the $(-2)$-curves not meeting the section and then the
  section gives a Gorenstein del Pezzo surface with Du Val
  singularities and $\rk\Pic Z=1$. The finitely many choices of a
  section differ by the action of the Mordell-Weil group of the
  elliptic fibration, and hence give isomorphic $Z$'s.

  Hence, the classification of Gorenstein log del Pezzo surfaces of
  degree 1 and rank 1 is equivalent to the classification of extremal
  rational elliptic fibrations with a section. The latter was done by
  Miranda and Persson in \cite{MirandaPersson}, and we just need to
  translate it to del Pezzo surfaces.

  On the level of graphs of exceptional curves, the transition from
  $\wZ$ to $\wZ'$ consists of inserting an extra $(-1)$-curve and
  changing $(-1)$-curves through $P$ to $(-2)$-curves. The graphs
  $A_n, D_n, E_n$ turn into the corresponding extended Dynkin graphs
  $\wA_n, \wD_n, \wE_n$. In addition, $A_1$ and $A_2$ can turn into
  graphs $*\wA_1$, $*\wA_2$ respectively. The elliptic fibres $\wA_0$
  and $*\wA_0$ are not seen in the graphs of $Z$.

  According to \cite[Thm 4.1]{MirandaPersson} there are 16 types of
  elliptic fibrations. The four special types $*\wA_0\wE_8$,
  $*\wA_1\wE_7$, $*\wA_2\wE_6$ and $2\wD_4$ are distinguished by the
  fact that the induced modular $j$-function $j:\bP^1\to \bP^1_j$ is
  constant and there are exactly two singular fibres.

  The subgroup $\Aut_j Y$ of automorphisms commuting with $j$ is always
  finite. Hence, in the four exceptional cases $\Aut Y$
  has dimension one and contains $\bC^*$.  In all other cases
  $j$-map is surjective, and hence the automorphism group is finite.

  By \cite[Thm 5.4]{MirandaPersson}, in fifteen of the sixteen cases the
  elliptic surface is unique. In the case $2\wD_4$, there are
  infinitely many isomorphism classes, one for each value
  $j\in \bA^1_j$.
\end{proof}


\medskip

If we consider log del Pezzo surfaces $Z$ of index $\le 2$ and without
Du Val singularities, we get an opposite case to the previous one.
Thus, any singularity of $Z$ must have index $2$.  This case includes
and is surprisingly similar to the classical case of non-singular del
Pezzo surfaces when there are no singularities at all.  Applying
Theorem \ref{thm4.1.1} we get the following

\begin{theorem}\label{thm4.1.4}
Up to deformation, there exist exactly 50 types of log del Pezzo surfaces $Z$
with singularities of index exactly $2$ (if a singularity
does exist). The DPN surfaces $Y$ of their right resolution
of singularities have empty Du Val part $\Duv \Gamma(Y)$,
zero root invariants, and are defined by their main invariants
$(r,a,\delta)$ (equivalently $(k,g,\delta)$), up to deformation.
The diagram $\Gamma(Y)$ can be obtained from the
diagram $\Gamma$ of cases 1 --- 50 of Table 3 (with the same main invariants)
as follows: $\Gamma(Y)$ consists of $\Log \Gamma(Y)=\Log \Gamma$ and
\begin{equation}\label{4.1.4}
\Var \Gamma(Y)=W(\Var \Gamma)
\end{equation}
where $W$ is generated by reflections in all vertices of $\Duv \Gamma$
(i. e. one should take $D=\emptyset$ in Theorem \ref{thm3.5.2}).
In cases 7, 8, 9, 10, 20 one can consider only diagrams $\Gamma$ of
cases 7a, 8a, 9a, 10a and 20a (diagrams 7b, 8b,c, 9b---f,
10b---m, 20b---d give the same).

The type of Dynkin diagram $\Duv \Gamma$ can be considered as analogous
to the type of root system which one usually associates to non-singular
del Pezzo surfaces. Its actual meaning is to give the type of
the Weyl group $W$ describing the varying part $\Var (\Gamma (Y))$
by \eqref{4.1.4}. In cases 7 --- 10, 20, one should (or can) take graphs
$\Gamma$ of cases 7a---10a, 20a.
\end{theorem}

\begin{proof} This case corresponds to $Y$ with empty $D\subset \Duv \Gamma$
of Theorem \ref{thm3.5.2}. Then the root invariant is 0. Thus, all
cases 7, 8, 9, 10 or 20 give the isomorphic root invariants and the
same diagrams, and we can consider only the corresponding cases
7a, 8a, 9a, 10a and 20a to calculate the diagrams.

Let us show that moduli spaces of DPN surfaces $Y$ with the
same main invariants $(r,a,\delta)$ and zero root invariant (i. e. $D=0$)
are connected.

It is enough to show connectedness of the moduli space of
the corresponding right DPN
pairs $(Y,C)$ where $C\in |-2K_Y|$ is non-singular. Taking double
covering $\pi:X\to Y$ ramified in $C$, it is enough to prove
connectedness of moduli $\cM_{(r,a,\delta)}$ of K3 surfaces with
non-symplectic involutions $(X,\theta)$ and $(S_X)_+=S$ where $S$
has invariants $(r,a,\delta)$. General such pairs have zero root
invariant, as we want, because general $(X,\theta)$ have $S_X=S$.
Connectedness of $\cM_{(r,a,\delta)}$ had been discussed in Section
\ref{subsec2.2} with proofs in Appendix: Section
\ref{subsec:maininv}.
\end{proof}

We remark that the result equivalent to Theorem \ref{thm4.1.4} was
first obtained in \cite{10}.

We remark that cases 1 --- 10 of Theorem \ref{thm4.1.4} give classical
non-singular del Pezzo surfaces. Therefore, Theorem \ref{thm4.1.4} and
all results of this work show that log del Pezzo surfaces of index
$\le 2$ are very similar to classical non-singular del Pezzo surfaces.

\subsection{An example: Enumeration of all possible types for $N=20$.}
\label{subsecN=20en} Let us consider enumeration of all types of
singularities and graphs of exceptional curves of log del Pezzo
surfaces of index 2 of type $N=20$, i. e. with the main invariants
$(r,a,\delta)=(10,8,1)$.

From Theorem \ref{thm4.1.1} and Table 3, cases 20a---d, we obtain
that all of them have one singularity $K_1$ of index 2 and Du Val
singularities which correspond to a subgraph of one of graphs
$D_8$, $D_62A_1$, $D_5A_3$ and $2D_4$. It follows that their Du
Val singularities are exactly of one of 52 types listed below:

\medskip

 \noindent $2D_4$;

\noindent $D_8$;

\noindent $D_7$;

\noindent $D_62A_1$, $D_6A_1$, $D_6$;

\noindent $D_5A_3$, $D_5 A_2$, $D_5 2A_1$, $D_5 A_1$, $D_5$;

\noindent $D_4A_3$, $D_4A_2$, $D_4 3A_1$, $D_4 2A_1$, $D_4 A_1$,
$D_4$;

\noindent $A_7$;

\noindent $A_6$;

\noindent $A_52A_1$, $A_5A_1$, $A_5$;

\noindent $A_4A_3$, $A_4A_2$, $A_42A_1$, $A_4A_1$, $A_4$;

\noindent $2A_3A_1$, $2A_3$, $A_3A_22A_1$, $A_3A_2A_1$, $A_3A_2$,
$A_34A_1$, $A_33A_1$, $A_32A_1$, $A_3A_1$, $A_3$;

\noindent $2A_22A_1$, $2A_2A_1$, $2A_2$, $A_24A_1$, $A_23A_1$,
$A_22A_1$, $A_2A_1$, $A_2$;

\noindent $6A_1$, $5A_1$, $4A_1$, $3A_1$; $2A_1$, $A_1$;

\noindent $\emptyset$.

Using calculations of root invariants of Lemma \ref{lemma3.4.4},
it is easy to calculate the root invariant for any of the
subgraphs. One can see that it is defined uniquely by the type of
Du Val singularities except the following 15 types of Du Val parts
of singularities for which we show all differences in their root
invariants.

$D_4A_3$: There are exactly {\it two possibilities} for the root
invariant (and then for the graph of exceptional curves). The
first one can be obtained by taking $D_4A_3$ as a subdiagram in
$D_8$ (case 20a), and the second by taking $D_4A_3$ as a
subdiagram in $D_5A_3$ (case 20c). In the second case, the
characteristic element can be written using elements of the
component $A_3$, and this is impossible in the first case.

$D_42A_1$: There are exactly {\it two possibilities} for the root
invariant. The first one can be obtained by taking $D_42A_1$ as a
subdiagram in $D_8$ (case 20a), and the second one by taking
$D_42A_1$ as a subdiagram in $D_62A_1$ (case 20b). In the second
case, the characteristic element can be written using elements of
the components $2A_1$, and it is impossible in the first case.

$A_7$: There are exactly {\it two possibilities} for the root
invariant. The group $H=\{0\}$ or $H\cong \bZ/2$. Both cases can
be obtained by taking subdiagrams in $D_8$ (case 20a).

$A_5A_1$: There are exactly {\it two possibilities} for the root
invariant: the group $H=\{0\}$ or $H\cong \bZ/2$. Both cases can
be obtained by taking subdiagrams in $D_8$ (case 20a).

$2A_3$: There are exactly {\it four possibilities} for the root
invariant. For the group $H=\{0\}$ the characteristic element can
be written using elements either of one component $A_3$ or only by
both components $A_3$. For the group $H\cong \bZ/2$ either
$\alpha=1$ or $\alpha=0$. Three of these cases can be obtained by
taking subdiagrams in $D_8$ (case 20a). The remaining case $H\cong
\bZ/2$ and $\alpha=0$ can be obtained by taking a subdiagram in
$D_5A_3$ (case 20c).

$A_3A_2$: There are exactly {\it two possibilities} for the root
invariant: $\alpha=0$ or $\alpha=1$. Both cases can be obtained by
taking subdiagrams in $D_8$ (case 20a).

$A_33A_1$: There are exactly {\it two possibilities} for the root
invariant: In the first case the characteristic element cannot be
written using elements of the components $A_3$ (it can be obtained
by taking a subdiagram in $D_8$, i. e. for the case 20a). For the
second case it can be written using elements of the component
$A_3$ (it can be obtained by taking a subdiagram in $D_62A_1$, i.
e. for the case 20b).

$A_32A_1$: There are exactly {\it five possibilities} for the root
invariant. For the group $H=\{0\}$ the characteristic element can
be written using elements either of one component $A_3$, or by
components $2A_1$, or using all three components $A_32A_1$. For
the group $H\cong \bZ/2$ either $\alpha=1$ or $\alpha=0$. Four of
these cases can be obtained by taking subdiagrams in $D_8$ (case
20a). The remaining case $H\cong \bZ/2$ and $\alpha=0$ can be
obtained considering a subdiagram in $D_62A_1$ (case 20b).

$A_3A_1$: There are exactly {\it two possibilities} for the root
invariant: $\alpha=0$ or $\alpha=1$. Both cases can be obtained by
taking subdiagrams in $D_8$ (case 20a).

$A_3$: There are exactly {\it two possibilities} for the root
invariant: $\alpha=0$ or $\alpha=1$. Both cases can be obtained by
taking subdiagrams in $D_8$ (case 20a).

$A_22A_1$: There are exactly {\it two possibilities} for the root
invariant: $\alpha=0$ or $\alpha=1$. Both cases can be obtained by
taking subdiagrams in $D_8$ (case 20a).

$5A_1$: There are exactly {\it two possibilities} for the root
invariant. For the first one the characteristic element can be
written using two pairs of components of $5A_1$.  For the second
one the characteristic element can be written using only one pair
of components of $5A_1$.

$4A_1$: There are exactly {\it four possibilities} for the root
invariant. For the group $H=\{0\}$ the characteristic element can
be written using elements either of two components $A_1$ or by
only four components $A_1$. For the group $H\cong \bZ/2$ either
$\alpha=1$ or $\alpha=0$. All four cases can be obtained by taking
subdiagrams in $2D_4$ (case 20c).

$3A_1$: There are exactly {\it two possibilities} for the root
invariant: $\alpha=0$ or $\alpha=1$. Both cases can be obtained by
considering subdiagrams in $D_8$ (case 20a).

$2A_1$: There are exactly {\it two possibilities} for the root
invariant: $\alpha=0$ or $\alpha=1$. Both cases can be obtained by
considering subdiagrams in $D_8$ (case 20a).

\medskip

Thus, for the types of Du Val singularities shown above (together
with the singularity $K_1$ of index two) we obtain the number shown above
of different types of log del Pezzo surfaces: their right
resolution of singularities can have that number of different
graphs of exceptional curves. By taking the corresponding sequence
of contractions of $-1$ curves, one can further investigate these
surfaces in details; in particular, one can enumerate connected
components of their moduli.

Thus, there are exactly $52+12+3\cdot 2+4=74$ different graphs of
exceptional curves on the right resolution of singularities of log
del Pezzo surfaces of index 2 with the main invariants
$(r,a,\delta)=(10, 8, 1)$ (i. e. N=20).

Of course, similar calculations can be done for all 50 types of
main invariants of log del Pezzo surfaces of index $\le 2$. The
considered case $N=20$ is one of the richest and most complicated.

\subsection{Application: Minimal projective compactifications of
af\-fine surfaces in $\bP^2$ by relatively minimal
  log del Pezzo surfaces of index $\le 2$.}
\label{subsec4.2} This is similar to \cite{15} in the Gorenstein case.

Let us consider one of the 45 relatively minimal surfaces of Theorem
\ref{thm4.1.2b}, \ref{thm4.1.2a} which are different from $\bP^2$ (i. e.
except the case 1). Let $\sigma:Y\to Z$ be its right resolution of
singularities, and $v_1,\dots,v_{r-1}$ a sequence of vertices of
$\Gamma (Y)$ such that the corresponding exceptional curves on $Y$
give a contraction of the sequence of curves of the 1st kind
$\tau:Y\to \bP^2$. Let $C\subset \bP^2$ be the union of images by
$\tau$ of all exceptional curves $E_v$, $v\in V(\Duv(\Gamma(Y)))\cup
V(\Log(\Gamma(Y)))$.  Then, the embedding $f=(\sigma\tau^{-1}):W\to Z$
gives a compactification of the affine surface $W=\bP^2-C$ of $\bP^2$.
The morphism $f$ is minimal in the sense that $f$ cannot be extended
through components of $C$ (see \cite{15} for details).  The
description of all such affine surfaces $W$ and such their
compactifications is then reduced to the description of subdiagrams of
$\Gamma(Y)$ (defined by $v_1,\dots,v_{r-1})$) which were described by
their connected components in Section \ref{subsec3.6}.

\subsection{Dimension of the moduli space}\label{subsec4.3.1}

For each triple of invariants
\begin{equation}
\left(k,\,g,\,\delta,\ \ \text{the root invariant}\right),
\end{equation}
equivalently,
\begin{equation}
\left(r,\,a,\,\delta,\ \
\text{the dual diagram of exceptional curves }
 \Gamma(Y)\right)
\end{equation}
one has the moduli space of pairs $\cM_{(Z,C)}$ of log del Pezzo
surfaces together with a smooth curve $C\in |-2K_Z|$. We have
established the equivalence between pairs $(Z,C)$ and K3 surfaces
$(X,\theta)$ with a non-symplectic involution. Hence instead of
moduli $\cM_{(Z,C)}$ of pairs $(Z,C)$ we can consider moduli
$\cM_{(X,\theta)}$ of pairs $(X,\theta)$.

By \eqref{dimmoduliSrootinv},
\begin{equation}
 \dim \cM_{(X,\theta)}=
 20-r-\#V(\Duv(\Gamma))=9+g-k-\#V(\Duv(\Gamma)).
\end{equation}
Moreover,
\begin{equation}
\dim|-2K_Z|=\dim|C_g|=3g-3
\end{equation}
(see Sections \ref{subsec1.4} and \ref{subsec1.5}). It follows that
the dimension of the parameter space $\cM_{(r,a,\delta),\Gamma(Y)}$
 of \emph{generic} surfaces $Z$ of type
$(r,a,\delta)$ and with the graph $\Gamma (Y)$ of exceptional
curves on the right resolution $Y$ of singularities (or with the
corresponding root invariant) is equal to

\begin{equation}
\begin{split}
\dim \cM_{(r,a,\delta),\Gamma(Y)} =
&12-2g-k-\#V(\Duv(\Gamma(Y))+\dim \Aut Z=\\
&\frac{r+3a}{2}-10-\#V(\Duv(\Gamma(Y))+ \dim \Aut Z
\end{split}
\label{4.3.1}
\end{equation}

Note that this formula may fail for non-generic surfaces. For example,
by Theorem \ref{thm4.1.2a} there are exactly two isomorphism classes
of Gorenstein surfaces with a single $E_8$-singularity. The formula
above gives
\begin{eqnarray*}
  \dim \cM_{(r,a,\delta),\Gamma(Y)} = \dim \Aut Z
\end{eqnarray*}
which is true for the generic surface that has trivial isomorphism
group and fails for the second surface which has $\Aut Z=\bC^*$.

\subsection{Some open questions} \label{subsec4.4}
\subsubsection{Finite characteristic} \label{subsubsec4.3.3}
It would be very interesting to generalize results of this work to
finite characteristic. As we had mentioned in Remark
\ref{rem3.4.7}, it seems, the main problem is to generalize
Theorem \ref{thm1.4.1}. We think that our results are valid in
characteristic $\ge 3$. As we have seen, in characteristic $2$ the
number of cases increases.

\subsubsection{Arithmetic of log del Pezzo surfaces of index $\le 2$}
\label{subsubsec4.3.4} There are many results (e. g.
see \cite{6}, \cite{7} and \cite{17}) where the arithmetic of
classical non-singular del Pezzo surfaces is studied. What is the
arithmetic of log del Pezzo surfaces of index $\le 2$ and
equivalent DPN surfaces of elliptic type?

\section{Appendix} \label{sec:appendix}
Here we add some important results and calculations which had been
used in the main part of the work and which are important by themselves.

\subsection{Integral symmetric bilinear forms. Elements of the discriminant forms
technique} \label{subsec:discrforms} Here, for readers'
convenience, we review results about integral symmetric bilinear
forms (lattices) which we used. We follow \cite{9}.

\subsubsection{Lattices}
\label{subsubsec:lattices} Everywhere in the sequel, by a {\it
lattice} we mean a free  $\bZ$-module of finite rank, with a
nondegenerate symmetric bilinear form with values in the ring
$\bZ$ of rational integers (thus, ``lattice'' replaces the phrase
``nondegenerate integral symmetric bilinear form'').

A lattice $M$ is called {\it even } if $x^2=x\cdot x$ is even for
each $x\in S$, and {\it odd} otherwise (here we denote by $x\cdot
y$ the value of the bilinear form of $M$ at the pair $(x,y)$). By
$M_1\oplus M_2$ we denote the orthogonal direct sum of lattices
$M_1$ and $M_2$. If $M$ is a lattice, we denote by $M(a)$ the
lattice obtained from $M$ by multiplying the form of $M$ by the
rational number $a\not=0$, assuming that $M(a)$ is also integral.

\subsubsection{Finite symmetric bilinear and quadratic forms}
\label{subsubsec:finiteforms} By a {\it finite symmetric bilinear
form} we mean a symmetric bilinear form $b:{\gA}\times {\gA}\to
\bQ/\bZ$ defined on a finite Abelian group $\gA$.

By a {\it finite quadratic form} we mean a map $q:\gA\to \bQ/2\bZ$
satisfying the following conditions:

1) $q(na)=n^2q(a)$ for all $n\in \bZ$ and $a\in \gA$.

2)$q(a+a^\prime)-q(a)-q(a^\prime)\equiv 2b(a,a^\prime)\pmod{2}$,
where $b:\gA\times \gA\to \bQ/\bZ$ is a finite symmetric bilinear
form, which we call {\it the bilinear form} of~$q$.

A finite quadratic form $q$ is nondegenerate when $b$ is
nondegenerate. In the usual way, we introduce the notion of
orthogonality ($\perp$) and of orthogonal sum ($\oplus$) of finite
symmetric bilinear and quadratic forms.

\subsubsection{The discriminant form of a lattice}
\label{subsubsec:discforms} The bilinear form of a lattice $M$
determines a canonical embedding $M\subset
M^\ast=\text{Hom}(M,\bZ)$. The factor group $\gA_M=M^\ast /M$ is
finite and Abelian, and its order is equal to $|det( M)|$. We
remind that the determinant $det(M)$ of $M$ equals
$det(e_i\cdot e_j)$ for some basis $e_i$ of the lattice $M$.  A
lattice $L$ is called {\it unimodular} if $det(L)=\pm 1$.

We extend the bilinear form of $M$ to one on $M^\ast$, taking values in $\bQ$. We put
$$
b_M(t_1+M,t_2+M)=t_1\cdot  t_2+\bZ
$$
where $t_1,t_2\in M^\ast$, and
$$
q_M(t+M)=t^2+2\bZ,
$$
if $M$ is even, where $t\in M^\ast$.

We obtain {\it the discriminant bilinear form }
$b_M:{\gA}_M\times {\gA}_M\to \bQ/\bZ$ and
{\it the discriminant quadratic form} $q_M:{\gA}_M\to \bQ/2\bZ$
(if $M$ is even) of the lattice $M$.
They are nondegenerate.

Similarly, one can define discriminant forms of $p$-adic lattices
over the ring $\bZ_p$
of $p$-adic integers for a prime $p$.
The decomposition of ${\gA}_M$ as a sum of its $p$-components
$({\gA}_M)_p$ defines
the decomposition  of $b_M$ and $q_M$ as the orthogonal sum of its
$p$-components
$(b_M)_p$ and $(q_M)_p$. They are equal to the discriminant forms of
the corresponding
$p$-adic lattices $M\otimes \bZ_p$.

We denote by:

$K_\theta^{(p)}(p^k)$ the 1-dimensional $p$-adic lattice determined by the
matrix $\langle \theta p^k\rangle$,
where $k\ge 1$ and $\theta\in \bZ^\ast_p$ (taken
$\mod (\bZ_p^\ast)^2$);

$U^{(2)}(2^k)$ and $V^{(2)}(2^k)$ the 2-dimensional $2$-adic lattices
determined by the matrices
$$
\left(
\begin{array}{cc}
0 & 2^k\\
2^k & 0
\end{array}\right),\ \ \
\left(
\begin{array}{cc}
2^{k+1} & 2^k\\
2^k & 2^{k+1}
\end{array}\right)
$$
respectively;

$q_\theta^{(p)}(p^k)$,  $u_+^{(2)}(2^k)$ and $v_+^{(2)}(2^k)$, the
discriminant quadratic forms of $K_\theta^{(p)}(p^k)$,
$U^{(2)}(2^k)$ and $V^{(2)}(2^k)$ respectively;

$b_\theta^{(p)}(p^k)$,  $u_-^{(2)}(2^k)$ and $v_-^{(2)}(2^k)$,
the bilinear forms of
$q_\theta^{(p)}(p^k)$,  $u_+^{(2)}(2^k)$ and $v_+^{(2)}(2^k)$ respectively.

These $p$-adic lattices and finite quadratic and bilinear forms
are called {\it elementary}.  Any $p$-adic lattice (respectively
finite non-degene\-rate quadratic, bilinear form) is an orthogonal
sum of elementary ones {\it (Jordan decomposition).}

 \subsubsection{Existence of an even lattice with a given discriminant
quadratic form}
The signature of a lattice $M$ is equal to
$sign\,M=t_{(+)}-t_{(-)}$ where $t_{(+)}$ and $t_{(-)}$ are numbers
of positive
and negative squares of the corresponding real form $M\otimes \bR$.
The formula
$$
sign\,q_M \mod 8=sign\, M\mod 8=t_{(+)}-t_{(-)}\mod 8
$$
where $M$ is an even lattice, correctly defines the signature
$\mod 8$ for non-degenerate finite quadratic forms. For elementary
finite quadratic forms we obtain respectively
$$
sign\,q_\theta^{(p)}(p^k)\equiv k^2(1-p)+4k\eta\mod 8
$$
where $p$ is odd and
$\left(\frac{\,\theta\,}{\,p\,}\right)=(-1)^\eta$, where we use
Legendre symbol;
$$
sign\,q_\theta^{(2)}(2^k)\equiv\theta+4\omega(\theta)k\mod 8
$$
where $\omega(\theta)\equiv (\theta^2-1)/8\mod 2$;
$$
sign\,v_+^{(2)}(2^k)\equiv 4k\mod 8;
$$
$$
sign\,u_+^{(2)}(2^k)\equiv 0\mod 8.
$$
In particular, $sign\,L\equiv 0\mod 8$ if $L$ is an even unimodular lattice.

We denote by $l(\gA)$ the minimal number of generators of a finite
Abelian group $\gA$. We consider an even lattice $M$ with the
invariants $(t_{(+)}, t_{(-)}, q)$ where $t_{(+)}$, $t_{(-)}$ are
its numbers of positive and negative squares, and $q\cong q_M$; we
denote by ${\gA}_q$ the group where $q$ is defined. {\it The
invariants $(t_{(+)}, t_{(-)}, q)$ are equivalent to the {\it
genus} of $M$} (see Corollary 1.9.4 in \cite{9}) .
 Thus, they define the isomorphism classes of the $p$-adic lattices
$M\otimes \bZ_p$ for all prime $p$, and of $M\otimes \bR$.

We have (see Theorem 1.10.1 in \cite{9}):

\begin{theorem} An even lattice with invariants $(t_{(+)},t_{(-)},q)$
exists if and only if the following conditions are simultaneously
satisfied:

1) $t_{(+)}-t_{(-)}\equiv sign\,q\mod 8$.

2) $t_{(+)}\ge 0$, $t_{(-)}\ge 0$, $t_{(+)}+t_{(-)}\ge l({\gA}_q)$.

3) $(-1)^{t_{(-)}}|{\gA}_q|\equiv det( K(q_p))\mod (\bZ_p^\ast)^2$
for all odd primes $p$ for which $t_{(+)}+t_{(-)}=l({\gA}_{q_p})$
(here $K(q_p)$ is the unique $p$-adic lattice with the
discriminant quadratic form $q_p$ and the rank $l({\gA}_{q_p})$).

4) $|{\gA}_q|\equiv \pm det(K(q_2))\mod (\bZ_2^\ast)^2$ if
$t_{(+)}+t_{(-)}=l({\gA}_{q_2})$ and
$q_2\not=q_\theta^{(2)}(2)\oplus q_2^\prime$ (here $K(q_2)$ is the
unique $2$-adic lattice with the discriminant quadratic form $q_2$
and the rank $l({\gA}_{q_2})$). \label{discrexistence}
\end{theorem}

 From $l({\gA}_q)=\max_{p}{l({\gA}_{q_p})}$, we obtain the important corollary.

\begin{corollary} An even lattice with invariants
$(t_{(+)},t_{(-)},q)$ exists if the following conditions are
simultaneously satisfied:

1) $t_{(+)}-t_{(-)}\equiv sign\,q\mod 8$.

2) $t_{(+)}\ge 0$, $t_{(-)}\ge 0$, $t_{(+)}+t_{(-)}> l({\gA}_q)$.
\label{discrexistence2}
\end{corollary}

Theorem 1.16.5 and Corollary 1.16.6 in \cite{9} give similar results
for odd lattices.

\subsubsection{Primitive embeddings into even unimodular lattices}
\label{subsubsec:primemb} We have a simple statement (see
Proposition 1.4.1 in \cite{9}).

\begin{proposition}
Let $M$ be an even lattice. Its overlattice $M\subset N$ of  finite
index is equivalent to the isotropic subgroup $H=N/M\subset {\gA}_M$
with respect to
$q_M$. Moreover, we have $q_N=q_M|(H^\perp)/H$.
\label{overlattice}
\end{proposition}

An embedding of lattices $M\subset L$ is primitive if $L/M$ is a
free $\bZ$-module.

Let $L$ be an even unimodular lattice, $M\subset L$ its primitive sublattice,
and $T=M^\perp_L$. Then $M\oplus T\subset L$ is an overlattice of a
finite index.
Applying Proposition \ref{overlattice}, we obtain that $H=L/(M\oplus T)$ is
the graph of an isomorphism $\gamma:q_M\cong -q_T$, and this is equivalent to a
primitive embedding $M\subset L$ into an even unimodular lattice with
$T=M^\perp_L$. Thus we have (see Proposition 1.6.1 in \cite{9})

\begin{proposition}
A primitive embedding of an even lattice $M$ into an even unimodular lattice,
in which the orthogonal complement is isomorphic to $T$, is
determined by an isomorphism $\gamma:q_M\cong
-q_T$.

Two such isomorphisms $\gamma$ and $\gamma^\prime$ determine
isomorphic primitive embeddings if and only if they are conjugate
via an automorphism of $T$. \label{primembedd1}
\end{proposition}

 From Theorem  \ref{discrexistence}  and Corollary \ref{discrexistence2} we
then obtain (Theorem 1.12.2 and Corollary 1.12.3 in \cite{9})

\begin{theorem}
The following properties are equivalent:

a) There exists a primitive embedding of an even lattice $M$ with invariants
$(t_{(+)},t_{(-)}, q)$
into some even unimodular lattice of signature $(l_{(+)},l_{(-)})$.

b) There exists an even lattice with invariants
$(l_{( +)}-t_{(+)},l_{(-)}-t_{(-)},-q)$

c) There exists an even lattice with invariants
$(l_{(-)}-t_{(-)},l_{(+)}-t_{(+)},q)$.

d) The following conditions are simultaneously satisfied:

1) $l_{(+)}-l_{(-)}\equiv 0\mod 8$.

2) $l_{(+)}-t_{(+)}\ge 0$, $l_{(-)}-t_{(-)}\ge 0$,
$l_{(+)}+l_{(-)}-t_{(+)}-t_{(-)}\ge l({\gA}_q)$.

3) $(-1)^{l_{(+)}-t_{(+)}}|{\gA}_q|\equiv det( K(q_p))\mod
(\bZ_p^\ast)^2$ for all odd primes $p$ for which
$l_{(+)}+l_{(-)}-t_{(+)}-t_{(-)}=l({\gA}_{q_p})$ (here $K(q_p)$ is
the unique $p$-adic lattice with the discriminant quadratic form
$q_p$ and the rank $l({\gA}_{q_p})$).

4) $|{\gA}_q|\equiv \pm det(K(q_2))\mod (\bZ_2^\ast)^2$ if
$l_{(+)}+l_{(-)}-t_{(+)}-t_{(-)}=l({\gA}_{q_2})$ and
$q_2\not=q_\theta^{(2)}(2)\oplus q_2^\prime$ (here $K(q_2)$ is the
unique $2$-adic lattice with the discriminant quadratic form $q_2$
and the rank $l({\gA}_{q_2})$). \label{primembedd2}
\end{theorem}

\begin{corollary} There exists a primitive embedding of an even lattice $M$
with invariants
$(t_{(+)},t_{(-)}, q)$ into some even unimodular lattice of signature
$(l_{(+)},l_{(-)})$ if the
following conditions are simultaneously satisfied:

1) $l_{(+)}-l_{(-)}\equiv 0\mod 8$.

2) $l_{(-)}-t_{(-)}\ge 0$, $l_{(+)}-t_{(+)}\ge 0$,
$l_{(+)}+l_{(-)}-t_{(+)}-t_{(-)}> l({\gA}_q)$.
\label{primembedd3}
\end{corollary}

It is well-known that an even unimodular lattice of signature
\linebreak $(l_{(+)},l_{(-)})$  is unique if it is indefinite (e.
g. see \cite{Serre}). The same is valid if $l_{(+)}+l_{(-)}\le 8$.
Thus, Theorem \ref{primembedd2} and Corollary \ref{primembedd3}
give existence of a primitive embedding of $M$ into these
unimodular lattices.

\subsubsection{Uniqueness}
\label{subsubsec:uniqueness} We restrict ourselves to the
following uniqueness result (see Theorem 1.14.2 in \cite{9} and
Theorem 1.2$^\prime$ in \cite{8}). We note that this is based on
fundamental results about spinor genus of indefinite lattices of
the rank $\ge 3$ due to M. Eichler and M. Kneser.

\begin{theorem} Let $T$ be an even indefinite lattice with the
invariants $(t_{(+)}, t_{(-)}, q)$
satisfying the following conditions:

a) $\rk T\ge l({\gA}_{q_p})+2$ for all $p\not=2$.

b) If $\rk T=l({\gA}_{q_2})$, then $q_2\cong u_+^{(2)}(2)\oplus q^\prime$ or
$q_2\cong v_+^{(2)}(2)\oplus q^\prime$.

Then the lattice $T$ is unique (up to isomorphisms), and
the canonical homomorphism
$O(T)\to O(q_T)$ is surjective.
\label{uniqueness1}
\end{theorem}

Applying additionally  Proposition  \ref{primembedd1}, we obtain
the following Analogue of Witt's Theorem for
primitive embeddings into even unimodular lattices
(see Theorem 1.14.4. in \cite{9}):

\begin{theorem} Let $M$ be an even lattice of signature
$(t_{(+)},t_{(-)})$, and let $L$ be an even unimodular
lattice of signature \linebreak $(l_{(+)},l_{(-)})$. Then a
primitive embedding of $M$ into $L$ is unique (up to
isomorphisms), provided the following conditions hold:

1) $l_{(+)}-t_{(+)}>0$ and $l_{(-)}-t_{(-)}>0$.

2) $l_{(+)}+l_{(-)}-t_{(+)}-t_{(-)}\ge 2+l({\gA}_{M_p})$ for all $p\not=2$.

3) If $l_{(+)}+l_{(-)}-t_{(+)}-t_{(-)}=l({\gA}_{M_2})$,
then $q_M\cong u_+^{(2)}(2)\oplus q^\prime$ or
$q_M\cong v_+^{(2)}(2)\oplus q^\prime$.

\label{Witt}
\end{theorem}


\subsection{Classification of main invariants and their
geometric interpretation}
\label{subsec:maininv} Here we apply results of Section
\ref{subsec:discrforms} to classify main invariants $S$ and
$(r,a,\delta)$ of non-symplectic involutions of K3 surfaces and
equivalent right DPN surfaces; moreover, we give their geometric
interpretation (types: elliptic, parabolic, hyperbolic, and
invariants $(k,g,\delta)$) (see Section \ref{subsec2.2}). We also
give proofs of results of Section \ref{subsec2.2} which were only
cited there. All these results had been obtained in
\cite{8,9,N1,10,31} and are well-known. We follow these papers.

We follow notation and considerations of Section \ref{subsec2.2}.

According to Section \ref{subsec2.2}, the set of main invariants $S$
of K3 surfaces with non-symplectic involution is exactly
{\it the set of isomorphism classes of 2-elementary even
hyperbolic lattices $S$ having a primitive embedding $S\subset L$
where $L\cong L_{K3}$ is an even unimodular lattice of signature
$(3,19)$.} Further we denote $L=L_{K3}$.

Since the lattice $T=S^\perp_L$ is also 2-elementary, let us more
generally consider all even 2-elementary lattices $M$. We denote
by $(t_{(+)},t_{(-)})$ the numbers of positive and negative
squares of $M$. Since $M$ is 2-elementary, the discriminant group
${\gA}_M\cong (\bZ/2\bZ)^a$ is a $2$-elementary group where $2^a$
is its order. We get the important {\it invariant $a$} of the
discriminant group ${\gA}_M$ (and $M$ itself). We have $a\in \bZ$
and $a\ge 0$.

The discriminant form $q_M$ of $M$ is a non-degenerate finite
quadratic form on the 2-elementary group ${\gA}_M$ (we call such a
form as 2-elementary either). By Jordan decomposition
(see Section \ref{subsubsec:discforms}),
the form $q_M$ is orthogonal sum of elementary finite quadratic forms
$u_+^{(2)}(2)$, $v_+^{(2)}(2)$ and $q_{\pm 1}^{(2)}(2)$ with the
signature $0\mod 8$, $4\mod 8$ and $\pm 1\mod 8$ respectively. If
$q_M$ is sum of only elementary forms $u_+^{(2)}(2)$ and
$v_+^{(2)}(2)$, then $q_M$ is {\it even:} it takes values only in
$\bZ/2\bZ$. Otherwise it is {\it odd:} at least one of its values
belongs to $\{-1/2,1/2\}\mod 2$. Therefore, we introduce an
important {\it invariant $\delta\in \{0,1\}$} of $q_M$ (and of $M$
itself). The $\delta=0$ if $q_M$ is even, and $\delta=1$ if $q_M$
is odd.

We have the important relations between elementary forms:
\linebreak  $2u_+^{(2)}(2)\cong 2v_+^{(2)}(2)$, $3q_{\pm
1}^{(2)}\cong q_{\mp 1}^{(2)}(2)\oplus v_+^{(2)}(2)$,
$q_1^{(2)}(2)\oplus q_{-1}^{(2)}(2)\oplus q_{\pm 1}^{(2)}(2)\cong
u_+^{(2)}(2)\oplus q_{\pm 1}^{(2)}(2)$. It follows that $q_M$ can
be written in the {\it canonical form} depending on its invariants
$a$, $\delta$ and $\sigma \mod 8=\text{sign}\,q_M\mod 8 \equiv
t_{(+)}-t_{(-)}\mod 8$. We have several cases:

\medskip

\noindent $\delta=0$: then $a\equiv 0\mod 2$, $\sigma\equiv 0\mod
4$, and $\sigma\equiv 0\mod 8$ if $a=0$. We have
$$
q_M\cong sv_+^{(2)}(2)\oplus (a/2-s)u_+^{(2)}(2)
$$
where $s=0$ or $1$ and $\sigma\equiv 4s\mod 8$.

\medskip

\noindent $\delta=1$: then $a\ge 1$, $\sigma\equiv a\mod 2$,
$\sigma\equiv \pm 1\mod 8$ if $a=1$, and $\sigma\not\equiv 4\mod
8$ if $a=2$. We have
$$
q_M\cong q_{\pm 1}^{(2)}(2)\oplus
\left((a-1)/2\right)u_+^{(2)}(2)\ if\ \sigma\equiv \pm 1\mod 8;
$$
$$
q_M\cong 2q_{\pm 1}\oplus (a/2-1)u_+^{(2)}(2) \ if\ \sigma\equiv
\pm 2\mod 8;
$$
$$
q_M\cong q_1^{(2)}(2)\oplus q_{-1}^{(2)}(2)\oplus
(a/2-1)u_+^{(2)}(2)\ if\  \sigma \equiv 0\mod 8;
$$
$$
q_M\cong q_1^{(2)}(2)\oplus q_{-1}^{(2)}(2)\oplus
v_+^{(2)}(2)\oplus (a/2-2)u_+^{(2)}(2)\ if\  \sigma \equiv 4\mod
8.
$$
Thus, the discriminant form $q_M$ is determined by its invariants
$(\sigma\equiv t_{(+)}-t_{(-)}\mod 8, a,\delta)$. Moreover, we
have listed above all conditions of existence of $q_M$ for the
given invariants $\sigma=t_{(+)}-t_{(-)}\mod 8$, $a\ge 0$ and
$\delta \in\{0,1\}$.

Assume that these conditions are satisfied. By Corollary
\ref{discrexistence2}, a 2-elementary lattice $M$ with invariants
$(t_{(+)},t_{(-)},q_M)$ exists if $t_{(+)}\ge 0$, $t_{(-)}\ge 0$
and $t_{(+)}+t_{(-)}>a=l({\gA}_M)$. The condition
$t_{(+)}+t_{(-)}\ge a$ is necessary for the existence. Assume that
$t_{(+)}+t_{(-)}=a$. If $\delta=1$, then $q_M\cong q_{\pm
1}^{(2)}(2)\oplus q^\prime$, and the lattice $M$ also does exist
by Theorem \ref{discrexistence}. If $\delta=0$, then $M(1/2)$ will
be an even unimodular lattice. It follows that the condition
$t_{(+)}-t_{(-)}\equiv 0\mod 8$ must be satisfied, and it is
sufficient for the existence of $M$ since an even unimodular
lattice $M(1/2)$ with the invariants $(t_{(+)},t_{(-)})$ does
exist under this condition. Thus, we finally listed all conditions
of existence of an even 2-elementary lattice $M$ with the
invariants $(t_{(+)},t_{(-)},a,\delta)$.

Moreover, we had proved that the invariants
$(t_{(+)},t_{(-)},a,\delta)$ define the discriminant quadratic
form $q_M$ of $M$. We have $l({\gA}_{M_p})=0$ if a prime $p$ is
odd. If $a\ge 3$, then $q_M\cong u_+^{(2)}(2)\oplus q^\prime$ or
$q_M\cong v_+^{(2)}(2)\oplus q^\prime$. By Theorem \ref{uniqueness1},
then the lattice $M$ is unique if it is indefinite and $\rk M\ge 3$.
Moreover, then the canonical homomorphism $O(M)\to O(q_M)$ is
epimorphic. If $\rk M\le 2$, then $M$ is one of lattices: $\langle
\pm 2\rangle$, $\langle \pm 2\rangle \oplus \langle \pm 2\rangle$,
$U$ or $U(2)$. One can easily check for them the same statements
directly.

Thus, finally we get the following classification result about
2-elemen\-tary even indefinite (except few exceptions) lattices.
It is Theorems 3.6.2 and 3.6.3 from \cite{9}.

\begin{theorem}
\label{2elementary} The genus of an even 2-elementary lattice $M$
is determined by the invariants $(t_{(+)},t_{(-)},a,\delta)$; and
if either $M$ is indefinite or $\rk M=2$, these invariants
determine the isomorphism class of $M$, and the canonical
homomorphism $O(M)\to O(q_M)$ is epimorphic.

An even 2-elementary lattice $M$ with invariants
$(t_{(+)},t_{(-)},a,\delta)$ exists if and only if all the
following conditions are satisfied (it being assumed that
$\delta=0$ or $1$, and that $a, t_{(+)},t_{(-)}\ge 0$):

1) $a\le t_{(+)}+t_{(-)}$;

2) $t_{(+)}+t_{(-)}+a\equiv 0\mod 2$;

3) $t_{(+)}-t_{(-)}\equiv 0\mod 4$ if $\delta=0$;

4) $(\delta=0,\ t_{(+)}-t_{(-)}\equiv 0\mod 8)$ if $a=0$;

5) $t_{(+)}-t_{(-)}\equiv \pm 1\mod 8$ if $a=1$;

6) $\delta=0$ if $(a=2,\ t_{(+)}-t_{(-)}\equiv 4\mod 8)$;

7) $t_{(+)}-t_{(-)}\equiv 0\mod 8$ if $(\delta=0,\
a=t_{(+)}+t_{(-)})$.
\end{theorem}

Let $S$ be a main invariant and $r=\rk S$. Since $S$ is
2-elementary even hyperbolic, by Theorem \ref{2elementary} it is
then determined by its invariants $(t_{(+)}=1, t_{(-)}=r-1,
a,\delta)$.

By Theorem \ref{primembedd2}, existence of a primitive embedding
$S\subset L_{K3}$ is equivalent to existence of a 2-elementary
even lattice $T=S^\perp$ with invariants
$(t_{(+)}=2,t_{(-)}=20-r,a,\delta)$ (indeed, $q_T\cong -q_S$ has
the same invariants $a$ and $\delta$).

Thus, {\it the set of main invariants $S$ is equal to the set of
$(r,a,\delta)$ such that both $(1,r-1,a,\delta)$ and
$(2,20-r,a,\delta)$ satisfy conditions 1) --- 7) of Theorem
\ref{2elementary}. It consists of exactly $(r,a,\delta)$
which are presented in Figure \ref{fig1}.}

By Theorem \ref{2elementary}, the orthogonal complement $T=S^\perp_L$
is uniquely determined by $(r,a,\delta)$, and the canonical
homomorphism $O(T)\to O(q_T)$ is epimorphic. By Proposition
\ref{primembedd1} the primitive embedding $S\subset L_{K3}$ is unique
up to automorphisms of $L_{K3}$.

Let us show that $O(S\subset L_{K3})$ contains an automorphism of
spinor norm $-1$ (i. e. it changes two connected components of the
quadric $\Omega_{S\subset L_{K3}}$, see \eqref{Omega1}). Using
Theorem \ref{2elementary}, it is easy to see that either
$T=\langle 2 \rangle\oplus \langle 2\rangle$, or $T= U\oplus
T^\prime$, or $T= U(2)\oplus T^\prime$, or $T=\langle 2 \rangle
\oplus \langle -2 \rangle \oplus T^\prime$. If $T\cong \langle 2
\rangle\oplus \langle 2\rangle$, we consider the automorphism
$\alpha$ of $T$ which is $+1$ on the first $\langle 2\rangle$ and
$-1$ on the second $\langle 2 \rangle$. In remaining cases we
consider an automorphism $\alpha$ of $T$ which is $-1$ on the first
2-dimensional hyperbolic summand of $T$, and which is $+1$ on
$T^\prime$. It is easy to see that $\alpha$ changes connected
components of the $\Omega_{S\subset L_{K3}}$. On the other hand,
$\alpha$ is identical on the $T^\ast/T$ and can be continued
identically on $S$. This extension gives an automorphism of
$L_{K3}$ which is identical on $S$ and has spinor norm $-1$.

It follows (see Section \ref{subsec2.2}) that for the fixed main
invariants $(r,a,\delta)$ the moduli space $\cM_{(r,a,\delta)}$ of
K3 surfaces with non-symplectic involution is connected.

\medskip

Now let us consider the geometric interpretation of main
invariants $(r,a,\delta)$ in terms of the set $C=X^\theta$ of the
fixed points. The set $C$ is non-singular. Indeed, if $x\in C$ a
singular point of $C$, then $\theta$ is the identity in the
tangent space $T_x$. Then $\theta^\ast (\omega_X)=\omega_X$ for
any $\omega_X\in H^{2,0}(X)$ and $\theta$ is symplectic. We get a
contradiction.

For a non-singular irreducible curve $C$ on a K3 surface $X$ we
have $g(C)=(C^2+C\cdot K_X)/2+1=C^2/2+1\ge 0$. It follows that
$C^2>0$, if $g(C)>1$. Since the Picard lattice $S_X$ is
hyperbolic, it follows that any two curves on $X$ of genus $\ge 2$
must intersect. It then follows that $X^\theta$ has one of types
A, B, C listed below:

Case A: $X^\theta=C_g+E_1+\cdots+E_k$ where $C_g$ is a
non-singular irreducible curve of genus $g\ge 0$  and
$C_g\not=\emptyset$, the curves $E_1,\dots,E_k$ are non-singular
irreducible  rational (i. e. $E_i^2=-2$). All curves $C_g$, $E_i$
are disjoint to each other.

Case B: $X^\theta=C_1^{(1)}+\cdots+C_1^{(k)}$ is disjoint union of
$k>1$ elliptic curves (we shall prove in a moment that actually
$k=2$).

Case C: $X^\theta=\emptyset$.

By Lefschetz formula, the Euler characteristics
$\chi(X^\theta)=2+r-(22-r)=2r-20$.

By Smith Theory (see \cite{Kha76}), the total Betti number
over $\bZ/2\bZ$ satisfies
$$
\dim H^*(X^\theta,\bZ/2\bZ)=\dim H^*(X,\bZ/2\bZ)-2a=24-2a\ \text{if}\
X^\theta\not=\emptyset.
$$

For any 2-dimensional cycle $Z\subset X$ one evidently has $Z\cdot
\theta(Z)\equiv Z\cdot X^\theta\mod 2$. Thus, $X^\theta\sim 0\mod
2$ in $H^2(X,\bZ)$ if and only if $(x,\theta^\ast x)\equiv 0\mod
2$ for any $x\in H^2(X,\bZ)$. Let us write $x\in L=H^2(X,\bZ)$ as
$x=x_++x_-$ where $x_+\in S^\ast$ and $x_-\in (S^\perp)^\ast$.
Then $x\cdot \theta^\ast(x)=x_+^2-x_-^2$. Moreover,
$x^2=x_+^2+x_-^2\equiv 0\mod 2$ because $L$ is even. Taking the sum,
we get $x\cdot \theta^\ast(x)\equiv 2x_+^2\mod 2$. Since
$H^2(X,\bZ)$ is unimodular,  any $x_+\in S^\ast$ appears in
this identity. It follows that $X^\theta\sim 0\mod 2$ in
$H^2(X,\bZ)$ if and only if $x_+^2\in \bZ$ for any $x_+\in
S^\ast$. Equivalently, the invariant $\delta =0$. Therefore
$$
\delta=0\ \text{if and only if}\ X^\theta\sim 0\mod 2\ \text{in}\
H^2(X,\bZ).
$$

In case B, elliptic curves $C_1^{(i)}$ belong to one elliptic
pencil $|C|$ of elliptic curves where it is known (see \cite{13}) that
$C$ is primitive in Picard lattice $S_X$. Assume that $k>2$.
Then $\theta$ is trivial on the base $\bP^1$ of the elliptic pencil.
Since it is also trivial on a fibre $C_1^{(i)}$ which is not multiple,
$\theta$ is symplectic, and we get a contradiction. Thus, in Case B
we have $k=2$ and $\delta=0$.

In case C, the quotient $Y=X/\{1, \theta\}$ is an Enriques surface.
It follows that $r=a=10$ and $\delta=0$ in this case.

Combining all these arguments, we obtain the geometric
interpretation of the invariants $(r,a,\delta)$ cited in Section
\ref{subsec2.2}

\medskip

\subsection{The analogue of Witt's theorem for 2-elementary finite forms}
\label{subsec:Witt'stheorem} Here we follow Section 1.9 in \cite{12}
to prove an important Lemma \ref{lemma2.4.3}.

We consider a 2-elementary finite bilinear forms $b:B\times B\to
\frac{1}{2}\bZ/\bZ$ and 2-elementary finite quadratic forms
$q:Q\to \frac{1}{2}\bZ/2\bZ$ on finite 2-elementary groups $B$,
$Q$.

In the previous section we gave classification of non-degenerate
2-elementary finite quadratic forms. Similarly one can classify
non-degenerate 2-elementary finite bilinear forms. They are
orthogonal sums of elementary forms $b_1^{(2)}(2)$ and
$u_-^{(2)}(2)$. The form $b_1^{(2)}(2)$ is the bilinear form of
quadratic forms $q_{\pm}^{(2)}(2)$, and the form  $u_-^{(2)}(2)$
is the bilinear form of quadratic forms $u_+^{(2)}(2)$ and
$v_+^{(2)}(2)$. We denote by $s_b$ the characteristic element of
$b$, i. e. $b(x,x)=b(s_b,x)$ for all $x\in B$. It is easy to see
that any non-degenerate 2-elementary finite bilinear form $b$ is
$$
b\cong m u_-^{(2)}(2)
$$
if $b(x,x)=0$ for all $x\in B$ (equivalently the characteristic
element $s_b=0$, these bilinear forms are the same as
skew-symmetric ones);
$$
b\cong b_1^{(2)}(2)\oplus mu_-^{(2)}(2),
$$
if $b(s_b,s_b)=\frac{1}{2}\mod 1$;
$$
b\cong 2b_1^{(2)}\oplus mu_-^{(2)}(2)
$$
if $s_b\not=0$ but $b(s_b,s_b)=0$.

 We prove (see Section 1.9 in \cite{12})

\begin{proposition}
\label{propWitt1} Let $b$ be a non-degenerate bilinear form on a
finite 2-elementary group $B$ and $\theta:H_1\to H_2$ be an
isomorphism of subgroups of $B$ which preserves the restrictions
$b|H_1$ and $b|H_2$ and that maps the characteristic element of
the form $b$ to itself (if, of course, it belongs to $H_1$). Then
$\theta$ extends to an automorphism of $b$.
\end{proposition}

\begin{proposition}
\label{propWitt2} Let $q$ be a quadratic form on a finite
2-elementary group $Q$ whose kernel is zero; that is
$$
\{x\in Q \suchthat x\perp Q\ and\ q(x)=0\}=0.
$$
Let $\theta:H_1\to H_2$ be an isomorphism of two subgroups of $Q$
that preserves the restrictions $q|H_1$ and $q|H_2$ and that maps
the elements of the kernel and the characteristic elements of the
bilinear form $q$ into the same sort of elements (of course, if
they belong to $H_1$). Then $\phi$ extends to an automorphism of
$q$ (an element $s\in Q$ is called characteristic if
$q(s,x)=q(x,x)$ for all $x\in Q$).
\end{proposition}

We shall prove the propositions by induction on the number of
generators of $B$ and $Q$. Let us begin with Proposition
\ref{propWitt1}. Suppose there exist $x_1\in H_1$ and
$x_2=\theta(x_1)\in H_2$ such that
$b(x_1,x_1)=b(x_2,x_2)=\frac{1}{2}\mod 1$. Write
$B_1=(x_1)^\perp_B$, $B_2=(x_2)^\perp_B$, $b_1=b|B_1$,
$b_2=b|B_2$, $H_1^\prime=(x_1)^\perp_{H_1}$,
$H_2^\prime=(x_2)^\perp_{H_2}$ and $\theta^\prime=\theta|H_1$.
Then the same conditions hold for the nondegenerate forms $b_1$
and $b_2$ defined on the subgroups $B_1$ and $B_2$, their
subgroups $H_1^\prime\subset B_1$ and $H_2^\prime\subset B_2$, and
an isomorphism $\theta^\prime:H_1^\prime \cong H_2^\prime$.
Everything reduces to extending $\theta^\prime$. Since $x_1$ and
$x_2$ are characteristic simultaneously, $b_1$ and $b_2$ are
isomorphic (this follows from classification of nondegenerate
bilinear forms). Therefore, the existence of an extension of
$\theta^\prime$ follows from the induction hypothesis. To complete
the proof it remains to consider the case when the function
$b(x,x)$ on $H_1$ and $H_2$ is zero. Denote by $s$ the
characteristic element of $B$. It is easy to check (using the
classification again) that the natural homomorphism
$$
O(b)\to O(s^\perp)
$$
is epimorphic (we always consider a subgroup with the restriction
of the form $b$ on the subgroup). In our case $H_1$ and $H_2$ lie
in $s^\perp$; therefore, it suffices to extend $\theta$ to an
automorphism of $s^\perp$. If $s\notin s^\perp$, this is obvious
since in this case $s^\perp$ is a nondegenerate skew-symmetric
form; for them the proposition is well-known and obvious. If $s\in
s^\perp$, then $[s]$ is the kernel of $s^\perp$, and, by the
hypothesis,
$$
[s]\cap H_1=\theta([s]\cap H_1)=[s]\cap H_2.
$$
Let
$$
\overline{\theta}:H_1/([s]\cap H_1)\cong H_2/([s]\cap H_2)
$$
be the isomorphism $\theta\mod[s]\cap H_1$. Then, because
$s^\perp/[s]$ is nondegenerate and skew-symmetric,
$\overline{\theta}$ extends to an automorphism $\overline{\psi}\in
O(s^\perp/[s])$. Let $\psi$ be a
lifting of $\overline{\psi}$ to an automorphism of $(s)^\perp$.
Then $\psi(x)-\theta(x)=\overline{g}(x)s$, if $x\in H_1$, where
$\overline{g}:H_1\to \bZ/2\bZ$ is a linear function. Extending
$\overline{g}$ to a linear function $g:s^\perp\to \bZ/2\bZ$, we
put
$$
\widetilde{\psi}(x)=\psi(x)+\widetilde{f}(x)s\ if\ x\in s^\perp.
$$
Evidently $\widetilde{\psi}\in O(s^\perp)$ is the desired
extension of $\theta$.

Let us prove Proposition \ref{propWitt2}. Assume that the bilinear
form of $q$ has a nonzero kernel. Then it is generated by an
element $r$, and $q(r)=1\mod 2$. Using Proposition
\ref{propWitt1}, we can extend $\theta$ to an automorphism $\psi$
of the bilinear form of $q$. The function $f(x)=q(x)-q(\psi(x))\in
\bZ/2\bZ$, where $x\in Q$, is linear and vanishes on $H_1$.
Evidently, $\widetilde{\psi}(x)=\psi(x)+f(x)r$, where $x\in Q$, is
the desired extension of $\theta$.

It remains to examine the case when the bilinear form $q$ is
nondegenerate. The case where there exist elements $x_1\in H_1$
and $x_2=\theta(x_1)\in H_2$ for which
$q(x_1)=q(x_2)=\frac{1}{2}\theta$ where $\theta\notin 2\bZ$, can
be examined similarly to the corresponding case of Proposition
\ref{propWitt1}. Therefore, we assume that $q(x)\equiv 0\mod 1$ if
$x\in H_1$. For a characteristic element $s\in Q$ the natural
homomorphism
$$
O(q)\to O(s^\perp)
$$
is epimorphic (this easily follows from the classification of
nondegenerate quadratic forms on 2-elementary groups given in Section
\ref{subsec:maininv}). In our case
$H_1$ and $H_2$ lie in $s^\perp$, and it suffices to extend
$\theta$ to an automorphism of $s^\perp$. If the bilinear form on $s^\perp$ is
non-degenerate, this follows from the classical Witt theorem over
the field with two elements. Suppose it is degenerate; then
the kernel of $s^\perp$ is generated by $s$. In the case $q(s)=0$, one
can pass to a form on $s^\perp/[s]$ and we argue in the same way as
in the proof of Proposition \ref{propWitt1}. But if $q(s)\equiv
1\mod 2$, then we pass to the first case, already treated.
Proposition \ref{propWitt2} is proved.

\subsection{Calculations of fundamental chambers of hyperbolic reflection
groups.}
\label{fundchambers}

Here we outline calculations of
fundamental chambers for hyperbolic reflection groups which had
been used in the main part of the work.

\subsubsection{Fundamental chambers $\cM^{(2,4)}$ of 2-elementary
even hyperbolic lattices of elliptic type (Table 1).}
\label{fundchambersofS} We consider all 50 types of 2-elementa\-ry
even hyperbolic lattices $S$ of elliptic type given by their full
invariants $(r,a,\delta)$. We outline the calculation of a
fundamental chamber $\cM^{(2,4)}$ (equivalently, the corresponding
Dynkin diagram $\Gamma (P(\cM^{(2,4)}))$) for the full reflection
group $W^{(2,4)}(S)=W(S)$. This is the group generated by
reflections $s_f$ in all roots $f$ of $S$. They are elements $f\in
S$ either with $f^2=-2$ or with $f^2=-4$ and $(f,S)\equiv 0\mod 2$.
The reflection $s_f\in O(S)$ is then given by
$$
x\mapsto x-\frac{2(x,f)f}{f^2},\ \ \ \forall x\in S.
$$

We use Vinberg's algorithm \cite{2} which we describe below. It can
be applied to any hyperbolic lattice $S$ and any of its reflection
subgroup $W\subset W(S)$ which is generated by reflections in some
precisely described subset $\Delta\subset S$ of primitive roots of
$S$ which is $W$-invariant.

First, we should choose a non-zero $H\in S$ with $H^2\ge 0$. Then
$H$ defines the half cone $V^+(S)$ such that $H\in
\overline{V^+(S)}$. We want to find a fundamental chamber
$\cM\subset \cL(S)=V^+(S)/\bR^+$ of $W$ containing $\bR^+h$.

Step 0. We consider the subset $\Delta_0$ of all roots from
$\Delta$ which are orthogonal to $H$. This set is either a finite
root system or affine root system. One should choose a bases $P_0$
in $\Delta_0$. For example, one can take another element $H_1\in
S$ such that $H_1^2\ge 0$, $(H,H_1)>0$ and $(H_1,\Delta_0)$ does
not contain zero. Then
$$
\Delta_0^+=\{f\in \Delta_0 \suchthat (f,H_1)>0\},
$$
and $P_0\subset \Delta_0^+$ consists of roots from $\Delta_0^+$
which are not non-trivial sums of others. For $f\in \Delta$ with
$(f,H)\ge 0$ we introduce the height
$$
h(f)=\frac{2(f,H)^2}{-f^2}.
$$
The height is equivalent to the hyperbolic distance between the
point $\bR^+H$ and the hyperplane $\cH_f$ which is orthogonal to
$f$. The set of all possible heights is a discrete ordered subset
\begin{equation}
h_0=0,\,h_1,\,h_2,\,\dots \,\,h_i,\,\dots, \label{pol1}
\end{equation}
of $\bR^+$. It is always a subset of non-negative integers, and
one can always take $\bZ^+$ as the set of possible heights.

The fundamental chamber $\cM\subset \cL(S)$ is defined by the set
$P(\cM)\subset \Delta$ of orthogonal roots to $\cM$ which is
\begin{equation}
P(\cM)=\bigcup_{0\le j}{P_j} \label{pol2}
\end{equation}
where $P_0$ is defined above and $P_j$ for $j>0$ consists of all
$f\in\Delta$ such that $(f,H)>0$, the height $h(f)=h_j$, and
\begin{equation}
\left(f,\bigcup_{0\le i\le  j-1}{P_i}\right)\ge 0. \label{pol3}
\end{equation}
Then
\begin{equation}
 \cM=\{\bR^+x\in \cL(S)\suchthat (x,P(\cM))\ge 0\}. \label{pol4}
\end{equation}

If $\cM$ has finite volume, the algorithm terminates after a
finite number $m$ of steps (i. e. all $P_j$ are empty for $j>m$,
and
\begin{equation}
P(\cM)=\bigcup_{0\le j\le  m}{P_j} \label{pol5}
\end{equation}
whenever \eqref{pol4} defines a polyhedron $\cM$ of finite volume
in $\cL(S)$ for $P(\cM)$ given by \eqref{pol5}.

\medskip

Below we apply this algorithm to 2-elementary hyperbolic lattices $S$ of
elliptic type and $W=W(S)=W^{(2,4)}(S)$. The set $\Delta\subset S$
consists of all $f\in S$ such that either $f^2=-2$ or
$f^2=-4$ and $(f,S)\equiv 0\mod 2$.

\medskip

{\bf Cases $S=\langle 2 \rangle \oplus lA_1$ where $0\le l\le 8$.}
Then $(r,a,\delta)=(1+l,1+l,1)$, $0\le l\le 8$. We use the
standard orthogonal basis $h$ for $\langle 2 \rangle$ where
$h^2=2$, and the standard orthogonal basis $v_1,\dots,v_l$ for
$lA_1$ where $v_1^2=\cdots =v_{l}^2=-2$.

We take $H=h$, $H_1=th+v_1+2v_2+\cdots + lv_{l}$ where $t>>0$.
Then $P(\cM^{(2,4)})$ consists of roots: $\beta_0=h-v_1-v_2$ if
$l=2$; $\beta_0=h-v_1-v_2-v_3$ if $l\ge 3$;
$\beta_1=v_{1}-v_2,\,\dots,\,\beta_{l-1}=v_{l-1}-v_{l}$ if $l\ge
2$; $\beta_{l}=v_{l}$ if $l\ge 1$.

For $2\le l\le 8$ the polyhedron $\cM^{(2,4)}$ is obviously a
simplex in $\cL(S)$ of finite volume. The Gram matrix
$(\beta_i,\beta_j)$ gives the Dynkin diagrams of cases $N=1$,
$N=3\ -\ 10$ of Table 1. Here we repeated calculations by Vinberg
in \cite{2}.

\medskip

{\bf Cases $S=U\oplus lA_1$, $0\le l\le 8$.} Then
$(r,a,\delta)=(2,0,0)$ if $l=0$, and $(r,a,\delta)=(2+l,l,1)$ if
$1\le l\le 8$. We use the standard basis $c_1,c_2$ for $U$ where
$c_1^2=c_2^2=0$ and $(c_1,c_2)=1$, and the standard orthogonal
basis $v_1,\dots,v_l$ for $lA_1$ as above.

We take $H=c_1$. We can take $P_0$ which consists of
$\beta_0=c_1-v_1$ if $l\ge 1$;
$\beta_1=v_1-v_2,\dots,\beta_{l-1}=v_{l-1}-v_l$ if $l\ge 2$;
$\beta_{l}=v_l$ if $l\ge 1$. Then $P(\cM^{(2,4)})$ consists of
$P_0$, $e=-c_1+c_2$ and additional elements
$\gamma_{1}=2c_1+2c_2-v_1-v_2-v_3-v_4-v_5$ if $l=5$, and
$\gamma_{1}=2c_1+2c_2-v_1-v_2-v_3-v_4-v_5-v_6$ if $l\ge 6$.

We shall discuss the finiteness of volume of $\cM^{(2,4)}$ later. We
obtain the diagrams of cases $N=11\ -\ 20$ excluding $N=15$ of Table
1.

\medskip

{\bf Cases $S=U\oplus D_4\oplus lA_1$, $0\le l\le 5$.} Then
$(r,a,\delta)=(6,2,0)$ if $l=0$, and $(r,a,\delta)=(6+l,2+l,1)$ if
$1\le l\le 5$. We use the standard bases $c_1,c_2$ for $U$ and
$v_1,\dots,v_l$ for $lA_1$ as above. We use the standard
orthogonal basis $\epsilon_1,\dots,\epsilon_m$ for $D_m\otimes
\bQ$ where $\epsilon_1^2=\cdots =\epsilon_m^2=-1$; the lattice
$D_m$ consists of all $x_1\epsilon_1+\cdots +x_m \epsilon_m$ where
$x_i\in \bZ$ and $x_1+\cdots +x_m\equiv 0\mod 2$.

We take $H=c_1$, and we can take $P_0$ which consists of
$\alpha_0=c_1-\epsilon_1-\epsilon_2$,
$\alpha_1=\epsilon_1-\epsilon_2$,
$\alpha_2=\epsilon_2-\epsilon_3$,
$\alpha_3=-\epsilon_1+\epsilon_2+\epsilon_3-\epsilon_4$,
$\alpha_4=2\epsilon_4$, $\beta_0=c_1-v_1$ if $l\ge 1$;
$\beta_1=v_1-v_2,\,\dots,\,\beta_{l-1}=v_{l-1}-v_l$ if $l\ge 2$;
$\beta_l=v_l$ if $l\ge 1$.

Then $P(\cM^{(2,4)})$ consists of $P_0$, $e=-c_1+c_2$,
$\gamma_1=2c_1+2c_2-\epsilon_1-\epsilon_2-\epsilon_3-\epsilon_4-v_1-v_2-v_3$
if $l=3$, and
$\gamma_1=2c_1+2c_2-\epsilon_1-\epsilon_2-\epsilon_3-
\epsilon_4-v_1-v_2-v_3-v_4$
if $l\ge 4$; $\gamma_2=2c_1+2c_2-v_1-v_2-v_3-v_4-v_5$ if $l=5$.

We obtain the diagrams of cases $N=21\ -\ 27$ excluding $N=25$ of
Table 1.

\medskip

{\bf Cases $S=U\oplus D_m\oplus lA_1$ where $m\equiv 0\mod 2$,
$m\ge 6$, $l\ge 0$, $m+2l\le 14$.} Then $r=2+m+l$, $a=l+2$.
Moreover $\delta=0$ if $l=0$ and $m\equiv 0\mod 4$, otherwise
$\delta=1$.

We use the standard bases $c_1,\,c_2$ for $U$, and
$\epsilon_1,\,\dots,\,\epsilon_m$ for $D_m\otimes \bQ$, and
$v_1,\,\dots,\,v_l$ for $lA_1$ as above.

We take $H=c_1$, and we can take $P_0$ which consists of
$\alpha_0=c_1-\epsilon_1-\epsilon_2$,
$\alpha_1=\epsilon_1-\epsilon_2,\,\dots,\,\alpha_{m-1}=
\epsilon_{m-1}-\epsilon_m$,
$\alpha_m=2\epsilon_m$; $\beta_0=c_1-v_1$ if $l\ge 1$;
$\beta_1=v_1-v_2,\,\dots,\,\beta_{l-1}=v_{l-1}-v_l$ if $l\ge 2$;
$\beta_l=v_l$ if $l\ge 1$.

Then $P(\cM^{(2,4)})$ consists of $P_0$, $e=-c_1+c_2$, and some
additional elements $\gamma_i$ depending on $m\ge 6$ and $l\ge 0$
where we always assume that $m+2l\le 14$ and $m\equiv 0\mod 2$.

If $m=6$ and $l=2$, one must add
$\gamma_1=2c_1+2c_2-\epsilon_1-\cdots -\epsilon_6-v_1-v_2$.

If $m=6$ and $l=3$, one must add
$\gamma_1=2c_1+2c_2-\epsilon_1-\cdots -\epsilon_6-v_1-v_2-v_3$,
$\gamma_2=2c_1+2c_2-2\epsilon_1-v_1-v_2-v_3$.

If $m=6$ and $l=4$, one must add
$\gamma_1=2c_1+2c_2-\epsilon_1-\cdots -\epsilon_6-v_1-v_2-v_3$,
$\gamma_2=2c_1+2c_2-2\epsilon_1-v_1-v_2-v_3-v_4$.

If $m=8$ and $l=1$, one must add
$\gamma_1=2c_1+2c_2-\epsilon_1-\cdots -\epsilon_8-v_1$.

If $m=8$ and $l=2$, one must add
$\gamma_1=2c_1+2c_2-\epsilon_1-\cdots -\epsilon_8-v_1-v_2$.

If $m=8$ and $l=3$, one must add
$\gamma_1=2c_1+2c_2-\epsilon_1-\cdots -\epsilon_8-v_1-v_2$,
$\gamma_2=2c_1+2c_2-\epsilon_1-\epsilon_2-\epsilon_3-\epsilon_4-v_1-v_2-v_3$,
$\gamma_3=2c_1+2c_2-2\epsilon_1-v_1-v_2-v_3$,

If $m=10$ and $l=0$, one must add
$\gamma_1=2c_1+2c_2-\epsilon_1-\cdots -\epsilon_{10}$.

If $m=10$ and $l=1$, one must add
$\gamma_1=2c_1+2c_2-\epsilon_1-\cdots -\epsilon_{10}-v_1$.

If $m=10$ and $l=2$, one must add
$\gamma_1=2c_1+2c_2-\epsilon_1-\cdots -\epsilon_{10}-v_1$,
$\gamma_2=2c_1+2c_2-\epsilon_1-\cdots-\epsilon_6-v_1-v_2$,
$\gamma_3=4c_1+4c_2-3\epsilon_1-\epsilon_2-\cdots
-\epsilon_{10}-2v_1-2v_2$.

If $m=12$ and $l=0$, one must add
$\gamma_1=2c_1+2c_2-\epsilon_1-\cdots -\epsilon_{12}$.

If $m=12$ and $l=1$, one must add
$\gamma_1=2c_1+2c_2-\epsilon_1-\cdots -\epsilon_{12}$,
$\gamma_2=2c_1+2c_2-\epsilon_1-\cdots -\epsilon_{8}-v_1$,
$\gamma_3=6c_1+
6c_2-4\epsilon_1-2\epsilon_2-2\epsilon_3-\cdots-2\epsilon_{11}-3v_1$.

We obtain the diagrams for $N=28\ - \  50$ of Table 1 except $N=30,\,
34,\, 40,\, 41,\, 44,\, 49,\, 50$ when either $\delta=0$ and
$a>2$, or $a\le 1$.

 \medskip

{\bf Cases $S=U(2)\oplus D_m$ where $m\equiv 0\mod 4$ and $0\le
m\le 12$.} Then $(r,a,\delta)=(2+m,4,0)$. We use the standard
bases $c_1,\,c_2$ for $U(2)$ where $c_1^2=c_2^2=0$ and
$(c_1,c_2)=2$, and the standard basis $\epsilon_1,\,\dots,\,
\epsilon_m$ for $D_m\otimes \bQ$ as above.

We use $H=c_1$ and denote $e=-c_1+c_2$ with $e^2=-4$.

If $m=0$, then $P_0=\emptyset$ and $P(\cM^{(2,4)})$ consists of
$e$.

If $m=4$, then $P_0$ consists of
$\alpha_0=c_1-\epsilon_1-\epsilon_2-\epsilon_3-\epsilon_4$,
$\alpha_1=\epsilon_1-\epsilon_2$,
$\alpha_2=\epsilon_2-\epsilon_3$,
$\alpha_3=-\epsilon_1+\epsilon_2+\epsilon_3-\epsilon_4$,
$\alpha_4=2\epsilon_4$. Then $P(\cM^{(2,4)})$ consists of $P_0$
and $e$.

If $m\ge 8$, then $P_0$ consists of $\alpha_0=c_1-2\epsilon_1$,
$\alpha_1=\epsilon_1-\epsilon_2,\,\dots,\,\alpha_{m-1}=\epsilon_{m-1}-
\epsilon_m$,
$\alpha_m=2\epsilon_m$.

If $m=8$, then $P(\cM^{(2,4)})$ consists of $P_0$, $e$ and
$\gamma_1=c_1+c_2-\epsilon_1-\cdots -\epsilon_8$.

If $m=12$, then $P(\cM^{(2,4)})$ consists of $P_0$, $e$ and
$\gamma_1=2c_1+c_2-\epsilon_1-\cdots - \epsilon_{12}$,
$\gamma_2=c_1+c_2-\epsilon_1-\cdots -\epsilon_6$.

We obtain the diagrams for $N=2,\,15,\,30,\,44$ of Table 1.

\medskip

{\bf Case $U(2)\oplus D_4\oplus D_4$.} Then
$(r,a,\delta)=(10,6,0)$. We use standard bases $c_1,c_2$ for
$U(2)$ and $\epsilon_1^{(1)},\,\dots,\, \epsilon_4^{(1)}$ for the
first $D_4$, and $\epsilon_1^{(2)},\,\dots,\,\epsilon_4^{(2)}$ for
the second $D_4$.

We take $H=c_1$ and $P_0$ which consists of
$\alpha_0^{(1)}=c_1-\epsilon_1^{(1)}-\epsilon_2^{(1)}-
\epsilon_3^{(1)}-\epsilon_4^{(1)}$,
$\alpha_1^{(1)}=\epsilon_1^{(1)}-\epsilon_2^{(1)}$,
$\alpha_2^{(1)}=\epsilon_2^{(1)}-\epsilon_3^{(1)}$,
$\alpha_3^{(1)}=-\epsilon_1^{(1)}+\epsilon_2^{(1)}+
\epsilon_3^{(1)}-\epsilon_4^{(1)}$,
$\alpha_4^{(1)}=2\epsilon_4^{(1)}$ and
$\alpha_0^{(2)}=c_1-\epsilon_1^{(2)}-\epsilon_2^{(2)}-\epsilon_3^{(2)}-
\epsilon_4^{(2)}$,
$\alpha_1^{(2)}=\epsilon_1^{(2)}-\epsilon_2^{(2)}$,
$\alpha_2^{(2)}=\epsilon_2^{(2)}-\epsilon_3^{(2)}$,
$\alpha_3^{(2)}=-\epsilon_1^{(2)}+\epsilon_2^{(2)}+\epsilon_3^{(2)}-
\epsilon_4^{(2)}$,
$\alpha_4^{(2)}=2\epsilon_4^{(2)}$.

Then $P(\cM^{(2,4)})$ consists of $P_0$ and $e=-c_1+c_2$.

We obtain the diagram for $N=25$ of Table 1.

\medskip

{\bf Cases $S=U\oplus E_7$, $U\oplus E_8$, $U\oplus E_8\oplus
A_1$, $U\oplus E_8\oplus E_7$, $U\oplus E_8\oplus E_8$.}
Respectively $(r,a,\delta)=(9,1,1),\ (10,0,0),\ (11,1,1),\
(17,1,1),\  (18,0,0)$. We use the standard basis $c_1,\,c_2$ for
$U$. For each irreducible root lattice $R_i=A_1,\,E_7,\, E_8$ of
the rank $t_i$ we use its standard basis
$r_1^{(i)},\,\dots,\,r_{t_i}^{(i)}$ of roots with the
corresponding Dynkin diagram. We denote by $r_{max}^{(i)}$ the
maximal root of $R_i$ corresponding to this basis.

For $S=U\oplus R$ where $R$ is the sum of irreducible root
lattices $R_i$ above, we take $H=c_1$ and $P_0$ which consists of
standard bases $r_1^{(i)},\,\dots,\,r_{t_i}^{(i)}$ of $R_i$ and
$r_0^{(i)}=c_1-r_{max}^{(i)}$.

Then $P(\cM^{(2,4)})$ consists of $P_0$ and $e=-c_1+c_2$, and one
additional element $\gamma_1$ if $S=U\oplus E_8\oplus E_7$. The
element $\gamma_1\in S$ is shown on the diagram $N=49$ of Table 1
as the right-most vertex. It can be easily computed using
pairings $(\gamma_1,\xi_i)$ prescribed by this diagram for basis
elements $\xi_i$ of $S$ given above.

We obtained the remaining diagrams of cases $N=34,\,40,\,41,\, 49,\,
50$ of Table 1.

\medskip

{\bf Finiteness of volume of polyhedra $\cM^{(2,4)}$ above.} To
prove finiteness of volume of the polyhedra $\cM^{(2,4)}$ defined
by the subsets $P=P(\cM^{(2,4)})\subset S$ calculated above with
the corresponding diagrams $\Gamma=\Gamma(P)$ of Table 1, one can
use methods developed by Vinberg in \cite{2}.

We remind that a subset $T\subset P$ is called {\it elliptic,
parabolic, hyperbolic,} if its Gram matrix is respectively
negative definite, semi-negative definite, hyperbolic. A
hyperbolic subset $T$ is called {\it Lann\'er} if each its proper
subset is elliptic. Dynkin diagrams of all Lann\'er subsets are
classified by Lann\'er, e.g. see Table 3 in \cite{2}. They have
at most 5 elements.

We exclude trivial cases $N=1,\,2,\,3,\,11$ when $\rk S\le 2$. In
all other cases, from our calculations, it easily follows that $P$
generates $S\otimes \bQ$, and $\Gamma (P)$ is connected. Moreover,
by the classification of affine Dynkin diagrams, one can check that
all connected components (for its Dynkin diagram) of any maximal
parabolic subset $T\subset P$ are also parabolic, and sum of their
ranks is $\rk S-2$. We remind that the rank of a connected
parabolic subset $T\subset P$ is equal to $\#T-1$.

From the classification of Lann\'er subsets, it easily follows that
the graph $\Gamma$ has no Lann\'er subgraphs if $N\not=45$, $47$.
By Proposition 1 in \cite{2}), then $\cM^{(2,4)}$ has finite
volume.

Assume that $N=45$ or $N=47$. Then the only Lann\'er subset
$L\subset P$ consists of two elements defining the broken edge (it
is the only one) of $\Gamma$. Finiteness of volume of $\cM^{(2,4)}$ is
then equivalent to $L^\perp\cap
\widetilde{\cM^{(2,4)}}=\emptyset$. Here
\begin{equation}
 \widetilde{\cM^{(2,4)}}=\{x\in S\otimes \bR\suchthat (x,P(\cM))\ge 0\}/\bR^+
 \label{pol4.1}
\end{equation}
is the natural extension of $\cM^{(2,4)}$. Let $K\subset P$
consists of all elements which are orthogonal to $L$. Looking at
the diagrams $\Gamma$ in Table 1, one can see that $K$ is elliptic
and has $\rk S-2$ elements. By Proposition 2 in \cite{2}, it is
enough to show that $(L\cup K)^\perp\cap
\widetilde{\cM^{(2,4)}}=\emptyset$ (it then implies that
$L^\perp\cap \widetilde{\cM^{(2,4)}}=\emptyset$). Since $\#K=\rk
S-2$, the $K^\perp\cap \widetilde{\cM^{(2,4)}}$ is the edge
(1-dimensional) $r_1$ of $\cM^{(2,4)}$. There are two more
elements $f_1\in P$ and $f_2\in P$ such that $K_1=K\cup \{f_1\}$
and $K_2=K\cup \{f_2\}$ are elliptic. It follows that the edge
$r_1$ terminates in two vertices $A_1$ and $A_2$ of $\cM^{(2,4)}$
which are orthogonal to $K_1$ and $K_2$ respectively. Any element
$\bR^+x\in r_1$ then has $x^2\ge 0$. It follows that $(x,L)\not=0$
because $L$ is a hyperbolic subset. It follows that $(L\cup
K)^\perp\cap \widetilde{\cM^{(2,4)}}=\emptyset$. Thus,
$\cM^{(2,4)}$ has finite volume for $N=45$, $47$ either.

\medskip

\subsubsection{Fundamental chambers $\cM^{(2,4)}_+$ of cases
$N=7$ (Table~2).} \label{fundchambersN=7} We use orthogonal basis
$h$, $v_1,\dots,\,v_6$ of $S\otimes \bQ$ where $h^2=8$,
$v_1^2=\cdots=v_6^2=-2$.

{\bf Case 7a.} As $P(\cM^{(2,4)})$, we can take $f_1=v_1-v_2$,
$f_2=v_2-v_3$, $f_3=v_3-v_4$, $f_4=v_4-v_5$, $f_5=v_4+v_5$ with
square $(-4)$ defining the root system $D_5$, and
$e=(-v_1-v_2-v_3-v_4+v_5)/2+h/4$ with square $(-2)$. They define
the diagram 7a. The Weyl group $W=W(D_5)$ (generated by
reflections in $f_1,\dots,f_5$) is the semi-direct product of
permutations of $v_i$, $1\le i\le 5$, and linear maps $v_i\to (\pm
)_iv_i$, $1\le i\le 5$, where $\prod_i{(\pm)_i}=1$, see \cite{1}.
It follows that $W(e)$ consists of
$$
e_{i_1i_2\dots i_k}=(\pm v_1\pm v_2\pm \cdots \pm v_5)/2+h/4
$$
where $1\le i_1<i_2<\cdots <i_k\le 5$ show where the signs
$(-)$ are placed, and $k\equiv 0\mod 2$ (their number is $16$ which is the
number of exceptional curves on non-singular del Pezzo surface of
degree $4$), e.g. we have $e=e_{1234}$.

{\bf Case 7b.} As the basis of the root subsystem $2A_1\oplus
A_3\subset D_5$, we take $f_1=v_1-v_2$, $f_6=v_1+v_2$,
$f_3=v_3-v_4$, $f_4=v_4-v_5$, $f_5=v_4+v_5$. Only
$e_{1234}=(-v_1-v_2-v_3-v_4+v_5)/2+h/4$,
$e_{1345}=(-v_1+v_2-v_3-v_4-v_5)/2+h/4$ (from the orbit $W(e)$)
have non-negative pairing with this basis. We obtain the diagram
7b of Table~2.

\subsubsection{Fundamental chambers $\cM^{(2,4)}_+$ of cases
$N=8$ (Table~2).} \label{fundchambersN=8} We use the orthogonal
basis $h$, $v_1,\dots,\,v_8$ over $\bQ$ with $h^2=6$,
$v_1^2=\cdots=v_8^2=-2$. As root system $E_6$ we can take (see
\cite{1}) all roots $\pm v_i\pm v_j$ ($1\le i<j\le 5$) and $\pm
\frac{1}{2}(\pm v_1\pm v_2 \pm v_3 \pm v_4 \pm v_5 -v_6-v_7+v_8)$
with even number of $(-)$. I. e. $E_6\subset E_8$ consists of all
roots in $E_8$ which are orthogonal to roots $v_6-v_7$ and
$v_7+v_8$ (they define $A_2$). We denote $W=W(E_6)$, the Weyl
group of $E_6$.

{\bf Case 8a.} As $P(\cM^{(2,4)})$, we can take
$f_1=(v_1-v_2-v_3-v_4-v_5-v_6-v_7+v_8)/2$, $f_2=v_1+v_2$,
$f_3=-v_1+v_2$, $f_4=-v_2+v_3$, $f_5=-v_3+v_4$, $f_6=-v_4+v_5$
(with square $-4$) defining the basis of the root system $E_6$,
and
$$
e=-v_5+\frac{1}{3}v_6+\frac{1}{3}v_7-\frac{1}{3}v_8+\frac{1}{3}h=
-v_5+v_6-\frac{2}{3}(v_6-v_7)-\frac{1}{3}(v_7+v_8)+\frac{1}{3}h
$$
(with square $-2$). They define the diagram 8a.

We have
$$
W(e)=W(-v_5+v_6)-\frac{2}{3}(v_6-v_7)-\frac{1}{3}(v_7+v_8)+\frac{1}{3}h
$$
where $W(-v_5+v_6)$ consists of all roots $\alpha$ of $E_8$ with
the properties: $(\alpha, v_6-v_7)=-2$ and $(\alpha, v_7+v_8)=0$.
Thus, $W(e)$ consists of all elements
$$
e_{\pm i}=\pm
v_i+v_6-\frac{2}{3}(v_6-v_7)-\frac{1}{3}(v_7+v_8)+\frac{1}{3}h,\
1\le i\le 5;
$$
$$
e_{i_1i_2\dots i_k}=
$$
$$\frac{1}{2}(\pm v_1\pm v_2\pm v_3\pm
v_4\pm
v_5+v_6-v_7+v_8)-\frac{2}{3}(v_6-v_7)-\frac{1}{3}(v_7+v_8)+\frac{1}{3}h
$$
where $1\le i_1<\cdots<i_k\le 5$ show where are $(-)$, and
$k\equiv 1\mod 2$;
$$
e_{78}=-v_7+v_8-\frac{2}{3}(v_6-v_7)-\frac{1}{3}(v_7+v_8)+\frac{1}{3}h.
$$
(Their number is 27, the number of lines on a non-singular cubic.)

{\bf Case 8b.} As a basis of $A_5\oplus A_1\subset E_6$ we can
take $f_1=(v_1-v_2-v_3-v_4-v_5-v_6-v_7+v_8)/2$, $f_3=-v_1+v_2$,
$f_4=-v_2+v_3$, $f_5=-v_3+v_4$, $f_6=-v_4+v_5$ and
$f_7=(-v_1-v_2-v_3-v_4-v_5+v_6+v_7-v_8)/2$. Only $e_{+1}=
v_1+v_6-\frac{2}{3}(v_6-v_7)-\frac{1}{3}(v_7+v_8)+\frac{1}{3}h$
and $e_5=(v_1+v_2+v_3+v_4-v_5+v_6-v_7+v_8)/2
-\frac{2}{3}(v_6-v_7)-\frac{1}{3}(v_7+v_8)+\frac{1}{3}h$ have
non-negative pairing with this basis. They define the diagram 8b
of Table~2.

{\bf Case 8c.} As a basis of $3A_2\subset E_6$ we can take
$f_1=(v_1-v_2-v_3-v_4-v_5-v_6-v_7+v_8)/2$, $f_3=-v_1+v_2$;
$f_5=-v_3+v_4$, $f_6=-v_4+v_5$; $f_2=v_1+v_2$,
$f_7=(-v_1-v_2-v_3-v_4-v_5+v_6+v_7-v_8)/2$. Only $e_{+3}=
v_3+v_6-\frac{2}{3}(v_6-v_7)-\frac{1}{3}(v_7+v_8)+\frac{1}{3}h$,
$e_{125}=(-v_1-v_2+v_3+v_4-v_5+v_6-v_7+v_8)/2
-\frac{2}{3}(v_6-v_7)-\frac{1}{3}(v_7+v_8)+\frac{1}{3}h$,
$e_{2}=(v_1-v_2+v_3+v_4+v_5+v_6-v_7+v_8)/2
-\frac{2}{3}(v_6-v_7)-\frac{1}{3}(v_7+v_8)+\frac{1}{3}h$ have
non-negative pairing with this basis. We obtain the diagram 8c of
Table~2.

\subsubsection{Fundamental chambers $\cM^{(2,4)}_+$ of cases
$N=9$} \label{fundchambersN=9} We use the orthogonal basis $h$,
$v_1,\dots,\,v_8$ over $\bQ$ with $h^2=4$,
$v_1^2=\cdots=v_8^2=-2$. As a root system $E_7$ we can take (see
\cite{1}) all roots $\pm v_i\pm v_j$ ($1\le i<j\le 6$), $\pm
(v_7-v_8)$, and $\pm \frac{1}{2}(\pm v_1\pm v_2 \pm v_3 \pm v_4
\pm v_5 \pm v_6+v_7-v_8)$ with even number of $(-)$. I. e.
$E_7\subset E_8$ consists of all roots in $E_8$ which are
orthogonal to the root $v_7+v_8$. We denote $W=W(E_7)$, the Weyl
group of $E_7$.

{\bf Case 9a.} As $P(\cM^{(2,4)})$ we can take
$f_1=(v_1-v_2-v_3-v_4-v_5-v_6-v_7+v_8)/2$, $f_2=v_1+v_2$,
$f_3=-v_1+v_2$, $f_4=-v_2+v_3$, $f_5=-v_3+v_4$, $f_6=-v_4+v_5$,
$f_7=-v_5+v_6$ (with square $-4$) defining the basis of the root
system $E_7$, and
$$
e=-v_6+\frac{1}{2}v_7-\frac{1}{2}v_8+\frac{1}{2}h=
-v_6+v_7-\frac{1}{2}(v_7+v_8)+\frac{1}{2}h
$$
(with square $-2$). They define the diagram 9a of Table~2.

The orbit $W(e)=W(-v_6+v_7)-\frac{1}{2}(v_7+v_8)+\frac{1}{2}h$
where $W(-v_6+v_7)$ consists of all roots $\alpha$ in $E_8$ with
the property $(\alpha,v_7+v_8)=-2$. It follows that $W(e)$
consists of
$$
e_{\pm i7}=\pm v_i+v_7-\frac{1}{2}(v_7+v_8)+\frac{1}{2}h,\ 1\le
i\le 6;
$$
$$
e_{\pm i8}=\pm v_i+v_8-\frac{1}{2}(v_7+v_8)+\frac{1}{2}h,\ 1\le
i\le 6;
$$
$$
e_{i_1\dots i_k}=\frac{1}{2}(\pm v_1\pm v_2\pm \cdots \pm
v_6+v_7+v_8)-\frac{1}{2}(v_7+v_8)+\frac{1}{2}h
$$
where $1\le i_1<\cdots<i_k\le 6$ show $(-)$ and $k\equiv 0\mod 2$.
Their number is 56, the number of exceptional curves on a
non-singular del Pezzo surface of degree 2.

{\bf Case 9b.} As a basis of $A_7\subset E_7$ we can take
$f_8=v_7-v_8,\,f_1,\,f_3,\,f_4$, $f_5$, $f_6,\,f_7$. Only
$e_0=\frac{1}{2}(v_1+v_2+ \cdots +
v_6+v_7+v_8)-\frac{1}{2}(v_7+v_8)+\frac{1}{2}h$ and $e_{56}$ have
non-negative pairing with the basis. We obtain the diagram 9b of
Table~2.

{\bf Case 9c.} As a basis of $A_2\oplus A_5\subset E_7$ we can
take $f_8=v_7-v_8,\,f_1$ and $f_2,\,f_4,\,f_5,\,f_6,\,f_7$. Only
$e_{-18}$, $e_{16}$, $e_{1456}$ have non-negative pairing with
this basis. We obtain the diagram 9c of Table~2.

{\bf Case 9d.} As a basis of $A_3\oplus A_1\oplus A_3\subset E_7$
we can take $f_8=v_7-v_8,\,f_1\,f_3$, and $f_2$, and
$f_5,\,f_6,\,f_7$. Only $e_{+38}$, $e_{26}$, $e_{12}$, $e_{1256}$
have non-negative pairing with this basis. We obtain the diagram
9d of Table~2.

{\bf Case 9e.} As a basis of $A_1\oplus D_6\subset E_7$ we can
take $f_8=v_7-v_8$ and $f_2,\,f_3,\,f_4,\,f_5,\,f_6,\,f_7$. Only
$e_{-68}$ and $e_{123456}$ have non-negative pairing with the
basis. We obtain the diagram 9e of Table~2.

{\bf Case 9f.} As a basis of $D_4\oplus 3A_1\subset E_7$ we can
take $f_2,\,f_3,\,f_4,\,f_5$ and $f_7$, $f_8=v_7-v_8$,
$f_9=-v_5-v_6$. Only $e_{-48}$, $e_{+58}$, $e_{1234}$, $e_{2346}$
have non-negative pairing with the basis. We obtain the diagram 9f
of Table~2.

{\bf Case 9g.} As a basis of $7A_1\subset E_7$ we can take
$u_1=v_1+v_2$, $u_2=-v_1+v_2$, $u_3=v_3+v_4$, $u_4=-v_3+v_4$,
$u_5=v_5+v_6$, $u_6=-v_5+v_6$, $u_7=v_7-v_8$. Only $e_{-28}$,
$e_{-48}$, $e_{-68}$, $e_{2456}$, $e_{2346}$, $e_{1246}$,
$e_{123456}$ have non-negative pairing with the basis. We obtain
the digram 9g described in Section \ref{subsubsec3.4.4}.

\subsubsection{Fundamental chambers $\cM^{(2,4)}_+$ of cases
$N=10$} \label{fundchambersN=10} We use the orthogonal basis $h$,
$v_1,\dots,\,v_8$ of $S\otimes \bQ$ with $h^2=2$,
$v_1^2=\cdots=v_8^2=-2$. As a root system $E_8$ we can take (see
\cite{1}) all roots $\pm v_i\pm v_j$ ($1\le i<j\le 8$) and
$\frac{1}{2}(\pm v_1\pm v_2\pm\cdots\pm v_8)$ with even number of
$(-)$. We denote $W=W(E_8)$, the Weyl group of $E_8$.

{\bf Case 10a.} As $P(\cM^{(2,4)})$ we can take
$f_1=(v_1-v_2-v_3-v_4-v_5-v_6-v_7+v_8)/2$, $f_2=v_1+v_2$,
$f_3=-v_1+v_2$, $f_4=-v_2+v_3$, $f_5=-v_3+v_4$, $f_6=-v_4+v_5$,
$f_7=-v_5+v_6$, $f_8=-v_6+v_7$ (with square $-4$) defining the
basis of the root system $E_8$, and $e=-v_7-v_8+h$. They define
the diagram 10a of Table~2.

The orbit $W(e)=W(-v_7-v_8)+h$ where $W(-v_7-v_8)$ consists of all
roots $\alpha$ in $E_8$. It follows that $W(e)$ consists of
$$
e_{\pm i,\pm j}=\pm v_i\pm v_j+h,\ 1\le i<j\le 8;
$$
$$
e_{i_1\dots i_k}=\frac{1}{2}(\pm v_1\pm v_2\pm \cdots \pm v_8)+h
$$
where $1\le i_1<\cdots<i_k\le 8$ show $(-)$ and $k\equiv 0\mod 2$.
Their number is 240, the number of exceptional curves on a
non-singular del Pezzo surface of degree 1.

{\bf Case 10b.} As a basis of $A_8\subset E_8$ we can take
$f_1,\,f_3,\,f_4,\,f_5,\,f_6,$ $f_7,\,f_8,\,f_9=-v_7-v_8$. Only
$e_{+1,+2}$, $e_{0}=(v_1+\cdots+v_8)/2+h$, $e_{67}$ have
non-negative pairing with the basis. We obtain the diagram 10b of
Table~2.

{\bf Case 10c.} As a basis of $A_1\oplus A_7\subset E_8$ we can
take $f_1$ and $f_2$, $f_4$, $f_5$, $f_6$, $f_7$, $f_8$,
$f_9=-v_7-v_8$. Only $e_{-1,+2}$, $e_{-1,+8}$, $e_{18}$, $e_{17}$,
$e_{1567}$ have non-negative pairing with the basis. We obtain the
diagram 10c of Table~2.

{\bf Case 10d.} As a basis of $A_2\oplus A_1\oplus A_5\subset E_8$
we can take $f_1,\,f_3$, and $f_2$, and
$f_5,\,f_6,\,f_7,\,f_8,\,f_9=-v_7-v_8$. Only $e_{-1,-2}$,
$e_{-2,+3}$, $e_{+3,+4}$, $e_{+3,+8}$, $e_{12}$, $e_{27}$,
$e_{28}$, $e_{1267}$ have non-negative pairing with the basis. We
obtain the diagram 10d of Table~2.

{\bf Case 10e.} As a basis of $A_4\oplus A_4\subset E_8$ we can
take $f_1,\,f_3,\,f_4,\,f_2$ and $f_6,\,f_7,\,f_8,\,f_9=-v_7-v_8$.
Only $e_{-3,+4}$, $e_{+4,+5}$, $e_{+4,+8}$, $e_{1237}$,
$e_{1238}$, $e_{23}$ have non-negative pairing with the basis. We
obtain the diagram 10d of Table~2.

{\bf Case 10f.} As a basis of $D_8\subset E_8$ we can take
$f_2,\,f_3,\,f_4,\,f_5,f_6,\,f_7,\,f_8$, $f_9=-v_7-v_8$. Only
$e_{-7,+8}$, $e_{234567}$ have non-negative pairing with the
basis. We obtain the diagram 10f of Table~2.

{\bf Case 10g.} As a basis of $D_5\oplus A_3\subset E_8$ we can
take $f_1\,f_2,\,f_3,\,f_4,\,f_5$ and $f_7,\,f_8$, $f_9=-v_7-v_8$.
Only $e_{-4,+5}$, $e_{+5,+6}$, $e_{+5,+8}$, $e_{1234}$, $e_{2348}$
have non-negative pairing with the basis. We obtain the diagram
10g of Table~2.

{\bf Case 10h.} As a basis of $E_6\oplus A_2\subset E_8$ we can
take $f_1$, $f_2$, $f_3$, $f_4$, $f_5$, $f_6$, and $f_8$, and
$f_9=-v_7-v_8$. Only $e_{-5,+6}$, $e_{+6,+7}$, $e_{+6,+8}$,
$e_{123458}$ have non-negative pairing with the basis. We obtain
the diagram 10h of Table~2.

{\bf Case 10i.} As a basis of $E_7\oplus A_1\subset E_8$ we can
take $f_1$, $f_2$, $f_3$, $f_4$, $f_5$, $f_6$, $f_7$, and
$f_9=-v_7-v_8$. Only $e_{-6,+7}$, $e_{+7,-8}$, $e_{+7,+8}$ have
non-negative pairing with the basis. We obtain the diagram 10i of
Table~2.

{\bf Case 10j.} As a basis of $2A_1\oplus D_6\subset E_8$ we can
take $f_2$ and $f_3$, and $f_5$, $f_6$, $f_7$, $f_8$,
$f_9=-v_7-v_8$, $f_{10}=-v_7+v_8$. Only $e_{-1,-2}$, $e_{+1,-2}$,
$e_{-2,+3}$, $e_{+3,+4}$, $e_{12}$, $e_{28}$ have non-negative
pairing with the basis. We obtain the diagram 10j of Table~2.

{\bf Case 10k.} As a basis of $2D_4\subset E_8$ we can take $f_2$,
$f_3$, $f_4$, $f_5$ and $f_7$, $f_8$, $f_9=-v_7-v_8$,
$f_{10}=-v_7+v_8$. Only $e_{-3,-4}$, $e_{-4,+5}$, $e_{+5,+6}$,
$e_{1234}$, $e_{2348}$ have non-negative pairing with the basis.
We obtain the diagram 10k of Table~2.

{\bf Case 10l.} As a basis of $2A_1\oplus 2A_3\subset E_8$ we can
take $f_2$; $f_3$; $f_5$, $f_6$, $f_{11}=v_3+v_4$; $f_9=-v_7-v_8$,
$f_8$, $f_{10}=-v_7+v_8$. Only $e_{-1,-2}$, $e_{+1,-2}$,
$e_{-2,-5}$, $e_{-2,+6}$, $e_{-4,-5}$, $e_{-5,+6}$, $e_{+6,+7}$,
$e_{123458}$, $e_{1245}$, $e_{2345}$, $e_{2458}$ have non-negative
pairing with the basis. We obtain the diagram 10l of Table~2.

{\bf Case 10m.} As a basis of $4A_2\subset E_8$ we can take $f_1$,
$f_3$; $f_2$,
$f_{10}=\frac{1}{2}(-v_1-v_2-v_3-v_4-v_5+v_6+v_7-v_8)$; $f_5$,
$f_6$; $f_8$, $f_9=-v_7-v_8$. Only $e_{-2,+3}$, $e_{+3,+4}$,
$e_{+3,-5}$, $e_{+3,+6}$, $e_{+3,+8}$, $e_{+6,+8}$, $e_{12}$,
$e_{1257}$, $e_{1267}$, $e_{25}$, $e_{27}$, $e_{28}$ have
non-negative pairing with the basis. We obtain the diagram 10m of
Table~2.

{\bf Case 10n.} As a basis of $D_4\oplus 4A_1\subset E_8$ we can
take $u_1=v_1+v_2$, $u_2=-v_1+v_2$, $u_3=-v_3+v_4$,
$u_4=-v_2+v_3$, $u_5=v_5+v_6$, $u_6=-v_5+v_6$, $u_7=-v_7+v_8$,
$u_8=v_7+v_8$ (this basis agrees with the one used in case 10n of
Lemma \ref{lemma3.4.4} if one replaces $f_i$ by $u_i$). Only
$e_{-3,-4}$, $e_{-4,-6}$,  $e_{-4,-8}$, $e_{-5,-6}$, $e_{+5,-6}$,
$e_{-6,-8}$, $e_{-7,-8}$, $e_{+7,-8}$, $e_{12345678}$,
$e_{123468}$, $e_{234568}$, $e_{234678}$ have non-negative pairing
with the basis. We obtain the diagram 10n of Figure \ref{fig3}.

{\bf Case 10o.} As a basis of $8A_1\subset E_8$ we can take
$u_1=-v_1+v_2$, $u_2=v_1+v_2$, $u_3=-v_3+v_4$, $u_4=v_3+v_4$,
$u_5=-v_5+v_6$, $u_6=v_5+v_6$, $u_7=-v_7+v_8$, $u_8=v_7+v_8$. The
set of indices $I=1,\,\dots,\,8$ has the structure of a
3-dimensional affine space over $F_2$ with  (affine) planes
$J\subset I$ determined by the property $\frac{1}{2}\sum_{j\subset
J}{u_j}\in E_8$. It is the same as the one used in case 10o of
Lemma \ref{lemma3.4.4}. For $i\in I$ we set $w_i=e_{+i,-(i+1)}$ if
$i$ is odd, and $w_i=e_{-(i-1),-i}$ if $i$ is even. For a plane
$\pi\subset I$ we set $w_\pi=-\frac{1}{2}\sum_{i\in \pi}{u_i}+h$.
The introduced elements $w_i$, $i\in I$, and  $w_\pi$, $\pi\subset
I$ is a plane, are the only elements (from the orbit $W(e)$) which
have non-negative pairing with the basis. We obtain the graph 10o
described in Section \ref{subsubsec3.4.4}.

{\bf Case of $7A_1\subset E_8$.} As a basis of $7A_1\subset E_8$
we can take $u_2=v_1+v_2$, $u_3=-v_3+v_4$, $u_4=v_3+v_4$,
$u_5=-v_5+v_6$, $u_6=v_5+v_6$, $u_7=-v_7+v_8$, $u_8=v_7+v_8$. We
denote $u_1=-v_1+v_2$ (the roots $\pm u_1$ are the only roots of
$E_8$ which are orthogonal to $7A_1\subset E_8$. The set of
indices $I=1,\,\dots,\,8$ has the structure of a 3-dimensional
affine space over $F_2$ with (affine) planes $J\subset I$
determined by the property $\frac{1}{2}\sum_{j\in J}{u_j}\in E_8$.
Taking $1\in I$ as an origin, makes the set $I$ to be a
3-dimensional vector space over $F_2$. As in the previous case, we
define $w_i$ for $i\in I$, and $w_\pi$ for an affine plane $\pi
\subset I$. We set $w_1^{(+)}=w_1=e_{+1,-2}$ and
$w_1^{(-)}=e_{-1,+2}$. If $1\in \pi$, we set $w_\pi^{(+)}=w_\pi$,
and $w_\pi^{(-)}=w_\pi+u_1$. The introduced elements $w_i$, $i\in
I-\{1\}$, $w_0^{(+)}$, $w_0^{(-)}$, $w_\pi$ for planes $\pi\subset
I-\{1\}$, and $w_\pi^{(+)}$, $w_\pi^{(-)}$ for planes $\pi\subset
I$ containing $1$ are the only elements (from the orbit $W(e)$)
which have non-negative pairing with the basis. We obtain the
graph described and used in Section \ref{subsubsec3.4.7} (cases 10n
and 10o).

\subsubsection{Fundamental chambers $\cM^{(2,4)}_+$ of cases
$N=20$} \label{fundchambersN=20} This case had been partly
described (including cases 20a and 20b below) at the end of Section
\ref{subsubsec3.4.4}; here we add further details of calculations
for readers' convenience. We use the orthogonal basis $h$,
$\alpha$, $v_1,\dots,\,v_8$ of $S\otimes \bQ$ with $h^2=2$,
$\alpha^2=v_1^2=\cdots=v_8^2=-2$. As a root system $D_8$ we can
take (see \cite{1}) all roots $\pm v_i\pm v_j$ ($1\le i<j\le 8$).
We denote $W=W(D_8)$, the Weyl group of $D_8$.

{\bf Case 20a.} As $P(\cM^{(2,4)})$ we can take $f_1=v_1-v_2$,
$f_2=v_2-v_3$, $f_3=v_3-v_4$, $f_4=v_4-v_5$, $f_5=v_5-v_6$,
$f_6=v_6-v_7$, $f_7=v_7-v_8$, $f_8=v_7+v_8$ (with square $-4$)
defining the basis of the root system $D_8$, and $\alpha$,
$b=\frac{h}{2}-\frac{\alpha}{2}-v_1$,
$c=h-\frac{1}{2}(v_1+v_2+\cdots +v_8)$ (with square $-2$). They
define the diagram 20a of Table~2.

The orbit $W(\alpha)$ consists of only $\alpha$; the orbit $W(b)$
consists of all
$$
b_{\pm i}=\frac{h}{2}-\frac{\alpha}{2}\pm v_i,\ 1\le i\le 8;
$$
the orbit $W(c)$ consists of all
$$
c_{i_1\dots i_k}=h+\frac{1}{2}(\pm v_1\pm v_2\pm\cdots \pm v_8)
$$
where $1\le i_1<i_2<\cdots <i_k\le 8$ show where are $(-)$, and
$k\equiv 0\mod 2$. Thus, $P(\cM^{(2)})$ has $1+2\cdot 8+2^7=81$
elements. This is the number of exceptional curves on the right
DPN surface with the main invariants $(r,a,\delta)=(10,8,1)$ and
the zero root invariant. One of them (corresponding to $\alpha$)
has square $(-4)$, all other are $(-1)$-curves.

{\bf Case 20b.} As a basis of $2A_1\oplus D_6\subset D_8$ we can
take $f_1$; $f_9=-v_1-v_2$; $f_3$, $f_4$, $f_5$, $f_6$, $f_7$,
$f_8$. Only $\alpha$, $b_{+2}$, $b_{-3}$, $c_{134567}$,
$c_{345678}$ have non-negative pairing with the basis. We obtain
the diagram 20b of Table~2.

{\bf Case 20c.} As a basis of $A_3\oplus D_5\subset D_8$ we can
take $f_1$, $f_2$, $f_9=-v_1-v_2$ and $f_4$, $f_5$, $f_6$, $f_7$,
$f_8$. Only $\alpha$, $b_{+3}$, $b_{-4}$, $c_{145678}$, $c_{4567}$
have non-negative pairing with the basis. We obtain the diagram
20c of Table~2.

{\bf Case 20d.} As a basis of $2D_4\subset D_8$ we can take $f_1$,
$f_2$, $f_3$, $f_9=-v_1-v_2$ and $f_5$, $f_6$, $f_7$, $f_8$. Only
$\alpha$, $b_{+4}$, $b_{-5}$, $c_{1567}$, $c_{5678}$ have
non-negative pairing with the basis. We obtain the diagram 20d of
Table~2.

{\bf Case 20e.} As a basis of $2A_1\oplus 2A_3\subset D_8$ we can
take $f_1$; $f_9=-v_1-v_2$; $f_3$, $f_4$, $f_{10}=-v_3-v_4$;
$f_6$, $f_7$, $f_8$. Only $\alpha$, $b_{+2}$, $b_{+5}$, $b_{-6}$,
$c_{1367}$, $c_{1678}$, $c_{3678}$, $c_{67}$ have non-negative
pairing with the basis. We obtain the diagram 20e of Figure
\ref{fig4}.

{\bf Case 20f.} As a basis of $4A_1\oplus D_4\subset D_8$ we can
take $f_1$; $f_9=v_1+v_2$; $f_3$; $f_{10}=v_3+v_4$; $f_5$, $f_6$,
$f_7$, $f_8$. Only $\alpha$, $b_{-1}$, $b_{-3}$, $b_{-5}$,
$c_{12345678}$, $c_{123567}$, $c_{134567}$, {$c_{135678}$ have
non-negative pairing with the basis. We obtain the diagram 20f of
Figure \ref{fig5}.

{\bf Case 20g.} As a basis of $8A_1\subset D_8$ we can take
$u_1=v_1-v_2$; $u_2=v_1+v_2$; $u_3=v_3-v_4$; $u_4=v_3+v_4$;
$u_5=v_5-v_6$; $u_6=v_5+v_6$; $u_7=v_7-v_8$; $u_8=v_7+v_8$. Only
$\alpha$, $b_{-1}$, $b_{-3}$, $b_{-5}$, $b_{-7}$, $c_{1357}$,
$c_{135678}$, $c_{134578}$, $c_{134567}$, $c_{123578}$,
$c_{123567}$, $c_{123457}$, $c_{12345678}$ have non-negative
pairing with the basis. We obtain the diagram 20g of Figure
\ref{fig6}.

{\bf Case $4A_1\oplus A_3\subset D_8$.} As a basis of $4A_1\oplus
A_3\subset D_8$ we can take $u_1=v_1-v_2$; $u_2=-v_1-v_2$;
$u_3=v_3-v_4$; $u_4=-v_3-v_4$; $u_5=v_5-v_6$, $u_6=-v_5-v_6$,
$u_7=v_6+v_7$. Only $\alpha$, $b_{+2}$, $b_{+4}$, $b_{-7}$,
$b_{-8}$, $b_{+8}$ and
$$
c_{(\mu_1,\mu_3,\mu_5)}= \frac{1}{2}\left((-1)^{\mu_1}
v_1+v_2+(-1)^{\mu_3} v_3+v_4+(-1)^{\mu_5}v_5+v_6-\right.
$$
$$
\left. -v_7+(-1)^{(\mu_1+\mu_3+\mu_5+1)}v_8\right)+h,
$$
$(\mu_1,\mu_3,\mu_5)\in (\bZ/2\bZ)^3$, have non-negative pairing
with the basis. We obtain the diagram which had been described in
Case 20e of Section \ref{subsubsec3.4.7}.

In a usual way, we identify $(\bZ/2\bZ)^3$ with the set of
vertices $V(K)$ of a 3-dimensional cube $K$. Thus the last set is
$c_v$, $v\in V(K)$. Each $u_i$, $1\le i \le 6$, defines a
2-dimensional face $\gamma_i$ of the cube $K$ which consists of
$c_v$ such that $(u_i,c_v)=2$. Therefore, we further write
$u_i=u_{\gamma_i}$, where $\gamma_i$ belongs to the set $\gamma(K)$
of 2-dimensional faces of $K$. The $u_7$ defines two distinguished
opposite faces $\gamma_5,\,\gamma_6\in \gamma(K)$ characterized by
$(u_5,u_7)=(u_6,u_7)=2$ (i. e. $u_5$, $u_7$, $u_6$ define the
component $A_3$).

We identify $b_{-7}$ with the pair $\{\gamma_5, \gamma_6\}\subset
\gamma(K)$ of distinguished opposite faces of $K$. We identify
$b_{+2}$ and $b_{+4}$ with the pairs $\{\gamma_1,\gamma_2\}$ and
$\{\gamma_3,\gamma_4\}$ of opposite faces of $K$. Thus, we further
numerate them as $b_{\overline{\gamma}}$, where $\overline{\gamma}$
belongs to the set $\overline{\gamma(K)}$ of pairs of opposite
2-dimensional faces of $K$. We have
$(u_{\gamma^\prime},b_{\overline{\gamma}})=2$ if $\gamma^\prime
\in \overline{\gamma}$ and $\overline{\gamma}\in
\overline{\gamma(K)}$ is different from the pair of two
distinguished opposite faces of $K$; otherwise it is~$0$.

We can identify $b_{-8}$ (respectively $b_{+8})$ with the four
vertices $c_{(\mu_1,\mu_3,\mu_5)}$ where $\mu_1+\mu_3+\mu_5+1\equiv
0\mod 2$ (respectively $\equiv 1\mod 2$). Each of these four vertices
contains one vertex from any two opposite vertices of $K$, no
three of its vertices belong to a face of $K$. Thus, we can further
denote $b_{\pm 8}$ by $b_t$ where $t$ belongs to the set
$\overline{V(K)}$ of these two fours. We have $(c_v,b_t)=2$ if $v\in
t\in \overline{V(K)}$. Otherwise, it is $0$.

{\bf Case $7A_1\subset D_8$.} As a basis of $7A_1\subset D_8$ we
can take $u_1=v_1-v_2$; $u_2=-v_1-v_2$; $u_3=v_3-v_4$;
$u_4=-v_3-v_4$; $u_5=v_5-v_6$; $u_6=-v_5-v_6$; $u_7=v_7-v_8$. Only
$\alpha$, $b_{+2}$, $b_{+4}$, $b_{+6}$, $b_{-7}$, $b_{+8}$ and
$$
c_{j_1j_2j_3(\pm 2)}= \frac{1}{2}\left((-1)^{j_1}
v_1+v_2+(-1)^{j_2} v_3+v_4+(-1)^{j_3}v_5+\right.
$$
$$
\left. +v_6\pm (v_7+v_8)\right)+h
$$
where $j_i=1$ or $2$, $1\le i \le3$, and $j_1+j_2+j_3+(\pm
2)\equiv 0\mod 2$, and
$$
c_{j_1j_2j_3 1}= \frac{1}{2}\left((-1)^{j_1} v_1+v_2+(-1)^{j_2}
v_3+v_4+(-1)^{j_3}v_5+\right.
$$
$$
\left. +v_6-v_7+v_8\right)+h
$$
where $j_i=1$ or $2$, $1\le i \le3$, and $j_1+j_2+j_3+1\equiv
0\mod 2$, have non-negative pairing with the basis. We obtain the
diagram which had been described in Case 20f,g of Section
\ref{subsubsec3.4.7}.

One can denote $f_{11}=u_1$, $f_{12}=u_2$, $f_{21}=u_3$,
$f_{22}=u_4$, $f_{31}=u_5$, $f_{32}=u_6$, $f_{41}=u_7$,
$b_1=b_{+2}$, $b_2=b_{+4}$, $b_3=b_{+6}$, $b_{4(-)}=b_{-7}$,
$b_{4(+)}=b_{+8}$


\providecommand{\bysame}{\leavevmode\hbox to3em{\hrulefill}\thinspace}
\providecommand{\MR}{\relax\ifhmode\unskip\space\fi MR }
\providecommand{\MRhref}[2]{%
  \href{http://www.ams.org/mathscinet-getitem?mr=#1}{#2}
}
\providecommand{\href}[2]{#2}

\end{document}